\documentclass[10pt,oneside,leqno]{amsart}

\makeatletter
\@namedef{subjclassname@2020}{\textup{2020} Mathematics Subject Classification}
\makeatother

\usepackage[top=2.5 cm,bottom=2 cm,left=2 cm,right=3 cm]{geometry}
\usepackage[utf8]{inputenc}
\usepackage{color}
\usepackage{amssymb}

\usepackage[english]{babel}

\usepackage{todonotes}

\usepackage[backref=page]{hyperref}

\usepackage{esint}
\usepackage{graphicx}
\usepackage{hyperref}
\usepackage{bm}
\usepackage{xcolor}
\usepackage[skins,breakable]{tcolorbox}

\usepackage{enumerate}

\usepackage{tikz}
\usetikzlibrary{arrows.meta, positioning, shapes.geometric}

\newcommand*{\mailto}[1]{\href{mailto:#1}{\nolinkurl{#1}}}

\numberwithin{equation}{section}

\newtheorem{example}{Example}[section]
\newtheorem{definition}[example]{Definition}
\newtheorem{theorem}[example]{Theorem}
 \newtheorem{proposition}[example]{Proposition}
\newtheorem{lemma}[example]{Lemma}
 \newtheorem{corollary}[example]{Corollary}
\newtheorem{remark}[example]{Remark}
\newtheorem*{maintheorem*}{Main Theorem}
\allowdisplaybreaks
\numberwithin{equation}{section}

\newcommand\so{
	\mathchoice
	{{\scriptstyle\mathcal{O}}}
	{{\scriptstyle\mathcal{O}}}
	{{\scriptscriptstyle\mathcal{O}}}
	{\scalebox{0.6}{$\scriptscriptstyle\mathcal{O}$}}
}

\renewcommand{\i}{\ifmmode\mathit{\mathchar"7010 }\else\char"10 \fi}
\renewcommand{\j}{\ifmmode\mathit{\mathchar"7011 }\else\char"11 \fi}
\newcommand{\R}{\mathbb{R}}
\newcommand{\N}{\mathbb{N}}
\newcommand{\Z}{\mathbb{Z}}
\DeclareMathOperator*{\sign}{sign}

\newcommand{\norm}[1]{\left\|#1\right\|}

\newcommand{\px}{\partial_x}
\newcommand{\py}{\partial_y}
\newcommand{\pt}{\partial_t}
\newcommand{\pxi}{\partial_\xi}

\newcommand{\pve}{\partial_\varepsilon}

\newcommand{\ve}{\varepsilon}

{%

\begin{enumerate}}%
{\end{enumerate}}

{%

\begin{enumerate}}%
{\end{enumerate}}

\begin{document}
\sloppy

\title[Regularity and Lipschitz metric for a Novikov system]
{Generic regularity and Lipschitz metric \\ for 
a two-component Novikov system}

\author[Karlsen]{K. H. Karlsen}
\address[Kenneth H. Karlsen]
{Department of Mathematics, University of Oslo, 
P.O.\ Box~1053, Blindern, NO-0316 Oslo, Norway}
\email{\mailto{kennethk@math.uio.no}}

\author[Rybalko]{Ya. Rybalko$^*$}
\address[Yan Rybalko]
{Department of Mathematics, University of Oslo, 
P.O.\ Box~1053, Blindern, NO-0316 Oslo, Norway}
\address{B.Verkin Institute for Low Temperature Physics and Engineering
of the National Academy of Sciences of Ukraine, 
47 Nauky Ave., Kharkiv, 61103, Ukraine}
\email{\mailto{rybalkoyan@gmail.com}}



\thanks{The work of Yan Rybalko was supported 
by the European Union’s Horizon Europe research 
and innovation programme under the 
Marie Sk\l{}odowska-Curie grant agreement No 101058830.
The work of Kenneth H. Karlsen was funded by the 
Research Council of Norway under project 351123 (NASTRAN)}

\date{\today}

\begin{abstract}
We investigate the Cauchy problem 
for a two-component generalization of the Novikov 
equation with cubic nonlinearity---an integrable system whose 
solutions may develop strong nonlinear 
phenomena such as gradient blow-up and 
interactions between peakon-like structures. 
Our study has two main objectives: first, to analyze the generic regularity 
of global conservative solutions; and second, to construct a new metric 
that guarantees the Lipschitz continuity of the flow. 
Building on the geometric framework developed 
by Bressan and Chen for quasilinear second-order wave equations, 
we prove that the solution retains $C^k$ regularity 
away from a finite number of piecewise $C^{k-1}$ characteristic curves.
Furthermore, we provide a description 
of the solution behavior in the vicinity of these curves. 
By introducing a Finsler norm on tangent vectors in the space of solutions, 
expressed in the transformed Bressan-Constantin variables, 
we introduce a Lipschitz metric representing 
the minimal energy transportation cost 
between two solutions.
\end{abstract}

\maketitle

{\small \tableofcontents}

\section{Introduction and main results}

In this paper we consider the following 
two-component Novikov system:
\begin{align}
	\nonumber
	\label{t-c-N}
	&m_t+(uvm)_x+u_xvm=0,
	&& m=m(t,x),\,u=u(t,x),\,v=v(t,x),\\
	&n_t+(uvn)_x+uv_xn=0,
	&& n=n(t,x),\\
	\nonumber
	&m=u-u_{xx},\,\,n=v-v_{xx},
	&& u,v\in\R,\quad
	t,x\in\R,
\end{align}
which was introduced in \cite{L19}
as an integrable vector generalization of the 
scalar Novikov equation \cite{N09}, having the 
form \eqref{t-c-N} with $v(t,x)=u(t,x)$:
\begin{equation}\label{N}
	m_t+(u^2m)_x+uu_xm=0,\quad m=u-u_{xx}.
\end{equation}
A substantial body of work has been devoted to 
various aspects of the scalar Novikov equation; 
see \cite{LS22} for a comprehensive overview of results 
related to \eqref{N}. Another important symmetry 
reduction of \eqref{t-c-N} is obtained by setting 
$v(t,x) = u(-t,-x)$, which leads to the following 
two-place (nonlocal) variant of the 
Novikov equation \cite{KR25}:
\begin{equation}\label{nN}
	m_t(t,x)
	+(u(t,x)u(-t,-x)m(t,x))_x
	+u_x(t,x)u(-t,-x)m(t,x)=0,\quad
	t,x\in\R.
\end{equation}

The two-place reductions of integrable systems were 
first studied by Ablowitz and Musslimani in \cite{AM13}, 
where they introduced, in particular, the 
nonlocal nonlinear Schr\"odinger (NNLS) equation
\begin{equation}
	\label{NNLS}
	\mathrm{i}q_t(t,x)+q_{xx}(t,x)
	+2\sigma q^{2}(t,x)\bar{q}(t,-x)=0,
	\quad q\in\mathbb{C},\quad
	\mathrm{i}^2=-1,\quad\sigma=\pm1,
\end{equation}
with $\bar{q}$ denoting the complex conjugate of $q$. 
Notice that both \eqref{nN} and \eqref{NNLS} 
possess the parity-time (PT) symmetry property, 
meaning that if $u(t,x)$ (respectively $q(t,x)$) is 
a solution of the two-place Novikov equation 
(respectively the NNLS equation), then so 
is $u(-t,-x)$ (respectively $\bar{q}(-t,-x)$).  
Two-place equations are particularly relevant 
because their dynamics are determined not only 
by local interactions but also by values of the solution 
at non-adjacent points. Such nonlocal dependence makes 
them well suited for modeling phenomena that involve 
strong spatial correlations or entanglement between 
events occurring at distant locations. Moreover, if the initial 
data $(u_0,v_0)(x)$ for \eqref{t-c-N} satisfy 
the symmetry condition $v_0(x)=u_0(-x)$, then 
the corresponding solution fulfills $v(t,x)=u(-t,-x)$ for 
all $t,x\in\R$, provided uniqueness holds.  For additional 
discussions and examples of two-place equations, 
see \cite{AM19,F16,LH17,LQ17} and the references therein.

Returning to the two-component Novikov system, 
it possesses a bi-Hamiltonian structure, a Lax pair, 
and an infinite hierarchy of conservation laws \cite{L19}. 
By means of inverse spectral methods, \cite{CS23} investigates 
the dynamics of multipeakon solutions of \eqref{t-c-N} 
(see also \cite{HLQ23, HLLQ24} for results on the 
stability of single-peakon solutions). 
In \cite{MM20}, the authors establish local 
well-posedness of \eqref{t-c-N} in suitable Besov spaces 
and show that blow-up can occur only 
through wave breaking (see also \cite{ZY25}). 
Assuming that the initial data $(m_0,n_0)$ satisfy 
certain sign assumptions, \cite{HQ20,ZQ21} obtain global 
solutions of the Novikov system. 
Using the method of characteristics and applying 
the Bressan-Constantin approach \cite{BC07}, the authors of 
\cite{HQ21} construct a global conservative weak 
solution of \eqref{t-c-N} subject to arbitrary initial data
$(u_0,v_0)\in \left(H^1\cap W^{1,4}\right)^2$.
In a recent refinement, \cite{KR25} revisits 
this method and establishes the existence of 
a global semigroup of conservative solutions 
to the two-component Novikov system, valid 
for a larger class of initial data 
$(u_0,v_0)\in\Sigma$, where $(\Sigma,d_\Sigma)$ 
denotes the complete---though not linear---metric 
space given by 
\begin{equation}\label{Sig}
	\Sigma=
	\left\{
	(f,g): f,g\in H^1(\mathbb{R})\mbox{ and }
	\int_{-\infty}^{\infty}
	f_x^2\,g_x^2\,dx<\infty
	\right\},
\end{equation}
and
\begin{equation}\label{d^2}
	d^2_\Sigma((f_1,g_1), (f_2,g_2))
	=\|f_1-f_2\|^2_{H^1}+\|g_1-g_2\|^2_{H^1}
	+\|f_{1x}\,g_{1x}-f_{2x}\,g_{2x}\|^2_{L^2}.
\end{equation}
Moreover, it is shown that in the 
Novikov system the measures $u_x^2\,dx$, $v_x\,dx$, 
and $u_x^2v_x^2\,dx$ may develop concentrations, 
whereas in the scalar case only $u_x^4\,dx$ 
exhibits such behavior \cite{CCL18}.  A summary of the 
results on the global semigroup of conservative weak solutions 
to \eqref{t-c-N} is provided in Section \ref{mr}.

\medskip

The main results of this work are grouped into 
three parts, which we now briefly outline before 
discussing them in detail in the following subsections.
All of these results are, to the best of our knowledge, 
new for the two-component Novikov system, 
and, in several aspects, they also sharpen or extend 
previously known results for the scalar equations. 
First, we establish the generic regularity of global 
conservative solutions, showing that the solution 
remains $C^k$ smooth away from a finite number of 
piecewise $C^{k-1}$ characteristic curves. 
Next, we analyze the behavior of these solutions 
near the characteristic curves, identifying how the 
spatial derivatives $u_x$ and $v_x$ may become singular 
and determining the rate of this loss of regularity. 
Finally, by introducing a Finsler norm on tangent vectors 
in the space of solutions, expressed in the 
transformed Bressan-Constantin variables, 
we define a Lipschitz metric that quantifies the 
minimal energy required to transport one 
solution into another along admissible deformation paths. 
This metric provides a geometric framework 
for measuring stability, ensures uniqueness, 
and establishes (Lipschitz) continuous dependence of 
solutions on the initial data within 
the conservative class.

\subsection{Generic regularity}
\begin{figure}
\centering{\includegraphics[scale=0.5]{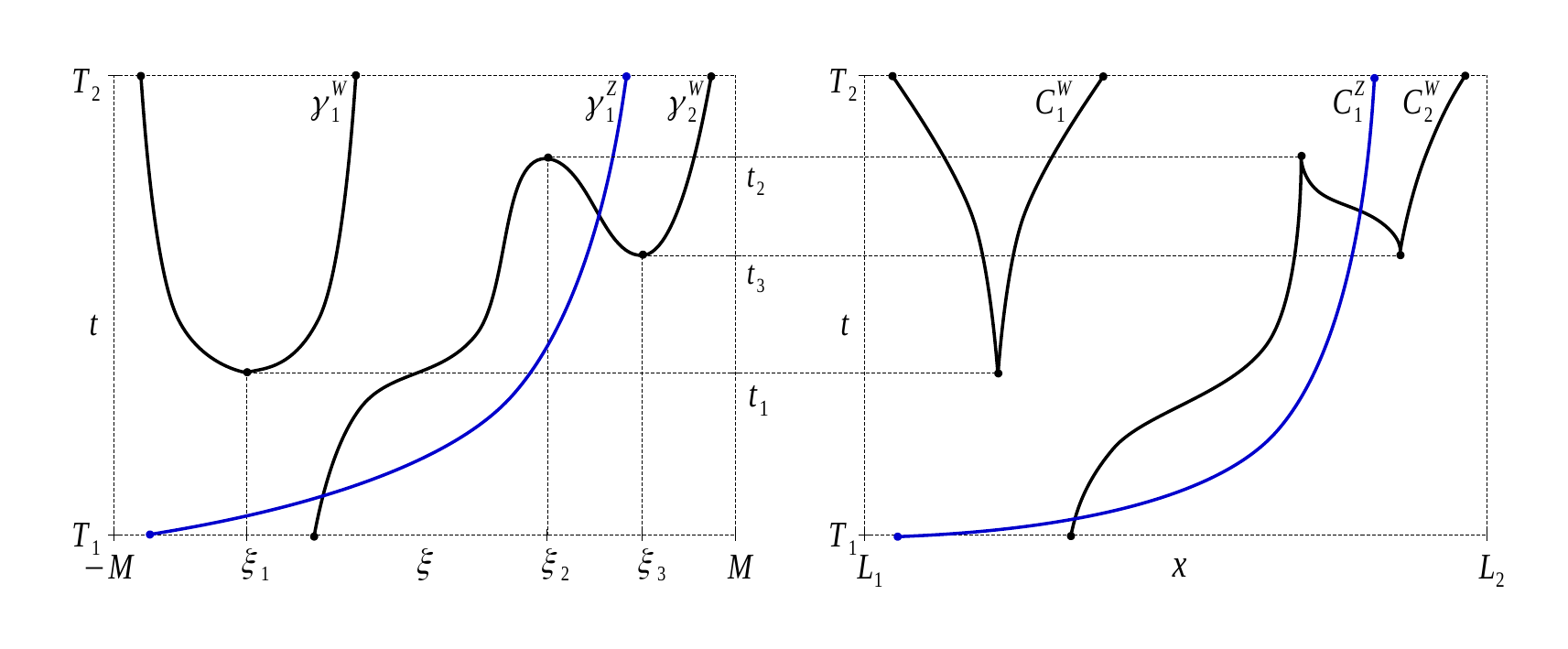}}
\caption{Illustration of 
the image of $\Gamma^W=\gamma_1^W\cup\gamma_2^W$ and 
$\Gamma^Z=\gamma_1^Z$ (see
\eqref{GammaWZ} and \eqref{gammaWZ} with $N_1=2$ and $N_2=1$)
under the characteristics map $(t,\xi)\mapsto(t,y(t,\xi))$.
Here $C_i^W$ is the image of $\gamma_i^W$, $i=1,2$, and
$C_1^Z$ is the image of $\gamma_1^Z$.
Notice that the characteristic curves $C_1^W$ and $C_2^W$ 
have discontinuous derivatives only at the images of the 
points $(t_i,\xi_i)$ where $W_\xi(t_i,\xi_i)=0$, $i=1,2,3$.}
\label{char-curves}
\end{figure}

The first objective of this paper is to establish the generic 
regularity properties of global conservative solutions of 
the Novikov system (see Theorem \ref{Thmgr}).
Our approach builds on the work of Bressan and Chen \cite{BC17},
who studied generic regularity for a 
quasilinear second-order wave equation.
Their method relies on transforming solutions from 
Eulerian variables to a new set of variables 
along characteristics, in which all possible singularities of 
the original solutions are resolved.
This change of variables follows closely 
the procedure developed by Bressan 
and Constantin \cite{BC07} for the Camassa–Holm 
equation, and we therefore refer to 
these as the Bressan-Constantin variables. 
In the transformed setting, the dynamics are governed by 
a semilinear system of ODEs in a Banach space, which 
admits a unique global solution.
Exploiting the regularity of the global solution in the $(t,\xi)$ 
plane---where $x = y(t,\xi)$ and $y(t,\cdot)$ denotes 
a characteristic---one can 
carry out a detailed analysis of the associated 
characteristic curves, which are defined by a 
corresponding system of ODEs 
(and can be viewed as level sets of certain functions). 
When these curves are mapped from the $(t,\xi)$ 
plane back to the Eulerian $(t,x)$ plane 
(see Figure~\ref{char-curves}), one observes that the loss of 
regularity in the solution occurs precisely along these curves, 
where the spatial derivative becomes unbounded 
(see \cite{BC17} and Section \ref{chv} for details).

The Bressan-Chen approach \cite{BC17} has been 
successfully applied to several peakon-type equations 
\cite{CCCS18,CT18,LZ17,TMQ21,ZLDW25,S-Y20,YZ25}. 
In particular, the generic regularity of the 
scalar Novikov equation \eqref{N} 
was analyzed in \cite[Theorem 1.2]{CCCS18}. 
The study \cite[Theorem 1.1]{CT18} further investigated the regularity 
of a modified two-component Camassa-Holm system, 
showing that its two components share the same regularity 
in the $(t,x)$ plane. In contrast, the two-component 
Novikov system exhibits 
a substantially more intricate nonlinear coupling between 
its components, leading to distinct regularity behaviors 
for $u$ and $v$.

More in detail, considering the metric space 
$(\Upsilon^k,d_{\Upsilon^k})$, where
\begin{equation}\label{Ups}
	\Upsilon^k=\left(C^k(\R)\right)^2\cap\Sigma,
	\quad k\in\N,
\end{equation}
and
\begin{equation}\label{d-Ups}
	d_{\Upsilon^k}((f_1,g_1), (f_2,g_2))
	=d_\Sigma((f_1,g_1), (f_2,g_2))
	+\|f_1-f_2\|_{C^k}+\|g_1-g_2\|_{C^k},
\end{equation}
we show that there exists a subset
$\mathcal{M}_T\subset\Upsilon^k$, $k\geq5$,
that is open in $\Upsilon^k$ and dense in
$\Upsilon^{k-1}\cap C^k_{\mathrm{loc}}(\R)$,
and two sets of characteristic curves, say
$\left\{\mathbf{C}_i^W\right\}_{i=1}^{N_1}, 
\left\{\mathbf{C}_j^Z\right\}_{j=1}^{N_2}
\subset[-T,T]\times\R$ 
for some $N_1,N_2\in\N\cup\{0\}$ (see \eqref{empts}),
which satisfy the following properties. Both components 
$u$ and $v$ are $C^k$–regular away from the characteristic sets 
$\mathbf{C}_i^W$ and $\mathbf{C}_j^Z$; 
see item (1) of Theorem \ref{Thmgr}. 
Moreover, $u$ is of class $C^1$ along each $\mathbf{C}_j^Z$ 
and $v$ is of class $C^1$ along each $\mathbf{C}_i^W$, 
as stated in item~(2) of Theorem \ref{Thmgr}.  

In contrast to the scalar Novikov case studied in~\cite{CCCS18}, 
where the solution exhibits a lower regularity than the initial data, 
we show that this loss is not optimal: 
the solution actually preserves its regularity 
away from the characteristic curves. 
This is consistent with known results on generic regularity 
for systems of conservation laws \cite[Theorem~5.1]{DG91} 
and scalar conservation laws \cite{Sch73,D85}.  

Finally, we prove that the solution $(u,v)(t)$ corresponding to 
the initial data $(u_0,v_0)\in\mathcal{M}_T$ 
does not develop concentration of energy 
and remains conservative for all $t\in\R$; 
see item~(3) of Theorem~\ref{Thmgr}.

The behavior of $u$ and $v$ near the characteristic curves 
depends crucially on the location of the point along the curve.  
A detailed analysis, given in Theorem~\ref{Thmchc}, 
reveals eight distinct types of points on the characteristic curves, 
corresponding to qualitatively different behaviors 
of the two-component solution $(u,v)(t,x)$ 
as $(t,x)$ approaches them.  
These cases are further discussed in 
Corollary \ref{Nchars} and Remark \ref{Nchc} 
for the specific setting of the scalar Novikov equation.  

The generic regularity of global conservative weak solutions 
to the two-component Novikov system is summarized 
in the following theorem.

\begin{theorem}
[Generic regularity of conservative solutions]\label{Thmgr}
Consider the Cauchy problem \eqref{t-c-N-n}--\eqref{iid} 
with $(u_0,v_0)\in\Upsilon^k$, $k\geq5$,
see \eqref{Ups}--\eqref{d-Ups}.
Then for any $T>0$ there exists a subset
$\mathcal{M}_T\subset\Upsilon^k$ that is open in $\Upsilon^k$
and dense in $\Upsilon^{k-1}\cap C^k_{\mathrm{loc}}(\R)$
such that 
for any initial data $(u_0,v_0)\in\mathcal{M}_T$,
the initial measure $\mu_0$
having a zero singular part, i.e., 
\begin{equation}\label{mu0i}
d\mu_0=\left((\px u_0)^2+(\px v_0)^2+(\px u_0)^2(\px v_0)^2
\right)dx,
\end{equation}
and $D_{W,0}=D_{Z,0}=\emptyset$ (see Remark \ref{Rmz}),
the corresponding global conservative solution 
$$
\left(u(t),v(t),\mu_{(t)}; D_W(t), D_Z(t)\right),
$$
given in Theorem \ref{Thm}
satisfies the following regularity properties.
We have a finite number of piecewise $C^{k-1}$ 
characteristic curves 
$\left\{\mathbf{C}_i^W\right\}_{i=1}^{N_1}, 
\left\{\mathbf{C}_j^Z\right\}_{j=1}^{N_2}
\subset[-T,T]\times\R$ 
for some $N_1,N_2\in\N\cup\{0\}$ (see \eqref{empts}) such that
\begin{enumerate}
	\item $u$ and $v$ are $k$ times continuously differentiable in
	$\left([-T,T]\times\R\right)\setminus
	\left(\bigcup\limits_{i=1}^{N_1}\mathbf{C}_i^W
	\cup
	\bigcup\limits_{j=1}^{N_2}\mathbf{C}_j^Z
	\right)$;
	
	\item $u$ is continuously differentiable in
	$\left([-T,T]\times\R\right)\setminus
	\left(\bigcup\limits_{i=1}^{N_1}\mathbf{C}_i^W
	\right)$, while $v$ is continuously differentiable in
	$\left([-T,T]\times\R\right)\setminus
	\left(\bigcup\limits_{j=1}^{N_2}\mathbf{C}_j^Z
	\right)$.
	
	\item the singular part of $\mu_{(t)}$
	is zero for all $t\in[-T,T]$, i.e.,
	\begin{equation}\label{mu-t-ac}
		d\mu_{(t)}=\left(u_x^2+v_x^2+u_x^2v_x^2\right)
		(t,x)\,dx,\quad 			
		\text{for all }\,
		t\in[-T,T].
	\end{equation}
	In particular, we have that
	$\mathrm{meas}(D_W(t))=\mathrm{meas}(D_Z(t))=0$,
	see \eqref{DN}, for all $t\in[-T,T]$.
	\end{enumerate}
Here the curves $\mathbf{C}_i^W$, $i=1,\dots, N_1$ and
$\mathbf{C}_j^Z$, $j=1,\dots,N_2$
have the following properties (see Figure \ref{char-curves}):
\begin{enumerate}[i)]
	\item $\bigcup\limits_{i=1}^{N_1}\mathbf{C}_i^W$ and
	$\bigcup\limits_{j=1}^{N_2}\mathbf{C}_j^Z$
	are the images of the map 
	$(t,\xi)\mapsto (t,y(t,\xi))$ of the level sets
	$\Gamma^W$ and
	$\Gamma^Z$ given by 
	(see also Figure \ref{char-curves})
	\begin{equation}\label{GammaWZ}
		\Gamma^W=\left\{
		(t,\xi)\in[-T,T]\times\R:
		W(t,\xi)=\pi\right\},\quad
		\Gamma^Z=\left\{
		(t,\xi)\in[-T,T]\times\R:
		Z(t,\xi)=\pi\right\},
	\end{equation}
	where $(U,V,W,Z,q)(t,\xi)$ is a solution of the associated ODE 
	system \eqref{ODE}--\eqref{id} given in Theorem \ref{Ckgl} 
	and the characteristic $y(t,\xi)$ is defined in \eqref{char1}.
	Moreover, we have
	\begin{equation}\label{sder}
	\begin{split}
		&(W_t,W_\xi)(t,\xi)\neq(0,0),\quad
		(W_\xi,W_{\xi\xi})(t,\xi)\neq(0,0),\quad
		\text{for all }(t,\xi)\in\Gamma^W,\\
		&(Z_t,Z_\xi)(t,\xi)\neq(0,0),\quad
		(Z_\xi,Z_{\xi\xi})(t,\xi)\neq(0,0),\quad
		\text{for all }(t,\xi)\in\Gamma^Z.
	\end{split}
	\end{equation}
	
	\item $\mathbf{C}_i^W\cap\mathbf{C}_j^W=\emptyset$ for all
	$i,j=1,\dots,N_1$, $i\neq j$ and
	$\mathbf{C}_i^Z\cap\mathbf{C}_j^Z=\emptyset$ for all
	$i,j=1,\dots,N_2$, $i\neq j$;
	
	\item $\mathbf{C}_i^W,\mathbf{C}_j^Z\subset[-T,T]\times\R$ 
	are bounded and have no self 
	intersections for all $i=1,\dots,N_1$ 
	and $j=1,\dots,N_2$;
	
	\item $\mathbf{C}_i^W$ and $\mathbf{C}_j^Z$ are either 
	closed curves or have their endpoints at 
	$\{\pm T\}\times\R$ for all $i=1,\dots,N_1$ and
	$j=1,\dots,N_2$;
\end{enumerate}
\end{theorem}

This theorem reveals a delicate balance of 
regularity in the two-component solution $(u,v)$.
Away from both families of characteristic curves, the solution 
is as smooth as one could expect.
Each component, however, loses smoothness along 
its own family of characteristics---$u$ across the $W$-curves 
and $v$ across the $Z$-curves. Across the other family, the 
situation is considerably better: $u$ remains differentiable 
across the $Z$-curves, and $v$ across the $W$-curves.
Thus, while each field exhibits singular behavior along its 
own characteristics, the cross-interaction preserves a 
notable degree of regularity.

Theorem \ref{Thmgr} contains the scalar 
Novikov equation as a special case.
\begin{corollary}
	Setting $u_0=v_0$ in Theorem \ref{Thmgr}, we obtain 
	the generic regularity of the global conservative solution 
	of the scalar Novikov equation with 
	$u_0\in C^k\cap H^1\cap W^{1,4}$, $k\geq 5$.
	In this case we have $u=v$, $W=Z$, $N_1=N_2$,
	$C_i^W=C_i^Z$, $i=1,\dots,N_1$, and thus for generic 
	initial data $u_0\in C^k\cap H^1\cap W^{1,4}$ the solution $u$ is $k$
	times continuously differentiable in
	$\left([-T,T]\times\R\right)\setminus
	\left(\bigcup\limits_{i=1}^{N_1}\mathbf{C}_i^W\right)$.
\end{corollary}

\begin{remark}[Concentration of energy]
The global conservative solutions established in 
Theorem~\ref{Thm} may, in general, exhibit finite-time 
energy concentration, a phenomenon also known from 
other peakon equations such as the Camassa-Holm 
equation \cite{BC07} and the scalar Novikov equation \cite{CCL18}. 
Nevertheless, item (3) of Theorem \ref{Thmgr} ensures 
that for generic initial data $(u_0,v_0)\in\mathcal{M}_T$, 
the corresponding global solution does not experience 
any concentration of energy at any time $t\in\R$. 
In particular, for generic initial data one has 
$N_W=N_Z=\emptyset$ (see \eqref{DN}), and the energies 
$E_{u_0}$, $E_{v_0}$, and $H_0$ 
remain conserved for all $t$, cf.~item~(4) in 
Definition~\ref{defsc}.
\end{remark}

\begin{remark}[Regularity rate]
	For the scalar Novikov equation~\eqref{N}, 
	it was shown in~\cite[Theorem~1.2]{CCCS18} that, for initial data in $C^3$, 
	the corresponding global solution is twice continuously differentiable
	outside a finite number of characteristic curves 
	(cf.~also~\cite[Theorem~3.6]{LZ17}, 
	\cite[Theorem~1.1]{TMQ21}, \cite[Theorem~1.1]{CT18}, 
	\cite[Theorem~1.2]{ZLDW25}, and \cite[Theorem~1.1]{S-Y20}).
		Here we impose stronger regularity assumptions on
		the initial data, which allow us to apply the parametric
		transversality theorem
	in Lemma \ref{L3},
	see Remark \ref{k-reg} for a detailed explanation.
		Moreover, we show that the two-component solution 
	retains its full regularity for all $(t,x)$ away 
	from the characteristic curves.
\end{remark}

\begin{remark}[Density]
		Notice that the open set $\mathcal{M}_T$ of initial data
		is dense in the topology
		$\Upsilon^{k-1}\cap C^k_{\mathrm{loc}}(\R)$,
		which is slightly weaker than the topology of $\Upsilon^k$.
	Because the two-component Novikov system considered here is nonlocal
		(see \eqref{t-c-N-n} below), even a compactly supported
		perturbation of the initial data immediately affects the solution
		for all $x\in\R$ and $t\neq0$.
		Therefore, in Step 5 of the proof of Theorem \ref{Thmgr},
	see Section \ref{prTh1}, we cannot use cutoff functions in the same way
	as for the quasilinear second order wave equation
	\cite[p.~350, step 4]{BC17}.
	Due to this technical limitation,
		our density result is weaker than the one in \cite{BC17};
	see Remark \ref{R-dens} for more details.
\end{remark}

\begin{remark}
Notice that the level sets $\Gamma^W$ and $\Gamma^Z$ 
(see~\eqref{GammaWZ}) are determined by the initial 
data $(u_0,v_0)$ through the unique solution of 
the associated Cauchy problem~\eqref{ODE}--\eqref{id} 
for the corresponding ODE system.
\end{remark}

\subsection{Behavior near characteristics}
Theorem~\ref{Thmgr} does not, however, describe the 
detailed behavior of the solution $(u,v)(t,x)$ as $(t,x)$ 
approaches the characteristic curves. 
In the next theorem,
using a characteristic-coordinate and transversality
approach related to that of \cite{BHY15},
we identify how the singularity rate 
of $u_x$ and $v_x$ depends on the local properties of 
points along these curves. Bressan, Huang and Yu
\cite{BHY15} study structurally stable singularities for a scalar nonlinear wave
equation. Here this framework is adapted to a nonlocal two-component
Novikov system, and we additionally construct a Lipschitz metric for
its global conservative flow. Recall that, for the class of 
initial data considered in Theorem \ref{Thmgr}, the Radon 
measure $\mu_{(t)}$ has no singular part and is 
therefore completely determined by $(u,v)(t,\cdot)$.

\begin{figure}
\centering{\includegraphics[scale=0.5]{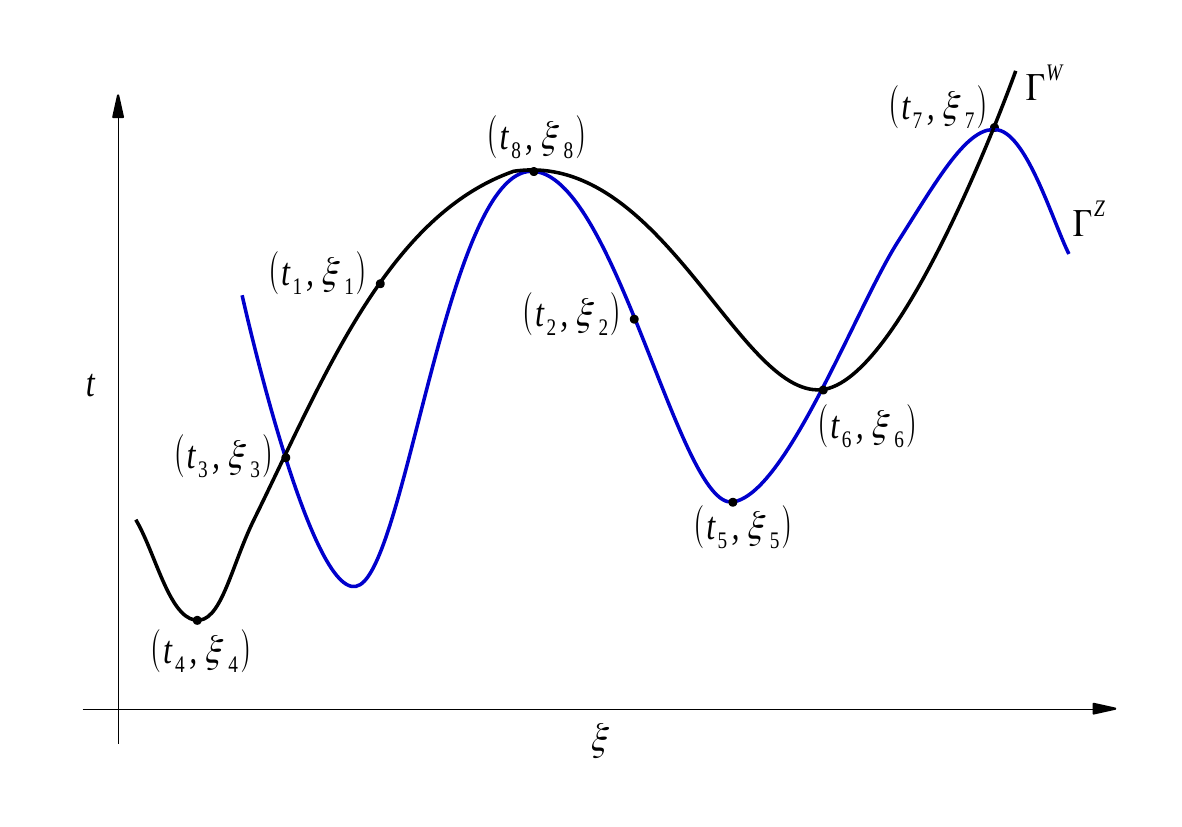}}
\caption{Eight types of points $(t_i,\xi_i)$, $i=1,\dots,8$ 
on the characteristic curves $\Gamma^W\cup\Gamma^Z$, 
see \eqref{GammaWZ}, discussed in Theorem \ref{Thmchc}.}
\label{char-WZ}
\end{figure}

\begin{theorem}[Solution behavior near the characteristic curves]
\label{Thmchc}
Consider the generic global conservative solution 
${(u,v)(t,x)}$ as given in Theorem \ref{Thmgr},
which corresponds to initial data
$(u_0,v_0)\in\mathcal{M}_T\subset\Upsilon^k$, $k\geq8$,
and a Radon measure $\mu_0$ having a zero singular part, see \eqref{mu0i}.
Consider also $(U,V,W,Z,q)(t,\xi)$ as given in Theorem \ref{Ckgl}
and $y(t,\xi)$ defined by \eqref{char1}.
Introduce the function (here $\so(\cdot)$ 
is the usual little-o notation)
\begin{equation}\label{elld}
	\ell_i(t,x)=x-x_i-(uv)(t_i,x_i)(t-t_i)+\so(t-t_i),
\end{equation}
which serves as a first order approximation of the monomial
$(\xi-\xi_i)^n$ in the Taylor expansion 
of the characteristic curve $x=y(t,\xi)$ at the points 
$t_i,x_i\in\R$, $i\in\{1,\dots,8\}$, 
that are specified below.

Then $u$ and $v$ have the following behavior 
as $(t,x)$ approaches eight different types 
of points on the characteristic curves 
(see Figure \ref{char-WZ}):
\begin{enumerate}
	\item if $W(t_1,\xi_1)=\pi$,
	$W_\xi(t_1,\xi_1)\neq0$, and 
	$Z(t_1,\xi_1)\neq\pi$, then
	\begin{equation}\label{uvch1}
	\begin{split}
		&u(t,x)=u(t_1,x_1)+a_{1,1}\ell_1^{\,2/3}(t,x)
		+\mathcal{O}\left(\ell_1(t,x)\right),\\
		&v(t,x)=v(t_1,x_1)+b_{1,1}\ell_1(t,x)
		+b_{1,2}(t-t_1)
		+\mathcal{O}\left(\ell_1^{\,4/3}(t,x)\right),
	\end{split}
	\end{equation}
	where $x_1=y(t_1,\xi_1)$, and $a_{1,1}\in\R\setminus\{0\}$,
	$b_{1,1},b_{1,2}\in\R$ are some constants;
	
	\item if $W(t_2,\xi_2)\neq\pi$, and
	$Z(t_2,\xi_2)=\pi$,
	$Z_\xi(t_2,\xi_2)\neq0$, then
	\begin{equation}
		\label{uvch2}
		\begin{split}
			&u(t,x)=u(t_2,x_2)+a_{2,1}\ell_2(t,x)
			+a_{2,2}(t-t_2)
			+\mathcal{O}\left(\ell_2^{\,4/3}(t,x)\right),\\
			&v(t,x)=v(t_2,x_2)+b_{2,1}\ell_2^{\,2/3}(t,x)
			+\mathcal{O}\left(\ell_2(t,x)\right),
		\end{split}
	\end{equation}
	where $x_2=y(t_2,\xi_2)$, and $b_{2,1}\in\R\setminus\{0\}$,
	$a_{2,1}, a_{2,2}\in\R$ are some constants;
	
	\item if $W(t_3,\xi_3)=\pi$, $W_\xi(t_3,\xi_3)\neq0$, and
	$Z(t_3,\xi_3)=\pi$, $Z_\xi(t_3,\xi_3)\neq0$, then
	\begin{equation}
		\label{uvch3}
		\begin{split}
			&u(t,x)=u(t_3,x_3)+a_{3,1}\ell_3^{\,4/5}(t,x)
			+\mathcal{O}\left(\ell_3(t,x)\right),\\
			&v(t,x)=v(t_3,x_3)+b_{3,1}\ell_3^{\,4/5}(t,x)
			+\mathcal{O}\left(\ell_3(t,x)\right),
		\end{split}
	\end{equation}
	where $x_3=y(t_3,\xi_3)$, and 
	$a_{3,1},b_{3,1}\in\R\setminus\{0\}$
	are some constants;
	
	\item if $W(t_4,\xi_4)=\pi$, $W_\xi(t_4,\xi_4)=0$, and
	$Z(t_4,\xi_4)\neq\pi$, then
	\begin{equation}
		\label{uvch4}
		\begin{split}
			&u(t,x)=u(t_4,x_4)+a_{4,1}\ell_4^{\,3/5}(t,x)
			+\mathcal{O}
			\left(\ell_4^{\,4/5}(t,x)\right),\\
			&v(t,x)=v(t_4,x_4)+b_{4,1}\ell_4(t,x)
			+b_{4,2}(t-t_4)
			+\mathcal{O}
			\left(\ell_4^{\,6/5}(t,x)\right),
		\end{split}
	\end{equation}
	where $x_4=y(t_4,\xi_4)$, and $a_{4,1}\in\R\setminus\{0\}$,
	$b_{4,1},b_{4,2}\in\R$
	are some constants;
	
	\item if $W(t_5,\xi_5)\neq\pi$, and 
	$Z(t_5,\xi_5)=\pi$, $Z_\xi(t_5,\xi_5)=0$, then
	\begin{equation}
		\label{uvch5}
		\begin{split}
			&u(t,x)=u(t_5,x_5)+a_{5,1}\ell_5(t,x)
			+a_{5,2}(t-t_5)
			+\mathcal{O}\left(
			\ell_5^{\,6/5}(t,x)\right),\\
			&v(t,x)=v(t_5,x_5)+b_{5,1}\ell_5^{\,3/5}(t,x)
			+\mathcal{O}\left(\ell_5^{\,4/5}(t,x)\right),
		\end{split}
	\end{equation}
	where $x_5=y(t_5,\xi_5)$, and $b_{5,1}\in\R\setminus\{0\}$,
	$a_{5,1},a_{5,2}\in\R$
	are some constants;
	
	\item if $W(t_6,\xi_6)=\pi$, $W_\xi(t_6,\xi_6)=0$, and
	$Z(t_6,\xi_6)=\pi$, $Z_\xi(t_6,\xi_6)\neq0$, then
	\begin{equation}
		\label{uvch6}
		\begin{split}
			&u(t,x)=u(t_6,x_6)+a_{6,1}\ell_6^{\,5/7}(t,x)
			+\mathcal{O}\left(
			\ell_6^{\,6/7}(t,x)\right),\\
			&v(t,x)=v(t_6,x_6)+b_{6,1}\ell_6^{\,6/7}(t,x)
			+\mathcal{O}\left(\ell_6(t,x)\right),
		\end{split}
	\end{equation}
	where $x_6=y(t_6,\xi_6)$, and 
	$a_{6,1}, b_{6,1}\in\R\setminus\{0\}$
	are some constants.
	
	\item if $W(t_7,\xi_7)=\pi$, $W_\xi(t_7,\xi_7)\neq0$, and
	$Z(t_7,\xi_7)=\pi$, $Z_\xi(t_7,\xi_7)=0$, then
	\begin{equation}
		\label{uvch7}
		\begin{split}
			&u(t,x)=u(t_7,x_7)+a_{7,1}\ell_7^{\,6/7}(t,x)
			+\mathcal{O}\left(\ell_7(t,x)\right),\\
			&v(t,x)=v(t_7,x_7)+b_{7,1}\ell_7^{\,5/7}(t,x)
			+\mathcal{O}\left(
			\ell_7^{\,6/7}(t,x)\right),
		\end{split}
	\end{equation}
	where $x_7=y(t_7,\xi_7)$, and 
	$a_{7,1}, b_{7,1}\in\R\setminus\{0\}$
	are some constants.
	
	\item if $W(t_8,\xi_8)=\pi$, $W_\xi(t_8,\xi_8)=0$, and
	$Z(t_8,\xi_8)=\pi$, $Z_\xi(t_8,\xi_8)=0$, then
	\begin{equation}
		\label{uvch8}
		\begin{split}
			&u(t,x)=u(t_8,x_8)+a_{8,1}\ell_8^{\,7/9}(t,x)
			+\mathcal{O}\left(
			\ell_8^{\,8/9}(t,x)\right),\\
			&v(t,x)=v(t_8,x_8)+b_{8,1}\ell_8^{\,7/9}(t,x)
			+\mathcal{O}\left(
			\ell_8^{\,8/9}(t,x)\right),
		\end{split}
	\end{equation}
	where $x_8=y(t_8,\xi_8)$, and 
	$a_{8,1}, b_{8,1}\in\R\setminus\{0\}$
	are some constants.
\end{enumerate}
\end{theorem}

The behavior of the solution $u$ to 
the scalar Novikov equation, which represents a 
special case of the two-component system \eqref{t-c-N}, 
is described in the following corollary.

\begin{corollary}\label{Nchars}
Taking $u=v$ in Theorem \ref{Thmchc}, we obtain 
the behavior of $u$ as $(t,x)$ approaches the characteristic 
curves in the case of the scalar Novikov equation \eqref{N}
(recall \eqref{elld} and that here $W=Z$):
\begin{enumerate}
	\item if $W(t_3,\xi_3)=\pi$, $W_\xi(t_3,\xi_3)\neq0$, then
	\begin{equation*}
		u(t,x)=u(t_3,x_3)+a_{3,1}\ell_3^{\,4/5}(t,x)
		+\mathcal{O}\left(\ell_3(t,x)\right),
	\end{equation*}
	where $x_3=y(t_3,\xi_3)$, and $a_{3,1}\in\R\setminus\{0\}$
	is some constant;
	
	\item if $W(t_8,\xi_8)=\pi$, $W_\xi(t_8,\xi_8)=0$, then
	\begin{equation*}
		u(t,x)=u(t_8,x_8)+a_{8,1}\ell_8^{\,7/9}(t,x)
		+\mathcal{O}\left(
		\ell_8^{\,8/9}(t,x)\right),\\
	\end{equation*}
	where $x_8=y(t_8,\xi_8)$, and 
	$a_{8,1}\in\R\setminus\{0\}$
	is some constant.
\end{enumerate}
\end{corollary}

\begin{remark}\label{Nchc}
The behavior of global solutions to the scalar 
Novikov equation \eqref{N} near characteristic 
curves was previously analyzed in~\cite{HLY24} 
(see also \cite[Theorem 1.3]{ZLDW25}). 
That work considers the same type of points 
as $(t_8,\xi_8)$ in Corollary~\ref{Nchars} 
(see \cite[Theorem 1.1, item 2]{HLY24}), 
and also treats points analogous to $(t_3,\xi_3)$, 
though under the additional assumption 
$W_{\xi\xi}(t_3,\xi_3)=0$ (cf.~\cite[Theorem 1.1, item 1]{HLY24}). 
We note that the term $(x-x_8)^{3/4}$ appears in 
the expansion \cite[Equation (1.2)]{HLY24}, 
even though $(x-x_8)$ may take negative values. 
Moreover, from \cite[Equations~(1.2)--(1.3)]{HLY24} it follows that 
$u(t_3,x_3)=u(t_8,x_8)=0$, which does not hold in general. 
Finally, in the Taylor expansion of the characteristic 
given in \cite[Equation~(2.17)]{HLY24}, 
one would naturally expect the presence of a term of the form 
$u^2(t_8,x_8)(t-t_8)$; similarly, in \cite[Equation~(2.25)]{HLY24}, 
a term $u^2(t_3,x_3)(t-t_3)$ would appear more consistent. 
These remarks suggest that certain expressions in \cite{HLY24} 
might merit further verification.
\end{remark}

Notice that for general initial data 
$(u_0,v_0,\mu_0;D_{W,0}, D_{Z,0})\in\mathcal{D}$ 
for the two-component Novikov equation 
(see \eqref{D-set} below), the solution can be 
guaranteed to be only H\"older continuous with exponent $1/2$; 
see item (1) in Theorem \ref{Thm} and item (3) in Definition \ref{defs}. 
For comparison, the global conservative 
solution of the scalar Novikov equation 
constructed in \cite{CCL18} is H\"older continuous with exponent $3/4$. 
In contrast, Theorem \ref{Thmchc} shows that the solution $(u,v)$ 
corresponding to more regular, generic initial data 
enjoys improved H\"older continuity properties.

\begin{corollary}[H\"older continuity]
Theorems \ref{Thmgr} and \ref{Thmchc} 
imply that for generic initial data 
$(u_0,v_0)\in\mathcal{M}_T$ the corresponding global solution
$(u,v)(t,x)$ is H\"older continuous 
on $[-T,T]\times\R$ with exponent $3/5$ for any 
fixed $T>0$, that is, for all $t_1,t_2,\in[-T,T]$ 
and $x_1,x_2\in\R$ we have
\begin{equation*}
	|u(t_1,x_1)-u(t_2,x_2)|,
	|v(t_1,x_1)-v(t_2,x_2)|\leq
	C\left(|t_1-t_2|^{3/5}
	+|x_1-x_2|^{3/5}\right),
\end{equation*}
for some $C>0$.
Moreover, in the case of the scalar Novikov equation \eqref{N}, 
the generic global conservative solution is H\"older continuous 
with exponent $7/9$, see Corollary \ref{Nchars}:
\begin{equation*}
	|u(t_1,x_1)-u(t_2,x_2)|\leq
	C\left(|t_1-t_2|^{7/9}
	+|x_1-x_2|^{7/9}\right),
\end{equation*}
for all $t_1,t_2,\in[-T,T]$, $x_1,x_2\in\R$, for
some $C>0$.
\end{corollary}

\begin{remark}
For comparison, we recall that the 
generic global conservative solution 
of the Camassa-Holm equation is H\"older continuous 
with exponent $3/5$, as shown in \cite[Theorem~3.7]{LZ17}, 
which is smaller than the exponent obtained for 
the scalar Novikov equation.
\end{remark}

\begin{corollary}
Since either $u_x$ or $v_x$ becomes singular 
along the characteristic curves, 
while the initial data are $C^k$-regular with $k\ge 3$, 
it follows that the curves $\mathbf{C}_i^W$ and $\mathbf{C}_j^Z$ 
in Theorem \ref{Thmgr} are separated from the initial line $\{t=0\}$ 
for all $i=1,\dots,N_1$ and $j=1,\dots,N_2$.
\end{corollary}

\subsection{Lipschitz metric}
The final main result of this paper addresses the 
Lipschitz metric for the two-component Novikov system.
To construct such a metric for the global semigroup of 
conservative weak solutions of \eqref{t-c-N}, we follow the 
approach of Bressan and Chen \cite{BC17A}, who introduced 
a Lipschitz metric for the quasilinear second-order wave equation
(see also \cite{BF05, CGH20, GHR11, GHR13} for related 
constructions of Lipschitz metrics for Camassa-Holm type 
equations using different methods).
The distance between two solutions is defined as the 
minimal cost required to transport the associated energy 
measure from one solution to the other.
Concretely, this is realized by taking the infimum over 
all admissible paths connecting the two (sufficiently regular) 
solutions of the length of these paths.

To define the path length in a way that ensures the 
Lipschitz property, one introduces a Finsler norm 
for the tangent vector $(r^\theta,s^\theta)(t)$ 
of the path $(u^\theta,v^\theta)(t)$, $\theta\in[0,1]$, 
which (formally) satisfies
\begin{equation}\label{Lippre}
	\|(r^\theta,s^\theta)(t)\|\leq C\|(r^\theta,s^\theta)(0)\|,
	\quad t\in[-T,T],\quad C=C(T)>0,
\end{equation}
see Theorem \ref{Lipn}. 
The task is to verify that \eqref{Lippre} holds 
in the transformed Bressan-Constantin variables 
for a specific class of paths (see Theorem \ref{egammat}).
The key advantage of computing path lengths in 
these transformed variables is that it resolves the singularities 
in the $x$-derivatives of the original solution $(u,v)$, which 
may develop in finite time even when starting from 
smooth initial data (see Theorem \ref{Thmgr}).

The main challenges in implementing 
the Bressan-Chen approach for our problem are:
(a) providing an appropriate definition of the tangent 
vector that satisfies \eqref{Lippre}, and  
(b) proving that the class of paths for which \eqref{Lippre} 
holds in the transformed variables is 
dense in a suitable space of paths.  

Point (a) was addressed in \cite{CCCS18} 
for the scalar Novikov equation \eqref{N}.
In the general case $u\neq v$, however, the analysis
becomes significantly more intricate.
Moreover, we slightly refine the definition of this norm by 
employing the weight function $e^{-\alpha|x|}$ with 
$\alpha\in(0,1)$, instead of $e^{-|x|}$ as in 
\cite[Equation (3.9)]{CCCS18}, which allows us to 
establish the necessary estimates in 
Theorem \ref{Lipn} and Proposition \ref{COm} (see Lemma \ref{LID} 
and Remark \ref{expLip} for details).

Concerning (b), we show that the appropriate class of paths, 
referred to as ``piecewise regular paths
under the ODE system''
(see Definition \ref{PrpT} below), is dense in a certain 
functional space of piecewise regular paths
(see Theorem \ref{dpath}).
The proof of Theorem \ref{dpath} revisits the 
methodology of Bressan and Chen \cite[Section 6]{BC17} 
and incorporates ideas from the analysis of generic regularity 
in Sections \ref{gODE} and \ref{prTh1}, in particular the 
use of the Thom's transversality theorem.

In \cite{BC17A, CCCS18}, the authors 
introduced ``piecewise regular paths’’ $\mathbf{u}^\theta(t)$ 
connecting solutions in the Eulerian variables.  
In contrast, our analysis is formulated entirely in 
the transformed variables, where the Bressan-Constantin 
coordinates eliminate the possible 
singularities of $u_x$ and $v_x$.  
Working with regular paths $\mathbf{U}^\theta(t)$ in this 
setting provides a more transparent framework for 
tracking the evolution between solutions, and 
builds naturally on the developments 
in \cite{BC17A,CCCS18,CCS17,CT18,CCS24, TMQ21, YZ25}.

A key ingredient is the density result for regular 
paths established in Theorem \ref{dpath}, which provides 
a natural class of well-behaved paths for 
the associated ODE system (see Definition \ref{PrpT}).
Notice that we impose on these paths
non-degeneracy conditions that differ from those in \cite{BC17, BC17A};
see \eqref{p-ins} and Section \ref{Sp-ins} for details.
By expressing the Finsler norm of tangent vectors 
in the transformed variables, we obtain effective 
control of path-length growth (Theorem \ref{egammat}), 
leading to the construction of a Lipschitz metric 
in this coordinate system (Section \ref{GdO}).  
Since this metric is Lipschitz continuous with 
respect to the ODE flow
(Theorem \ref{LRO}, Proposition \ref{COm}),
it overcomes the lack of a global bijection between the Eulerian 
and transformed coordinates (cf.~Figure \ref{path-m}).  
Transporting the metric back to the Eulerian variables then 
yields a Lipschitz metric for global conservative 
solutions of the two-component Novikov system 
(Definition \ref{Ldist}), culminating 
in the theorem below.

\begin{theorem}[Lipschitz metric in $\mathcal{D}$]
\label{ThLip}
Consider initial data 
$\mathbf{u}_0,\hat{\mathbf{u}}_0\in\mathcal{D}$, where 
\begin{equation*}
	\mathbf{u}_0=\left(u_0,v_0,\mu_0;
	D_{W,0},D_{Z,0}\right),\quad
	\hat{\mathbf{u}}_0=\left(\hat{u}_0,\hat{v}_0,\hat{\mu}_0;
	\hat{D}_{W,0},\hat{D}_{Z,0}\right),
\end{equation*}
and assume that
$u_0,v_0,\hat{u}_0,\hat{v}_0\in W^{1,1}(\R)$.
Then the corresponding global conservative solutions
$\mathbf{u}(t)$ and $\hat{\mathbf{u}}(t)$ given by Theorem
\ref{Thm} satisfy
$
u(t),v(t),\hat{u}(t),\hat{v}(t)\in W^{1,1}(\R),
$
$t\in[-T,T]$, for any $T>0$
and there exists a metric $d_{\mathcal{D},K}(\cdot,\cdot)$, 
defined in Definition \ref{Ldist}, that ensures the 
solutions satisfy the following Lipschitz property:
\begin{equation}\label{Lipgs}
	d_{\mathcal{D},K}(\mathbf{u}(t),\hat{\mathbf{u}}(t))
	\leq C
	d_{\mathcal{D},K}(\mathbf{u}_0,\hat{\mathbf{u}}_0),\quad
	C=C(T,K,\alpha)>0,\quad
	t\in[-T,T].
\end{equation}
Here $\alpha\in(0,1)$ is the parameter in the weight
$e^{-\alpha|\cdot|}$, see \eqref{Rtr}, while
the energy constant $K>0$ is taken as in \eqref{K-unb}
where the endpoints $\mathbf{U}^0$ and $\mathbf{U}^1$ correspond to
$\mathbf{u}_0$ and $\hat{\mathbf{u}}_0$, respectively,
in the transformed variables.
\end{theorem}
\begin{remark}[$L^1(\R)$ decay]
	\label{L1dec}
		The additional $W^{1,1}(\R)$ assumption on
		$(u(t),v(t))$ is needed to justify rigorously the
	completion of the space of regular functions equipped with
	the geodesic Lipschitz metric
	$d_{\Omega_{\mathrm{tcN}}^1,K}(\cdot,\cdot)$ defined in the transformed 
	variables, see Section \ref{GdO} and
	Remark \ref{compL1} below for details.
\end{remark}

The article is organized as follows.
In Section \ref{mr} we recall the results 
on the global semigroup of conservative 
weak solutions of \eqref{t-c-N}, and we 
introduce the notations and mathematical 
facts used throughout the paper.
Section \ref{cvss} provides a brief overview of 
the Bressan-Constantin approach to constructing global 
solutions of the two-component Novikov 
system (see \cite{HQ21, KR25} for details).
Sections \ref{gODE}–\ref{prTh2} contain the analysis 
of the generic regularity of \eqref{t-c-N} and 
the proofs of Theorems \ref{Thmgr} and \ref{Thmchc}.
Finally, Section \ref{LM} is devoted to the 
construction of a Lipschitz metric for the 
global conservative solutions of \eqref{t-c-N}, 
culminating in the proof of Theorem \ref{ThLip}.

\section{Preliminaries}\label{mr}
In this section, we recall the existence 
of a global semigroup of conservative weak solutions 
to the two-component Novikov system, as obtained 
in \cite{KR25}. We also provide some useful 
estimates and notations that will be 
used throughout the rest of the work.

Applying the operator $(1-\px^2)^{-1}$ 
to both sides of \eqref{t-c-N}, 
we obtain the following nonlocal system:
\begin{align}\label{t-c-N-n}
	\begin{split}
		&u_t+uvu_x+\px P_1+P_2=0,\quad 
		u=u(t,x),\,\,v=v(t,x),\,\, P_j=P_j(t,x),\,\,j=1,2,\\
		&v_t+uvv_x+\px S_1+S_2=0,\quad
		S_j=S_j(t,x),\,\,j=1,2,
	\end{split}
	\\
	\label{iid}
	&u(0,x)=u_0(x),\quad v(0,x)=v_0(x),
\end{align}
where 
\begin{equation*}
	\begin{split}
		&P_1(t,x)=(1-\px^2)^{-1}\left(
		u^2v+uu_xv_x+\frac{1}{2}vu_x^2
		\right)(t,x),\\
		&P_2(t,x)=\frac{1}{2}(1-\px^2)^{-1}
		\left(u_x^2v_x\right)(t,x),
	\end{split}
\end{equation*}
and
\begin{equation*}
	\begin{split}
		&S_1(t,x)=(1-\px^2)^{-1}\left(
		uv^2+vu_xv_x+\frac{1}{2}uv_x^2
		\right)(t,x),\\
		&S_2(t,x)=\frac{1}{2}(1-\px^2)^{-1}
		\left(u_xv_x^2\right)(t,x).
	\end{split}
\end{equation*}
Differentiating \eqref{t-c-N-n} with respect to $x$, 
we obtain
\begin{equation}\label{tcNnd}
	\begin{split}
		&u_{xt}+uv u_{xx}-u^2v+\frac{1}{2}vu_x^2
		+P_1+\px P_2=0,\\
		&v_{xt}+uvv_{xx}-uv^2+\frac{1}{2}uv_x^2
		+S_1+\px S_2=0.
	\end{split}
\end{equation}
This system will be employed to define 
weak solutions in the 
Definition \ref{defs} that follows.

The two-component Novikov system possess the 
following conservation laws \cite{HQ20, HLQ23}:
\begin{equation}\label{consq}
	\begin{split}
		&E_u(t)=\int_{-\infty}^{\infty}
		\left(u^2+u_x^2\right)(t,x)\,dx
		=:E_{u_0},\quad
		E_v(t)=\int_{-\infty}^{\infty}
		\left(v^2+v_x^2\right)(t,x)\,dx
		=:E_{v_0},\\
		&G(t)=\int_{-\infty}^{\infty}
		\left(uv+u_xv_x\right)(t,x)\,dx
		=:G_{0},\\
		&H(t)=\int_{-\infty}^{\infty}
		\left(3u^2v^2+u^2v_x^2
		+u_x^2v^2
		+4uu_xvv_x
		-u_x^2v_x^2\right)
		(t,x)\,dx=:H_0,
	\end{split}
\end{equation}
for all $t\in\R$. Notice that the
conservation laws \eqref{consq} imply the following bound
for $\|u_xv_x\|_{L^2}$ (notice that 
$7E_{u_0}E_{v_0}-H_0\geq 0$ \cite{KR25}):
\begin{equation*}
	\|u_xv_x\|_{L^2}^2\leq7E_{u_0}E_{v_0}-H_0.
\end{equation*}

We are now in a position to define what we 
mean by weak and conservative weak solutions 
of \eqref{t-c-N-n}, beginning with the former. 
For a more detailed discussion, we refer 
to \cite{KR25} and the references therein.

\begin{definition}[\cite{KR25}, Global weak solution 
of \eqref{t-c-N-n}--\eqref{iid}]\label{defs}
Suppose that $(u_0,v_0)\in\Sigma$, where 
$\Sigma$ is defined in \eqref{Sig}.
We say that a vector function $(u,v)(t,x)$ is a 
global weak solution of the Cauchy 
problem \eqref{t-c-N-n}-\eqref{iid} 
on $\mathbb{R}$, if
$(u,v)(0,x)=(u_0,v_0)(x)$ for all 
$x\in\mathbb{R}$ and $(u,v)$ satisfies 
the equations
\begin{equation*}
	\begin{split}
		&\int_{-T}^T\int_{-\infty}^{\infty}
		\biggl((\px u)(\pt\phi_u+uv\px\phi_u)\\
		&\left.\qquad\qquad\,\,\,\,
		+\left(u^2v+u(\px u)\px v+\frac{1}{2}v(\px u)^2
		-P_1-\px P_2\right)\phi_u\right)\,dx\,dt=0,\\
		&\int_{-T}^T\int_{-\infty}^{\infty}
		\biggl((\px v)(\pt\phi_v+uv\px\phi_v)\\
		&\left.\qquad\qquad\,\,\,\,+\left(uv^2+v(\px u)\px v
		+\frac{1}{2}u(\px v)^2-S_1-\px S_2\right)\phi_v\right)\,dx\,dt=0,
		\end{split}
	\end{equation*}
for all test functions 
$\phi_u,\phi_v\in C^\infty((-T,T)\times\mathbb{R})$ 
with compact support and arbitrary $T>0$.
Moreover, $u$ and $v$ have the following properties:
\begin{enumerate}
	\item $(u,v)(t,\cdot)$ belongs to 
	$\Sigma$
	for any fixed $t\in\mathbb{R}$;
	
	\item $u(t,\cdot)$, $v(t,\cdot)$ are Lipschitz 
	continuous with values in $L^2$, that is, for 
	all $t_1,t_2\in[-T,T]$, for any fixed $T>0$, 
	the functions $u$ and $v$ satisfy
	\begin{equation*}
		\|u(t_1,\cdot)-u(t_2,\cdot)\|_{L^2},
		\|v(t_1,\cdot)-v(t_2,\cdot)\|_{L^2}
		\leq C|t_1-t_2|,
	\end{equation*}
	for some $C=C(E_{u_0},E_{v_0},H_0,T)>0$, 
	where $E_{u_0},E_{v_0}$ 
	and $H_0$ are defined in \eqref{consq};
		
	\item $u(t,x)$ and $v(t,x)$ are H\"older continuous 
	on $[-T,T]\times\mathbb{R}$ with exponent $1/2$ for any 
	fixed $T>0$, that is, for all $t_1,t_2,\in[-T,T]$ 
	and $x_1,x_2\in\mathbb{R}$ we have
	\begin{equation*}
		|u(t_1,x_1)-u(t_2,x_2)|,
		|v(t_1,x_1)-v(t_2,x_2)|\leq
		C\left(|t_1-t_2|^{1/2}
		+|x_1-x_2|^{1/2}\right),
	\end{equation*}
	for some $C=C(E_{u_0},E_{v_0},H_0,T)>0$.
\end{enumerate}
\end{definition}

Next, we introduce the concept of a conservative weak 
solution.

\begin{definition}
[\cite{KR25}, Global conservative weak solution 
of \eqref{t-c-N-n}--\eqref{iid}]\label{defsc}
Suppose that $(u_0,v_0)\in\Sigma$, where $\Sigma$ 
is defined in \eqref{Sig}. 
We define a vector function $(u,v)(t,x)$ to be a 
global conservative weak solution to the Cauchy 
problem \eqref{t-c-N-n}-\eqref{iid} if it meets 
two criteria. First, $(u,v)$ must be a global weak solution 
in accordance with Definition \ref{defs}, and secondly, it 
must adhere to the following five conditions:
\begin{enumerate}
	\item $\int_{-\infty}^{\infty}
	\left(uv+u_xv_x\right)(t,x)\,dx=G_{0}$, 
	for any $t\in\mathbb{R}$;
		
	\item there exist positive Radon measures
	$\lambda_t^{(u)}$,
	$\lambda_t^{(v)}$, and $\lambda_{t}^{(uv)}$
	on $\R$
	such that
	\begin{enumerate}
		\item the absolutely continuous parts of
		$\lambda_t^{(u)}$, $\lambda_t^{(v)}$, 
		and $\lambda_t^{(uv)}$ with respect to the 
		Lebesgue measure on $\R$ have the following form:
		\begin{equation*}
			d\lambda_t^{(u,ac)}=u_x^2\,dx,
			\quad
			d\lambda_t^{(v,ac)}=v_x^2\,dx, \quad 
			d\lambda_t^{(uv,ac)}
			=u_x^2v_x^2\,dx,
		\end{equation*}
		while the nonzero singular parts of
		$\lambda_t^{(u)}$, $\lambda_t^{(v)}$, 
		and $\lambda_t^{(uv)}$ are supported, 
		for a.e.\,\,$t\in\R$, on the sets where
		$v(t,\cdot)=0$, $u(t,\cdot)=0$, and 
		$(uv)(t,\cdot)=0$ respectively;
			
		\item the following conservation 
		laws hold for any $t\in\mathbb{R}$:
		\begin{equation*}
		\begin{split}
			&\int_{-\infty}^{\infty}u^2(t,x)\,dx
			+\lambda_t^{(u)}(\R)
			=E_{u_0},\quad
			\int_{-\infty}^{\infty}v^2(t,x)\,dx
			+\lambda_t^{(v)}(\R)
			=E_{v_0},\\
			&\int_{-\infty}^{\infty}
			\left(3u^2v^2
			+4uvu_xv_x\right)(t,x)\,dx
			+\int_{-\infty}^{\infty}u^2(t,x)
			\,d\lambda_t^{(v)}
			+\int_{-\infty}^{\infty}v^2(t,x)
			\,d\lambda_t^{(u)}
			-\lambda_t^{(uv)}(\R)
			=H_0;
		\end{split}
		\end{equation*}
	
		\item the following inequality
		holds for any $t\in\R$:
		\begin{equation*}
			\int_{-\infty}^{\infty}
			\left(3u^2v^2
			+4uvu_xv_x-u_x^2v_x^2\right)
			(t,x)\,dx
			+\int_{-\infty}^{\infty}u^2(t,x)
			\,d\lambda_t^{(v)}
			+\int_{-\infty}^{\infty}v^2(t,x)
			\,d\lambda_t^{(u)}
			\geq H_0;
		\end{equation*}
		\end{enumerate}
		
		\item 
		the following inequalities 
		hold for any $t\in\R$:
		\begin{align*}
			&\|u(t,\cdot)\|_{H^1}^2\leq E_{u_0},\quad
			\|v(t,\cdot)\|_{H^1}^2\leq E_{v_0}.
		\end{align*}
	\end{enumerate}
	Moreover, introducing the sets
	\begin{subequations}
	\label{DN}
	\begin{align}
		\label{DN-a}
		&D_W(t)=\left\{
		\xi\in\mathbb{R}:\cos\frac{W(t,\xi)}{2}=0\right\},
		\quad
		D_Z(t)=\left\{
		\xi\in\mathbb{R}:\cos\frac{Z(t,\xi)}{2}=0\right\},
		\\
		&N_W=\left\{
		t\in\mathbb{R}:\mathrm{meas}(D_W(t))> 0\right\},
		\quad
		N_Z=\left\{
		t\in\mathbb{R}:\mathrm{meas}(D_Z(t))> 0\right\},
	\end{align}
	\end{subequations}
	where the solution $(U,V,W,Z,q)(t,\xi)$ 
	of the associated ODE system \eqref{ODE} depends 
	solely on the initial conditions, 
	as detailed in Section \ref{cvss} below, 
	and $\mathrm{meas}(\cdot)$ denotes 
	the Lebesgue measure on $\mathbb{R}$, we have
	\begin{enumerate}\setcounter{enumi}{3}
		\item $\|u(t,\cdot)\|_{H^1}^2=E_{u_0}$, 
		for any $t\in\mathbb{R}\setminus N_W$;
		\quad $\|v(t,\cdot)\|_{H^1}^2=E_{v_0}$,
		for any $t\in\mathbb{R}\setminus N_Z$;
		$$
		\int_{-\infty}^{\infty}
		\left(3u^2v^2+u^2v_x^2+v^2u_x^2
		+4uvu_xv_x-u_x^2v_x^2\right)
		(t,x)\,dx=H_0,
		$$
		\qquad\qquad\qquad\qquad\qquad\qquad
		\qquad\qquad\qquad\qquad\qquad\qquad 
		for any $t\in\mathbb{R}\setminus (N_W\cup N_Z)$;
		\item 
		\begin{itemize}
			\item if $\mathrm{meas}(N_W)=0$, 
			then $\lambda_t^{(u)}$ is a measure-valued 
			solution $w_u$ of the following continuity 
			equation with source term:
			$$
			\pt w_u+\px(uvw_u)
			=2u_x(u^2v-P_1-\px P_2)+uu_x^2v_x;
			$$
			
			\item if $\mathrm{meas}(N_Z)=0$, then
			$\lambda_t^{(v)}$ is a measure-valued solution 
			$w_v$ of the following continuity equation with source term:
			$$
			\pt w_v+\px(uvw_v)
			=2v_x(uv^2-S_1-\px S_2)+vu_xv_x^2;
			$$
			
			\item if $\mathrm{meas}(N_W)$, 
			$\mathrm{meas}(N_Z)=0$, then
			$\lambda_t^{(uv)}$ is a measure-valued 
			solution $w_{uv}$ of the following 
			continuity equation with source term:
			\begin{equation*}
				\pt w_{uv}+\px(uvw_{uv})
				=2u_xv_x
				\left((u^2v-P_1-\px P_2)v_x
				+(uv^2-S_1-\px S_2)u_x
				\right);
			\end{equation*}
		\end{itemize}
	\end{enumerate}
\end{definition}

Introduce a positive Radon measure $\mu$
on $\R$ such that
\begin{equation}\label{mu}
	\mu=\mu^{ac}+\mu^s,\quad
	d\mu^{ac}=
	\left(
	u_x^2+v_x^2+u_x^2v_x^2
	\right)dx,\quad
	(u,v)\in\Sigma,
\end{equation}
where the metric space $(\Sigma,d_\Sigma)$ is 
defined by \eqref{Sig}--\eqref{d^2}, and
$\mu^{ac}$ and $\mu^{s}$ respectively denote the absolutely 
continuous and singular parts of $\mu$ with 
respect to the Lebesgue measure on $\R$.
Given a measure $\mu$, we consider the following function 
(cf.~\cite[Section 6]{BC07}, \cite[Equation (7.1)]{KR25}):
\begin{equation}\label{y-def}
	y(\xi)=
	\begin{cases}
		\sup\left\{
		y\in\mathbb{R}:
		y+\mu\left([0,y]\right)
		\leq\xi\right\},&\xi\geq 0,\\
		\inf\left\{
		y\in\mathbb{R}:
		|y|+\mu\left([-y,0)\right)
		\leq-\xi\right\},&\xi< 0,
	\end{cases}
\end{equation}
where $\mu$ satisfies \eqref{mu} 
and introduce the (measurable) sets 
$D_W, D_Z\subset\R$ such that
(notice that $y_\xi\leq 1$ for a.e. $\xi\in\R$,
see \cite[Section 7]{KR25})
\begin{equation}\label{pi-s}
	D_W\cup D_Z=\left\{
	\xi\in\R:y_\xi(\xi)=0
	\right\},\quad
	\text{up to a set of Lebesgue measure zero.}
\end{equation}
Taking into account that $\mathrm{meas}\left(\left\{
\xi\in\R:y_\xi(\xi)=0
\right\}\right)=\mu^s(\R)$, we 
conclude from \eqref{pi-s} that
\begin{equation}\label{mzer}
	\mathrm{meas}(D_W),\mathrm{meas}(D_Z)
	\leq\mu^s(\R).
\end{equation}

Now, we define the following set, as per 
Equation (2.15) in \cite{KR25}:
\begin{equation}
	\begin{split}
		\label{D-set}
		\mathcal{D}
		=&\left\{\left(u,v,\mu;D_W,D_Z\right)
		:(u,v)\in\Sigma,\,\mu
		\mbox{ is a positive Radon measure 
		which satisfies }(\ref{mu}),\right.\\ 
		&\left.\qquad\qquad\qquad\qquad\qquad\quad
		\qquad\qquad
		\mbox{and }D_W,D_Z\mbox{ satisfy } 
		(\ref{pi-s})\mbox{ with }y(\xi)
		\mbox{ given in }(\ref{y-def})
		\right\}.
	\end{split}
\end{equation}
Then we have the following well-posedness result 
for the two-component Novikov system \eqref{t-c-N-n}:

\begin{theorem}
[Global semigroup of conservative solutions]
\label{Thm}
Consider initial data 
$\left(u_0,v_0,\mu_0;D_{W,0},D_{Z,0}\right)\in\mathcal{D}$ 
and define a flow map 
$\Psi_t:\mathbb{R}\times\mathcal{D}
\to\mathcal{D}$ as follows:
\begin{equation}\label{semiPsi}
\Psi_t\left(u_0,v_0,\mu_0;D_{W,0},D_{Z,0}\right)
=\left(u(t),v(t),\mu_{(t)};D_W(t),D_Z(t)\right),
\end{equation}
where $\left(u(t),v(t),\mu_{(t)}\right)$
and $(D_W(t),D_Z(t))$
are defined by \eqref{uvdef}--\eqref{mut} 
and \eqref{DN-a}, respectively, in terms of the solution
$(U,V,W,Z,q)$ of the ODE system \eqref{ODE}--\eqref{id2} 
given by Theorem \ref{gwp}.
Then $\left(u(t),v(t),\mu_{(t)};D_W(t),D_Z(t)\right)$ 
satisfies the following properties:
\begin{enumerate}
	\item $(u,v)$ is a global conservative weak solution 
	of \eqref{t-c-N-n}--\eqref{iid} in the 
	sense of Definition \ref{defsc};
	
	\item $\Psi_t$ satisfies the semigroup property, that 
	is (i) $\Psi_0=\mathrm{id}$ and 
	(ii) $\Psi_{t+\tau}=\Psi_{t}\circ\Psi_{\tau}$;
	
	\item assuming that 
	\begin{enumerate}[a)]
		\item $d_\Sigma\left((u_{0,n},v_{0,n}),
		(u_0,v_0)\right)\to0$ as $n\to\infty$,
		
		\item 
		$\mu_{0,n}
		\overset{\ast}{\rightharpoonup} \mu_0$ 
		(weakly-$\ast$) as $n\to\infty$,

		\item $D_{W,0}=\bigcap\limits_{0<\ve<1}
		D_{W,0}^{\ve}$
		and 
		$D_{Z,0}=\bigcap\limits_{0<\ve<1}
		D_{Z,0}^\ve$
		up to a set of Lebesgue measure zero, 
		where (cf.~\eqref{DN-a})
		\begin{equation*}
			\begin{split}
				&D_{W,0}^\ve=
				\bigcup\limits_{m=1}^\infty
				\bigcap\limits_{n=m}^\infty
				\left\{
				\xi\in\R:
				\left|\cos\frac{W_{0,n}(\xi)}{2}\right|<\ve
				\right\},\\
				&D_{Z,0}^\ve=
				\bigcup\limits_{m=1}^\infty
				\bigcap\limits_{n=m}^\infty
				\left\{
				\xi\in\R:
				\left|\cos\frac{Z_{0,n}(\xi)}{2}\right|<\ve
				\right\},
			\end{split}
		\end{equation*}
		and $W_{0,n}$, $Z_{0,n}$ are defined by \eqref{id2} with
		$\left(u_{0,n}, v_{0,n},\left(D_{W,0}\right)_n,
		\left(D_{Z,0}\right)_n\right)$ instead of
		$\left(u_{0}, v_{0},D_{W,0},D_{Z,0}\right)$,
	\end{enumerate}
	we have $\forall T>0$,
	\begin{equation*}
		\|(u_n-u)\|
		_{L^\infty([-T,T]\times\mathbb{R})}
		+\|(v_n-v)\|
		_{L^\infty([-T,T]\times\mathbb{R})}\to 0,
		\quad n\to\infty.
	\end{equation*}
	Here, $(u_n,v_n)$ are the solutions of 
	\eqref{t-c-N-n} that correspond to the 
	initial data 
	$\left(u_{0,n},v_{0,n},\mu_{0,n};
	\left(D_{W,0}\right)_n,\left(D_{Z,0}\right)_n
	\right)$;
	
	\item $\mu_{(t)}(\mathbb{R})\leq 
	E_{u_0}+E_{v_0}+7E_{u_0}E_{v_0}-H_0$ for any 
	$t\in\mathbb{R}$;
	
	\item if $\mathrm{meas}(N_W)$, $\mathrm{meas}(N_Z)=0$, 
	then $\mu_{(t)}$ is a measure-valued solution $w$ 
	of the following continuity equation with source term:
	\begin{equation*}
	\begin{split}
		w_t+(uvw)_x
		=&2u_x(1+v_x^2)(u^2v-P_1-\px P_2)
		+2v_x(1+u_x^2)(uv^2-S_1-\px S_2)
		\\&+uu_x^2v_x+vu_xv_x^2.
	\end{split}
	\end{equation*}
\end{enumerate}
\end{theorem}

\begin{remark}\label{RemKR}
In \cite[Theorem 2.5]{KR25}, a flow map 
$\Psi_t(u_0,v_0,\mu_0)=(u(t),v(t),\mu_{(t)})$ is introduced. 
However, to ensure the semigroup property, it 
is necessary to also include the sets $D_W(t)$ and $D_Z(t)$, 
which appear in the definition of the initial data $W_0, Z_0$, 
see \eqref{id2}. The proof of items (1), (2), (4), and (5) 
in Theorem \ref{Thm} follows from \cite[Section 7]{KR25}.
Notice that the definition of the initial data $W_0,Z_0$ 
given in \eqref{id2} allows one to establish 
that $\tilde{W}_0(\sigma)=W(\tau,\xi)$ and
$\tilde{Z}_0(\sigma)=Z(\tau,\xi)$ for all $\tau\in\R$, where 
(see the proof of item (2) in \cite[Section 7]{KR25})
\begin{equation}\label{sigd}
	\sigma=\sigma(\xi)=\int_{\xi_0}^\xi q(\tau,\eta)\,d\eta,
\end{equation}
and $\tilde{W}_0,\tilde{Z}_0$ are 
defined by \eqref{id2} in terms of
\[
\left(\tilde{u}_0,\tilde{v}_0,\tilde{\mu}_0;
\tilde{D}_{W,0},\tilde{D}_{Z,0}\right)
=\Psi_\tau\left(u_0,v_0,\mu_0;D_{W,0},D_{Z,0}\right).
\]

To show item (3), we use the same arguments 
as in the proof of item (3) in \cite[Section 7]{KR25}, 
paying additional attention to the limits
\begin{equation}\label{lim-L2}
	\|W_{0,n}-W_0\|_{L^2(\R)}\to0,\quad
	\|Z_{0,n}-Z_0\|_{L^2(\R)}\to0,\quad n\to\infty.
\end{equation}
To prove \eqref{lim-L2},
we must show, in addition to \cite{KR25}, that
$\|W_{0,n}-W_0\|_{L^2(D_{W,0})}\to0$ 
and $\|Z_{0,n}-Z_0\|_{L^2(D_{Z,0})}\to0$ as $n\to\infty$.
Since $D_{W,0}\subset D_{W,0}^\ve$ for any small 
$\ve>0$, up to a set of Lebesgue measure zero 
(see assumption (3,c)), we conclude that 
for almost every $\xi\in D_{W,0}$ there exists 
$m\in\N$, with $m=m(\xi)$, such that
\begin{equation}\label{aexi}
	\xi\in\bigcap\limits_{n=m}^\infty
	\left\{
	\eta\in\R:\left|\cos\frac{W_{0,n}(\eta)}{2}\right|<\ve
	\right\}.
\end{equation}
Assuming, without loss of generality, that
$W_0,W_{0,n}\in(-\pi+\ve_0,\pi+\ve_0]$, for some 
$0<\ve\ll\ve_0<1$ (see Remark \ref{WZval} below), we obtain
the following estimate for all $n\geq m$ from \eqref{aexi}:
\begin{equation*}
	\left|W_{0,n}(\xi)-W_{0}(\xi)\right|
	\leq\left|W_{0,n}(\xi)-\pi\right|
	+\left|\pi-W_{0}(\xi)\right|
	=\left|W_{0,n}(\xi)-\pi\right|\leq C\ve,
\end{equation*}
for some $C>1$ and almost every $\xi\in D_{W,0}$.
Since $\ve>0$ is arbitrary, we conclude that
$W_{0,n}\to W_{0}$ almost everywhere in $D_{W,0}$.
Finally, taking into account that 
$|W_{0,n}|$ and $|W_0|$ are uniformly 
bounded on $\R$ and that $\mathrm{meas}(D_{W,0})<\infty$, 
see \eqref{mzer}, we conclude that
$\|W_{0,n}-W_0\|_{L^2(D_{W,0})}\to0$.
The proof that $\|Z_{0,n}-Z_0\|_{L^2(D_{Z,0})}\to0$ 
is similar.
\end{remark}

\begin{remark}
For the two-component Novikov system, 
concentration of the three energies $u_x^2\,dx$, 
$v_x^2\,dx$, and $u_x^2 v_x^2\,dx$ may 
occur (see item (2) in Definition \ref{defsc}).
To track these energies for all $t$, we introduce 
the positive Radon measures $\lambda_t^{(u)}$, 
$\lambda_t^{(v)}$, and $\lambda_t^{(uv)}$, which 
can be expressed in terms of $(U,V,W,Z,q)$ as 
follows (see \cite[Equation (6.6)]{KR25}):
\begin{equation}
\label{rm}
\begin{split}
	&\lambda^{(u)}_t\left([a,b]\right)
	=\int_{\{\xi:y(t,\xi)\in[a,b]\}}
	\left(q\sin^2\frac{W}{2}
	\cos^2\frac{Z}{2}\right)(t,\xi)\,d\xi,
	\\
	&\lambda^{(v)}_t\left([a,b]\right)
	=\int_{\{\xi:y(t,\xi)\in[a,b]\}}
	\left(q\cos^2\frac{W}{2}
	\sin^2\frac{Z}{2}\right)(t,\xi)\,d\xi,
	\\
	&\lambda^{(uv)}_t\left([a,b]\right)
	=\int_{\{\xi:y(t,\xi)\in[a,b]\}}
	\left(q\sin^2\frac{W}{2}
	\sin^2\frac{Z}{2}\right)(t,\xi)
	\,d\xi,
\end{split}
\end{equation}

Notice, however, that the map 
\begin{equation*}
	\widetilde\Psi_t
	\left(u_0,v_0,\lambda_0^{(u)},\lambda_0^{(v)},
	\lambda_0^{(uv)}\right)=
	\left(u(t),v(t),\lambda_t^{(u)},\lambda_t^{(v)},
	\lambda_t^{(uv)}\right)
\end{equation*}
is not a semigroup, since one cannot retrieve the sets 
$D_W(t)$ and $D_Z(t)$ from the measures
$\lambda_t^{(u)}$, $\lambda_t^{(v)}$, and
$\lambda_t^{(uv)}$ and thus 
define the initial data \eqref{id2}. 
Indeed, consider, for example, the characteristic 
$y(\tau,\xi)$ such that $y_\xi=0$ on the 
maximal interval $\xi\in(\xi_1,\xi_4)$, 
and assume that $W=\pi$ 
for $\xi\in(\xi_1,\xi_2)$ and $Z=\pi$ for
$\xi\in(\xi_3,\xi_4)$, $\xi_3<\xi_2$.
Introducing $\hat{y}=y(\tau,\xi)$, 
$\xi\in(\xi_1,\xi_4)$, we have
\[
\left(\lambda^{(u)}_\tau+\lambda^{(uv)}_\tau\right)
\left(\left\{\hat{y}\right\}\right)=
\int_{\xi_1}^{\xi_2}q(\tau,\eta)\,d\eta=
\sigma(\xi_2)-\sigma(\xi_1),
\]
with $\sigma(\xi)$ given by \eqref{sigd}.
Thus, once the measures \eqref{rm} are fixed, 
we retain only the Lebesgue measure of 
the set on which $\tilde{W}_0(\sigma)
=W(\tau,\xi(\sigma))=\pi$, but 
not its precise location.
\end{remark}

\begin{remark}[$L^1(\R)$ solutions of \eqref{t-c-N-n}--\eqref{iid}]
	If the initial data satisfy the additional assumptions
		$u_0,v_0\in W^{1,1}(\R)$, then the corresponding global
		solution $(u(t),v(t))$ in \eqref{semiPsi} retains this
		regularity and decay. More precisely, by Remark \ref{L1ODE},
	\eqref{uvdef}, and \eqref{pxuv} below,
	\begin{equation*}
		u(t,\cdot),\,v(t,\cdot)\in W^{1,1}(\R),\quad t\in[-T,T].
	\end{equation*}
	Indeed, \eqref{uvdef} implies that
	(see \eqref{DN-a} and \eqref{pxiy}; we omit $t$ here)
	\begin{equation*}
		\|u\|_{L^1(\R)}=\int_{\xi\in\R\setminus\{D_W\cup D_Z\}}
		|U(\xi)|y_\xi(\xi)\,d\xi
		=\int_{-\infty}^\infty|U(\xi)|y_\xi(\xi)\,d\xi
		\leq\|q\|_{L^\infty}\|U\|_{L^1}.
	\end{equation*}
		It then follows from \eqref{pxuv} that
	\begin{equation*}
	\begin{split}
		\|u_x\|_{L^1(\R)}&=\int_{\xi\in\R\setminus\{D_W\cup D_Z\}}
		\left|\tan\frac{W(\xi)}{2}\right|y_\xi(\xi)\,d\xi
		=\int_{-\infty}^\infty
		q\left|\sin\frac{W(\xi)}{2}\cos\frac{W(\xi)}{2}\right|
		\cos^2\frac{Z(\xi)}{2}\,d\xi\\
		&\leq\|q\|_{L^\infty}\|W\|_{L^1}.
	\end{split}
	\end{equation*}
	Employing Remark \ref{L1ODE}, we conclude that
	the norms $\|u\|_{L^1}$ and $\|u_x\|_{L^1}$ are uniformly bounded.
	The same argument applies to $v$; hence the $W^{1,1}$ norms
	of both components are uniformly bounded.
\end{remark}

\medskip

\noindent\textbf{Notations.} Introduce the function
\begin{equation}\label{D}
	D(t,x)=\left(
	\left(1+u_x^2\right)
	\left(1+v_x^2\right)\right)(t,x).
\end{equation}
Consider the linear operator 
$\mathcal{I}_\alpha$ defined by
\begin{equation}\label{I-op}
	\mathcal{I}_\alpha(f)
	=\int_{-\infty}^{\infty}
	f(x)e^{-\alpha|x|}\,dx,\quad \alpha\in(0,1),
\end{equation}
where $f\in L^1_{\mathrm{loc}}(\R)\cap
L^\infty(\R\setminus[-R,R])$ with some $R>0$.
Notice that $\mathcal{I}_\alpha$ is monotone, i.e.,
\begin{equation}\label{mon}
	\mathcal{I}_\alpha(f)\leq
	\mathcal{I_\alpha}(g),\quad
	\text{if }f\leq g.
\end{equation}
Also we adopt notation
\begin{equation}\label{empts}
	\{a_i\}_{i=n_1}^{n_2}=\emptyset
	\quad\text{and}\quad
	\bigcup\limits_{i=n_1}^{n_2}a_i=\emptyset,\quad
	\text{if }n_2<n_1.
\end{equation}

We use the following Banach space:
\begin{equation}
	\label{EB}
	E_{q_\infty}=\left(H^{1}(\R)
	\cap W^{1,4}(\R)\right)^2\times
	\left(L^2(\R)
	\cap L^{\infty}(\R)\right)^2\times
	\left(L^{\infty}(\R)\cap L^1_{q_\infty}(\R)\right),
	\quad q_{\infty}\in\R,
\end{equation}
where
\begin{equation*}
	L^1_{q_{\infty}}(\R)=\left\{q\in L^1_{\mathrm{loc}}(\R):
	q(x)=p(x)+q_\infty,\,\,\text{with some }\,p\in L^1(\R)\right\},
	\quad q_{\infty}\in\R,
\end{equation*}
equipped with the norm
\begin{equation}\label{nE}
	\|(U,V,W,Z,q)\|_{E_{q_\infty}}=\|U\|_{H^1\cap W^{1,4}}
	+\|V\|_{H^1\cap W^{1,4}}
	+\|W\|_{L^2\cap L^{\infty}}
	+\|Z\|_{L^2\cap L^{\infty}}
	+\|q\|_{L^{\infty}}+\|q-q_\infty\|_{L^1}.
\end{equation}
Define a closed subset 
$\Omega\subset E_{q_\infty}$ with $q_\infty=1$ as follows
(see also Remark \ref{rEB} below):
\begin{equation}\label{Om}
\begin{split}
	\Omega= & \Bigl\{
	(U,V,W,Z,q)\in E_1:
	\|U\|_{H^1\cap W^{1,4}},\|V\|_{H^1\cap W^{1,4}}
	\leq R_1,
	\,\,
	\|W\|_{L^2\cap L^\infty},\|Z\|_{L^2\cap L^\infty}\leq R_2,\\
	&
	\qquad\qquad\qquad\qquad\quad\quad
	q^-\leq q(\xi)\leq q^+,
	\,\,\xi\in\mathbb{R}
	\Bigr\},
\end{split}
\end{equation}
for some $R_1,R_2,q^-,q^+>0$.
In the context of the global solutions
of the two-component Novikov system \eqref{t-c-N-n},
we consider functions $(U,V,W,Z,q)\in\Omega$,
which satisfy the following compatibility conditions (see Section \ref{cvss} below):
\begin{equation}\label{ccUV}
	U_\xi(\xi)=\left(\frac{q}{2}\sin W\cos^2\frac{Z}{2}\right)(\xi),
	\quad
	V_\xi(\xi)=\left(\frac{q}{2}\cos^2\frac{W}{2}\sin Z\right)(\xi),
	\quad \xi\in\R.
\end{equation}
This motivates us to define the following subset of $\Omega$, see \eqref{Om}
(here, the subscript $\mathrm{tcN}$ refers to the two-component Novikov system):
\begin{equation}\label{tcNOm}
	\Omega_{\mathrm{tcN}}=\left\{
	(U,V,W,Z,q)\in\Omega,\text{ which satisfy }\eqref{ccUV}
	\right\}.
\end{equation}
\begin{remark}
	The set $\Omega_{\mathrm{tcN}}$ defined in \eqref{tcNOm}
	(see also \eqref{Om} and \eqref{EB}) is a closed subset of $E_1$.
\end{remark}
\begin{proof}
	Consider any fundamental sequence
	$(U_n,V_n,W_n,Z_n,q_n)\in\Omega_{\mathrm{tcN}}$,
	which converges to $(U,V,W,Z,q)$ in $E_1$.
	Since $\Omega_{\mathrm{tcN}}\subset\Omega$
	and $\Omega$ is a closed subset of $E_1$, we have that
	$(U,V,W,Z,q)$ belongs to $\Omega$.
	Thus, the task is to show that
	it satisfies conditions in \eqref{ccUV}.
	We have that
	\begin{equation}\label{WnZnqn}
		(W_n,Z_n,q_n)(\xi)\to(W,Z,q)(\xi),\quad\text{for a.e. }\xi\in\R.
	\end{equation}
	Moreover, after passing to a
	common subsequence if necessary, $\pxi U_n\to\pxi U$ and
	$\pxi V_n\to\pxi V$ a.e.~on $\R$.
	This, together with \eqref{WnZnqn} and that $(U_n,V_n,W_n,Z_n,q_n)$
	satisfies \eqref{ccUV}, yields \eqref{ccUV} for $(U,V,W,Z,q)$.
\end{proof}

We will also 
use the following Banach spaces (see \eqref{EB} with $q_{\infty}=0$):
\begin{equation}\label{E-0}
	\mathcal{C}^k_{E_0}=
	\left(\left(C^{k}(\R)\right)^2\times
	\left(C^{k-1}(\R)\right)^3\right)
	\cap E_0,\quad k\in\N,
\end{equation}
and (see \eqref{tcNOm})

\begin{equation}\label{CktcN}
	\mathcal{C}^k_{\Omega_{\mathrm{tcN}}}=
	\left(\left(C^{k}(\R)\right)^2\times
	\left(C^{k-1}(\R)\right)^3\right)
	\cap\Omega_{\mathrm{tcN}},\quad k\in\N.
\end{equation}
In Section \ref{LM}, we will use the following spaces,
which impose $L^1(\R)$ decay on $(U,V,W,Z)$:
\begin{equation}\label{OmL1tcN}
\begin{split}
	&E_0^1=\left\{(U,V,W,Z,q)\in E_0:
	U,V\in W^{1,1}(\R),\,\, W,Z\in L^1(\R)
	\right\},\\
	&\Omega_{\mathrm{tcN}}^1=
	\left\{(U,V,W,Z,q)\in\Omega_{\mathrm{tcN}}:
	U,V\in W^{1,1}(\R),\,\, W,Z\in L^1(\R)
	\right\},
\end{split}
\end{equation}
and
\begin{equation}\label{P-pth}
	\mathcal{P}^{k,1}_{\mathrm{tcN}}[\alpha,\beta]=
	C\left([\alpha,\beta],
	\mathcal{C}^{k}_{\Omega_{\mathrm{tcN}}^1}\times
	\mathcal{C}^{k}_{E_0^1}\right),
	\quad k\in\N,\quad \alpha<\beta,
\end{equation}
where (cf.~\eqref{E-0} and \eqref{CktcN})
\begin{equation}\label{C1ktcN}
\begin{split}
	&\mathcal{C}^k_{E_0^1}=
	\left(\left(C^{k}(\R)\right)^2\times
	\left(C^{k-1}(\R)\right)^3\right)
	\cap E_0^1,\quad k\in\N,\\
	&\mathcal{C}^k_{\Omega_{\mathrm{tcN}}^1}=
	\left(\left(C^{k}(\R)\right)^2\times
	\left(C^{k-1}(\R)\right)^3\right)
	\cap \Omega_{\mathrm{tcN}}^1,\quad k\in\N.
\end{split}
\end{equation}
The norms in the spaces mentioned above 
are induced by the direct sum of the 
Banach spaces, similar to \eqref{nE}. 
Finally, we adopt the notations (see \eqref{D-set})
\begin{equation}\label{u-bf}
	\mathbf{u}=\left(u,v,\mu;
	D_{W},D_{Z}\right),\quad
	\mathbf{u}\in\mathcal{D},
\end{equation}
and (cf.~\eqref{Om})
\begin{equation}\label{U-bf}
	\mathbf{U}^\theta
	=\left(U^\theta,V^\theta,W^\theta,Z^\theta,q^\theta
	\right),\quad \mathbf{U}^\theta\in
	\Omega_{\mathrm{tcN}}^1\,
	\text{ for all fixed }\,\theta\in[0,1].
\end{equation}

\medskip

\noindent\textbf{Basic facts.}
We will employ the following 
estimates throughout the paper:
\begin{align}
	\label{Cin}
	&2|a|\leq 1+a^2,\quad a\in\R,\\
	\label{Sob}
	&\|f\|_{L^\infty(\R)}\leq\|f\|_{H^1(\R)}\quad
	\text{(Sobolev inequality)},\\
	\label{sinx}
	&|\sin x|\leq|x|,\quad x\in\R.
\end{align}

Let us recall several notions concerning 
$C^k$ maps between smooth manifolds.

A regular value of such a map is a point in the 
target manifold at which the differential is surjective 
at every point in its preimage; in this case, the preimage 
forms a smooth submanifold (see, 
for example, \cite[Definition 2.9]{MN20}).

\begin{definition}[Regular value]\label{RV}
Let $f:\mathcal{X}\to\mathcal{Y}$ be 
a $C^1$ mapping between manifolds.
Then $y\in\mathcal{Y}$ is called a 
regular value of $f$, if for 
all $x\in f^{-1}(y)$ one has
$$
\left.(df)\right|_x(T_x\mathcal{X})
=T_{f(x)}\mathcal{Y}.
$$
In particular, if $f^{-1}(y)=\emptyset$, 
then $y$ is a regular value.
\end{definition}

We next introduce the more general notion of 
transversality, which describes how a regular function 
meets a submanifold in its codomain.
In this setting, one requires that the intersection 
occur cleanly and without tangencies, in a way that remains 
stable under small perturbations 
of the map (see \cite[Definition 2.17]{MN20}).

\begin{definition}[Transversality]
Let $f:\mathcal{X}\to\mathcal{Y}$ be a $C^1$ 
mapping between manifolds and $\mathcal{W}\subset\mathcal{Y}$ 
a submanifold. We say that $f$ is transverse 
to $\mathcal{W}$ at a point 
$x\in\mathcal{X}$, if $f(x)\in\mathcal{W}$ and 
$$
\left.(df)\right|_x(T_x\mathcal{X})
+T_{f(x)}\mathcal{W}
=T_{f(x)}\mathcal{Y}.
$$
We say that $f$ is transverse to $\mathcal{W}$ 
if it is either transverse at every point 
$x\in\mathcal{X}$ such that $f(x)\in\mathcal{W}$ or 
$f^{-1}(\mathcal{W})=\emptyset$.
\end{definition}

\begin{remark}\label{tr=rv}
If $\mathcal{W}=\{y\}$, then $f$ is transverse to 
$\mathcal{W}$ if and only if $y$ 
is a regular value of $f$.
\end{remark}

The transversality theorem, a crucial component 
in the proof of the generic regularity 
result presented in Theorem \ref{Thmgr}, is 
as follows; see \cite[Chapter 3, Theorem 2.7]{H76}.

\begin{theorem}[Parametric transversality theorem]
\label{BTL}
Let $F:\mathcal{X}\times\mathcal{N}\to\mathcal{Y}$ 
be a mapping between manifolds and 
$\mathcal{W}\subset\mathcal{Y}$ a submanifold.
Assume that
\begin{equation}\label{F-reg}
	F\in C^k(\mathcal{X}\times\mathcal{N}),
	\quad k>\max\{0,\dim\mathcal{Y}+\dim\mathcal{W}
	-\dim\mathcal{X}\}.
\end{equation}
Denote $f^\nu(x)=F(x,\nu)$, where 
$x\in\mathcal{X}$ and $\nu\in\mathcal{N}$. 
If $F$ is transverse to $\mathcal{W}$, then
there exists a dense subset 
$\tilde{\mathcal{N}}\subset\mathcal{N}$ such that
$f^\nu$ is transverse to $\mathcal{W}$ for all 
$\nu\in\tilde{\mathcal{N}}$.
\end{theorem}

Finally, we recall the regular value theorem (see 
\cite[Chapter 1, Theorem 3.2]{H76} and 
\cite[Theorem 2.3]{MN20}).

\begin{theorem}[Regular value theorem]\label{RVT}
Let $f:\mathcal{X}\to\mathcal{Y}$ be a $C^k$ 
mapping between manifolds, $k\geq 1$ and 
$\dim(\mathcal{Y})\leq\dim(\mathcal{X})$.
If $y\in f(\mathcal{X})$ is a regular value, then
$f^{-1}(y)$ is a $C^k$ submanifold of 
$\mathcal{X}$ of dimension 
$\dim(\mathcal{X})-\dim(\mathcal{Y})$.
\end{theorem}

In Section \ref{Tvss}, we apply the operator 
$\mathcal{I}_\alpha$ to functions of the form 
$D(1-\partial_x^2)^{-1}(f)$, where $D$ is defined 
in \eqref{D} and $f$ is a given function (see, for 
instance, \eqref{P4D} below).
To obtain the expression appearing on the right-hand 
side of estimate \eqref{clin}, we make use of the 
inequality stated in the following lemma.

\begin{lemma}\label{LID}
	Assume that $D$ has the form $D(t,x)=1+\tilde{D}(t,x)$, where
	\begin{equation}\label{Dtilu}
		\left\|\tilde{D}(t,\cdot)\right\|_{L^1(\R)}\leq C,\quad
		\forall t\in[-T,T],\quad\forall T>0,\quad
		\text{for some }\,C>0.
	\end{equation}
	Consider $f\in L^1_{\mathrm{loc}}(\R)\cap
	L^\infty(\R\setminus[-R,R])$ for some $R>0$.
	Then for each $t$ we have, see \eqref{I-op},
	\begin{equation}\label{ID}
		\mathcal{I}_\alpha\left(
		\left|
		\int_{-\infty}^{\infty}e^{-|x-y|}f(y)\,dy
		\right|D(t,x)
		\right)
		\leq C(\alpha)
		\mathcal{I}_\alpha\left(|f|\right),\quad
		\text{for some }\,C(\alpha)>0.
	\end{equation}
\end{lemma}

\begin{proof}
	Recalling the definition \eqref{I-op} 
	of $\mathcal{I}_\alpha$ 
	and changing the order of integration, 
	we have (the argument $t$ is dropped here)
	\begin{equation}
		\label{ID-1}
		\mathcal{I}_\alpha\left(
		\left|
		\int_{-\infty}^{\infty}e^{-|x-y|}f(y)\,dy
		\right|D(x)
		\right)
		\leq
		\int_{-\infty}^{\infty}
		\left(
		\int_{-\infty}^{\infty}
		e^{-|x-y|-\alpha|x|+\alpha|y|}D(x)\,dx\right)
		|f(y)|e^{-\alpha|y|}\,dy.
	\end{equation}
	For $y\geq0$, we obtain the following estimate:
	\begin{equation}\label{ID-i}
		\begin{split}
			\int_{-\infty}^{\infty}
			e^{-|x-y|-\alpha|x|+\alpha|y|}D(x)\,dx&=
			e^{(\alpha-1)y}
			\left(
			\int_{-\infty}^0 e^{(1+\alpha)x}
			D(x)\,dx
			+\int_{0}^y e^{(1-\alpha)x}
			D(x)\,dx
			\right)\\
			&\quad
			+e^{(1+\alpha)y}\int_y^\infty 
			e^{-(1+\alpha)x}D(x)\,dx.
		\end{split}	
	\end{equation}
	Recalling that $D = 1 + \tilde{D}$ and the 
	uniform bound \eqref{Dtilu}, we 
	conclude from \eqref{ID-i} that
	\begin{equation*}
		\int_{-\infty}^{\infty}
		e^{-|x-y|-\alpha|x|+\alpha|y|}D(x)\,dx
		\leq\frac{1+e^{(\alpha-1)y}}{1+\alpha}+
		\frac{1-e^{(\alpha-1)y}}{1-\alpha}
		+\left(2+e^{(\alpha-1)y}\right)
		\left\|\tilde{D}(t,\cdot)
		\right\|_{L^1(\R)}\leq C(\alpha),
	\end{equation*}
	where $y\geq0$. Arguing in the 
	same way for $y<0$, we obtain
	\begin{equation}\label{ID-2}
		\int_{-\infty}^{\infty}
		e^{-|x-y|-\alpha|x|+\alpha|y|}D(x)\,dx
		\leq C(\alpha),
		\quad y\in\R.
	\end{equation}
	Combining \eqref{ID-1} and \eqref{ID-2}, we 
	arrive at \eqref{ID}.
\end{proof}

\begin{remark}\label{expLip}
	The condition $\alpha<1$ is essential 
	for the derivation of \eqref{ID}.  
	Indeed, if one formally sets $\alpha=1$ 
	and proceeds as in \eqref{ID-i},  
	one arrives at the estimate (cf.\,\eqref{ID-2})
	\begin{equation*}
		\int_{-\infty}^{\infty}
		e^{-|x-y|-|x|+|y|}D(x)\,dx
		\leq C(E_{u_0},E_{v_0},H_0)(1+|y|),
		\quad y\in\R,
	\end{equation*}
	which is not sufficient to recover \eqref{ID} 
	in the case $\alpha=1$.   In light of this observation, 
	it appears natural that in \cite[Equation (3.9)]{CCCS18},  
	where the norm of the tangent vector is defined 
	for the Novikov equation,  the factor $e^{-|x|}$ could be 
	replaced by $e^{-\alpha|x|}$ with $\alpha\in(0,1)$.
	We also note that the weight function is important in the
	estimates presented in equations
	\eqref{3D-Om-est}--\eqref{Rnmthest}, though
	the norm \eqref{rstr} is finite
	even with $\alpha=0$ for solutions in $W^{1,1}$.
\end{remark}

\section{Generic regularity}\label{chv}

In this section, we establish the generic regularity 
result presented in Theorem \ref{Thmgr}. We employ the 
term ``generic" here because the properties (1)--(3) of 
the solution in Theorem \ref{Thmgr} hold for all 
initial data $(u_0,v_0)$ originating 
from the subset
$\mathcal{M}_T\subset\Upsilon^k$, which is open in $\Upsilon^k$
and dense in the topology specified in Theorem \ref{Thmgr}.

\subsection{Change of variables and a semilinear system}
\label{cvss}
In this subsection, we briefly recall 
the Bressan-Constantin 
formalism for reducing \eqref{t-c-N-n}, with 
initial data $\left(u_0,v_0,\mu_0;
D_{W,0},D_{Z,0}\right)\in\mathcal{D}$ (see \eqref{D-set}), 
to an equivalent ODE system. Further details can 
be found in \cite{HQ21, KR25}.

Define the characteristic $y(t, \xi)$ 
as the solution to the ODE:
\begin{equation}\label{char}
	\begin{split}
		y_t(t,\xi)=u(t,y(t,\xi))v(t,y(t,\xi)),
		\quad y(0,\xi)=y_0(\xi),
	\end{split}
\end{equation}
where the initial data $y_0(\xi)$ is a 
monotone increasing function given by
\eqref{y-def} with $\mu_0$ instead of $\mu$.
Introduce the following new 
unknowns (see \cite{HQ21, KR25}):
\begin{equation}\label{UV}
	U(t,\xi)=u(t,y(t,\xi)),\quad
	V(t,\xi)=v(t,y(t,\xi)),
\end{equation}
\begin{equation}\label{WZ}
	W(t,\xi)=2\arctan u_x(t,y(t,\xi)),\quad
	Z(t,\xi)=2\arctan v_x(t,y(t,\xi)),
\end{equation}
and (recall \eqref{D})
\begin{equation}\label{q}
	q(t,\xi)=\left(\left(1+u_x^2\right)
	\left(1+v_x^2\right)\right)(t,y(t,\xi))y_\xi(t,\xi)
	=D(t,y(t,\xi))y_\xi(t,\xi).
\end{equation}

It can be shown \cite{HQ21,KR25} that if $(u,v)$ 
satisfies \eqref{t-c-N}, then $(U,V,W,Z,q)(t,\xi)$ 
formally solves the following system of ODEs in a Banach space
(we drop the arguments $t,\xi$ for simplicity):
\begin{equation}\label{ODE}
	\begin{split}
		&U_t=-\px P_1-P_2,\\
		&V_t=-\px S_1-S_2,\\
		&W_t=
		2U^2V\cos^2\frac{W}{2}-V\sin^2\frac{W}{2}
		-2\left(P_1+\px P_2\right)\cos^2\frac{W}{2},\\
		&Z_t=
		2UV^2\cos^2\frac{Z}{2}-U\sin^2\frac{Z}{2}
		-2\left(S_1+\px S_2\right)\cos^2\frac{Z}{2},\\
		&q_t=
		q\left(U^2V+\frac{1}{2}V-P_1-\px P_2\right)
		\sin W
		+q\left(UV^2+\frac{1}{2}U
		-S_1-\px S_2\right)\sin Z,
	\end{split}
\end{equation}
subject to the initial data (see 
\eqref{UV}, \eqref{WZ} and \eqref{q})
\begin{equation}\label{id2}
	\begin{split}
		&U_0(\xi):= U(0,\xi)=u_0(y_0(\xi)),\quad
		V_0(\xi):= V(0,\xi)=v_0(y_0(\xi)),\\
		&W_0(\xi):= W(0,\xi)=
		\begin{cases}
			\pi&\mbox{if }\xi\in D_{W,0},\\
			2\arctan(\px u_0)(y_0(\xi)),&
			\mbox{otherwise},\\
		\end{cases}\\
		&Z_0(\xi):= Z(0,\xi)=
		\begin{cases}
			\pi&\mbox{if }\xi\in D_{Z,0},\\
			2\arctan(\px v_0)(y_0(\xi)),&
			\mbox{otherwise},\\
		\end{cases}\\
		& q_0(\xi):= q(0,\xi)=1.
	\end{split}
\end{equation}
Here (we slightly abuse notations by writing, for example, 
$P_1(t,\xi)$ instead of $P_1(t,y(t,\xi))$)
\begin{align}
	\label{P-i-d}
	&P_1(t,\xi)=\frac{1}{2}
	\int_{-\infty}^{\infty}
	\mathcal{E}(t,\xi,\eta)
	p_1(t,\eta)\,d\eta,\quad
	P_2(t,\xi)=\frac{1}{8}
	\int_{-\infty}^{\infty}
	\mathcal{E}(t,\xi,\eta)
	p_2(t,\eta)\,d\eta,\\
	\nonumber
	&(\px P_1)(t,\xi)=\frac{1}{2}
	\left(
	\int_{\xi}^{\infty}
	-\int_{-\infty}^{\xi}
	\right)
	\mathcal{E}(t,\xi,\eta)
	p_1(t,\eta)\,d\eta,\\
	\label{P2-i-d}
	&(\px P_2)(t,\xi)=\frac{1}{8}
	\left(
	\int_{\xi}^{\infty}
	-\int_{-\infty}^{\xi}
	\right)
	\mathcal{E}(t,\xi,\eta)
	p_2(t,\eta)\,d\eta,
\end{align}
and
\begin{align*}
	&S_1(t,\xi)=\frac{1}{2}\int_{-\infty}^{\infty}
	\mathcal{E}(t,\xi,\eta)
	s_1(t,\eta)\,d\eta,\quad
	S_2(t,\xi)=\frac{1}{8}
	\int_{-\infty}^{\infty}
	\mathcal{E}(t,\xi,\eta)
	s_2(t,\eta)\,d\eta,\\
	&(\px S_1)(t,\xi)=\frac{1}{2}
	\left(
	\int_{\xi}^{\infty}
	-\int_{-\infty}^{\xi}
	\right)
	\mathcal{E}(t,\xi,\eta)
	s_1(t,\eta)\,d\eta,\\
	&(\px S_2)(t,\xi)=\frac{1}{8}
	\left(
	\int_{\xi}^{\infty}
	-\int_{-\infty}^{\xi}
	\right)
	\mathcal{E}(t,\xi,\eta)
	s_2(t,\eta)\,d\eta,
\end{align*}
where
\begin{equation}
	\label{E}
	\mathcal{E}(t,\xi,\eta)
	=\exp\left(-\left|\int_\eta^\xi
	\left(q\cos^2\frac{W}{2}
	\cos^2\frac{Z}{2}\right)(t,s)\,ds\right|\right),
	\quad t,\xi,\eta\in\mathbb{R},
\end{equation}
\begin{equation}
	\label{p1q1}
	\begin{split}
		&p_1(t,\xi)=q(t,\xi)\left(
		U^2V
		\cos^2\frac{W}{2}
		\cos^2\frac{Z}{2}
		+\frac{1}{4}U\sin W\sin Z
		+\frac{1}{2}V\sin^2\frac{W}{2}
		\cos^2\frac{Z}{2}
		\right)(t,\xi),
		\\ & 
		s_1(t,\xi)=q(t,\xi)
		\left(UV^2
		\cos^2\frac{W}{2}
		\cos^2\frac{Z}{2}
		+\frac{1}{2}U\cos^2\frac{W}{2}
		\sin^2\frac{Z}{2}
		+\frac{1}{4}V\sin W\sin Z
		\right)(t,\xi),
	\end{split}
\end{equation}
and
\begin{equation}
	\label{p2q2}
	\begin{split}
		p_2(t,\xi)=\left(q\sin^2\frac{W}{2}
		\sin Z\right)(t,\xi),\quad
		s_2(t,\xi)=\left(q\sin W\sin^2\frac{Z}{2}
		\right)(t,\xi).
	\end{split}
\end{equation}

\begin{remark}\label{Rmz}
Recalling \eqref{id2} and \eqref{EB}, we observe that in 
the definition of $W_0$ and $Z_0$ the sets
$D_{W,0}$ and $D_{Z,0}$ may be chosen up to 
sets of Lebesgue measure zero.
In particular, if $\mu_0^s=0$, then 
by \eqref{mzer} one may take 
$D_{W,0}=D_{Z,0}=\emptyset$.
\end{remark}

\begin{remark}
If the initial measure $\mu_0$ has no singular part, then
$y_0$ can be determined as follows:
\begin{equation}\label{y0}
	\int_0^{y_0(\xi)}\left(1+(\px u_0)^2(x)\right)
	\left(1+(\px v_0)^2(x)\right)dx=\xi.
\end{equation}
Moreover, in this case the initial data \eqref{id2} 
are equivalent to (see Remark \ref{Rmz})
\begin{equation}\label{id}
	\begin{split}
		&U_0(\xi)=u_0(y_0(\xi)),\quad
		V_0(\xi)=v_0(y_0(\xi)),\\
		&W_0(\xi)
		=2\arctan(\px u_0)(y_0(\xi)),\quad
		Z_0(\xi)
		=2\arctan(\px v_0)(y_0(\xi)),\\
		& q_0(\xi)=1.
	\end{split}
\end{equation}
\end{remark}
Given initial data
$\left(u_0,v_0,\mu_0;D_{W,0},D_{Z,0}\right)
\in\mathcal{D}$ (see \eqref{D-set}),
direct computation shows that the functions 
$U_0$, $V_0$, $W_0$, $Z_0$, and $q_0$ defined in
\eqref{id2} satisfy the compatibility conditions \eqref{ccUV}.
Therefore,
we consider the Cauchy problem 
\eqref{ODE}--\eqref{id2} in the closed subset 
$\Omega_{\mathrm{tcN}}\subset E_1$ defined by 
\eqref{EB} and \eqref{tcNOm}.
Then, arguing along the same lines as in \cite[Sections IV--V]{HQ21} 
(see also \cite[Sections 4--5]{KR25}), we obtain the
following global well-posedness result for the ODE 
system (see, however, Remarks \ref{Uin} and \ref{rEB} below).
\begin{theorem}\label{gwp}
Consider initial data $(U_0,V_0,W_0,Z_0,q_0)\in\Omega_{\mathrm{tcN}}$,
where $\Omega_{\mathrm{tcN}}$ is defined in \eqref{tcNOm}.
Then there exists a unique global 
solution $(U,V,W,Z,q)(t,\xi)$ of
\eqref{ODE} with initial data
$(U_0,V_0,W_0,Z_0,q_0)(\xi)$ satisfying
\begin{equation}
	\label{glODE}
	(U,V,W,Z,q)\in C\left(
	[-T,T], \Omega_{\mathrm{tcN}}
	\right),\quad
	\mbox{for any $T>0$}.
\end{equation}

Moreover, we have the following 
conservation laws
(see Remark \ref{cluU} below):
\begin{subequations}\label{consqU}
\begin{equation}
\begin{split}
	&E_u(t)=\int_{-\infty}^{\infty}
	\left(U^2\cos^2\frac{W}{2}+
	\sin^2\frac{W}{2}\right)(t,\xi)
	\left(
	q\cos^2\frac{Z}{2}
	\right)(t,\xi)\,d\xi = E_{u_0},\\
	&E_{v}(t)=\int_{-\infty}^{\infty}
	\left(V^2\cos^2\frac{Z}{2}+
	\sin^2\frac{Z}{2}\right)(t,\xi)
	\left(
	q\cos^2\frac{W}{2}
	\right)(t,\xi)\,d\xi = E_{v_0},
\end{split}
\end{equation}
and
\begin{equation}
	G(t)=\int_{-\infty}^{\infty}\left(
	qUV\cos^2\frac{W}{2}\cos^2\frac{Z}{2}
	+\frac{q}{4}\sin W\sin Z
	\right)(t,\xi)\,d\xi=G_0,
\end{equation}
as well as
\begin{equation}
\begin{split}
	H(t)=\int_{-\infty}^{\infty}
	&\left(3U^2V^2\cos^2\frac{W}{2}\cos^2\frac{Z}{2}
	+U^2\cos^2\frac{W}{2}\sin^2\frac{Z}{2}
	+V^2\sin^2\frac{W}{2}\cos^2\frac{Z}{2}\right.\\
	&\left.\quad
	+UV\sin W\sin Z-\sin^2\frac{W}{2}\sin^2\frac{Z}{2}\right)
	(t,\xi)q(t,\xi)\,d\xi= H_0,
\end{split}
\end{equation}
\end{subequations}
for any $t\in[-T,T]$.
\end{theorem}

\begin{remark}[Compatibility condition \eqref{ccUV}]
	\label{Uin}
	Assuming that the initial data
		$(U_0,V_0,W_0,Z_0,q_0)$ satisfy the compatibility conditions
	\eqref{ccUV}, it is shown in \cite[Equations (5.1)--(5.2)]{HQ21}
	that these conditions are preserved
	by the ODE system \eqref{ODE}:~the corresponding unique global solution
		satisfies \eqref{ccUV} for every fixed $t\in\R$.
	Recalling that $(U_0,V_0,W_0,Z_0,q_0)$, defined in \eqref{id2} 
	in terms of the initial data of the Cauchy problem for
		the two-component Novikov system, satisfies the compatibility
	conditions \eqref{ccUV}, we are interested only in the solutions
		of the ODE system \eqref{ODE} that satisfy \eqref{ccUV},
	as considered in Theorem \ref{gwp}.
	Note that if the initial data do not satisfy \eqref{ccUV},
		the conservation laws \eqref{consqU} need not hold, and the
		global well-posedness argument for the ODE system \eqref{ODE}
		no longer applies.
\end{remark}

\begin{remark}[Subset $\Omega$]
	\label{rEB}
	Definition \ref{Om} slightly differs with 
	that given in \cite[Equation (4.3)]{KR25}
	(see also the subset $\Lambda$ in \cite[Section IV]{HQ21}).
	Namely, here we have an additional restriction 
	$(q(\cdot)-1)\in L^1(\R)$, which will be useful for 
	defining the norm of the tangent vector of a 
	regular path, see Section \ref{LM}, particularly \eqref{y^ve} below.
	Notice that this property is automatically 
	satisfied for the unique global solution of \eqref{ODE}--\eqref{id2} 
	obtained in \cite{HQ21,KR25}.
	Indeed, having in hand a
	global solution $(U,V,W,Z,q)$, the
	differential equation for $q$ in \eqref{ODE} 
	yields the following estimate 
	(we omit the argument $t$ for brevity;
	here $L^p=L^p(\R)$, $p=1,2,\infty$)
	\begin{equation*}
		\begin{split}
			\|q(\cdot)-1\|_{L^1}\leq&\|q_0(\cdot)-1\|_{L^1}
			+T\|q\|_{L^\infty}\|W\|_{L^2}
			\left\|U^2V+\frac{1}{2}V-P_1-\px P_2\right\|_{L^2}\\
			&+T\|q\|_{L^\infty}\|Z\|_{L^2}
			\left\|UV^2+\frac{1}{2}U-S_1-\px S_2\right\|_{L^2},
		\end{split}
	\end{equation*}
	where we have used the Cauchy-Schwarz 
	inequality and \eqref{sinx} for $W\in\R$.
	Taking into account that the norms 
	$\|U\|_{L^2\cap L^\infty}$, $\|V\|_{L^2\cap L^\infty}$,
	$\|W\|_{L^2}$, $\|Z\|_{L^2}$, $\|q\|_{L^\infty}$,
	$\|P_1\|_{L^2}$, $\|\px P_2\|_{L^2}$,
	$\|S_1\|_{L^2}$, and $\|\px S_2\|_{L^2}$ 
	are uniformly bounded, see \cite[Section V]{HQ21},
	we conclude that the norm $\|q(t,\cdot)-1\|_{L^1}$ is
	uniformly bounded as well, provided that the initial data
	$q_0$ satisfies $\|q_0(\cdot)-1\|_{L^1}<\infty$.
\end{remark}

\begin{remark}\label{cluU}
If $y(t,\cdot)$ is strictly monotone, then 
the conservation laws \eqref{consq}
are equivalent to their ODE counterparts \eqref{consqU} 
expressed in the variables $(U,V,W,Z,q)$.  
For this reason, and with a slight abuse of 
notation, we use the same symbols 
to refer to both versions.
\end{remark}

\begin{remark}\label{WZval}
	Observing that the right-hand side of \eqref{ODE} is invariant under
	the addition of multiples of $2\pi$ to either $W$ or $Z$, and using the uniform-in-$t$ bounds for $\|W(t,\cdot)\|_{L^\infty}$ and
	$\|Z(t,\cdot)\|_{L^\infty}$, we may regard
	the values of $W$ and $Z$ as lying on the circle $\R/2\pi\Z$. 
	Thus, without loss of generality, we assume 
	throughout the paper that
	\begin{equation}\label{perWZ}
		W(t,\xi), Z(t,\xi)\in(-\pi,\pi],\quad t,\xi\in\R.
	\end{equation}
\end{remark}

\begin{remark}[$L^1(\R)$ solutions of \eqref{ODE}--\eqref{id2}]
\label{L1ODE}
	Suppose that the initial data
	$(U_0,V_0,W_0,Z_0,q_0)\in\Omega_{\mathrm{tcN}}^1$;
	see \eqref{OmL1tcN}. Then the corresponding global solution \eqref{glODE}
	belongs to $\Omega_{\mathrm{tcN}}^1$ as well for all $t\in[-T,T]$.
		These additional decay assumptions on the solutions
	will be used in Section \ref{GdO} to rigorously justify
	the completion arguments, see Proposition \ref{COm}
	and Remark \ref{compL1}.
	
	We now derive the required $L^1(\R)$ a priori bounds from \eqref{ODE}.
	The first equation in \eqref{ODE}
	as well as the estimate (see \cite[Lemma 4.1]{KR25} and \eqref{Om})
	\[
	\mathcal{E}(t,\xi,\eta)\leq\Gamma(\xi-\eta),
	\quad\xi,\eta\in\R,\quad
	\Gamma(\xi)=\exp\left(\frac{q^-}{4}\left(R_2^2-|\xi|\right)\right),
	\]
	yield (recall \eqref{P-i-d};
	here $L^p=L^p(\R)$, $p=1,2,\infty$ and we omit $t$ for simplicity)
	\begin{equation}\label{UL1}
		\|U\|_{L^1}\leq \|U_0\|_{L^1}+T\|\Gamma\|_{L^1}
		\left(\|p_1\|_{L^1}+\|p_2\|_{L^1}\right).
	\end{equation}
	Using \eqref{p1q1}, \eqref{Sob}, \eqref{sinx}, 
	the Cauchy-Schwarz inequality and the 
	uniform bounds for $\|U\|_{L^\infty\cap L^2}$, $\|V\|_{L^2}$,
	$\|W\|_{L^2}$, and $\|q\|_{L^\infty}$ (see \cite[Section 5]{KR25}),
	we obtain that
	\begin{equation*}
	\begin{split}
		\|p_1\|_{L^1}\leq \|q\|_{L^\infty}\left(\|U\|_{L^\infty}\|UV\|_{L^1}
		+\frac{1}{4}\|UW\|_{L^1}+\frac{1}{4}\|VW\|_{L^1}\right)
		\leq C(E_{u_0},E_{v_0}, H_0),
	\end{split}
	\end{equation*}
	where $E_{u_0}$, $E_{v_0}$, and $H_0$ are given in \eqref{consqU}.
	The function $p_2$, see \eqref{p2q2}, can be estimated similarly,
	and thus we have the needed estimate for
	$\|U\|_{L^1}$ from \eqref{UL1}.
	The norms $\|V\|_{L^1}$, $\|W\|_{L^1}$, and $\|Z\|_{L^1}$
	can be treated similarly.
	Finally, the a priori bounds for $\|U_\xi\|_{L^1}$ and
	$\|V_\xi\|_{L^1}$ can be derived from the compatibility conditions \eqref{ccUV}.
\end{remark}

\begin{figure}
\centering{\includegraphics[scale=0.5]{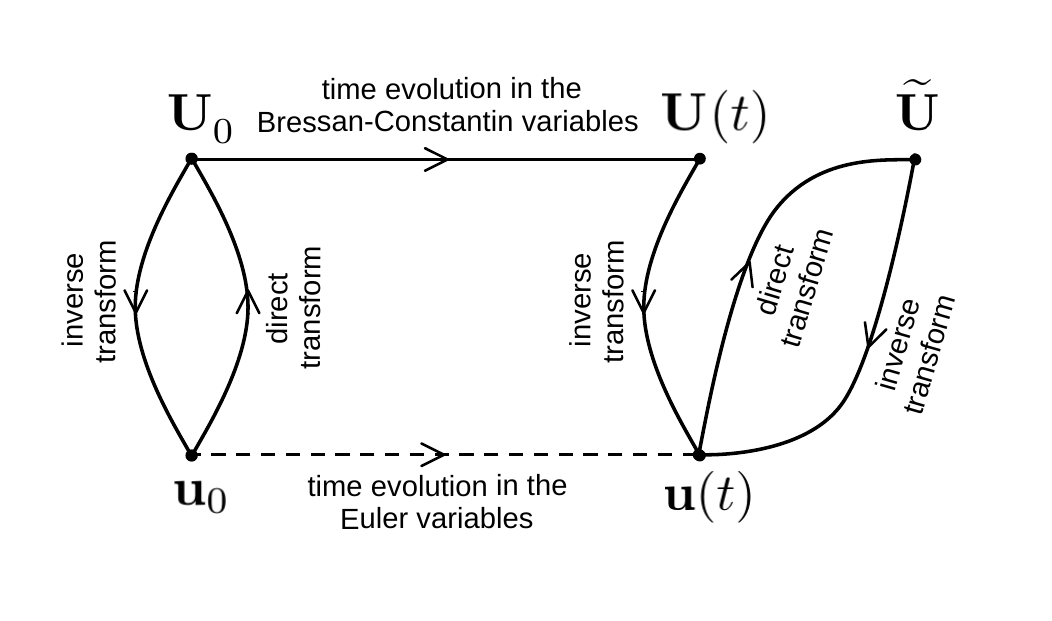}}
\caption{A diagram of the Bressan-Constantin 
approach applied to the two-component Novikov equation. 
Starting from the initial data 
$\mathbf{u}_0=\left(u_0,v_0,\mu_0;D_{W,0},D_{Z,0}\right)$, 
the direct transform \eqref{id2} yields 
$\mathbf{U}_0=(U_0,V_0,W_0,Z_0,q_0)$, 
which serves as the initial data 
for the associated ODE system \eqref{ODE}. 
The initial Eulerian data $\mathbf{u}_0$ 
can be recovered from $\mathbf{U}_0$ via the inverse 
transform \eqref{uvdef}--\eqref{mut} 
and \eqref{DN-a} evaluated at $t=0$.
Given the unique global solution 
$U(t)=(U,V,W,Z,q)(t)$ of the ODE system, 
Theorem~\ref{Thm} provides a corresponding 
global conservative solution 
$\mathbf{u}(t)=(u(t),v(t),\mu_{(t)};D_W(t),D_Z(t))$ 
of the Novikov system. We emphasize that 
applying the direct transform \eqref{id2} 
to $\mathbf{u}(t)$ produces a solution 
$\widetilde{\mathbf{U}}(t)$ 
which, in general, does \emph{not} coincide with $\mathbf{U}(t)$ 
(in particular, one has $\tilde{q}\equiv 1$ in $\widetilde{\mathbf{U}}$).  
Consequently, the Bressan-Constantin and Eulerian variables 
do not form a bijective correspondence 
(cf.\,\cite[Section~7]{KR25}, \cite[Proof of Theorem~3]{BC07}, 
and \cite[Section~8]{BC07d}).}
\label{inv-d-t}
\end{figure}

Let us show how to define the 
global solution $(u(t),v(t),\mu_{(t)};D_W(t),D_Z(t))$ 
given in Theorem \ref{Thm} of the two-component 
Novikov system in terms of the global solution 
$(U, V, W, Z, q)(t)$ of the associated ODE system 
(see Figure \ref{inv-d-t} for the illustration). 

The representation of the characteristic $y(t,\xi)$ 
in terms of $(U,V)$ reads, as per \eqref{char},
\begin{equation}\label{char1}
	y(t,\xi)=y_0(\xi)
	+\int_0^t\left(UV\right)(\tau,\xi)\,d\tau,
	\quad t,\xi\in\mathbb{R},
\end{equation}
where $y_0(\xi)$ is given by \eqref{y-def} 
with $\mu_0$ instead of $\mu$.
Since $(U,V,W,Z,q)(t,\cdot)$ satisfies the
compatibility conditions \eqref{ccUV} for every fixed
$t\in\R$, the characteristic $y(t,\xi)$ satisfies the following
differential equation (see \cite[Equation (3.2)]{HQ21}):
\begin{equation}\label{pxiy}
	y_\xi(t,\xi)=\left(q\cos^2\frac{W}{2}\cos^2\frac{Z}{2}\right)(t,\xi),
	\quad t,\xi\in\R.
\end{equation}
Therefore, recalling that $y_0(0)=0$, we can define
the initial characteristic $y_0(\xi)$ only in terms 
of the initial data of the Cauchy problem of the ODE system:
\begin{equation}\label{y0ode}
	y_0(\xi)=\int_0^\xi
	\left(q\cos^2\frac{W}{2}\cos^2\frac{Z}{2}\right)(0,\eta)
	\,d\eta.
\end{equation}

Then we define the flow $\left(u(t),v(t),\mu_{(t)};
D_W(t),D_Z(t)\right)$ in Theorem \ref{Thm} as follows:
\begin{align}\label{uvdef}
	u(t,x)=U(t,\xi),\quad v(t,x)=V(t,\xi),\quad
	\mbox{if }x=y(t,\xi),\quad t\in\R,
\end{align}
and (see \cite[Equation (7.3)]{KR25})
\begin{equation}\label{mut}
	\mu_{(t)}([a,b])
	=\int\limits_{\{\xi:y(t,\xi)\in[a,b]\}}
	\left(
	\cos^2\frac{W}{2}\sin^2\frac{Z}{2}
	+\sin^2\frac{W}{2}\cos^2\frac{Z}{2}
	+\sin^2\frac{W}{2}\sin^2\frac{Z}{2}
	\right)(t,\xi)q(t,\xi)\,d\xi,
\end{equation}
while $D_W(t)$ and $D_Z(t)$ are given in \eqref{DN-a}. 
Also notice that $(u_x,v_x)$ can be found by
(see \cite[Equation (6.4)]{KR25})
\begin{equation}\label{pxuv}
	\begin{split}
	&u_x(t,x)=\tan\frac{W(t,\xi)}{2},\quad x=y(t,\xi),
	\quad\xi\in\R\setminus D_W(t),\\
	&v_x(t,x)=\tan\frac{Z(t,\xi)}{2},\quad\,\,\, x=y(t,\xi),
	\quad\xi\in\R\setminus D_Z(t).
	\end{split}
\end{equation}
\begin{remark}[Inverse transform]
		Given initial data in $\mathcal{D}$ for the Cauchy problem
		\eqref{t-c-N-n},
	equations \eqref{char1}, \eqref{uvdef}, \eqref{mut}, and \eqref{DN-a}
	construct the corresponding global conservative weak solution
		described in Theorem \ref{Thm}.
	On the other hand, if we start with a global solution
	of the ODE system in $\Omega_{\mathrm{tcN}}$,
	see Theorem \ref{gwp}, then, thanks to
	compatibility conditions \eqref{ccUV},
	these equations also
	define a global conservative weak solution
		of the Cauchy problem for the two-component Novikov equation.
		Here the initial data in $\mathcal{D}$ are defined through the
		characteristic $y_0(\xi)$ in \eqref{y0ode},
	in terms of the initial data of the problem for the ODE system.
		Notice that the conditions \eqref{ccUV} are required for
	global well-posedness of the ODE Cauchy problem
	(see Remark \ref{Uin}), for the conservation laws \eqref{consqU}
	(and, consequently, \eqref{consq}),
	and for the identity \eqref{pxiy}, which is essential
	for proving that $(u,v)$ is a weak solution of \eqref{t-c-N-n}--\eqref{iid}
	(see \cite[Section VI, Step 7]{HQ21}).
\end{remark}

We will also require the global well-posedness 
of regular solutions $(U, V, W, Z, q)$ to \eqref{ODE}.
\begin{theorem}\label{Ckgl}
Consider initial data
$(U_0,V_0,W_0,Z_0,q_0)\in
\mathcal{C}^k_{\Omega_{\mathrm{tcN}}}$, $k\in\N$,
with $\mathcal{C}^k_{\Omega_{\mathrm{tcN}}}$ given in \eqref{CktcN}.
Then there exists a unique global 
solution $(U,V,W,Z,q)(t,\xi)$ of
\eqref{ODE} with initial data
$(U_0,V_0,W_0,Z_0,q_0)(\xi)$ such that  
\begin{equation*}
	(U,V,W,Z,q)\in C\left(
	[-T,T],\mathcal{C}^k_{\Omega_{\mathrm{tcN}}}
	\right),\quad
	\mbox{for any $T>0$}.
\end{equation*}
\end{theorem}

\begin{proof}[Sketch of the proof]	
The argument follows the structure of 
the proof of Theorem~\ref{gwp} in \cite{HQ21}.  
For completeness, we outline the main steps 
and refer the reader to \cite{HQ21} for full details.

Using the same estimates as in \cite[Lemma~4.1]{HQ21}  
(see also \cite[Appendix~A]{KR25}), one 
shows that the right-hand side of \eqref{ODE} 
is Lipschitz continuous in 
$\mathcal{C}^k_{\Omega_{\mathrm{tcN}}}$
for any $k\in\N$. The arguments in \cite[Section~V]{HQ21} 
then yield a global solution of  
\eqref{ODE}--\eqref{id2} in 
$\Omega_{\mathrm{tcN}}$.

To verify that this global solution actually belongs to  
$(C^{k})^{2}\times (C^{k-1})^{3}$, we use a priori bounds.  
Uniform bounds on $|U|$, $|V|$, and $q$ 
provide uniform control of the nonlocal terms  
$|\pxi P_j|$, $|\pxi S_j|$,  
$|\pxi(\px P_j)|$, and $|\pxi(\px S_j)|$  
(see \cite[Section~V]{HQ21}).  
From the first two equations in \eqref{ODE}, we 
then obtain uniform estimates for  
$|U_\xi|$ and $|V_\xi|$.

For $k=2$, applying Gronwall's inequality to the 
last three equations of \eqref{ODE}  
yields uniform bounds for $|W_\xi|$, $|Z_\xi|$, and $|q_\xi|$.  
Substituting these into the first two equations in \eqref{ODE}  
gives uniform estimates for $|U_{\xi\xi}|$ and $|V_{\xi\xi}|$.  
Iterating this procedure produces the a priori bounds
\[
\left|\pxi^{i} U\right|,\quad\left|\pxi^{i} V\right|,\quad
\left|\pxi^{\,i-1} W\right|,\quad
\left|\pxi^{\,i-1} Z\right|,\quad
\left|\pxi^{\,i-1} q\right|,
\qquad i=1,\dots,k,
\]
as required.
\end{proof}

\subsection{Generic solutions of the ODE system}
\label{gODE}
We begin with Lemma \ref{L2}, which provides the perturbations
used in Lemma \ref{L3} to show that, generically within
$\Omega_{\mathrm{tcN}}$, the relevant triples formed from a solution
of \eqref{ODE} and its derivatives never equal $(\pi,0,0)$.
This allows us to eliminate degenerate elements 
from the solution set $(u,v)$, thereby 
facilitating the proof of Theorem \ref{Thmgr}, 
see Section \ref{prTh1} below.

In Lemma \ref{L3} we show that, generically, 
the value $(\pi,0,0)$ is never attained by 
the three-dimensional maps 
$(W,W_\xi,W_{\xi\xi})$, 
$(W,W_t,W_{\xi})$, $(Z,Z_\xi,Z_{\xi\xi})$, 
and $(Z,Z_t,Z_{\xi})$. To establish this, we consider 
three-parameter perturbations of, for instance, 
$(W,W_\xi,W_{\xi\xi})$ at points $(t,\xi)$ where
$(W,W_\xi,W_{\xi\xi})(t,\xi)=(\pi,0,0)$. 
Ensuring that $(\pi,0,0)$ becomes a regular 
value of the perturbed map requires that 
its Jacobian have full rank $3$ (recall 
Definition \ref{RV} and see \eqref{r1}, \eqref{r2}).
Following the approach in \cite{BC17} (see 
also \cite{CCCS18, LZ17}), we construct 
such perturbations in the next lemma,
emphasizing that they satisfy the compatibility conditions \eqref{ccUV}.
The proof given here differs substantially from that in
\cite{CCCS18}; see Remark \ref{psol} for more details.

\begin{lemma}\label{L2}
Consider the global solution 
$(U,V,W,Z,q)\in C\left(
[-T,T], \mathcal{C}^k_{\Omega_{\mathrm{tcN}}}
\right)$,
$k\geq3$, of the Cauchy problem for \eqref{ODE}
given in Theorem \ref{Ckgl}. Then, for any 
$(t_0, \xi_0) \in \R^2$, we have
\begin{enumerate}[1)]
	\item there exists a three-parameter family of
	global solutions
	\begin{equation}\label{3pars}
		\left(\tilde{U},\tilde{V},\tilde{W},
		\tilde{Z},\tilde{q}\right)
		(\cdot,\cdot;\nu)\in
		C\left(
		[-T,T], \mathcal{C}^k_{\Omega_{\mathrm{tcN}}}
		\right),\quad \nu=(\nu_1,\nu_2,\nu_3),
	\end{equation}
	of \eqref{ODE} for any $T>0$,
	with $|\nu|<\ve$, $\ve>0$, such that:
	\begin{enumerate}[a)]
		\item at $\nu=0$ we have 
		$\left(\tilde{U},\tilde{V},\tilde{W},
		\tilde{Z},\tilde{q}\right)
		(t,\xi;0)=(U,V,W,Z,q)(t,\xi)$ for all $t,\xi$;
		\item at $(t,\xi;\nu)=(t_0,\xi_0;0)$ we have
		either $\left(\tilde{W},
		\tilde{W}_\xi,\tilde{W}_{\xi\xi}\right)(t_0,\xi_0;0)
		\neq(\pi,0,0)$ or
		\begin{equation}\label{r1}
			\mathrm{rank}\left(\left.
			D_\nu\left(
			\tilde{W},\tilde{W}_\xi,\tilde{W}_{\xi\xi}\right)
			\right|_{(t,\xi;\nu)=(t_0,\xi_0;0)}
			\right)=3;
		\end{equation}
	\end{enumerate}
	\item there exists a three-parameter family \eqref{3pars} 
	of global solutions of \eqref{ODE} such that 1a) holds and
	either $\left(\tilde{W},
	\tilde{W}_t,\tilde{W}_{\xi}\right)(t_0,\xi_0;0)
	\neq(\pi,0,0)$ or
	\begin{equation}\label{r2}
		\mathrm{rank}\left(\left.
		D_\nu\left(
		\tilde{W},\tilde{W}_t,\tilde{W}_{\xi}\right)
		\right|_{(t,\xi;\nu)=(t_0,\xi_0;0)}
		\right)=3;
	\end{equation}
	\item claims 1)--2) hold with the 
	function $Z$ in place of $W$ 
	in \eqref{r1} and \eqref{r2}.
\end{enumerate}
\end{lemma}

\begin{proof}
\textbf{Proof of item 1.}
If $(W,W_\xi,W_{\xi\xi})(t_0,\xi_0)\neq (\pi,0,0)$,
then we can take a trivial perturbation 
\[
\left(\tilde{U},\tilde{V},\tilde{W},
\tilde{Z},\tilde{q}\right)
(t,\xi;\nu)=(U,V,W,Z,q)(t,\xi),\quad\text{for any }\,t,\xi,\nu,
\]
and item 1 follows.
It remains to consider the nontrivial case
$(W,W_\xi,W_{\xi\xi})(t_0,\xi_0)=(\pi,0,0)$.
\medskip

	\textbf{Step 1.}
	Define the perturbations of $W(t_0,\xi)$ and $Z(t_0,\xi)$
	as follows:
	\begin{equation}\label{WZ-pert}
	\begin{split}
		&\tilde{W}(t_0,\xi;\nu,a)=W(t_0,\xi)
		+\sum\limits_{j=1}^3\nu_jf_j(\xi)+a\psi(\xi),\\
		&\tilde{Z}(t_0,\xi;b)=Z(t_0,\xi)+b\chi(\xi).
	\end{split}
	\end{equation}
	Here the smooth functions $f_j$, $j=1,2,3$, $\psi$, and $\chi$ satisfy the following properties:
	\begin{equation}\label{fj-psi}
	\begin{split}
		&f_j(\xi)=0,\quad\text{for }\,|\xi-\xi_0|\geq\delta_0,
		\quad\text{with some }\,\delta_0>0,\quad j=1,2,3,\\
		&\psi^2(\xi)+\chi^2(\xi)=0,\quad\text{for }\,
		|\xi-\xi_1|\geq\delta_1,
		\quad\text{with some }\,\delta_1>0\text{ and }
		\xi_1\gg\xi_0.
	\end{split}
	\end{equation}
The additional scalar parameters $a$ and $b$
in \eqref{WZ-pert}
will be determined as smooth functions of $\nu$ from the system
of equations \eqref{Hnuab}.
These equations will ensure that
the perturbations $\tilde{U}$ and $\tilde{V}$, see \eqref{pUV} below,
are compact perturbations of, respectively, $U$ and $V$.

\medskip

\textbf{Step 2.}
We leave $q$ unperturbed and set
(here and below, $\nu$, $a$, and $b$ lie in a neighborhood of zero)
\begin{equation}\label{q-pert}
	\tilde{q}(t_0,\xi;\nu)=q(t_0,\xi),\quad\text{for any }\,\nu,
\end{equation}
The compatibility conditions for $\tilde{U}$ and $\tilde{V}$
therefore read as follows; see \eqref{ccUV}. Here and below,
we suppress the argument $t_0$:
\begin{equation}\label{pccUV}
\begin{split}
	&\tilde{U}_\xi(\xi;\nu,a,b)=\frac{q(\xi)}{2}
	\sin\tilde{W}(\xi;\nu,a)
	\cos^2\frac{\tilde{Z}(\xi;b)}{2},\\
	&\tilde{V}_\xi(\xi;\nu,a,b)=\frac{q(\xi)}{2}
	\cos^2\frac{\tilde{W}(\xi;\nu,a)}{2}
	\sin\tilde{Z}(\xi;b).
\end{split}
\end{equation}
Using \eqref{ccUV}, we obtain the following representations for
$\tilde{U}$ and $\tilde{V}$ from \eqref{pccUV}:
\begin{equation}\label{pUV}
\begin{split}
	&\tilde{U}(\xi;\nu,a,b)=U(\xi)
	+\tilde{H}_1(\xi;\nu,a,b),\\
	&\tilde{V}(\xi;\nu,a,b)=V(\xi)
	+\tilde{H}_2(\xi;\nu,a,b),
\end{split}
\end{equation}
where $U(\xi)=U(t_0,\xi)$, $V(\xi)=V(t_0,\xi)$, and
\begin{equation}\label{Hj}
\begin{split}
	&\tilde{H}_1(\xi;\nu,a,b)=
	\frac{1}{2}\int_{-\infty}^\xi q(\xi^\prime)
	\left(\sin\tilde{W}(\xi^\prime;\nu,a)
	\cos^2\frac{\tilde{Z}(\xi^\prime;b)}{2}
	-\sin W(\xi^\prime)\cos^2\frac{Z(\xi^\prime)}{2}
	\right)d\xi^\prime,\\
	&\tilde{H}_2(\xi;\nu,a,b)=
	\frac{1}{2}\int_{-\infty}^\xi q(\xi^\prime)
	\left(\cos^2\frac{\tilde{W}(\xi^\prime;\nu,a)}{2}
	\sin\tilde{Z}(\xi^\prime;b)
	-\cos^2\frac{W(\xi^\prime)}{2}\sin Z(\xi^\prime)
	\right)d\xi^\prime.
\end{split}
\end{equation}

By \eqref{fj-psi} and \eqref{pUV},
the perturbations $\tilde{U}$ and $\tilde{V}$
are equal to the initial functions $U$ and $V$,
respectively, for $\xi\leq\xi_0-\delta_0$:
\begin{equation}\label{comppUV}
	\tilde{U}(\xi;\nu,a,b)=U(\xi),\quad\text{and}\quad
	\tilde{V}(\xi;\nu,a,b)=V(\xi),\quad\text{for any }\,
	\xi\leq\xi_0-\delta_0,\text{ and any }\,\nu,\, a,\, b.
\end{equation}
Taking into account that (see \eqref{WZ-pert} and \eqref{fj-psi})
\begin{equation*}
	\tilde{W}(\xi;\nu,a)=W(\xi),\quad\text{and}\quad
	\tilde{Z}(\xi;b)=Z(\xi),\quad\text{for any }\,
	\xi\geq\xi_1+\delta_1,\text{ and any }\,\nu,\, a,\, b,
\end{equation*}
\eqref{Hj} gives
\begin{equation}\label{tHjc}
	\tilde{H}_j(\xi;\nu,a,b)=\tilde{H}_j(\xi_1+\delta_1;\nu,a,b),
	\quad\text{for any }\,
	\xi\geq\xi_1+\delta_1,\text{ and any }\,\nu,\, a,\, b,
	\quad j=1,2.
\end{equation}
Therefore, combining \eqref{pUV} and \eqref{tHjc},
we conclude that the following condition is sufficient:
\begin{equation}\label{Hnuab}
	\tilde{H}(\nu,a,b):=\begin{pmatrix}
		\tilde{H}_1(\infty;\nu,a,b)\\
		\tilde{H}_2(\infty;\nu,a,b)
	\end{pmatrix}
	=\begin{pmatrix}
		0\\0
	\end{pmatrix}.
\end{equation}
If \eqref{Hnuab} holds, then $\tilde{U}$ and $\tilde{V}$
are equal to $U$ and $V$,
respectively, for all sufficiently large $\xi$
(cf.~\eqref{comppUV}):
\begin{equation*}
	\tilde{U}(\xi;\nu,a,b)=U(\xi),\quad\text{and}\quad
	\tilde{V}(\xi;\nu,a,b)=V(\xi),\quad\text{for any }\,
	\xi\geq\xi_1+\delta_1,\text{ and any }\,\nu,\, a,\, b.
\end{equation*}
In other words, choosing $a$ and $b$ so that \eqref{Hnuab}
holds makes $\tilde{U}$ and $\tilde{V}$ compact perturbations
of $U$ and $V$, respectively, while preserving the compatibility
conditions \eqref{pccUV}.

\medskip

\textbf{Step 3.}
It remains to determine $a$ and $b$ as smooth functions of $\nu$
from the system of equations \eqref{Hnuab}.
To this end, it suffices to show that the determinant of the Jacobian
of $\tilde{H}$ with respect to $(a,b)$ is nonzero
(see \eqref{Hnuab}):
\begin{equation}\label{JacH}
\det\left.D_{(a,b)}\tilde{H}\right|_{(\nu,a,b)=(0,0,0)}
=\det\begin{pmatrix}
	\partial_a \tilde{H}_1&\partial_b \tilde{H}_1\\
	\partial_a \tilde{H}_2&\partial_b \tilde{H}_2
\end{pmatrix}(\infty;0,0,0)\neq0.
\end{equation}
Direct calculations show that (see \eqref{Hj} and \eqref{WZ-pert})
\begin{equation}\label{paH1}
	\partial_a\tilde{H}_1(\infty;0,0,0)=
	\frac{1}{2}\int_{-\infty}^\infty
	\left(q\cos W\cos^2\frac{Z}{2}\right)(\xi)\psi(\xi)\,d\xi
	=\frac{1}{2}\int_{\xi_1-\delta_1}^{\xi_1+\delta_1}
	\left(q\cos W\cos^2\frac{Z}{2}\right)(\xi)\psi(\xi)\,d\xi,
\end{equation}
where in the last equality we have used that $\psi(\xi)=0$
for all $|\xi-\xi_1|\geq\delta_1$, see \eqref{fj-psi}.
Employing that
\[
W,Z\in C^1(\R)\cap L^2(\R),\quad\text{and}\quad(q-1)\in C^1(\R)\cap L^1(\R),
\]
we have the following limits:
\begin{equation}\label{WZ-lim}
	W^2(\xi)+Z^2(\xi)+(q(\xi)-1)^2\to 0,\quad\xi\to\pm\infty.
\end{equation}
Therefore, \eqref{paH1} yields the following 
representation for $\partial_a\tilde{H}_1(\infty;0,0,0)$:
\begin{equation}\label{aspaH1}
	\partial_a\tilde{H}_1(\infty;0,0,0)=\frac{1}{2}
	\int_{\xi_1-\delta_1}^{\xi_1+\delta_1}\psi(\xi)\,d\xi+\so(1),
\end{equation}
where $\so(1)$ is a function of $\xi_1$ vanishing as $\xi_1\to\infty$.
Normalizing $\psi(\xi)$ as follows:
\begin{equation}\label{npsi}
\int_{-\infty}^\infty\psi(\xi)\,d\xi
=\int_{\xi_1-\delta_1}^{\xi_1+\delta_1}\psi(\xi)\,d\xi=2,
\end{equation}
we obtain from \eqref{aspaH1} that
\begin{equation}\label{lpaH1}
	\partial_a\tilde{H}_1(\infty;0,0,0)=1+\so(1),
	\quad\xi_1\to\infty.
\end{equation}

Arguing similarly, we obtain the following representations for the
remaining entries of the Jacobian in \eqref{JacH}:
\begin{equation}\label{lpaHj}
\begin{split}
	&\partial_b\tilde{H}_1(\infty;0,0,0)=\so(1),
	\quad
	\partial_a\tilde{H}_2(\infty;0,0,0)=\so(1),
	\quad\xi_1\to\infty\\
	&\partial_b\tilde{H}_2(\infty;0,0,0)=
	\frac{1}{2}
	\int_{\xi_1-\delta_1}^{\xi_1+\delta_1}\chi(\xi)\,d\xi+\so(1),
	\quad\xi_1\to\infty.
\end{split}
\end{equation}
Imposing the following normalization condition for $\chi(\xi)$ (cf.~\eqref{npsi}):
\begin{equation*}
	\int_{-\infty}^\infty\chi(\xi)\,d\xi
	=\int_{\xi_1-\delta_1}^{\xi_1+\delta_1}\chi(\xi)\,d\xi=2,
\end{equation*}
we obtain from \eqref{lpaH1} and \eqref{lpaHj} that
\begin{equation*}
\det\begin{pmatrix}
	\partial_a\tilde{H}_1&\partial_b\tilde{H}_1\\
	\partial_a\tilde{H}_2&\partial_b\tilde{H}_2
\end{pmatrix}(\infty;0,0,0)
=\det\begin{pmatrix}
	1&0\\
	0&1
\end{pmatrix}+\so(1)=1+\so(1),\quad\xi_1\to\infty,
\end{equation*}
which yields \eqref{JacH} for $\xi_1$ large enough.
Since $\tilde{H}(\nu,a,b)$ is smooth in $(\nu,a,b)$ and
$\tilde{H}(0,0,0)=(0,0)^T$, the implicit function
theorem yields smooth functions $a(\nu)$ and $b(\nu)$, defined for
$\nu$ in a neighborhood of zero, such that
\eqref{Hnuab} is fulfilled:
\[
\tilde{H}(\nu,a(\nu),b(\nu))=0,\quad|\nu|<\delta,
\quad\text{for some sufficiently small }\,\delta>0.
\]
Moreover, the smooth functions 
$a(\nu)$, $b(\nu)$ vanish as $\nu\to0$ and satisfy the following estimate:
\begin{equation*}
	\left|a(\nu)\right|+\left|b(\nu)\right|\leq C|\nu|,
	\quad |\nu|<\delta,\quad C>0.
\end{equation*}

\medskip

\textbf{Step 4.}
Applying Theorem \ref{Ckgl} to the Cauchy problem for the ODE
system \eqref{ODE} with initial data
$\left(\tilde{U},\tilde{V},\tilde{W},
\tilde{Z},\tilde{q}\right)(\xi;\nu)$ 
taken at $t=t_0$, we obtain the family of global solutions
\eqref{3pars} satisfying item 1.a.

We are now in a position to verify the rank condition \eqref{r1}.
Using \eqref{WZ-pert} and recalling that $\psi(\xi_0)=0$ (see \eqref{fj-psi}),
direct calculations show that
\begin{equation}\label{Dr1}
	\left.D_\nu\left(
	\tilde{W},\tilde{W}_\xi,\tilde{W}_{\xi\xi}\right)
	\right|_{(t,\xi,\nu)=(t_0,\xi_0;0)}
	=
	\renewcommand{\arraystretch}{1.2}{
	\begin{pmatrix}
		f_1& \pxi f_1&\pxi^2 f_1\\
		f_2& \pxi f_2&\pxi^2 f_2\\
		f_3& \pxi f_3&\pxi^2 f_3\\
	\end{pmatrix}}(\xi_0).
\end{equation}
Taking $f_j(\xi)$, $j=1,2,3$, such that the matrix
on the right-hand side of \eqref{Dr1} has nonzero determinant,
we arrive at \eqref{r1}, and item 1 is proved.

\medskip

\textbf{Proof of item 2.}
The proof of item 2 follows the same reasoning as that of item 1.
The only difference is that the Jacobian in \eqref{r2} has the
following form (cf.~\eqref{Dr1}):
\begin{equation}\label{Dr2}
	\left.D_\nu\left(
	\tilde{W},\tilde{W}_t,\tilde{W}_{\xi}\right)
	\right|_{(t,\xi,\nu)=(t_0,\xi_0;0)}
	=
	\renewcommand{\arraystretch}{1.2}{
		\begin{pmatrix}
			f_1(\xi_0)& \partial_{\nu_1}\tilde{W}_t(t_0,\xi_0;0)
			&\pxi f_1(\xi_0)\\
			f_2(\xi_0)& \partial_{\nu_2}\tilde{W}_t(t_0,\xi_0;0)
			&\pxi f_2(\xi_0)\\
			f_3(\xi_0)& \partial_{\nu_3}\tilde{W}_t(t_0,\xi_0;0)
			&\pxi f_3(\xi_0)\\
	\end{pmatrix}}.
\end{equation}
Taking $f_j(\xi)$ such that
\begin{equation*}
	\begin{pmatrix}
		f_1& \pxi f_1\\
		f_2& \pxi f_2\\
		f_3& \pxi f_3
	\end{pmatrix}(\xi_0)=
	\begin{pmatrix}
		1& 0\\
		0& 0\\
		0& 1
	\end{pmatrix},
\end{equation*}
we conclude that the determinant of the matrix on the right-hand side of
\eqref{Dr2} is nonzero if and only if
\begin{equation}\label{Pnu2Wt}
\partial_{\nu_2}\tilde{W}_t(t_0,\xi_0;0)\neq0.
\end{equation}

Using the following 
differential equation for $\tilde{W}_t$ (see \eqref{ODE}):
\begin{equation*}
	\tilde{W}_t(t_0,\xi_0;\nu)=
	\left(2\tilde{U}^2\tilde{V}\cos^2\frac{\tilde{W}}{2}
	-\tilde{V}\sin^2\frac{\tilde{W}}{2}
	-2\left(\tilde{P}_1
	+\px \tilde{P}_2\right)\cos^2\frac{\tilde{W}}{2}\right)(t_0,\xi_0;\nu),
\end{equation*}
where the nonlocal terms $\tilde{P}_1$ and $\px \tilde{P}_2$
are defined as in \eqref{P-i-d}, \eqref{P2-i-d}
with $\left(\tilde{U},\tilde{V},\tilde{W},
\tilde{Z},\tilde{q}\right)$ instead of $(U,V,W,Z,q)$
and using that $\tilde{W}(t_0,\xi_0;0)=W(t_0,\xi_0)=\pi$,
we obtain that (see \eqref{WZ-pert} and \eqref{pUV})
\begin{equation}\label{pnu2pW}
	\partial_{\nu_2}\tilde{W}_t(t_0,\xi_0;0)
	=-\partial_{\nu_2}\tilde{V}(t_0,\xi_0;0)
	=\frac{1}{4}\int_{\xi_0-\delta_0}^{\xi_0}
	\left(q\sin W\sin Z\right)(t_0,\xi)f_2(\xi)\,d\xi.
\end{equation}

If there exists $\delta_0>0$ for which
$(\sin W\sin Z)(t_0,\xi)\neq0$
for some $\xi\in[\xi_0-\delta_0,\xi_0]$, then it suffices to choose
$f_2(\xi)$ not orthogonal
to the function $\left(q\sin W\sin Z\right)(t_0,\cdot)$
in $L^2(\xi_0-\delta_0,\xi_0)$.

If, for every $\delta_0>0$,
the function
$(\sin W\sin Z)(t_0,\cdot)\equiv0$
on the interval $\xi\in[\xi_0-\delta_0,\xi_0]$,
we must slightly modify the perturbation of $Z$, cf.~\eqref{WZ-pert}:
\begin{equation*}
	\tilde{Z}(t_0,\xi;\nu_2,b)=Z(t_0,\xi)+\nu_2g_2(\xi)+b\chi(\xi),
\end{equation*}
where the smooth function $g_2(\xi)$ satisfies the following condition:
\begin{equation*}
	g_2(\xi)=0,
	\quad\text{for }\,
	|\xi-\xi_2|\geq\delta_2,
	\quad\text{with some }\,\delta_2>0\text{ and }
	\xi_2\ll\xi_0.
\end{equation*}
Then we have the following expression for
$\partial_{\nu_2}\tilde{W}_t$, cf.~\eqref{pnu2pW}
\begin{equation*}
\begin{split}
	\partial_{\nu_2}\tilde{W}_t(t_0,\xi_0;0)
	&=-\partial_{\nu_2}\tilde{V}(t_0,\xi_0;0)
	=\frac{1}{2}\int_{\xi_2-\delta_2}^{\xi_2+\delta_2}
	\left(q\cos^2\frac{W}{2}\cos Z\right)
	(t_0,\xi)g_2(\xi)\,d\xi\\
	&=\frac{1}{2}\int_{\xi_2-\delta_2}^{\xi_2+\delta_2}
	g_2(\xi)\,d\xi+\so(1),\quad\xi_2\to-\infty.
\end{split}
\end{equation*}
Taking $\xi_2$ sufficiently negative and normalizing $g_2(\xi)$ as follows:
\[
\int_{\xi_2-\delta_2}^{\xi_2+\delta_2}
g_2(\xi)\,d\xi=1,
\]
we obtain \eqref{Pnu2Wt} and therefore item 2 is established.

\medskip

\textbf{Proof of item 3.} Item 3 follows by the same arguments as
items 1 and 2.

\end{proof}

\begin{remark}\label{psol}
	The construction in \cite[Section 4.2]{CCCS18},
	where the authors consider the Novikov equation
	with $U=V$ and $W=Z$,
	suggests defining
	a three-parameter perturbation 
	$\left(\tilde{U},\tilde{W},\tilde{q}\right)(\cdot,\cdot;\nu)$
	of the solution $(U,W,q)(\cdot,\cdot)$
	by perturbing its initial data at $t=0$ as follows
	(see \cite[Equation (4.12)]{CCCS18}):
	\begin{equation}
	\begin{split}
		&\tilde{U}_0(\xi;\nu)=
		U_0(\xi)+\sum\limits_{j=1}^3\nu_j U_{0}^{(j)}(\xi),
		\quad
		\tilde{W}_0(\xi;\nu)=
		W_0(\xi)+\sum\limits_{j=1}^3\nu_j W_{0}^{(j)}(\xi),\\
		&\tilde{q}_0(\xi;\nu)=q_0(\xi)
		+\sum\limits_{j=1}^3\nu_j q_{0}^{(j)}(\xi),
	\end{split}
	\end{equation}
	with smooth and compactly supported functions
	$U_{0}^{(j)}$, $W_{0}^{(j)}$, and $q_{0}^{(j)}$, $j=1,2,3$.
	To satisfy the compatibility condition \eqref{ccUV},
	one must have (cf.~\eqref{pccUV})
	\begin{equation}\label{ttpxiU}
		\pxi\tilde{U}_0(\xi;\nu)=
		\left(\frac{\tilde{q}_0}{2}\sin \tilde{W}_0\cos^2\frac{\tilde{W}_0}{2}\right)(\xi;\nu),
		\quad\xi\in\R,
	\end{equation}
	for all $\nu$ in a neighborhood of zero.
	Observe that the Taylor expansion in $\nu$ of the
	right-hand side of \eqref{ttpxiU} involves an
	infinite number of terms of the form 
	$\nu_1^{n_1}\nu_2^{n_2}\nu_3^{n_3}\tilde{F}_{n_1,n_2,n_3}(\xi)$
	with some $n_j\in\N\cup\{0\}$, $j=1,2,3$, and
	$\tilde{F}_{n_1,n_2,n_3}(\xi)$ defined through
	$W_0^{(j)}(\xi)$ and $q_0^{(j)}(\xi)$, $j=1,2,3$.
	Thus, one must consider nonzero
	functions $W_0^{(j)}(\xi)$ and $q_0^{(j)}(\xi)$, $j=1,2,3$
	such that they satisfy the following infinite family of constraints:
	$$
	\tilde{F}_{n_1,n_2,n_3}\equiv0,\quad n_1+n_2+n_3\geq2.
	$$
	Proving the existence of such functions appears difficult.
	
	In this paper, we construct the perturbed solutions
	in a different way.
	Starting with perturbations of $W$ and $Z$ as in \eqref{WZ-pert},
	we define $\tilde{U}$ and $\tilde{V}$ from the compatibility
	conditions \eqref{pccUV}; see \eqref{pUV}.
	To remove possibly nonzero defects
	$\tilde{H}_j(\infty;\nu,a,b)$ in the limits of
	$\tilde{U}(\xi)$ and $\tilde{V}(\xi)$ as $\xi\to\infty$,
	we use the two far-field amplitudes $a=a(\nu)$ and $b=b(\nu)$,
	see \eqref{WZ-pert}.
	They provide two independent degrees of freedom with which to
	satisfy \eqref{Hnuab}, thereby ensuring that
	$\tilde{U}$ and $\tilde{V}$ defined in \eqref{pUV} are compact 
	perturbations of $U$ and $V$, respectively.
\end{remark}

We will apply Lemma \ref{L2} in the 
proof of the next result at points $(t_0,\xi_0)$ 
where one of the four vectors $(W,W_\xi,W_{\xi\xi})$, 
$(W,W_t,W_\xi)$, $(Z,Z_\xi,Z_{\xi\xi})$, 
or $(Z,Z_t,Z_\xi)$ attains the value $(\pi,0,0)$.
For instance, if $(W,W_\xi,W_{\xi\xi})(t_0,\xi_0)=(\pi,0,0)$, 
we use the perturbation of the solution 
described in item 1) of Lemma \ref{L2}.
This strategy was originally developed by Bressan 
and Chen in their analysis of a quasilinear 
second order wave equation \cite[Lemma 5]{BC17}.

\begin{lemma}\label{L3}
Let $\mathcal{K}$ denote the family of all
solutions of \eqref{ODE} for which (recall \eqref{CktcN})
\begin{equation}\label{sol-cl}
(U,V,W,Z,q)\in C\left(
[-T,T], \mathcal{C}^k_{\Omega_{\mathrm{tcN}}}
\right),\quad k\geq5.
\end{equation}
Consider a subfamily 
$\mathcal{K}^\prime\subset\mathcal{K}$
such that any solution
$(U,V,W,Z,q)\in\mathcal{K}^\prime$ satisfies
the following properties:
\begin{equation}\label{ins}
	\begin{split}
		&(W,W_\xi,W_{\xi\xi})(t,\xi)\neq(\pi,0,0),\quad
		(W,W_t,W_{\xi})(t,\xi)\neq(\pi,0,0),\\
		&(Z,Z_\xi,Z_{\xi\xi})(t,\xi)\neq(\pi,0,0),\quad
		(Z,Z_t,Z_{\xi})(t,\xi)\neq(\pi,0,0),
	\end{split}
\end{equation}
for all $(t,\xi)\in[-T,T]\times\R$.

Then $\mathcal{K}^\prime$ is 
a relatively open and dense subset of
$\mathcal{K}$ in the topology of
$C\left([-T,T],
\mathcal{C}^k_{\Omega_{\mathrm{tcN}}} \right)$.
\end{lemma}

\begin{proof}
Inspired by the proof of \cite[Lemma 5]{BC17},
we will use the representation
\begin{equation}\label{Kpr-int}
\mathcal{K}^\prime=\bigcap\limits_{j=1}^4
\mathcal{K}^\prime_j,
\end{equation}
where the subfamilies 
$\mathcal{K}^\prime_j\subset\mathcal{K}^\prime$,
$j=1,\dots,4$, are such that one 
of the four conditions in \eqref{ins} is satisfied.
For example, $\mathcal{K}^\prime_1$ 
consists of all solutions $(U,V,W,Z,q)$ such that
$(W,W_\xi,W_{\xi\xi})(t,\xi)\neq(\pi,0,0)$
for all $(t,\xi)\in\Lambda_{T,M}$.

Since \eqref{Kpr-int} is a finite intersection,
it is sufficient to show that each subfamily
$\mathcal{K}^\prime_j$, $j=1,\dots,4$, is a
relatively open and dense subset of $\mathcal{K}$.
We provide the detailed proof for
$\mathcal{K}^\prime_1$, as the remaining 
subfamilies can be handled in a similar manner.

\medskip
\textbf{Step 1.}
Consider any $(U,V,W,Z,q)\in\mathcal{K}^\prime_1$.
Since $W,Z\in C\left([-T,T],L^2\cap C^1\right)$,
there exists $M=M(T)>0$ such that
\begin{equation}\label{WZ1}
	|W(t,\xi)|, |Z(t,\xi)|<1,\quad
	\text{for all }t\in[-T,T],\,\,
	|\xi|\geq M.
\end{equation}
Notice that the solutions of the ODE system \eqref{ODE}
corresponding to small perturbations
of the initial data $(U,V,W,Z,q)(0,\xi)$ in
$\mathcal{C}^k_{\Omega_{\mathrm{tcN}}}$
remain close to the original solution $(U,V,W,Z,q)(t,\xi)$
in the topology $C\left([-T,T],
\mathcal{C}^k_{\Omega_{\mathrm{tcN}}} \right)$.
Therefore, \eqref{WZ1} implies that 
$\mathcal{K}^\prime_1$ is a relatively open 
subset of $\mathcal{K}$ in this topology.

\medskip

\textbf{Step 2.}
In view of Step 1, it remains to prove that
$\mathcal{K}^\prime_1$ is a dense subset
of $\mathcal{K}$.
For any solution 
$(U,V,W,Z,q)\in\mathcal{K}$ we define
\begin{equation}\label{Lambd_T}
	\Lambda_{T,M}=\left\{
	(t,\xi)\in\R^2:
	|t|\leq T,\,|\xi|\leq M
	\right\},\quad
	T,M>0,
\end{equation}
where $M=M(T)>0$ is the same as in \eqref{WZ1}.
We consider the set
$\Lambda_{T+\delta,M+\delta}$ with some $\delta>0$,
which allows us to work on
$\mathrm{int}(\Lambda_{T+\delta,M+\delta})$, a manifold without
boundary, while retaining
$\Lambda_{T,M}$ as a subset in Step 3 below.
For any point 
$(t_1,\xi_1)\in\Lambda_{T+\delta,M+\delta}$ the vector 
$(W,W_\xi,W_{\xi\xi})(t_1,\xi_1)$ 
can either equal or not equal 
to the fixed vector $(\pi,0,0)$.
Therefore, the following
two cases are possible:
\begin{enumerate}[(I)]
	\item $(W,W_\xi,W_{\xi\xi})(t_1,\xi_1)
	\neq(\pi,0,0)$, which implies that
	the $(W,W_\xi,W_{\xi\xi})(t,\xi)
	\neq(\pi,0,0)$
	for all $(t,\xi)$ in a closed neighborhood
	$\overline{\mathcal{U}_{(t_1,\xi_1)}}
	\cap\Lambda_{T+\delta,M+\delta}$;

	\item $(W,W_\xi,W_{\xi\xi})(t_1,\xi_1)
	=(\pi,0,0)$.
	Here we consider the perturbation 
	$\left(\tilde{U},\tilde{V},\tilde{W},
	\tilde{Z},\tilde{q}\right)(t,\xi;\nu)$,
	$\nu=(\nu_1,\nu_2,\nu_3)$, of the solution
	given in Lemma \ref{L2}, part 1), and
	a neighborhood 
	$\mathcal{U}_{(t_1,\xi_1)}$
	such that \eqref{r1} holds
	for all $(t,\xi)\in\overline{\mathcal{U}_{(t_1,\xi_1)}}
	\cap\Lambda_{T+\delta,M+\delta}$.	
\end{enumerate}
Next we take a point 
$(t_2,\xi_2)\in\Lambda_{T+\delta,M+\delta}\setminus
\mathcal{U}_{(t_1,\xi_1)}$ and argue in 
the same way as in (I) and (II) for $(t_1,\xi_1)$.

We explain how to superpose the perturbations when both
$(t_1,\xi_1)$ and $(t_2,\xi_2)$
fall into case (II).
Here we apply Lemma \ref{L2} with 
three new parameters $(\nu_4,\nu_5,\nu_6)$ in place of
$(\nu_1,\nu_2,\nu_3)$ and with
\[
\left(\tilde{U},\tilde{V},\tilde{W},
\tilde{Z},\tilde{q}\right)(t_2,\xi;\nu_1,\nu_2,\nu_3)
\quad\text{in place of}\quad
(U,V,W,Z,q)(t_0,\xi).
\]
For example, the new perturbation
$\check{W}$
has the following form (cf.~\eqref{WZ-pert})
\begin{equation}\label{ch-W}
	\check{W}(t_2,\xi;\nu_1,\dots,\nu_6)
	=\tilde{W}(t_2,\xi;\nu_1,\nu_2,\nu_3)
	+\sum\limits_{j=4}^6\nu_jf_j(\xi)+a(\nu_4,\nu_5,\nu_6)\psi(\xi).
\end{equation}
This
results in the perturbation
$\left(\check{U},\check{V},\check{W},
\check{Z},\check{q}\right)(t,\xi;\nu)$,
with $\nu=(\nu_1,\dots,\nu_6)$.
Lemma \ref{L2} yields
\begin{equation*}
	\mathrm{rank}\left(\left.
	D_{(\nu_4,\nu_5,\nu_6)}\left(
	\check{W},\check{W}_\xi,\check{W}_{\xi\xi}\right)
	\right|_{(t,\xi;\nu)=(t_2,\xi_2;0)}
	\right)=3.
\end{equation*}
Then, \eqref{ch-W} yields that
\begin{equation}\label{t2ch-W}
\check{W}(t_2,\xi;\nu_1,\nu_2,\nu_3,0,0,0)
=\tilde{W}(t_2,\xi;\nu_1,\nu_2,\nu_3),\quad
\text{for all }\,\xi,\nu_1,\nu_2,\nu_3.
\end{equation}
In view of the uniqueness of the solution of the ODE system \eqref{ODE}, we have \eqref{t2ch-W} for all $t\in[-T,T]$
and therefore the rank condition \eqref{r1} holds
for the new perturbation
$\check{W}$ at $t=t_1$ as well:
\begin{equation*}
	\mathrm{rank}\left(\left.
	D_{(\nu_1,\nu_2,\nu_3)}\left(
	\tilde{W},\tilde{W}_\xi,\tilde{W}_{\xi\xi}\right)
	\right|_{(t,\xi;\nu)=(t_1,\xi_1;0)}\right)=3.
\end{equation*}
	
Continuing this procedure, we obtain
the open cover of $\Lambda_{T+\delta,M+\delta}$ by the
neighborhoods $\mathcal{U}_{(t_\beta,\xi_\beta)}$,
$\beta\in\mathcal{B}$.
Taking into account that $\Lambda_{T+\delta,M+\delta}$
is a compact set, we have a finite cover 
of $\Lambda_{T+\delta,M+\delta}$, that is
\begin{equation}
	\label{Lcov}
	\Lambda_{T+\delta,M+\delta}\subset\bigcup
	\limits_{i=1}^{\tilde{N}}\mathcal{U}_{(t_i,\xi_i)},
	\quad
	\text{for some }\tilde{N}>0.
\end{equation}
In what follows, we consider 
a perturbed solution
$\left(\tilde{U},\tilde{V},\tilde{W},
\tilde{Z},\tilde{q}\right)(t,\xi;\nu)$
with $\nu=(\nu_1,\dots,\nu_{3\tilde{N}_1})$, 
$\tilde{N}_1\leq \tilde{N}$, obtained in the neighborhoods 
$\mathcal{U}_{(t_i,\xi_i)}$, $i=1,\dots, \tilde{N}$.
Notice that
\begin{equation}\label{nu=0}
	(U,V,W,Z,q)(t,\xi)=
	\left(\tilde{U},\tilde{V},\tilde{W},
	\tilde{Z},\tilde{q}\right)(t,\xi;0),\quad
	(t,\xi)\in\Lambda_{T+\delta,M+\delta},
\end{equation}
and for any 
$(t,\xi)\in\mathcal{U}_{(t_i,\xi_i)}$,
$i=1,\dots,\tilde{N}$, we have either
\begin{equation}\label{rp}
\begin{split}
	&\left(\tilde{W},\tilde{W}_\xi,
	\tilde{W}_{\xi\xi}\right)(t,\xi;0)
	\neq(\pi,0,0),\qquad\qquad\qquad\qquad\,\,\,
	\text{for all }(t,\xi)\in
	\overline{\mathcal{U}_{(t_i,\xi_i)}},\,\,\text{or}\\
	&\mathrm{rank}\left(\left.
	D_{(\nu_{3j-2},\nu_{3j-1},\nu_{3j})}\left(
	\tilde{W},\tilde{W}_\xi,\tilde{W}_{\xi\xi}\right)
	\right|_{(t,\xi;0)}
	\right)=3,\quad
	\text{for all }(t,\xi)\in
	\overline{\mathcal{U}_{(t_i,\xi_i)}},
\end{split}
\end{equation}
with $i=1,\dots,\tilde{N}$ and 
the corresponding $j=1,\dots,\tilde{N_1}$.

\medskip
\textbf{Step 3.}
Consider a map
\begin{equation}\label{tW-map}
\left(\tilde{W},\tilde{W}_\xi,
\tilde{W}_{\xi\xi}\right):
\mathrm{int}(\Lambda_{T+\delta,M+\delta})\times(-\ve,\ve)^{3\tilde{N}_1}
\mapsto\R^3,
\end{equation}
where a small $\ve>0$ is taken such that \eqref{rp} holds
for all $\nu\in(-\ve,\ve)^{3\tilde{N}_1}$ and for all
$i=1,\dots,\tilde{N}$, i.e., we have either
\begin{equation}\label{rpnu}
\begin{split}
	&\left(\tilde{W},\tilde{W}_\xi,
	\tilde{W}_{\xi\xi}\right)(t,\xi;\nu)
	\neq(\pi,0,0),\qquad\qquad\qquad\qquad\,\,\,
	\text{for all }(t,\xi)\in\mathcal{U}_{(t_i,\xi_i)},
	\,\,\nu\in(-\ve,\ve)^{3\tilde{N}_1},
	\,\,\text{or}\\
	&\mathrm{rank}\left(\left.
	D_{(\nu_{3j-2},\nu_{3j-1},\nu_{3j})}\left(
	\tilde{W},\tilde{W}_\xi,\tilde{W}_{\xi\xi}\right)
	\right|_{(t,\xi;\nu)}
	\right)=3,\quad
	\text{for all }(t,\xi)\in\mathcal{U}_{(t_i,\xi_i)},
	\,\,\nu\in(-\ve,\ve)^{3\tilde{N}_1},
\end{split}
\end{equation}
with $i=1,\dots,\tilde{N}$ and the 
corresponding $j=1,\dots,\tilde{N_1}$.

We apply the parametric transversality theorem (Theorem \ref{BTL})
for the three-dimensional map given in \eqref{tW-map}
and the zero-dimensional submanifold
$\mathcal{W}=\{(\pi,0,0)\}$. 
Observe that the manifolds
$\mathcal{X}=\mathrm{int}(\Lambda_{T+\delta,M+\delta})$,
$\mathcal{N}=(-\ve,\ve)^{3\tilde{N}_1}$,
and $\mathcal{Y}=\R^3$ are manifolds without boundary.
Taking into account that
\begin{equation}\label{reg-ptt}
	\max\{0,\dim\mathcal{Y}+\dim\mathcal{W}
	-\dim\mathcal{X}\}=\max\{0,3+0-2\}=1,
\end{equation}
the map $\left(\tilde{W},\tilde{W}_\xi,
\tilde{W}_{\xi\xi}\right)$
satisfies the regularity condition \eqref{F-reg}:
\begin{equation*}
	\left(\tilde{W},\tilde{W}_\xi,
	\tilde{W}_{\xi\xi}\right)\in
	C^2\left(
	\mathrm{int}(\Lambda_{T+\delta,M+\delta})\times(-\ve,\ve)^{3\tilde{N}_1}
	\right),
\end{equation*}
where we have used that $k\geq5$ in \eqref{sol-cl}
and the construction of the perturbations (see \eqref{3pars}
and Remark \ref{k-reg} below).
Taking into account \eqref{rpnu},
we have that $(\pi,0,0)$ is a regular value 
for the map 
$\left(\tilde{W},\tilde{W}_\xi,\tilde{W}_{\xi\xi}\right)$
(as per Definition \ref{RV}).
Thus, we can employ the parametric transversality result 
for this map with the
submanifold $\mathcal{W}=\{(\pi,0,0)\}$,
which yields that
$\left(\tilde{W},\tilde{W}_\xi,\tilde{W}_{\xi\xi}\right)
(\cdot,\cdot;\nu)$ is transverse to $\{(\pi,0,0)\}$
for $\nu$ from a dense subset of $(-\ve,\ve)^{3\tilde{N}_1}$.
Since the map is three-dimensional, while the 
domain
$\mathrm{int}(\Lambda_{T+\delta,M+\delta})$
is two-dimensional, we conclude 
that the transversality condition can only be 
satisfied when the preimage 
of $\{(\pi,0,0)\}$ is empty.

Therefore, there exists a sequence 
$\{\nu^n\}_{n=1}^\infty
\subset(-\ve,\ve)^{3\tilde{N}_1}$
such that $\nu^n\to 0$ as $n\to\infty$ and
\[
\left(\tilde{W},\tilde{W}_\xi,
\tilde{W}_{\xi\xi}\right)(t,\xi;\nu^n)
\neq(\pi,0,0),\quad
\text{for all } (t,\xi)\in\mathrm{int}(\Lambda_{T+\delta,M+\delta}),\,\,
n\in\N,
\]
and thus $\left(\tilde{U},\tilde{V},\tilde{W},
\tilde{Z},\tilde{q}\right)
(\cdot,\cdot;\nu^n)\in\mathcal{K}^\prime_1$
(note that 
	$\Lambda_{T,M}\subset\mathrm{int}(\Lambda_{T+\delta,M+\delta})$).
Finally, using \eqref{nu=0} and taking into account that
at each step of the construction of the perturbed solution
it remains in an arbitrarily small neighborhood of the original solution
in the topology $C\left([-T,T],
\mathcal{C}^k_{\Omega_{\mathrm{tcN}}} \right)$,
provided $\nu$ is sufficiently small (see the perturbations
\eqref{WZ-pert}, \eqref{q-pert}, and \eqref{pUV}
with $t_i$ in place of $t_0$),
we conclude that
\begin{equation*}
	\left(\tilde{U},\tilde{V},\tilde{W},
	\tilde{Z},\tilde{q}\right)
	(\cdot,\cdot;\nu^n)\to
	\left(U,V,W,Z,q\right)(\cdot,\cdot)
	\quad\text{as }n\to\infty,
\end{equation*}
in the topology $C\left([-T,T],
\mathcal{C}^k_{\Omega_{\mathrm{tcN}}} \right)$.
Thus, we have shown that $\mathcal{K}^\prime_1$ 
is dense in $\mathcal{K}^\prime$.
\end{proof}
The following remark explains the regularity requirement
on $k$ in \eqref{sol-cl}.
\begin{remark}\label{k-reg}
	The regularity assumption on $k$ in \eqref{sol-cl}
	is determined by the application of the parametric
	transversality theorem
	to the maps $\left(\tilde{W},\tilde{W}_\xi,
	\tilde{W}_{\xi\xi}\right)$,
	$\left(\tilde{W},\tilde{W}_t,
	\tilde{W}_{\xi}\right)$,
	$\left(\tilde{Z},\tilde{Z}_\xi,
	\tilde{Z}_{\xi\xi}\right)$,
	and 
	$\left(\tilde{Z},\tilde{Z}_t,
	\tilde{Z}_{\xi}\right)$,
	cf.~\eqref{tW-map}.
	Since Theorem \ref{BTL} is employed for the manifolds
	$\mathcal{X}=\mathrm{int}(\Lambda_{T+\delta,M+\delta})$,
	$\mathcal{N}=(-\ve,\ve)^{3\tilde{N}_1}$,
	and $\mathcal{Y}=\R^3$, as well as
	the zero-dimensional submanifold 
	$\mathcal{W}=\{(\pi,0,0)\}\subset\mathcal{Y}$,
	we have \eqref{reg-ptt},
	and therefore the above maps must be of class at least $C^2$,
	see \eqref{F-reg}.
	For $k=5$, we have $\tilde{W},\tilde{Z}\in C^4(\R)$, which yields that
	$\tilde{W}_{\xi\xi},\tilde{Z}_{\xi\xi}\in C^2(\R)$ and
	thus Theorem \ref{BTL} can be applied.

	Note that replacing the actual maps by smooth approximations
	does not, by itself, show that the original solution family has
	the desired transversality property;
	cf.~\cite[p.~348,~Step~2]{BC17},
	\cite[p.~2930,~Step~4]{LZ17}, and \cite[Lemma~4.4]{CCCS18}.
	Indeed, suppose, for example, that
	$\left(\tilde{F}_1,\tilde{F}_2,\tilde{F}_3\right)(t,\xi;\nu)$
	is a smooth approximation of $\left(\tilde{W},\tilde{W}_\xi,
	\tilde{W}_{\xi\xi}\right)(t,\xi;\nu)$. 
	There can exist points $(t^\prime,\xi^\prime;\nu^\prime)$ such that
	\[
	\left|\left(\tilde{F}_1,\tilde{F}_2,\tilde{F}_3\right)
	(t^\prime,\xi^\prime;\nu^\prime)-(\pi,0,0)\right|
	\ll
	\left|\left(\tilde{F}_1,\tilde{F}_2,\tilde{F}_3\right)
	(t^\prime,\xi^\prime;\nu^\prime)
	-\left(\tilde{W},\tilde{W}_\xi,
	\tilde{W}_{\xi\xi}\right)
	(t^\prime,\xi^\prime;\nu^\prime)\right|.
	\]
	Consequently, the original map
	$\left(\tilde{W},\tilde{W}_\xi,
	\tilde{W}_{\xi\xi}\right)$
	may still attain $(\pi,0,0)$ at
	$(t^\prime,\xi^\prime;\nu^\prime)$, even if it is arbitrarily close
	to $\left(\tilde{F}_1,\tilde{F}_2,\tilde{F}_3\right)$
	in a $C^r$ topology.
\end{remark}

\subsection{Proof of Theorem \ref{Thmgr}}\label{prTh1}

By combining Lemmas \ref{L2}, and \ref{L3}, 
we obtain the generic regularity of solutions to 
the two-component Novikov system, thereby 
completing the proof of Theorem \ref{Thmgr}.

\medskip
\textbf{Step 1.}
For any initial data 
$(\tilde{u}_0,\tilde{v}_0)\in\Upsilon^k$, 
$k\geq5$, we consider a neighborhood
(recall \eqref{Ups}--\eqref{d-Ups})
\begin{equation*}
	\mathcal{U}_{(\tilde{u}_0,\tilde{v}_0)}(\delta)=
	\{(u_0,v_0)\in\Upsilon^k:
	d_{\Upsilon^k}((u_0,v_0),(\tilde{u}_0,\tilde{v}_0))
	<\delta\},
\end{equation*}
with some small $\delta>0$.
To prove our theorem, it is enough to show that
there exists 
a subset that is open in $\Upsilon^k$ and dense
in $\Upsilon^{k-1}\cap C^k_{\mathrm{loc}}(\R)$
$\tilde{\mathcal{M}}_T\subset
\mathcal{U}_{(\tilde{u}_0,\tilde{v}_0)}(\delta)$
such that for any 
$(u_0,v_0)\in\tilde{\mathcal{M}}_T$ the global 
conservative solution $(u,v)(t,x)$ given 
in Theorem \ref{Thm} is locally of class $C^k$ 
in the complement of finitely 
many $C^{k-1}$ curves in $[-T,T]\times\R$.

Since $(u_0,v_0)\in\Upsilon^k$, 
Theorem \ref{Ckgl} implies that
\begin{equation}\label{ODEreg}
	(U,V,W,Z,q)\in C\left([-T,T],
	\mathcal{C}^k_{\Omega_{\mathrm{tcN}}}
	\right).
\end{equation}
Here and throughout the proof $(U,V,W,Z,q)$
denotes the solution of \eqref{ODE} subject to initial 
data \eqref{id} corresponding to $(u_0,v_0)$.

\medskip
\textbf{Step 2.}
Taking into account that
$W,Z\in C\left([-T,T],L^2\cap C^1\right)$,
there exists $\tilde{M}_T>0$ such that
\begin{equation}\label{WZl1}
	|W(t,\xi)|, |Z(t,\xi)|<1,\quad
	\text{for all }t\in[-T,T],\,\,
	|\xi|\geq\tilde{M}_T.
\end{equation}
From \eqref{WZl1} we have that
the value $\pi$ is attained neither by 
$W$ nor $Z$ for sufficiently large $\xi$.
Taking into account that $y(t,\cdot)$ 
is monotone increasing and using
\eqref{uvdef}, \eqref{pxiy}, \eqref{ODEreg}, 
we conclude that there exists
$$
Q_{T,\delta}=\left\{(t,x)\in\R^2:|t|\leq T,
\,|x|\leq R_{T,\delta}\right\},\quad
\text{for some }R_{T,\delta}>0,
$$
such that  $u$ and $v$ are $C^k$ in
$\left([-T,T]\times\R\right)
\setminus Q_{T,\delta}$ for all initial data in 
$\mathcal{U}_{(\tilde{u}_0,\tilde{v}_0)}(\delta)$.

Consequently, it remains to analyze 
the singularities
in the rectangle $Q_{T,\delta}$.

\medskip
\textbf{Step 3.}
For any fixed 
$(u_0,v_0)\in\mathcal{U}_{(\tilde{u}_0,\tilde{v}_0)}(\delta)$,
we define the set $\tilde{G}$ as the image of the map
$(t,y(t,\xi))$ defined on $\Lambda_{T,M}$ 
(recall \eqref{Lambd_T}):
\begin{equation}\label{tGu0v0}
	\tilde{G}_{(u_0,v_0)}=\left\{
	(t,x): x=y(t,\xi),\,\,(t,\xi)\in\Lambda_{T,M}
	\right\}.
\end{equation}
Taking $M\geq\tilde{M}_T$ large enough 
(recall \eqref{WZl1}), we obtain that
$Q_{T,\delta}\subset\tilde{G}_{(u_0,v_0)}$ 
for all $(u_0,v_0)\in\mathcal{U}_{(\tilde{u}_0,\tilde{v}_0)}$.
Then $\tilde{\mathcal{M}}_T\subset 
\mathcal{U}_{(\tilde{u}_0,\tilde{v}_0)}$
is defined as follows: we say that the 
initial data $(u_0,v_0)$ 
belongs to $\tilde{\mathcal{M}}_T$,
if for the corresponding solution
$(U,V,W,Z,q)$ of the ODE system we have
\eqref{ins} for all $(t,\xi)\in\Lambda_{T,M}$.
Below we show that $\tilde{\mathcal{M}}_T$
is open and dense subset of 
$\mathcal{U}_{(\tilde{u}_0,\tilde{v}_0)}(\delta)$.

\medskip
\textbf{Step 4.} Here we prove that $\tilde{\mathcal{M}}_T$
is an open in the topology of $\Upsilon^k$.
Suppose to the contrary that there exists
$(u_0,v_0)\in\tilde{\mathcal{M}}_T$
and a sequence $\{(u_{0,n},v_{0,n})\}_{n=1}^\infty$
such that $(u_{0,n},v_{0,n})\to(u_0,v_0)$
as $n\to\infty$ in $\Upsilon^k$ and 
$(u_{0,n},v_{0,n})\not\in\tilde{\mathcal{M}}_T$
for all $n\in\N$.
By the definition of $\tilde{\mathcal{M}}_T$, see step 3,
there exists a sequence 
$\{(t_n,\xi_n)\}_{n=1}^\infty\subset\Lambda_{T,M}$
such that one of the vectors in \eqref{ins}
corresponding to the initial data $(u_{0,n},v_{0,n})$
attains $(\pi,0,0)$.
Consequently, there exists a subsequence, also denoted 
by $\{(t_n,\xi_n)\}_{n=1}^\infty$, such that, for example,
$(W_n,\pxi W_n, \pxi^2 W_n)=(\pi,0,0)$ 
(the other vectors can be treated in the same way).

Since $\Lambda_{T,M}$ is a compact set, there exists
$(\tilde{t},\tilde{\xi})\in\Lambda_{T,M}$ and
a subsequence, which we again denote by 
$\{(t_n,\xi_n)\}_{n=1}^\infty$, such that
$\left(t_n,\xi_n\right)\to(\tilde{t},\tilde{\xi})$.
Using that $u_{0,n}\to u_0$ in $C^k(\R)$, we conclude that 
$W_n\to W$ in $C\left([-T,T],C^{k-1}(\R)\right)$
(see \eqref{id}), and therefore
$\left(W,W_\xi,W_{\xi\xi}\right)
(\tilde{t},\tilde{\xi})=(\pi,0,0)$, 
which is a contradiction.

\medskip
\textbf{Step 5.} 
Let us show that $\tilde{\mathcal{M}}_T$ is dense
in $\mathcal{U}_{(\tilde{u}_0,\tilde{v}_0)}(\delta)$, i.e.,
for any $(u_0,v_0)\in
\mathcal{U}_{(\tilde{u}_0,\tilde{v}_0)}(\delta)$
there exists 
$\{(u_{0,n},v_{0,n})\}_{n=1}^\infty
\subset\tilde{\mathcal{M}}_T$ such that
\begin{equation}\label{conv-un}
	(u_{0,n},v_{0,n})\to(u_0,v_0),\quad n\to\infty,
	\quad\text{in }\,\Upsilon^{k-1}\cap C^k(-R,R),\,\text{ for any }\, R>0.
\end{equation}
Recall that $(U,V,W,Z,q)(t,\xi)$ is the solution
of the ODE system \eqref{ODE}
corresponding to the initial data \eqref{y0}--\eqref{id} associated with
	$(u_0,v_0)$. In particular, $q(0,\xi)=1$.
Then Lemma \ref{L3} implies that there exists a sequence
$\{(U_n,V_n,W_n,Z_n,q_n)\}_{n=1}^\infty$ such that
\begin{enumerate}[(i)]
	\item $(U_n,V_n,W_n,Z_n,q_n)\to(U,V,W,Z,q)$
	in $C\left([-T,T],\mathcal{C}^k_{\Omega_{\mathrm{tcN}}}\right)$,
		\item \eqref{ins}
		holds for every $(U_n,V_n,W_n,Z_n,q_n)$, $n\in\N$,
	and $(t,\xi)\in[-T,T]\times\R$.
\end{enumerate}
Note that, in general, $q_n(0,\xi)\not\equiv1$.

Set $(u_{0,n},v_{0,n})(x)=(u_n,v_n)(0,x)$, where
$(u_n,v_n)(t,x)$ is obtained from
$(U_n,V_n,W_n,Z_n,q_n)(t,\xi)$ through \eqref{uvdef}:
\begin{equation}\label{un-def}
	u_n(t,y_n(t,\xi))=U_n(t,\xi),\quad
	v_n(t,y_n(t,\xi))=V_n(t,\xi),\quad t\in[-T,T],\,\,\xi\in\R,
	\quad n\in\N.
\end{equation}
Here $y_n(t,\xi)$ has the following form
(see \eqref{char1} and \eqref{y0ode}):
\begin{equation}\label{yny0n}
	y_n(t,\xi)=y_{0,n}(\xi)+\int_0^t(U_nV_n)(\tau,\xi)\,d\tau,
	\quad \text{with}\quad 
	y_{0,n}(\xi)=\int_0^\xi\left(q_n\cos^2\frac{W_n}{2}
	\cos^2\frac{Z_n}{2}\right)(0,\eta)\,d\eta.
\end{equation}
Note that the definition of $y_{0,n}(\xi)$ given in \eqref{yny0n}
is consistent with \eqref{y0}, where one should take
$(u_{0,n},v_{0,n})$ in place of $(u_0,v_0)$.

\medskip

\textbf{Step 5.1.} Here we prove \eqref{conv-un}. 
Since $W_n(0,\cdot)$, $Z_n(0,\cdot)$, and $q_n(0,\cdot)$
converge to, respectively, 
$W(0,\cdot)$, $Z(0,\cdot)$, and $q(0,\cdot)\equiv1$ in $C^{k-1}(\R)$,
we obtain from the second equation in \eqref{yny0n} that
\begin{equation}
	\label{y0y0n}
\|y_0-y_{0,n}\|_{C^r(\R)}\to0,\quad 1\leq r\leq k,
\quad n\to\infty.
\end{equation}
We now show that \eqref{y0y0n} also holds for $r=0$.
For all $\xi\in\R$ we have the following estimates:
\begin{equation}\label{|y0y0n|}
\begin{split}
	|y_0-y_{0,n}|(\xi)&\leq\int_{0}^\xi
	\left|(q-q_n)\cos^2\frac{W_n}{2}\cos^2\frac{Z_n}{2}\right|
	(0,\eta)\,d\eta\\
	&\quad+\int_{0}^\xi q(0,\eta)
	\left|\cos^2\frac{W}{2}\cos^2\frac{Z}{2}-
	\cos^2\frac{W_n}{2}\cos^2\frac{Z_n}{2}\right|
	(0,\eta)\,d\eta\\
	&\leq\|q_n(0,\cdot)-1\|_{L^1(\R)}+I_n,
\end{split}
\end{equation}
with
\begin{equation*}
	I_n=\int_{-\infty}^\infty\left|\cos^2\frac{W}{2}\cos^2\frac{Z}{2}-
	\cos^2\frac{W_n}{2}\cos^2\frac{Z_n}{2}\right|(0,\eta)d\eta.
\end{equation*}
The integral $I_n$ can be estimated as follows:
\begin{equation}\label{I_n}
\begin{split}
	I_n&\leq\int_{-\infty}^\infty
	\left|\cos^2\frac{W}{2}-
	\cos^2\frac{W_n}{2}\right|(0,\eta)\,d\eta
	+\int_{-\infty}^\infty
	\left|\cos^2\frac{Z}{2}-
	\cos^2\frac{Z_n}{2}\right|(0,\eta)\,d\eta\\
	&\leq4\int_{-\infty}^\infty
	\left|\sin\frac{W-W_n}{4}
	\sin\frac{W+W_n}{4}\right|(0,\eta)\,d\eta
	+4\int_{-\infty}^\infty
	\left|\sin\frac{Z-Z_n}{4}
	\sin\frac{Z+Z_n}{4}\right|(0,\eta)\,d\eta
\end{split}
\end{equation}
Using \eqref{sinx} and the Cauchy--Schwarz inequality, we obtain that
\begin{equation}\label{W0W0n}
	\int_{-\infty}^\infty
	\left|\sin\frac{W-W_n}{4}
	\sin\frac{W+W_n}{4}\right|(0,\eta)\,d\eta
	\leq C\|(W-W_n)(0,\cdot)\|_{L^2(\R)},
\end{equation}
where $C>0$ depends on the $L^2(\R)$ norms of $W(0,\cdot)$ and
$W_n(0,\cdot)$.
Applying the same argument to the second integral on the right-hand
side of \eqref{I_n}, which involves $Z$ and $Z_n$, and recalling
\eqref{|y0y0n|}, we obtain
$\|y_0-y_{0,n}\|_{C(\R)}\to0$ as $n\to\infty$.
Together with \eqref{y0y0n}, this yields
\begin{equation}
	\label{Cky0y0n}
	\|y_0-y_{0,n}\|_{C^k(\R)}\to0,
	\quad n\to\infty.
\end{equation}

Recalling \eqref{un-def}, we obtain the following estimates
(here we choose $\xi=\xi_n$ so that $x=y_{0,n}(\xi)$ and use the
notation $U_{0,n}(\xi)=U_n(0,\xi)$ and $U_0(\xi)=U(0,\xi)$):
\begin{equation}\label{u0n-est}
\begin{split}
	|u_{0,n}-u_0|(x)&=|u_{0,n}-u_0|(y_{0,n}(\xi))\\
	&\leq |u_{0,n}(y_{0,n}(\xi))-u_0(y_0(\xi))|
	+|u_0(y_0(\xi))-u_0(y_{0,n}(\xi))|\\
	&=\left|U_{0,n}-U_0\right|(\xi)+|u_0(y_0(\xi))-u_0(y_{0,n}(\xi))|\\
	&\leq\left\|U_{0,n}-U_0\right\|_{C(\R)}
	+\|u_0^\prime\|_{C(\R)}\|y_{0,n}-y_0\|_{C(\R)},
\end{split}
\end{equation}
which, together with
\eqref{Cky0y0n} and the limit
$\|U_{0,n}-U_0\|_{C^k(\R)}\to0$ as $n\to\infty$,
yields that $\|u_{0,n}-u_0\|_{C(\R)}\to0$, $n\to\infty$.
Since $(u_0,v_0)$ and $(u_{0,n},v_{0,n})$ are regular,
\eqref{Cky0y0n} also shows that the corresponding characteristics
have no critical points:
\begin{equation}\label{pxiyn}
	c\leq|\pxi y_0|\leq C,\quad c\leq|\pxi y_{0,n}|\leq C,\quad
	\text{for all }\,n\in\N,\quad c,C>0.
\end{equation}
Using \eqref{pxiyn} and that (see \eqref{un-def})
\begin{equation}\label{uprime}
u_0^\prime(y_0(\xi))=\frac{\pxi U_0(\xi)}{\pxi y_0(\xi)},
\quad
u_{0,n}^\prime(y_{0,n}(\xi))=
\frac{\pxi U_{0,n}(\xi)}{\pxi y_{0,n}(\xi)},
\end{equation}
we conclude, reasoning as in \eqref{u0n-est}, that
$\|u_{0,n}^\prime-u_0^\prime\|_{C(\R)}\to0$ as $n\to\infty$.
Arguing similarly for $v_{0,n}$ and
for the higher derivatives of 
$(u_{0,n},v_{0,n})$, we obtain that 
\begin{equation}\label{u0n-k-1}
	\|u_{0,n}-u_0\|_{C^{k-1}(\R)},\,
	\|v_{0,n}-v_0\|_{C^{k-1}(\R)}\to0,\quad n\to\infty.
\end{equation}
Note that we do not claim \eqref{u0n-k-1} in $C^k(\R)$ since,
in general, $\pxi^kU_0(\xi)$ is not uniformly continuous on $\R$,
cf. the last inequality in \eqref{u0n-est}.
Thus, we obtain $C^k$ convergence only on compact sets:
\begin{equation}\label{u0n-k}
	\|u_{0,n}-u_0\|_{C^{k}(-R,R)},\,
	\|v_{0,n}-v_0\|_{C^{k}(-R,R)}\to0,
	\quad\text{for all }\,R>0,
	\quad n\to\infty,
\end{equation}

It remains to show that
$(u_{0,n},v_{0,n})\to(u_0,v_0)$ in $\Sigma$	(see \eqref{Sig}--\eqref{d^2}).
Using \eqref{un-def} with $t=0$,
the third inequality in \eqref{u0n-est},
\eqref{pxiyn},
and taking $dx=\pxi y_{0,n}\,d\xi$,
we have that
\begin{equation*}
	\|u_{0,n}-u_0\|_{L^2}^2\leq C\|U_{0,n}-U_0\|_{L^2}^2
	+C\int_{-\infty}^\infty
	\left(U_0(\xi)-U_0(y_{0}^{-1}[y_{0,n}(\xi)])\right)^2\,d\xi.
\end{equation*}
Taking into account that
$y_0^{-1}\circ y_{0,n}\to\operatorname{id}$ uniformly on $\R$,
see \eqref{Cky0y0n} and \eqref{pxiyn},
as well as $U_{0,n}\to U_0$ in $L^2(\R)$,
we have that $u_{0,n}\to u_0$ in $L^2(\R)$.
The same argument, together with \eqref{uprime}, yields convergence
of $u_{0,n}$ in $L^4(\R)$ and of its derivative in
$L^2(\R)\cap L^4(\R)$ (recall that $U_{0,n}\to U_0$ in
$H^1\cap W^{1,4}$):
\begin{equation}\label{unH1}
	u_{0,n}\to u_0\quad \text{in }\,
	H^1(\R)\cap W^{1,4}(\R),\quad n\to\infty.
\end{equation}
Using that $V_{0,n}\to V_0$
in $H^1(\R)\cap W^{1,4}(\R)$, we obtain similar convergence for $\{v_{0,n}\}_{n=1}^\infty$:
\begin{equation}\label{vnH1}
	v_{0,n}\to v_0\quad \text{in }\,
	H^1(\R)\cap W^{1,4}(\R),\quad n\to\infty.
\end{equation}
Finally, combining \eqref{u0n-k-1} and \eqref{unH1},
we have that 
$$
\begin{aligned}
	\norm{u_{0,n}'v_{0,n}'-u_0'v_0'}_{L^2(\R)}
	&\leq\norm{u_{0,n}'-u_0'}_{L^4(\R)}
	\norm{v_{0,n}'}_{L^4(\R)}\\
	&\quad+\norm{u_0'}_{L^4(\R)}
	\norm{v_{0,n}'-v_0'}_{L^4(\R)}\to0,
\end{aligned}
$$
which yields convergence in $\Sigma$
(recall \eqref{d^2}):
\[
d_{\Sigma}\left((u_{0,n},v_{0,n}),(u_0,v_0)\right)\to0,
\quad n\to\infty.
\]

\textbf{Step 5.2.}
Now we must show that the solutions
\[
\left(\tilde{U}_n,\tilde{V}_n,\tilde{W}_n,\tilde{Z}_n,\tilde{q}_n\right)
\in C\left([-T,T],\mathcal{C}^k_{\Omega_{\mathrm{tcN}}}\right),\quad n\in\N,
\]
corresponding to the initial data $(u_{0,n},v_{0,n})$
defined in \eqref{un-def} with $t=0$,
satisfy \eqref{ins} for all $(t,\xi)\in\Lambda_{T,M}$ (recall Step 3).
Although, in general,
$\left(\tilde{U}_n,\tilde{V}_n,\tilde{W}_n,\tilde{Z}_n,\tilde{q}_n\right)$
does not coincide with
$(U_n,V_n,W_n,Z_n,q_n)$, since
$\tilde{q}(0,\xi)=1$ for all $\xi\in\R$, see \eqref{id},
the two solutions coincide after relabeling \cite{BC07,BC07d,KR25}.
More precisely, define
\begin{equation*}
	\sigma_n(\xi)=\int^{\xi}_0q_n(0,\eta)\,d\eta,
	\quad n\in\N.
\end{equation*}
Then
\begin{equation}\label{Un-tUn}
\begin{split}
	&\left(\tilde{U}_n,\tilde{V}_n,\tilde{W}_n,
	\tilde{Z}_n\right)(t,\sigma_n(\xi))
	=(U_n,V_n,W_n,Z_n)(t,\xi),
	\quad\xi\in\R,\quad n\in\N,\\
	&\tilde{q}_n(t,\sigma_n(\xi))=\frac{q_n(t,\xi)}
	{q_n(0,\xi)},
	\quad\xi\in\R,\quad n\in\N.
\end{split}
\end{equation}

The first identity in \eqref{Un-tUn} yields
\begin{equation*}
	\partial_{\sigma_n}(\tilde W_n)(t,\sigma_n(\xi))
	=\frac{(W_n)_\xi(t,\xi)}{q_n(0,\xi)},\quad
	\partial_{\sigma_n}^2(\tilde W_n)(t,\sigma_n(\xi))
	=\frac{(W_n)_{\xi\xi}(t,\xi)}{q_n^2(0,\xi)}
	-\frac{(W_n)_\xi(t,\xi)(q_n)_\xi(0,\xi)}{q_n^3(0,\xi)},
\end{equation*}
and, since $\sigma_n(\xi)$ is independent of $t$,
$$
(\tilde W_n)_t(t,\sigma_n(\xi))=(W_n)_t(t,\xi).
$$
Thus, at any point
where $W_n=\pi$ and $(W_n)_\xi=0$,
$\partial_{\sigma_n}^2(\tilde W_n)=0$ if and only if 
$(W_n)_{\xi\xi}=0$.
Applying the same computation to $\tilde{Z}_n$, we conclude that
\eqref{ins} holds for
$\left(\tilde{U}_n,\tilde{V}_n,\tilde{W}_n,
\tilde{Z}_n,\tilde{q}_n\right)$
for all $(t,\xi)\in[-T,T]\times\R$.

Since $\sigma_n^\prime(\xi)=q_n(0,\xi)$ and
$(U_n,V_n,W_n,Z_n)(t,\xi)$ satisfies the compatibility conditions
\eqref{ccUV}, it follows from
\eqref{Un-tUn} that $\left(\tilde{U}_n,\tilde{V}_n,\tilde{W}_n,
\tilde{Z}_n\right)(t,\sigma_n)$ satisfies \eqref{ccUV} as well.
Taking into account that $(U_n,V_n,W_n,Z_n)$ satisfies
the ODE system \eqref{ODE}, 
$d\sigma_n=q_n(0,\xi)\,d\xi$, and that the
ODE for $q_t$ in \eqref{ODE} is linear in $q$,
the function $\left(\tilde{U}_n,\tilde{V}_n,\tilde{W}_n,
\tilde{Z}_n\right)$ also satisfies \eqref{ODE}.
Changing the variables $\eta=\sigma_n(\zeta)$ in the integral representation of the initial characteristic $\tilde{y}_{0,n}$
related to $\left(\tilde{U}_n,\tilde{V}_n,\tilde{W}_n,
\tilde{Z}_n\right)$ (cf.~\eqref{y0} and \eqref{y0ode}), we conclude that
\begin{equation*}
\begin{split}
	\tilde{y}_{0,n}(\sigma_n(\xi))&=
	\int_0^{\sigma_n(\xi)}\left(\cos^2\frac{\tilde{W}_n}{2}
	\cos^2\frac{\tilde{Z}_n}{2}\right)(0,\eta)\,d\eta\\
	&=\int_0^{\xi}\left(q_n\cos^2\frac{W_n}{2}
	\cos^2\frac{Z_n}{2}\right)(0,\zeta)\,d\zeta=y_{0,n}(\xi).
\end{split}
\end{equation*}
Finally, recalling that $(U_n,V_n,W_n,Z_n)(0,\xi)$ and 
$\left(\tilde{U}_n,\tilde{V}_n,\tilde{W}_n,
\tilde{Z}_n\right)(0,\sigma_n)$ correspond to the same data
$(u_{0,n},v_{0,n})$ in the Eulerian variables, we have that
the initial data of the solutions $(U_n,V_n,W_n,Z_n)(t,\xi)$ and 
$\left(\tilde{U}_n,\tilde{V}_n,\tilde{W}_n,
\tilde{Z}_n\right)(t,\sigma_n(\xi))$
of the ODE system \eqref{ODE} are the same.

\medskip
\textbf{Step 6.}
In this step, we analyze the level sets where 
$W$ or $Z$ attains the critical value $\pi$ 
and, using the regular value theorem (which ensures that 
these level sets are $C^{k-1}$ one-dimensional submanifolds) 
together with compactness (which guarantees that only 
finitely many connected components can occur), show that each 
such set decomposes into finitely many $C^{k-1}$ one-dimensional 
curves $\mathbf{C}_i^W$ and $\mathbf{C}_j^Z$.
We then verify that these curves satisfy 
properties i)–iv) stated in Theorem\ref{Thmgr}.

Recalling \eqref{GammaWZ}, \eqref{WZl1} 
and that $M\geq\tilde{M}_T$, we conclude 
that $\Gamma^W,\Gamma^Z\subset\Lambda_{T,M}$. 
The conditions \eqref{ins} imply that
for all $(t,\xi)\in\Gamma^W$ we have 
$(W_t,W_\xi)(t,\xi)\neq(0,0)$,
while for all $(t,\xi)\in\Gamma^Z$ we have 
$(Z_t,Z_\xi)(t,\xi)\neq(0,0)$.
Then the regular value theorem (see Theorem \ref{RVT}) 
implies that both $\Gamma^W$ and $\Gamma^Z$ are one-dimensional 
submanifolds of $\Lambda_{T,M}$ of class $C^{k-1}$.
Since $\Lambda_{T,M}$ is a compact set, we conclude that
\begin{equation}\label{gammaWZ}
	\Gamma^W=\bigcup\limits_{i=1}^{N_1}\gamma_i^W,\quad
	\Gamma^Z=\bigcup\limits_{i=1}^{N_2}\gamma_i^Z,\quad
	\text{for some }N_1,\,N_2\in\N\cup\{0\},
\end{equation}
where $\gamma_i^W$, $i=1,\dots,N_1$ and 
$\gamma_j^Z$, $j=1,\dots,N_2$ are $C^{k-1}$ 
curves which satisfy
items ii)--iv) in Theorem \ref{Thmgr} with 
$\gamma_i^W$ and $\gamma_j^Z$ in place of 
$\mathbf{C}_i^W$ and $\mathbf{C}_j^Z$, respectively, for
$i=1,\dots,N_1$ and $j=1,\dots, N_2$.
Concerning item iv) we notice that neither 
$W$ nor $Z$ attains the value $\pi$ on the 
boundary $|\xi|=M$ and $|t|<T$ (recall \eqref{WZl1}).

Define $\mathbf{C}_i^W$, $i=1,\dots,N_1$, and 
$\mathbf{C}_j^Z$, $j=1,\dots,N_2$ as the images of the curves
$\gamma_i^W$, $i=1,\dots,N_1$, and $\gamma_j^Z$, $j=1,\dots,N_1$,
respectively, under the map $(t,\xi)\mapsto(t,y(t,\xi))$ 
(recall \eqref{char1} and \eqref{uvdef}):
\begin{equation*}
	\begin{split}
		&\mathbf{C}_i^W=\left\{(t,y(t,\xi)):
		(t,\xi)\in\gamma_i^W\right\},\quad i=1,\dots,N_1,\\
		&\mathbf{C}_j^Z=\left\{(t,y(t,\xi)):
		(t,\xi)\in\gamma_j^Z\right\},\quad j=1,\dots,N_2.
	\end{split}
\end{equation*}
Recalling \eqref{GammaWZ} and \eqref{gammaWZ}, 
we have established item i) of Theorem \ref{Thmgr}.

We provide a detailed analysis 
of the image of the curve $\gamma_1^W$
under the map $(t,\xi)\mapsto(t,y(t,\xi))$.
The remaining curves can be studied in a similar manner.
Consider an arbitrary point $(\tilde{t},\tilde{\xi})\in\gamma_1^W$.
If $W_\xi(\tilde{t},\tilde{\xi})\neq0$,
then in a neighborhood of $(\tilde{t},\tilde{\xi})$, the
curve $\gamma_1^W$ can be locally represented as 
$\xi=\xi(t)$ by the implicit function theorem.
Consequently, the image $(t,y(t,\xi(t)))$ possesses 
a nonzero tangent vector and is of class $C^{k-1}$.

In the case $W_\xi(\tilde{t},\tilde{\xi})=0$,
the curve $\gamma_1^W$ locally has the form $t=t(\xi)$
(recall that $W_t(\tilde{t},\tilde{\xi})\neq0$, 
see \eqref{ins} and \eqref{GammaWZ}), where 
$t^\prime(\tilde{\xi})=-\frac{W_\xi}{W_t}
(\tilde{t},\tilde{\xi})=0$. 
Thus, the image of $(t(\xi),y(t(\xi),\xi))$ 
has zero tangent vector at $\xi=\tilde{\xi}$ 
and the image of $\gamma_1^W$ loses 
the regularity at this point. 
Using that $(W_\xi,W_{\xi\xi})(\tilde{t},\tilde{\xi})\neq(0,0)$, 
we conclude that $t^{\prime\prime}(\tilde{\xi})
=-\frac{W_{\xi\xi}}{W_t}(\tilde{t},\tilde{\xi})\neq0$. 
Therefore, the critical point $(\tilde{t},\tilde{\xi})
\in\gamma_1^W$ is isolated, i.e., there exists a neighborhood 
$\mathcal{U}_{(\tilde{t},\tilde{\xi})}\subset\gamma_1^W$ 
of the point $(\tilde{t},\tilde{\xi})$ such that 
for all $(t,\xi)\in\mathcal{U}_{(\tilde{t},\tilde{\xi})}
\setminus(\tilde{t},\tilde{\xi})$, we have
$W_\xi(t,\xi)\neq0$. Since $\gamma_1^W$ is a 
compact set, there are a finite number
of points $(\tilde{t},\tilde{\xi})\in\gamma$ such that
$W_\xi(\tilde{t},\tilde{\xi})=0$. 
Thus, the image of $\gamma_1^W$ under the 
map $(t,\xi)\mapsto(t,y(t,\xi))$ is a continuous 
curve which consists of a finite number 
of $C^{k-1}$ curves, i.e., $\mathbf{C}_1^W$ 
is a piecewise $C^{k-1}$ curve.

Arguing in the same way for all $\gamma_i^W$ and $\gamma_j^Z$, 
we conclude that $\mathbf{C}_i^W$ and $\mathbf{C}_j^Z$ are 
piecewise $C^{k-1}$ curves, $i=1,\dots,N_1$, $j=1,\dots,N_2$.
Moreover, since $\gamma_i^W$ and $\gamma_j^Z$ satisfy properties 
ii)--iv) in Theorem \ref{Thmgr} and $y(t,\xi)$ is 
monotone (see \eqref{pxiy}), we conclude that 
$\mathbf{C}_i^W$ and $\mathbf{C}_j^Z$ 
also satisfy the same properties (see Figure \ref{char-curves} 
for an illustration). Consequently, items i)--iv) 
of Theorem \ref{Thmgr} are established.

\medskip
\textbf{Step 7.}
Now we are at the position to show that for all 
$(u_0,v_0)\in\tilde{\mathcal{M}}_T$ 
the corresponding solution $(u,v)$
satisfies properties (1)--(3) described 
in Theorem \ref{Thmgr}. 
Step 2 above implies that it is sufficient
to study the regularity properties of 
$u$ and $v$ in the rectangle $Q_{T,\delta}$.

Recall from step 3 that 
$Q_{T,\delta}\subset\tilde{G}_{(u_0,v_0)}$ for all 
$(u_0,v_0)\in\mathcal{U}_{(\tilde{u}_0,\tilde{v}_0)}$ 
(see \eqref{tGu0v0}). Then for all $(t,\xi)\in\Lambda_{T,M}
\setminus\left(\Gamma^W\cup\Gamma^Z\right)$
we have $W(t,\xi)\neq\pi$ and $Z(t,\xi)\neq\pi$.
Combining \eqref{uvdef}, \eqref{pxiy} and \eqref{gammaWZ}, we 
conclude that $u$ and $v$ are 
$C^k$ in a small neighborhood of $(t,y(t,\xi))$ for all such 
$(t,\xi)$ and therefore item (1) is established.
Then, using \eqref{pxuv}, we arrive at item (2).

Recalling that the critical points 
$(t,\xi)\in\Gamma^W$ with $W_\xi(t,\xi)=0$ 
and $(t,\xi)\in\Gamma^Z$ with $Z_\xi(t,\xi)=0$ 
are isolated  (see the discussion in Step~6 above), 
we infer from \eqref{pxiy} that, 
for each fixed $t$, the characteristic 
map $y(t,\cdot)$ is strictly monotone.  
Hence we may change variables via $x=y(t,\xi)$ 
in \eqref{mut} and, using \eqref{WZ} 
and \eqref{q}, obtain \eqref{mu-t-ac}.

\begin{remark}[Step 5]\label{R-dens}
		Our proof of the density of $\tilde{\mathcal{M}}_T$
		differs from those in \cite[p.\,350]{BC17},
	\cite[p.\,2932]{LZ17}, \cite[p.\,1117]{CCCS18}.
	Taking into account that the nonlocal peakon equations
	do not, in general, have finite propagation speed \cite{C05},
		solutions corresponding to different initial data that
		coincide on a compact subset of $\R$
	need not coincide on the domain
	$[-T,T]\times[-R,R]$ in the $(t,x)$-plane
	\cite{LP20}, cf.~\cite{CCCS18}, \cite{LZ17}.
	Here, applying Lemma \ref{L3},
	we define the sequence $(u_{0,n},v_{0,n})$,
		which converges to any prescribed $(u_0,v_0)$, directly from
	the corresponding solutions of the associated ODE system,
	see \eqref{un-def}--\eqref{yny0n}.
	Then, taking advantage of the convergence of these solutions
	in $C\left([-T,T],\mathcal{C}^k_{\Omega_{\mathrm{tcN}}}\right)$,
	see item (i) in Step 5,
	we eventually establish \eqref{conv-un} in Step 5.1
	as well as \eqref{ins} for the corresponding sequence
	of solutions of the ODE system in Step 5.2.
	Note that since $\px^k u_0$ and $\px^k v_0$ are
		continuous but need not be uniformly continuous on $\R$,
	\eqref{conv-un} involves convergence in $C^k$
		on compact sets only; cf.~\eqref{u0n-est}--\eqref{u0n-k}.
\end{remark}
\subsection{Proof of Theorem \ref{Thmchc}}
\label{prTh2}
In what follows we must calculate
every $\xi$-derivative of the three
basic functions $U$, $V$, and $y$, see 
\eqref{ccUV} and \eqref{pxiy}, up to ninth order, i.e.,
$3\times 9=27$ distinct derivatives.
Even a single one of these is already
unwieldy: after $n$ differentiations 
a triple product such as
$q\cos^2\frac{W}{2}\cos^2\frac{Z}{2}$ splits into
$\binom{n+2}{2}$ multinomial triples, so 
for $n=8$ one begins with
$45$ top-level terms.
Every factor inside those triples is itself a
composite of $W$ or $Z$,
and must next be expanded by the chain
rule. The $m^{\text{th}}$ outer derivative 
of a single-variable
function contributes as many monomials as 
there are ordered partitions
of $m$ (i.e., the $m$-th Bell number);
for $m=8$ this number is $4140$.
Thus, in the worst case, the
ninth-order derivative of $y(t,\xi)$ (for example) 
contains the following number of elementary monomials:
$$
\sum\limits_{m=0}^8\sum\limits_{i=0}^m B_iB_{m-i}=14924,\quad
\text{where }B_i\text{ is the $i$-th Bell number.}
$$
On the jump curves
$W=\pi$ and $Z=\pi$, however, the identities $\sin\pi=0$ and
$\cos\frac{\pi}{2}=0$ annihilate almost 
all of these terms; for instance
$\pxi^9y$ (with $W_\xi=Z_\xi=0$ as required by the compatibility
conditions along the curve) collapses to a single 
surviving monomial, see \eqref{Tayy2} below. Because the 
unreduced formulas would fill many pages while conveying no
additional insight, we record 
only the values taken at the special
angles and omit the intermediate algebra.
\medskip

\textbf{Case (1)}.
Using that $W(t_1,\xi_1)=\pi$, we obtain 
the following from \eqref{ccUV} and \eqref{pxiy}:
\begin{equation}\label{TayUV1}
	\begin{split}
		&U_\xi(t_1,\xi_1)=V_\xi(t_1,\xi_1)=
		V_{\xi\xi}(t_1,\xi_1)=y_\xi(t_1,\xi_1)
		=y_{\xi\xi}(t_1,\xi_1)=0,\\
		&U_{\xi\xi}(t_1,\xi_1)
		=-\frac{1}{2}\left(qW_\xi\cos^2\frac{Z}{2}\right)(t_1,\xi_1),
		\quad
		V_{\xi\xi\xi}(t_1,\xi_1)
		=\frac{1}{4}\left(qW_\xi^2\sin Z\right)(t_1,\xi_1),
	\end{split}
\end{equation}
as well as
\begin{equation}\label{Tayy1}
	y_{\xi\xi\xi}(t_1,\xi_1)
	=\frac{1}{2}\left(qW_\xi^2\cos^2\frac{Z}{2}\right)(t_1,\xi_1).
\end{equation}
Notice that since $W_\xi(t_1,\xi_1)\neq0$ and
$Z(t_1,\xi_1)\neq\pi$, we have $U_{\xi\xi}(t_1,\xi_1)\neq0$ and
$y_{\xi\xi\xi}(t_1,\xi_1)\neq0$ (recall \eqref{perWZ}).
From \eqref{TayUV1} we conclude that 
the Taylor expansions for $U$ and $V$ read as follows:
\begin{equation}\label{TayUV1e}
	\begin{split}
		&U(t,\xi)=U(t_1,\xi_1)+\tilde{a}_{1,1}(\xi-\xi_1)^2
		+\tilde{a}_{1,2}(t-t_1)
		+\mathcal{O}\left((\xi-\xi_1)^3\right)
		+\so(t-t_1),\\
		&V(t,\xi)=V(t_1,\xi_1)+\tilde{b}_{1,1}(\xi-\xi_1)^3
		+\tilde{b}_{1,2}(t-t_1)
		+\mathcal{O}\left((\xi-\xi_1)^4\right)
		+\so(t-t_1),
	\end{split}
\end{equation}
with some $\tilde{a}_{1,1}\in\R\setminus\{0\}$ and
$\tilde{a}_{1,2},\tilde{b}_{1,1},\tilde{b}_{1,2}\in\R$.
Notice that since $U_t$ and/or $V_t$ can vanish at the point
$(t_1,\xi_1)$, the constants $\tilde{a}_{1,2}$ and/or
$\tilde{b}_{1,2}$ can be zero in \eqref{TayUV1e}.
Recalling that $x_1=y(t_1,\xi_1)$, we 
obtain from \eqref{char1}, \eqref{TayUV1}, 
and \eqref{Tayy1} that
\begin{equation}\label{Tayy1e}
	x=x_1+\tilde{c}_{1,1}(\xi-\xi_1)^3
	+(uv)(t_1,x_1)(t-t_1)
	+\mathcal{O}\left((\xi-\xi_1)^4\right)
	+\so(t-t_1),\quad
	x=y(t,\xi),
\end{equation}
with $\tilde{c}_{1,1}\in\R\setminus\{0\}$.
Observe that $y_t(t_1,\xi_1)=(uv)(t_1,x_1)$ 
can be zero in \eqref{Tayy1e}.
Combining \eqref{elld} and \eqref{Tayy1e}, 
we conclude that
\begin{equation*}
	\xi-\xi_1=\hat{c}_{1,1}l_1^{1/3}(t,x), \quad
	\hat{c}_{1,1}=\tilde{c}_{1,1}^{-1/3}\neq0,
\end{equation*}
which, together with \eqref{uvdef} 
and \eqref{TayUV1e}, yields \eqref{uvch1}.
\medskip

\textbf{Case (2)}. 
The proof follows similar steps as in Case (1).
Here we use that (cf.\,\,\eqref{TayUV1}--\eqref{Tayy1})
\begin{equation*}
	\begin{split}
		&U_\xi(t_2,\xi_2)=U_{\xi\xi}(t_2,\xi_2)=
		V_{\xi}(t_2,\xi_2)=y_\xi(t_2,\xi_2)
		=y_{\xi\xi}(t_2,\xi_2)=0,\\
		&U_{\xi\xi\xi}(t_2,\xi_2)
		=\frac{1}{4}\left(qZ_\xi^2\sin W\right)(t_2,\xi_2),
		\quad V_{\xi\xi}(t_2,\xi_2)
		=-\frac{1}{2}\left(qZ_\xi
		\cos^2\frac{W}{2}\right)(t_2,\xi_2)\\
		&y_{\xi\xi\xi}(t_2,\xi_2)
		=\frac{1}{2}\left(qZ_\xi^2
		\cos^2\frac{W}{2}\right)(t_2,\xi_2),
	\end{split}
\end{equation*}
which imply (cf.\,\,\eqref{TayUV1e})
\begin{equation*}
	\begin{split}
		&U(t,\xi)=U(t_2,\xi_2)+\tilde{a}_{2,1}(\xi-\xi_2)^3
		+\tilde{a}_{2,2}(t-t_2)
		+\mathcal{O}\left((\xi-\xi_2)^4\right)
		+\so(t-t_2),\\
		&V(t,\xi)=V(t_2,\xi_2)+\tilde{b}_{2,1}(\xi-\xi_2)^2
		+\tilde{b}_{2,2}(t-t_2)
		+\mathcal{O}\left((\xi-\xi_2)^3\right)
		+\so(t-t_2),
	\end{split}
\end{equation*}
with some $\tilde{b}_{2,1}\in\R\setminus\{0\}$ and
$\tilde{a}_{2,1},\tilde{a}_{2,2},\tilde{b}_{2,2}\in\R$.
Using similar expansion as in \eqref{Tayy1e}, we arrive at 
\eqref{uvch2}
\medskip

\textbf{Case (3)}.
Taking into account that $W(t_3,\xi_3)=Z(t_3,\xi_3)=\pi$, 
we obtain from \eqref{ccUV} and \eqref{pxiy} that
\begin{equation}\label{TayUV3}
	\begin{split}
		&\pxi^iU(t_3,\xi_3)=\pxi^iV(t_3,\xi_3)=0,\quad i=1,2,3,\\
		&\pxi^4U(t_3,\xi_3)
		=-\frac{3}{4}\left(qW_\xi Z^2_{\xi}\right)(t_3,\xi_3),\quad
		\pxi^4V(t_3,\xi_3)
		=-\frac{3}{4}\left(qW^2_\xi Z_{\xi}\right)(t_3,\xi_3),
	\end{split}
\end{equation}
and
\begin{equation}\label{Tayy3}
	\pxi^i y(t_3,\xi_3)=0,\quad i=1,2,3,4,\quad
	\pxi^5 y(t_3,\xi_3)=\frac{3}{2}
	\left(qW_\xi^2Z_\xi^2\right)(t_3,\xi_3).
\end{equation}
Using that $W_\xi(t_3,\xi_3)\neq0$ 
and $Z_\xi(t_3,\xi_3)\neq0$, 
we have from \eqref{TayUV3} and \eqref{Tayy3} that
\begin{equation*}
	\begin{split}
		&U(t,\xi)=U(t_3,\xi_3)+\tilde{a}_{3,1}(\xi-\xi_3)^4
		+\tilde{a}_{3,2}(t-t_3)
		+\mathcal{O}\left((\xi-\xi_3)^5\right)
		+\so(t-t_3),\\
		&V(t,\xi)=V(t_3,\xi_3)+\tilde{b}_{3,1}(\xi-\xi_3)^4
		+\tilde{b}_{3,2}(t-t_3)
		+\mathcal{O}\left((\xi-\xi_3)^5\right)
		+\so(t-t_3),
	\end{split}
\end{equation*}
and
\begin{equation*}
	x=x_3+\tilde{c}_{3,1}(\xi-\xi_3)^5
	+(uv)(t_3,x_3)(t-t_3)
	+\mathcal{O}\left((\xi-\xi_3)^6\right)
	+\so(t-t_3),\quad
	x=y(t,\xi),
\end{equation*}
with some $\tilde{a}_{3,1},\tilde{b}_{3,1},\tilde{c}_{3,1}
\in\R\setminus\{0\}$ and
$\tilde{a}_{3,2},\tilde{b}_{3,2}\in\R$.
Then arguing as in Case (1), we arrive at \eqref{uvch3}.
\medskip

\textbf{Case (4)}.
Using that $W(t_4,\xi_4)=\pi$, $W_\xi(t_4,\xi_4)=0$, 
as well as \eqref{ccUV} and \eqref{pxiy}, we obtain
\begin{equation}\label{TayUV4}
	\begin{split}
		&\pxi^i U(t_4,\xi_4)=0,\quad i=1,2,\quad
		\pxi^j V(t_4,\xi_4)=0,\quad j=1,2,3,4,\\
		&U_{\xi\xi\xi}(t_4,\xi_4)
		=-\frac{1}{2}\left(qW_{\xi\xi}
		\cos^2\frac{Z}{2}\right)
		(t_4,\xi_4),\quad
		\pxi^5 V(t_4,\xi_4)=\frac{3}{4}
		\left(qW_{\xi\xi}^2\sin Z\right)(t_4,\xi_4),
	\end{split}
\end{equation}
and
\begin{equation}\label{Tayy4}
	\pxi^i y(t_4,\xi_4)=0,
	\quad i=1,2,3,4,\quad
	\pxi^5 y(t_4,\xi_4)=\frac{3}{2}
	\left(qW_{\xi\xi}^2\cos^2\frac{Z}{2}
	\right)(t_4,\xi_4).
\end{equation}
Taking into account that $Z(t_4,\xi_4)\neq\pi$ and
$W_{\xi\xi}(t_4,\xi_4)\neq0$ (see \eqref{sder}), 
we conclude from \eqref{TayUV4} and \eqref{Tayy4} that
\begin{equation}\label{TayUV4e}
	\begin{split}
		&U(t,\xi)=U(t_4,\xi_4)+\tilde{a}_{4,1}(\xi-\xi_4)^3
		+\tilde{a}_{4,2}(t-t_4)
		+\mathcal{O}\left((\xi-\xi_4)^4\right)
		+\so(t-t_4),\\
		&V(t,\xi)=V(t_4,\xi_4)+\tilde{b}_{4,1}(\xi-\xi_4)^5
		+\tilde{b}_{4,2}(t-t_4)
		+\mathcal{O}\left((\xi-\xi_4)^6\right)
		+\so(t-t_4),
	\end{split}
\end{equation}
and
\begin{equation}\label{Tayy4e}
	x=x_4+\tilde{c}_{4,1}(\xi-\xi_4)^5
	+(uv)(t_4,x_4)(t-t_4)
	+\mathcal{O}\left((\xi-\xi_4)^6\right)
	+\so(t-t_4),\quad
	x=y(t,\xi),
\end{equation}
with some $\tilde{a}_{4,1},\tilde{c}_{4,1}
\in\R\setminus\{0\}$ and
$\tilde{a}_{4,2},\tilde{b}_{4,1},
\tilde{b}_{4,2}\in\R$.
Then arguing as in Case (1), 
we obtain \eqref{uvch4}.
\medskip

\textbf{Case (5)}.
The proof follows a line of reasoning 
analogous to that used in Case (4) above.
Here we have (cf.\,\,\eqref{TayUV4} and \eqref{Tayy4})
\begin{equation*}
\begin{split}
	&\pxi^i U(t_5,\xi_5)=0,\quad i=1,2,3,4,\quad
	\pxi^j V(t_5,\xi_5)=0,\quad j=1,2,\\
	&\pxi^5 U(t_5,\xi_5)=\frac{3}{4}
	\left(qZ_{\xi\xi}^2\sin W\right)(t_5,\xi_5),\quad
	V_{\xi\xi\xi}(t_5,\xi_5)
	=-\frac{1}{2}\left(qZ_{\xi\xi}\cos^2\frac{W}{2}\right)
	(t_5,\xi_5),\\
	&\pxi^i y(t_5,\xi_5)=0,\quad i=1,2,3,4,\quad
	\pxi^5 y(t_5,\xi_5)=\frac{3}{2}
	\left(qZ_{\xi\xi}^2\cos^2\frac{W}{2}\right)(t_5,\xi_5),
	\end{split}
\end{equation*}
which yield
(cf.\,\,\eqref{TayUV4e})
\begin{equation*}
	\begin{split}
		&U(t,\xi)=U(t_5,\xi_5)+\tilde{a}_{5,1}(\xi-\xi_5)^5
		+\tilde{a}_{5,2}(t-t_5)
		+\mathcal{O}\left((\xi-\xi_5)^6\right)
		+\so(t-t_5),\\
		&V(t,\xi)=V(t_5,\xi_5)+\tilde{b}_{5,1}(\xi-\xi_5)^3
		+\tilde{b}_{5,2}(t-t_5)
		+\mathcal{O}\left((\xi-\xi_5)^4\right)
		+\so(t-t_5),
	\end{split}
\end{equation*}
with $\tilde{b}_{5,1}\in\R\setminus\{0\}$ 
and $\tilde{a}_{5,1},\tilde{a}_{5,2},\tilde{b}_{5,2}\in\R$. 
Using an expression analogous 
to \eqref{Tayy4e}, we 
obtain \eqref{uvch5}.
\medskip

\textbf{Case (6)}.
Using that $W(t_6,\xi_6)=Z(t_6,\xi_6)=\pi$, 
$W_\xi(t_6,\xi_6)=0$, we obtain the following expressions:
\begin{equation}\label{TayUV6}
	\begin{split}
		&\pxi^i U(t_6,\xi_6)=0,\quad i=1,\dots,4\quad
		\pxi^j V(t_6,\xi_6)=0,\quad j=1,\dots,5,\\
		&\pxi^5U(t_6,\xi_6)
		=-\frac{3}{2}\left(qW_{\xi\xi}Z_\xi^2\right)
		(t_6,\xi_6),\quad
		\pxi^6 V(t_6,\xi_6)=-\frac{15}{4}
		\left(qW_{\xi\xi}^2Z_\xi\right)(t_6,\xi_6),
	\end{split}
\end{equation}
and
\begin{equation}\label{Tayy6}
	\pxi^i y(t_6,\xi_6)=0,\quad i=1,\dots,6,\quad
	\pxi^7y(t_6,\xi_6)=\frac{45}{4}
	\left(qW^2_{\xi\xi}Z^2_\xi\right)(t_6,\xi_6).
\end{equation}
Recalling that $Z_\xi(t_6,\xi_6)\neq0$ and
$W_{\xi\xi}(t_6,\xi_6)\neq0$ (see \eqref{sder}), 
we conclude that $\pxi^5U(t_6,\xi_6)\neq0$, 
$\pxi^6V(t_6,\xi_6)\neq0$, and
$\pxi^7y(t_6,\xi_6)\neq0$.
Using the following Taylor expansions 
of $U$, $V$, and $y$
\begin{equation}\label{TayUV6e}
	\begin{split}
		&U(t,\xi)=U(t_6,\xi_6)+\tilde{a}_{6,1}(\xi-\xi_6)^5
		+\tilde{a}_{6,2}(t-t_6)
		+\mathcal{O}\left((\xi-\xi_6)^6\right)
		+\so(t-t_6),\\
		&V(t,\xi)=V(t_6,\xi_6)+\tilde{b}_{6,1}(\xi-\xi_6)^6
		+\tilde{b}_{6,2}(t-t_6)
		+\mathcal{O}\left((\xi-\xi_6)^7\right)
		+\so(t-t_6),
	\end{split}
\end{equation}
and
\begin{equation}\label{Tayy6e}
	x=x_6+\tilde{c}_{6,1}(\xi-\xi_6)^7
	+(uv)(t_6,x_6)(t-t_6)
	+\mathcal{O}\left((\xi-\xi_6)^8\right)
	+\so(t-t_6),\quad
	x=y(t,\xi),
\end{equation}
with some $\tilde{a}_{6,1},\tilde{b}_{6,1},\tilde{c}_{6,1}
\in\R\setminus\{0\}$ and
$\tilde{a}_{6,2},\tilde{b}_{6,2}\in\R$, 
we arrive at \eqref{uvch6}.
\medskip

\textbf{Case (7)}.
The proof proceeds similarly to Case (6). 
We have the following expressions 
(cf.\,\,\eqref{TayUV6} and \eqref{Tayy6}):
\begin{equation*}
	\begin{split}
		&\pxi^i U(t_7,\xi_7)=0,\quad i=1,\dots,5\quad
		\pxi^j V(t_7,\xi_7)=0,\quad j=1,\dots,4,\\
		&\pxi^6 U(t_7,\xi_7)=-\frac{15}{4}
		\left(qW_\xi Z_{\xi\xi}^2\right)(t_7,\xi_7),\quad
		\pxi^5V(t_7,\xi_7)
		=-\frac{3}{2}\left(qW_\xi^2 Z_{\xi\xi}\right)
		(t_7,\xi_7),\\
		&	\pxi^i y(t_7,\xi_7)=0,\quad i=1,\dots,6,\quad
		\pxi^7y(t_7,\xi_7)=\frac{45}{4}
		\left(qW^2_\xi Z^2_{\xi\xi}\right)(t_7,\xi_7),
	\end{split}
\end{equation*}
which imply (cf.\,\,\eqref{TayUV6e})
\begin{equation*}
	\begin{split}
		&U(t,\xi)=U(t_7,\xi_7)+\tilde{a}_{7,1}(\xi-\xi_7)^6
		+\tilde{a}_{7,2}(t-t_7)
		+\mathcal{O}\left((\xi-\xi_7)^7\right)
		+\so(t-t_7),\\
		&V(t,\xi)=V(t_7,\xi_7)+\tilde{b}_{7,1}(\xi-\xi_7)^5
		+\tilde{b}_{7,2}(t-t_7)
		+\mathcal{O}\left((\xi-\xi_7)^6\right)
		+\so(t-t_7),
	\end{split}
\end{equation*}
with $\tilde{a}_{7,1},\tilde{b}_{7,1}
\in\R\setminus\{0\}$ and
$\tilde{a}_{7,2},\tilde{b}_{7,2}\in\R$.
Using the Taylor expansion as 
in \eqref{Tayy6e}, we obtain \eqref{uvch7}.
\medskip

\textbf{Case (8)}.
Recalling that $W(t_8,\xi_8)=Z(t_8,\xi_8)=\pi$ and
$W_\xi(t_8,\xi_8)=Z_\xi(t_8,\xi_8)=0$, we obtain
\begin{equation}\label{TayUV2}
	\begin{split}
		&\pxi^i U(t_8,\xi_8)=0,\quad i=1,\dots,6\quad
		\pxi^j V(t_8,\xi_8)=0,\quad j=1,\dots,6,\\
		&\pxi^7U(t_8,\xi_8)
		=-\frac{45}{4}\left(qW_{\xi\xi}Z_{\xi\xi}^2\right)
		(t_8,\xi_8),\quad
		\pxi^7 V(t_8,\xi_8)=-\frac{45}{4}
		\left(qW_{\xi\xi}^2Z_{\xi\xi}\right)(t_8,\xi_8),
	\end{split}	
\end{equation}
and
\begin{equation}\label{Tayy2}
	\pxi^i y(t_8,\xi_8)=0,\quad i=1,\dots,8,\quad
	\pxi^9y(t_8,\xi_8)=\frac{315}{2}
	\left(qW_{\xi\xi}^2Z_{\xi\xi}^2\right)(t_8,\xi_8).
\end{equation}
Equations \eqref{TayUV2} and \eqref{Tayy2} 
imply that $U$, $V$ and $y$ have the 
following Taylor expansions as $\xi\to\xi_8$:
\begin{equation}\label{TayUV8e}
	\begin{split}
		&U(t,\xi)=U(t_8,\xi_8)
		+\tilde{a}_{8,1}(\xi-\xi_8)^7
		+\tilde{a}_{8,2}(t-t_8)
		+\mathcal{O}\left((\xi-\xi_8)^8\right)
		+\so(t-t_8),\\
		&V(t,\xi)=V(t_8,\xi_8)
		+\tilde{b}_{8,1}(\xi-\xi_8)^7
		+\tilde{b}_{8,2}(t-t_8)
		+\mathcal{O}\left((\xi-\xi_8)^8\right)
		+\so(t-t_8),
	\end{split}
\end{equation}
and
\begin{equation}
	\label{Tayy8e}
	x=x_8+\tilde{c}_{8,1}(\xi-\xi_8)^9
	+(uv)(t_8,x_8)(t-t_8)
	+\mathcal{O}\left((\xi-\xi_8)^{10}\right)
	+\so(t-t_8),\quad
	x=y(t,\xi),
\end{equation}
with some $\tilde{a}_{8,1},
\tilde{b}_{8,1},\tilde{c}_{8,1}
\in\R\setminus\{0\}$ and
$\tilde{a}_{8,2},\tilde{b}_{8,2}\in\R$.
Thus, \eqref{TayUV8e} and 
\eqref{Tayy8e} imply \eqref{uvch8}.
\hfill $\square$

\begin{remark}\label{multyz}
The proof also shows that
for $(t,\xi)\in\left(\Gamma^W\cap\Gamma^Z\right)$,
see \eqref{GammaWZ}, there exists $i\in\{2,\dots,9\}$
such that $\pxi^i y(t,\xi)\neq0$ provided that
$W$ and $Z$ satisfy \eqref{ins}.
\end{remark}

\section{Lipschitz metric for global solutions}\label{LM}

In this section, we introduce a metric 
$d_{\mathcal{D},K}(\cdot,\cdot)$ on $\mathcal{D}$ 
under which the global conservative solutions of the 
two-component Novikov system, as stated in 
Theorem~\ref{Thm}, satisfy a Lipschitz 
continuity property. Given any 
$\mathbf{u},\hat{\mathbf{u}}\in\mathcal{D}$ 
satisfying $u,v,\hat{u},\hat{v}\in W^{1,1}(\R)$,
we define $d_{\mathcal{D},K}(\mathbf{u},\hat{\mathbf{u}})$ 
as the geodesic distance between the corresponding 
functions $\mathbf{U}^0,\mathbf{U}^1\in
\Omega_{\mathrm{tcN}}^1$ obtained 
via the direct transform \eqref{id2} in 
the Bressan-Constantin variables 
(see \eqref{OmL1tcN}, \eqref{u-bf}, and \eqref{U-bf}).  
The geodesic distance between $\mathbf{U}^0$ 
and $\mathbf{U}^1$ is the infimum of the 
lengths of all paths connecting them
and satisfying the compatibility conditions \eqref{ccUV}
(see Definition~\ref{Gdist} below).  
We emphasize that the metric $d_{\mathcal{D},K}(\cdot,\cdot)$ 
is constructed in the transformed variables, where 
all potential singularities of solutions 
to the Novikov system are resolved.

We next introduce a suitable definition 
of the length of a path $\mathbf{U}^\theta$, 
$\theta\in[0,1]$, that satisfies a Lipschitz 
continuity property under the ODE system \eqref{ODE}.  
Following \cite{BC17A,CCCS18}, we first 
derive an appropriate norm of the tangent vectors 
for paths in the Eulerian variables $(u,v)$ 
(see \eqref{rstr} and Theorem~\ref{Lipn}).  
Rewriting this norm in the transformed 
variables (see \eqref{Rtr}), we obtain 
a Lipschitz estimate for the length 
$\|\mathbf{U}^\theta\|_{\mathcal{L}}$ 
(Definition~\ref{lPR}) along regular paths 
for the ODE system \eqref{ODE} 
(see Definition~\ref{PrpT} and Theorem~\ref{egammat}).  
Since any path can be approximated by 
regular paths under \eqref{ODE} 
(Theorem~\ref{dpath}), a completion argument 
shows that the resulting geodesic metric 
on $\Omega_{\mathrm{tcN}}^1$, and therefore the metric 
$d_{\mathcal{D},K}(\cdot,\cdot)$ on $\mathcal{D}$, 
satisfies the desired Lipschitz property.
We emphasize that the additional $L^1(\R)$ decay assumptions on
solutions in $\Omega_{\mathrm{tcN}}^1$
are used here to justify the completion argument rigorously
in Proposition \ref{COm} (see also Remark \ref{compL1}).

\subsection{Paths of solutions}\label{Sp-ins}

In Theorem \ref{dpath}, we show that 
any sufficiently regular path 
$\mathbf{U}_0^\theta(\xi)
=(U_0^\theta,V_0^\theta,W_0^\theta,Z_0^\theta,q_0^\theta)(\xi)$, 
$\theta\in[0,1]$, can be approximated 
by another regular path $\widehat{\mathbf{U}}_0^\theta(\xi)$ 
whose endpoints remain fixed
and evolution under the ODE system 
\eqref{ODE} satisfies the following non-degeneracy 
conditions for all $(t,\xi)\in[-T,T]\times\R$ and 
for all $\theta\in(0,1)$ (cf.~\eqref{ins}):
\begin{equation}\label{p-ins}
\begin{split}
	&\left(W^\theta,W^\theta_\xi,
	W^\theta_{\xi\xi}, W^\theta_{\xi\xi\xi}\right)
	(t,\xi)\neq(\pi,0,0,0),\\
	&\left(Z^\theta,Z^\theta_\xi,
	Z^\theta_{\xi\xi}, Z^\theta_{\xi\xi\xi}\right)
	(t,\xi)\neq(\pi,0,0,0).
\end{split}
\end{equation}
Because the endpoints are fixed, we cannot guarantee \eqref{p-ins}
for $\theta=0,1$.
We then show that the evolution of such piecewise
regular approximations
$\widehat{\mathbf{U}}_0^\theta$ under \eqref{ODE} 
enjoys the Lipschitz property 
with respect to the associated notion 
of path length (see Theorem \ref{egammat}), 
which explains the central role of such approximations. 
The proof of Theorem \ref{egammat} relies essentially
on the fact that there are only finitely many points $\xi\in\R$
at which $W^\theta(\xi)=\pi$ or $Z^\theta(\xi)=\pi$,
as guaranteed by \eqref{p-ins}.
Note that we do not need the $t$-derivatives here,
which are crucial in \eqref{ins} for establishing the regular characteristic curves for the generic regularity results.
In Definition \ref{PrpT} below, we 
refer to these paths as piecewise regular paths under 
the ODE system \eqref{ODE}.

The proof of Theorem \ref{dpath} proceeds 
along lines similar to those of
Lemma \ref{L3}. The major difference 
is that here we consider the
four-dimensional maps
$\left(W^\theta,W^\theta_\xi,W^\theta_{\xi\xi},
W^\theta_{\xi\xi\xi}\right)(t,\xi)$
and
$\left(Z^\theta,Z^\theta_\xi,Z^\theta_{\xi\xi},
Z^\theta_{\xi\xi\xi}\right)(t,\xi)$,
which depend on the three independent variables 
$(t,\xi,\theta)$, not only on $(t,\xi)$ as 
in Lemma \ref{L3}. Since the domain is three-dimensional, while the
codomain is
four-dimensional, transversality 
here implies that the preimage 
of $\{(\pi,0,0,0)\}$ for the generic maps is empty,
cf.~Lemma \ref{L3}.

The proof of Theorem \ref{dpath} is 
inspired by \cite[Theorem~2]{BC17}, where
Bressan and Chen established a similar
result for a quasilinear second order wave equation.
In that work, however, the authors consider three-dimensional
maps, analogous to \eqref{ins} for our problem.
Therefore, the perturbed path (i) has, in general, different endpoints,
and (ii) satisfies \eqref{ins} for all except possibly
a finite subset of $[0,1]$.
In view of (i), one can define only regular, not piecewise regular,
paths under the ODE system
\eqref{ODE}, cf.~Definition \ref{PrpT}.
Also, in the proof of Theorem \ref{LRO},
one would need to add connector paths from the original
$\mathbf{U}^i$, $i=0,1$,
to the regularized endpoints, which requires an explicit
short-path estimate in the constructed geodesic metric.
The latter is somewhat involved and can be obtained
by considering the ``perturbed'' straight-line path
constructed in step 2 of the proof of Proposition \ref{COm},
which satisfies the compatibility conditions \eqref{ccUV} for all
$\theta\in[0,1]$.

\begin{theorem}\label{dpath}
Consider a path
$\mathbf{U}_0^\theta\in\mathcal{C}^{k}_{\Omega_{\mathrm{tcN}}^1}$, $\theta\in[0,1]$ 
(see \eqref{C1ktcN} and \eqref{U-bf}), such that
\begin{equation}\label{id-path}
	\left(\mathbf{U}_0^\theta,\frac{d}{d\theta}
	\mathbf{U}_0^\theta
	\right)\in \mathcal{P}^{k,1}_{\mathrm{tcN}}[0,1],\quad
	k\geq 6,
\end{equation}
where 
$\mathcal{P}^{k,1}_{\mathrm{tcN}}[0,1]$
is defined in \eqref{P-pth}.
Then for any $\ve>0$ there exists a perturbed path
\[\left(\widehat{\mathbf{U}}_0^\theta,\frac{d}{d\theta}
\widehat{\mathbf{U}}_0^\theta\right)\in
\mathcal{P}^{k,1}_{\mathrm{tcN}}[0,1],
\]
such that
\begin{enumerate}
	\item the endpoints remain fixed:
	$\left.\left(\widehat{\mathbf{U}}_0^\theta,\frac{d}{d\theta}
	\widehat{\mathbf{U}}_0^\theta\right)\right|_{\theta=0,1}=
	\left.\left(\mathbf{U}_0^\theta,\frac{d}{d\theta}
	\mathbf{U}_0^\theta
	\right)\right|_{\theta=0,1}$;
	
	\item 
	$\left\|\left(\widehat{\mathbf{U}}_0^\theta-
	\mathbf{U}_0^\theta,
	\frac{d}{d\theta}\bigl(\widehat{\mathbf{U}}_0^\theta-
	\mathbf{U}_0^\theta\bigr)
	\right)\right\|_{\mathcal{P}^{k,1}_{\mathrm{tcN}}[0,1]}
	<\ve$;

	\item for the global conservative solutions
	(see Theorem \ref{Ckgl})
	\[
	\mathbf{U}^\theta(t,\xi)
	=\left(U^\theta,V^\theta,W^\theta,Z^\theta,q^\theta
	\right)(t,\xi) \quad\text{and}\quad
	\widehat{\mathbf{U}}^\theta(t,\xi)
	=\left(\widehat{U}^\theta,\widehat{V}^\theta,
	\widehat{W}^\theta,\widehat{Z}^\theta,
	\widehat{q}^{\,\theta}\right)(t,\xi),
	\]
	of the ODE system \eqref{ODE} subject 
	to the initial data
	$\mathbf{U}^\theta_0$ and 
	$\widehat{\mathbf{U}}^\theta_0$, 
	respectively, we have 
	(here we drop the arguments $(t,\xi)$ of the 
	solutions and write $\mathbf{U}^\theta$ and 
	$\widehat{\mathbf{U}}^\theta$, respectively)
	\begin{enumerate}
		\item $\left(
		\mathbf{U}^\theta,
		\frac{d}{d\theta}
		\left(\mathbf{U}^\theta\right)
		\right),
		\left(
		\widehat{\mathbf{U}}^\theta,
		\frac{d}{d\theta}
		\left(\widehat{\mathbf{U}}^\theta\right)
		\right)\in
		C\left([-T,T],\mathcal{P}^{k,1}_{\mathrm{tcN}}[0,1]\right)$,
		for any $T>0$;
		
		\item $\left\|\left(\widehat{\mathbf{U}}^\theta-
		\mathbf{U}^\theta,
		\frac{d}{d\theta}\left(\widehat{\mathbf{U}}^\theta-
		\mathbf{U}^\theta\right)
		\right)\right\|_{C\left([-T,T],
		\mathcal{P}^{k,1}_{\mathrm{tcN}}[0,1]\right)}<C\ve$, 
		for some $C=C(T)>0$;
	
		\item the functions $\widehat{W}^\theta(t,\xi)$ and
		$\widehat{Z}^\theta(t,\xi)$ satisfy the
		non-degeneracy conditions \eqref{p-ins}
		with $\widehat{W}^\theta$ and
		$\widehat{Z}^\theta$ in place of, respectively,
		$W^\theta$ and
		$Z^\theta$,
		for all $(t,\xi,\theta)\in[-T,T]\times\R\times(0,1)$.
	\end{enumerate}
\end{enumerate}
\end{theorem}

\begin{proof}
First, following the same argument as in Theorem \ref{Ckgl},
we can show that there exists a unique global solution
$\left(\mathbf{U}^\theta,\frac{d}{d\theta}\mathbf{U}^\theta
\right)\in C\left([-T,T],\mathcal{P}^{k,1}_{\mathrm{tcN}}[0,1]\right)$ 
of the ODE system subject to initial data \eqref{id-path}
(see also Remark \ref{L1ODE}).
Here $\pt\left(\mathbf{U}^\theta\right)$ 
satisfies \eqref{ODE} with $\mathbf{U}^\theta$ 
instead of $\mathbf{U}$, while
$\pt\left(\frac{d}{d\theta}\mathbf{U}^\theta\right)$ 
satisfies the following linear (in 
$\frac{d}{d\theta}\mathbf{U}^\theta$) ODE system:
\begin{equation*}
	\pt\left(\frac{d}{d\theta}\mathbf{U}^\theta\right)
	=D_{\mathbf{U}^\theta}\mathbf{F}(\mathbf{U}^\theta)
	\cdot\frac{d}{d\theta}\mathbf{U}^\theta,
\end{equation*}
where $\mathbf{F}$ equals the right-hand 
side of \eqref{ODE}. 
In $D_{\mathbf{U}^\theta}
\mathbf{F}(\mathbf{U}^\theta)$, 
	the partial derivative of, for example,
	the nonlocal term $P_2$ with respect to $U^\theta$
is a linear operator that acts on 
$\frac{d}{d\theta}U^\theta$ as 
follows (recall \eqref{P-i-d}):
\[
\partial_{U^\theta}P_2\cdot\frac{d}{d\theta}U^\theta=
\frac{1}{8}\int_{-\infty}^\infty
\partial_{U^\theta}\left(\mathcal{E}p_2\right)\cdot
\frac{d}{d\theta}U^\theta\,d\eta.
\]
Thus, we have item (3a) for 
$\left(\mathbf{U}^\theta,
\frac{d}{d\theta}\mathbf{U}^\theta
\right)$.

Using that $W^\theta,Z^\theta\in 
C\left([-T,T],L^2\cap C^1\right)$, 
for all fixed $\theta$, and that 
$W^\theta,Z^\theta\in C([-T,T]\times\R\times[0,1])$, 
we have (cf.~\eqref{WZl1})
\begin{equation*}
	\left|W^\theta(t,\xi)\right|,
	\left|Z^\theta(t,\xi)\right|<1,\quad
	\text{for all }t\in[-T,T],\,\,
	|\xi|\geq M,\,\,\theta\in[0,1],
\end{equation*}
and therefore $W^\theta(t,\xi)$ and $Z^\theta(t,\xi)$ 
satisfy \eqref{ins} for 
such $(t,\xi,\theta)$ (see (3c)).

Then we examine the quantities 
$\left(W^\theta,\,W^\theta_\xi,\,
W^\theta_{\xi\xi},W^\theta_{\xi\xi\xi}\right)(t,\xi)$ for all
$(t,\xi,\theta)\in
\mathrm{int}(\Lambda_{T+\delta,M+\delta})\times(0,1)$
with some small $\delta>0$
(recall $\eqref{Lambd_T}$) and some fixed $T>0$.
Following the same 
strategy as in step~2 of the proof of
Lemma~\ref{L3},
but using the four-parameter perturbations
in \eqref{4p-tW} below,
we cover $\Lambda_{T+\delta,M+\delta}\times[0,1]$ 
by a finite collection of open neighborhoods
\[
\mathcal U_{(t_i,\xi_i,\theta_i)},\qquad i=1,\dots,\tilde N,
\]
centered at points $(t_i,\xi_i,\theta_i)
\in\Lambda_{T+\delta,M+\delta}\times[0,1]$,
such that (cf.~$\eqref{Lcov}$)
\[
\Lambda_{T+\delta,M+\delta}\times[0,1]
\subset \bigcup_{i=1}^{\tilde N}
\mathcal U_{(t_i,\xi_i,\theta_i)}.
\]
Here we must take additional care with the parameter
$\theta$ in the construction of the perturbations.
In \eqref{WZ-pert} (see also \eqref{ch-W}) we keep $t=t_i$ fixed
and then define $\mathbf{U}^\theta$ in a 
neighborhood of $t=t_i$ as a unique solution of the ODE system \eqref{ODE}.
To preserve the $\theta$-dependence, we consider the following
four-parameter
perturbations of 
$W$ and $Z$, cf.~\eqref{WZ-pert} (we take here $t=t_1$):
\begin{equation}\label{4p-tW}
\begin{split}
	&\tilde{W}^\theta(t_1,\xi;\nu,a)=W^\theta(t_1,\xi)
	+\beta(\theta)
	\sum\limits_{j=1}^4\nu_jf_j(\xi)+a\psi(\xi),\\
	&\tilde{Z}^\theta(t_1,\xi;b)=Z^\theta(t_1,\xi)+b\chi(\xi),
\end{split}
\end{equation}
where $f_j$, $j=1,\dots,4$, $\psi$, and $\chi$ satisfy assumptions
analogous to those in \eqref{fj-psi}, and
\begin{equation*}
\beta\in C^\infty[0,1],\quad
0<\beta(\theta)<\,1\text{ for all }\,\theta\in(0,1),\quad
\left.
\left(\beta,\beta^\prime\right)(\theta)
\right|_{\theta=0,1}=0,
\end{equation*}
while $a=a(\theta,\nu)$ and $b=b(\theta,\nu)$
are yet to be determined.
We then proceed as in \eqref{q-pert}--\eqref{paH1}, with the
additional variable $\theta\in[0,1]$.
Since $\mathbf{U}^\theta(t_1,\xi)$ is $C^1$ in $\theta$,
the estimates \eqref{WZ-lim}, \eqref{aspaH1}, and \eqref{lpaHj}
hold uniformly for $\theta\in[0,1]$.
Moreover, the values of $a=a(\theta,\nu)$ and $b=b(\theta,\nu)$
are small for all $\theta\in [0,1]$ and $\nu$ in a neighborhood of 
zero, which implies that we can determine
$a(\theta,\nu)$ and $b(\theta,\nu)$
from the system of equations \eqref{Hnuab} for all $\theta\in[0,1]$.
Here both $a(\theta,\nu)$ and $b(\theta,\nu)$ are $C^1$ in $\theta\in[0,1]$,
$a(\theta,0)=b(\theta,0)=0$ for any $\theta\in[0,1]$,
and $a(\theta,\nu)=b(\theta,\nu)=0$ for $\theta=0,1$ and any $\nu$.

The derivatives in $\nu$, however, depend on $\beta(\theta)$,
and the corresponding Jacobian has the following form (cf.~\eqref{Dr1}):
\begin{equation*}
	\left.D_\nu\left(
	\tilde{W}^\theta,\tilde{W}^\theta_\xi,
	\tilde{W}^\theta_{\xi\xi},\tilde{W}^\theta_{\xi\xi\xi}\right)
	\right|_{(t,\xi,\nu)=(t_1,\xi_1;0)}=
	\beta(\theta)
	\renewcommand{\arraystretch}{1.2}{
	\begin{pmatrix}
		f_1& \pxi f_1& \pxi^2 f_1& \pxi^3 f_1\\
		f_2& \pxi f_2& \pxi^2 f_2& \pxi^3 f_2\\
		f_3& \pxi f_3& \pxi^2 f_3& \pxi^3 f_3\\
		f_4& \pxi f_4& \pxi^2 f_4& \pxi^3 f_4\\
	\end{pmatrix}}(\xi_1).
\end{equation*}
Following this procedure, as in Step 2 of Lemma \ref{L3}, we
construct a perturbed family
\[
\left(\widetilde{\mathbf{U}}^\theta,
\frac{d}{d\theta}\widetilde{\mathbf{U}}^\theta
\right)(t,\xi;\nu),\qquad
\nu=(\nu_1,\dots,\nu_{4\tilde N_1}),
\quad \tilde N_1\le \tilde N,
\]
that satisfies item (1) for every $\nu$, with
$\widetilde{\mathbf{U}}^\theta$ in place of
$\widehat{\mathbf{U}}^\theta$,
and such that (cf.~$\eqref{nu=0}$),
\begin{equation}\label{1nu=0}
\begin{split}
	&\left(\mathbf{U}^\theta,
	\frac{d}{d\theta}\mathbf{U}^\theta\right)(t,\xi)=
	\left(\widetilde{\mathbf{U}}^\theta,
	\frac{d}{d\theta}\widetilde{\mathbf{U}}^\theta
	\right)(t,\xi;0),\quad
	(t,\xi,\theta)\in\Lambda_{T+\delta,M+\delta}\times[0,1].
\end{split}
\end{equation}

Proceeding analogously to step 3 in Lemma \ref{L3}, 
we arrive at the following result (cf.~\eqref{rpnu}).
\begin{equation*}
	\begin{split}
		&\left(\tilde{W}^\theta,\tilde{W}^\theta_\xi,
		\tilde{W}^\theta_{\xi\xi},
		\tilde{W}^\theta_{\xi\xi\xi}\right)(t,\xi;\nu)
		\neq(\pi,0,0,0),\quad\text{or}\\
		&\mathrm{rank}\left(\left.
		D_{(\nu_{4j-3},\nu_{4j-2},\nu_{4j-1},\nu_{4j})}\left(
		\tilde{W}^\theta,\tilde{W}^\theta_\xi,
		\tilde{W}^\theta_{\xi\xi},
		\tilde{W}^\theta_{\xi\xi\xi}\right)
		\right|_{(t,\xi;\nu)}
		\right)=4,
	\end{split}
\end{equation*}
for $(t,\xi,\theta)\in
\mathcal{U}_{(t_i,\xi_i,\theta_i)}\cap
\left(\R^2\times(0,1)\right)$,
$\nu\in(-\ve_1,\ve_1)^{4\tilde{N}_1}$,
with some small $\ve_1>0$, $i=1,\dots,\tilde{N}$, 
and the corresponding $j=1,\dots,\tilde{N_1}$.
Recalling Definition \ref{RV}, we conclude that
$(\pi,0,0,0)$ is a regular value of the map
$(t,\xi,\theta;\nu)\mapsto\left(\tilde{W}^\theta,
\tilde{W}^\theta_\xi,
\tilde{W}^\theta_{\xi\xi},\tilde{W}^\theta_{\xi\xi\xi}\right)(t,\xi;\nu)$ for
$(t,\xi,\theta;\nu)\in\mathrm{int}(\Lambda_{T+\delta,M+\delta})
\times(0,1)\times(-\ve_1,\ve_1)^{3\tilde{N}_1}$.

We proceed similarly for the vector
$\left(\tilde{Z}^\theta,\tilde{Z}^\theta_\xi,
\tilde{Z}^\theta_{\xi\xi},\tilde{Z}^\theta_{\xi\xi\xi}\right)$
and obtain 
a perturbed solution 
$\left(\widetilde{\mathbf{U}}^\theta,
\frac{d}{d\theta}\widetilde{\mathbf{U}}^\theta
\right)(t,\xi;\nu)$,
$\nu=(\nu_1,\dots,\nu_{4\tilde{N}_2})$,
$\tilde{N}_1\leq \tilde{N}_2$, which satisfies \eqref{1nu=0}, and
for which $(\pi,0,0,0)$ is a regular value of
each of the following two maps:
\begin{equation}\label{4map}
	(t,\xi,\theta;\nu)\mapsto
	\begin{cases}
		\left(\tilde{W}^\theta,
		\tilde{W}^\theta_\xi,
		\tilde{W}^\theta_{\xi\xi},
		\tilde{W}^\theta_{\xi\xi\xi}\right)(t,\xi;\nu),\\
		\left(\tilde{Z}^\theta,\tilde{Z}^\theta_\xi,
		\tilde{Z}^\theta_{\xi\xi},
		\tilde{Z}^\theta_{\xi\xi\xi}\right)(t,\xi;\nu)
	\end{cases}	
	(t,\xi,\theta;\nu)\in
	\mathrm{int}(\Lambda_{T+\delta,M+\delta})
	\times(0,1)\times(-\ve_2,\ve_2)^{4\tilde{N}_2},
\end{equation}
for some small $\ve_2\leq\ve_1$.

Applying the parametric transversality theorem
(Theorem \ref{BTL}) to the map \eqref{4map}, defined on the
disjoint union of two copies of the manifold
\[
\mathrm{int}(\Lambda_{T+\delta,M+\delta})
\times(0,1)\times(-\ve_2,\ve_2)^{4\tilde{N}_2},
\]
we conclude that 
there exists a dense subset 
$\tilde{\mathcal{N}}\subset(-\ve_2,\ve_2)^{4\tilde{N}_2}$,
such that $(\pi,0,0,0)$ is a regular value
of the four-dimensional maps
$\left(\tilde{W}^\theta,
\tilde{W}^\theta_\xi,
\tilde{W}^\theta_{\xi\xi},
\tilde{W}^\theta_{\xi\xi\xi}\right)(\cdot,\cdot;\nu)$
and
$\left(\tilde{Z}^\theta,
\tilde{Z}^\theta_\xi,
\tilde{Z}^\theta_{\xi\xi},
\tilde{Z}^\theta_{\xi\xi\xi}\right)(\cdot,\cdot;\nu)$
for each fixed $\nu\in\tilde{\mathcal{N}}$.
Taking into account that the domains of these maps are three-dimensional,
we conclude that the value $(\pi,0,0,0)$ is not attained 
for such values of $\nu$:
\begin{equation*}
\begin{split}
&\left(\tilde{W}^\theta,
\tilde{W}^\theta_\xi,
\tilde{W}^\theta_{\xi\xi},
\tilde{W}^\theta_{\xi\xi\xi}\right)(t,\xi;\nu)\neq(\pi,0,0,0),\\
&\left(\tilde{Z}^\theta,
\tilde{Z}^\theta_\xi,
\tilde{Z}^\theta_{\xi\xi},
\tilde{Z}^\theta_{\xi\xi\xi}\right)(t,\xi;\nu)\neq(\pi,0,0,0),
\end{split}
\end{equation*}
for all $(t,\xi,\theta;\nu)\in
\mathrm{int}(\Lambda_{T+\delta,M+\delta})
\times(0,1)\times\tilde{\mathcal{N}}$.
Thus, the perturbation satisfies (3c).
Notice that here we have (cf.~\eqref{reg-ptt})
\[
\max\{0,\dim\mathcal{Y}+\dim\mathcal{W}
-\dim\mathcal{X}\}=\max\{0,4+0-3\}=1,
\]
and thus we take $k\geq6$ in \eqref{id-path} to employ the parametric transversality theorem, see \eqref{F-reg}.

Finally, recalling that 
$\left(\widetilde{\mathbf{U}}^\theta,
\frac{d}{d\theta}\widetilde{\mathbf{U}}^\theta\right)(0,\xi;\nu)$ 
is a compact-in-$\xi$ perturbation of 
$\left(\mathbf{U}_0^\theta,
\frac{d}{d\theta}\mathbf{U}_0^\theta\right)(\xi)$ 
(see \eqref{L2}),
and arguing as in Step~3 
of Lemma \ref{L3}, we obtain that 
$\left(\widetilde{\mathbf{U}}^\theta,
\frac{d}{d\theta}\widetilde{\mathbf{U}}^\theta\right)(t,\xi;\nu)$ 
satisfies conditions (2), (3a), and (3b), with 
$\widetilde{\mathbf{U}}^\theta$ in place of 
$\widehat{\mathbf{U}}^\theta$, provided 
$\nu$ is taken sufficiently small.
\end{proof}

We next define piecewise regular paths connecting
two prescribed points $\mathbf{U}^0,\mathbf{U}^1\in
\mathcal{C}^{k}_{\Omega_{\mathrm{tcN}}^1}$, $k\geq6$.
\begin{definition}[Piecewise regular path]\label{prp}
The path 
\[
\mathbf{U}^\theta\in 
C\left([0,1],\mathcal{C}^k_{\Omega_{\mathrm{tcN}}^1}\right),
\quad k\geq6,
\]
is called piecewise regular if
there exists a finite set $\{\theta^*_i\}_{i=0}^{N^*}$,
$0=\theta_0^*<\theta_1^*<\dots<\theta_{N^*-1}^*<\theta_{N^*}^*=1$,
such that (see \eqref{P-pth})
\begin{equation}\label{Uprp}
	\mathbf{U}^\theta\in\mathcal{P}^{k,1}_{\mathrm{tcN}}
	[\theta_i^*,\theta_{i+1}^*],\quad
	i=0,\dots,N^*-1.
\end{equation}
\end{definition}

Applying Theorem \ref{dpath} to each regular segment in
\eqref{Uprp}, we obtain a piecewise regular path satisfying
the non-degeneracy conditions \eqref{p-ins} for all $t,\xi$ and
$\theta\in[0,1]\setminus\{\theta_i^*\}_{i=0}^{N^*}$.
This motivates the following definition of a piecewise regular
path under the ODE system \eqref{ODE}.
\begin{definition}[Piecewise regular path under ODE system \eqref{ODE}]
\label{PrpT}
Consider a piecewise regular path of initial data $\mathbf{U}^\theta_0$
with some $k\geq6$,
see Definition \ref{prp}.
Let $\mathbf{U}^\theta(t)$, $t\in[-T,T]$, be the evolution of 
$\mathbf{U}^\theta_0$ under the ODE system 
\eqref{ODE}, which satisfies (see items (1), (3a) 
in Theorem \ref{dpath}, and \eqref{Uprp}):
\begin{equation}\label{n-g^t}
\begin{split}
	&\left(
	\mathbf{U}^\theta,
	\frac{d}{d\theta}
	\left(\mathbf{U}^\theta\right)
	\right)\in C\left([-T,T],
	\mathcal{P}^{k,1}_{\mathrm{tcN}}[\theta_i^*,\theta_{i+1}^*]\right),
	\quad k\geq6,\quad i=0,\dots, N^*-1,\\
	&\mathbf{U}^{\theta}(t)=
	\left(U^\theta,V^\theta,W^\theta,Z^\theta,q^\theta\right)(t),
	\quad \theta\in[0,1].
\end{split}
\end{equation}
We say that $\mathbf{U}^\theta_0$ is a piecewise regular
path under the ODE system \eqref{ODE} if
$W^\theta$ and $Z^\theta$ satisfy
the non-degeneracy conditions \eqref{p-ins} 
for all $(t,\xi)\in[-T,T]\times\R$ and
$\theta\in[0,1]\setminus\{\theta_i^*\}_{i=0}^{N^*}$.	
\end{definition}

\begin{remark}\label{PrpT-r}
Note that, for each fixed $t\in[-T,T]$
and $\theta\in[0,1]\setminus\{\theta_i^*\}_{i=0}^{N^*}$, the 
functions $W^\theta(t,\cdot)$ and 
$Z^\theta(t,\cdot)$ in \eqref{n-g^t} attain 
the value $\pi$ at only finitely many points. 
This follows from the non-degeneracy conditions 
\eqref{p-ins} for $W^\theta$ and $Z^\theta$, which 
ensure that any point $\xi\in\R$ satisfying
$W^\theta(t,\xi)=\pi$ or $Z^\theta(t,\xi)=\pi$ is isolated.
Indeed, if $W^\theta(t,\xi)=\pi$, then at least one of
$W^\theta_\xi$, $W^\theta_{\xi\xi}$, or
$W^\theta_{\xi\xi\xi}$ is nonzero at $(t,\xi)$; the same argument
applies to $Z^\theta$. The far-field bounds exclude the value $\pi$
for large $|\xi|$, so the set of such points is finite.
\end{remark}

\begin{remark}[Uniform in $\theta$ conservation laws]
For any 
$\mathbf{U}^\theta
\in C\left([0,1],\mathcal{C}^k_{\Omega_{\mathrm{tcN}}^1}\right)$,
$k\geq6$, there exists a non-negative constant 
\begin{equation}\label{Unccl}
	K=\sup\limits_{\theta\in[0,1]}
	\max\left\{E_u^\theta(0),E_v^\theta(0),
	H^\theta(0)\right\},\quad
	K<\infty,
\end{equation}
such that
\begin{equation*}
	\sup\limits_{\theta\in[0,1]}
	\max\left\{E_u^\theta(t),E_v^\theta(t),
	H^\theta(t)\right\}
	\leq K,
\end{equation*}
for any $t\in[-T,T]$.
Here $E_u^\theta(t)$, $E_v^\theta(t)$, and
$H^\theta(t)$
are defined as in \eqref{consqU} with
$(U^\theta,V^\theta,W^\theta,Z^\theta,q^\theta)(t)$ 
instead of $(U,V,W,Z,q)(t)$ (see \eqref{n-g^t}).
\end{remark}

Our goal is to define the norm 
$\|\cdot\|_{\mathbf{U}^\theta(t)}$
of the tangent vector 
$\mathbf{R}^\theta(t)=\frac{d}{d\theta}
\mathbf{U}^\theta(t)$, 
$t\in[-T,T]$, for any initial path 
$\mathbf{U}^\theta_0$ in such a way that
\[
\left\|\mathbf{R}^\theta(t)\right\|_{\mathbf{U}^\theta(t)}
\leq C\left\|\mathbf{R}^\theta_0\right\|_{\mathbf{U}^\theta_0},\quad
C=C(T, K)>0,
\]
for every $\theta\in[0,1]\setminus\{\theta_i^*\}_{i=0}^{N^*}$
($K$ is given in \eqref{Unccl}).
With this estimate at our disposal, we conclude that
\begin{equation*}
	\left\|\mathbf{U}^\theta(t)\right\|_{\mathcal{L}}\leq
	C\left\|\mathbf{U}^\theta_0\right\|_{\mathcal{L}},
	\quad \mbox{for all }t\in[-T,T],
\end{equation*}
with some $C=C(T,K)>0$. Here 
$\left\|\mathbf{U}^\theta(t)\right\|_{\mathcal{L}}$ 
is a length of the path $\mathbf{U}^\theta(t)$, 
as defined in \eqref{npath} below.

To introduce a suitable definition of the norm 
$\|\cdot\|_{\mathbf{U}^\theta(t)}$, 
it is helpful to first formulate
this norm in the original variables $(u,v)$. 
In the next subsection, we define the 
norm of the tangent vector associated
with a family of smooth perturbed 
solutions of the two-component Novikov system.

\subsection{Tangent vectors for smooth solutions}\label{Tvss}
Assuming that $(u,v)$ is sufficiently 
smooth solution of \eqref{t-c-N-n}, we introduce 
the following one-parameter family 
of perturbed solutions of \eqref{t-c-N-n}:
\begin{equation}\label{rs}
	\begin{split}
		&u^\ve(t,x)=u(t,x)+\ve r(t,x)+\so(\ve),
		\quad
		r=\left.\pve u^\ve\right|_{\ve=0},\\
		&v^\ve(t,x)=v(t,x)+\ve s(t,x)+\so(\ve),
		\quad
		s=\left.\pve v^\ve\right|_{\ve=0}.
	\end{split}
\end{equation}
It is evident that after taking the 
partial derivative in $\ve$ 
at $\ve=0$ of the generic equation
\begin{equation}\label{apde}
	u_t^\ve+H(u^\ve,u_x^\ve,v^\ve,v_x^\ve)=0,
\end{equation}
we obtain the following evolutionary PDE for $r$:
\begin{equation}\label{ar_t}
	r_t+H_ur+H_{u_x}r_x+H_vs+H_{v_x}s_x=0.
\end{equation}
In our problem $H$ involves nonlocal terms in the form
$\tilde{P}=(1-\px^2)(\tilde{p}(u,u_x,v,v_x))$,
where the partial derivative of $\tilde{P}$ in, for example, $u$
is a linear operator acting on $r$ as follows:
$$
\partial_u\tilde{P}\cdot r=\left(1-\px^2\right)
\left(\partial_u\tilde{p}(u,u_x,v,v_x)\cdot r\right).
$$
Combining \eqref{t-c-N-n}, \eqref{apde} and \eqref{ar_t}, we 
obtain the following equations for $r,s$:
\begin{equation}\label{rs-t}
	\begin{split}
		&r_t+uvr_x+u_xvr+uu_xs+\px P_3+P_4=0,\\
		&s_t+uvs_x+uv_xs+vv_xr+\px S_3+S_4=0,
	\end{split}
\end{equation}
where
\begin{equation}\label{P34}
	\begin{split}
		&P_3=(1-\px^2)^{-1}\left(
		u_xvr_x+uv_xr_x+u_xv_xr+2uvr+uu_xs_x
		+\frac{1}{2}u_x^2s
		+u^2s
		\right),\\
		&P_4=\frac{1}{2}(1-\px^2)^{-1}\left(
		2u_xv_xr_x+u_x^2s_x
		\right),
	\end{split}
\end{equation}
and
\begin{equation*}
	\begin{split}
		&S_3=(1-\px^2)^{-1}\left(
		u_xvs_x+uv_xs_x+u_xv_xs+2uvs+vv_xr_x
		+\frac{1}{2}v_x^2r+v^2r
		\right),\\
		&S_4=\frac{1}{2}(1-\px^2)^{-1}\left(
		2u_xv_xs_x+v_x^2r_x
		\right).
	\end{split}
\end{equation*}
Differentiating \eqref{rs-t} in $x$, we have
\begin{equation}\label{rsx-t}
	\begin{split}
		&r_{tx}+uvr_{xx}+u_xvr_x+(u_{xx}v-2uv)r
		+\left(
		uu_{xx}+\frac{1}{2}u_x^2-u^2
		\right)s+P_3+\px P_4=0,\\
		&s_{tx}+uvs_{xx}+uv_xs_x+(uv_{xx}-2uv)s
		+\left(
		vv_{xx}+\frac{1}{2}v_x^2-v^2
		\right)r+S_3+\px S_4=0.
	\end{split}
\end{equation}
Notice that \eqref{rs-t} and \eqref{rsx-t} are consistent 
with \cite[equation (3.2)]{CCCS18} 
and \cite[equation (3.3)]{CCCS18}, respectively.
Consider two characteristics 
$x^\ve(t)$ and $x(t)$ corresponding to 
solutions $(u^\ve,v^\ve)$ 
and $(u,v)$, respectively:
\begin{equation}\label{ch1}
	\frac{d}{dt}x^\ve(t)
	=(u^\ve v^\ve)(x^\ve),\quad
	\frac{d}{dt}x(t)
	=(uv)(x).
\end{equation}
Then the horizontal shift $h(t,x)$ 
defined along the characteristics by
\begin{equation}\label{xve}
	x^\ve(t)=x(t)+\ve h(t,x(t))+\so(\ve),
	\quad h=\left.\pve x^\ve\right|_{\ve=0},
\end{equation}
satisfies the following equation (see \eqref{ch1}):
\begin{equation}\label{ah-t}
	h_t+uvh_x=\left.\frac{d}{d\ve}\left(
	(u^\ve v^\ve)(x^\ve)\right)\right|_{\ve=0}.
\end{equation}
Using that
\begin{equation}\label{fve}
	\frac{d}{d\ve}f^\ve(x^\ve)=
	\pve f^\ve(x^\ve)+\px f^\ve(x^\ve)\pve x^\ve,
\end{equation}
for any function $f^\ve$,
we obtain the following linear equation 
for $h$ from \eqref{ah-t} and \eqref{rs}:
\begin{equation}\label{h-t}
	h_t+uvh_x=(u_xv+uv_x)h+us+vr.
\end{equation}
Taking derivative in $x$ of \eqref{h-t}, we obtain
\begin{equation}\label{h-tx}
	h_{tx}+uvh_{xx}=(u_{xx}v+2u_xv_x+uv_{xx})h
	+us_x+u_xs+vr_x+v_xr.
\end{equation}

Introduce the set
\[
\mathcal{A}=\left\{
\mbox{smooth and uniformly bounded solutions }
h(t,x) \mbox{ of } \eqref{h-t}
\right\}.
\]
For each fixed $t$, we define the norm 
of the tangent vector 
$(r,s)(t,\cdot)$ associated with 
the solution $(u,v)(t,\cdot)$ as follows:
\begin{equation}\label{rstr}
	\|(r,s)(t,\cdot)\|_{(u,v)(t,\cdot)}
	=\inf_{h\in\mathcal{A}}
	\left(\sum_{i=1}^{6}
	\mathcal{I}_\alpha\!\left(|f_i(t,\cdot)|\right)
	\right),
\end{equation}
where the linear operator $\mathcal{I}_\alpha$ 
is defined in \eqref{I-op}, and the functions 
$f_i=f_i(t,x)$, $i=1,\dots,6$, are given by the 
following formulas (recall \eqref{fve} together 
with \eqref{D}, \eqref{rs}, and \eqref{xve}; 
we suppress the dependence on $t$ for brevity):
\begin{equation}\label{fi}
	\begin{split}
		&f_1(x)=D(x)\left.\frac{d}{d\ve}\left(
		x^\ve\right)\right|_{\ve=0}
		=\left(hD\right)(x),\\
		&f_2(x)=D(x)\left.\frac{d}{d\ve}\left(
		u^\ve(x^\ve)\right)\right|_{\ve=0}
		=\left((r+u_xh)D\right)(x),\\
		&f_3(x)=D(x)\left.\frac{d}{d\ve}\left(
		v^\ve(x^\ve)\right)\right|_{\ve=0}
		=\left((s+v_xh)D\right)(x),\\
		&f_4(x)=D(x)\left.\frac{d}{d\ve}\left(
		\arctan u_x^\ve(x^\ve)\right)\right|_{\ve=0}
		=\left((r_x+u_{xx}h)
		(1+v_x^2)\right)(x),\\
		&f_5(x)=D(x)\left.\frac{d}{d\ve}\left(
		\arctan v_x^\ve(x^\ve)\right)\right|_{\ve=0}
		=\left((s_x+v_{xx}h)
		(1+u_x^2)\right)(x),\\
		&f_6(x)=\left.\frac{d}{d\ve}\left(
		D^\ve(x^\ve)\frac{dx^\ve}{dx}\right)\right|_{\ve=0}
		=\left(2u_x(1+v_x^2)(r_x+u_{xx}h)
		+2v_x(1+u_x^2)(s_x+v_{xx}h)
		+h_xD\right)(x).
	\end{split}
\end{equation}
Note that in the expression for $f_6$ we use 
$\frac{dx^\ve}{dx}=1+\ve h_x+\so(\ve)$, and 
$D^\ve(x^\ve)=\left(\left(1+(u_x^\ve)^2\right)
\left(1+(v_x^\ve)^2\right)\right)(x^\ve)$. 
The right-hand side of \eqref{rstr} 
measures the cost of transporting 
energy from $u(t,x)$ to its perturbation 
$u^\ve(t,x^\ve)$, where the six terms 
$f_i$ in \eqref{fi} quantify the discrepancies 
between the solution and its perturbation.

\begin{remark}
Observe that \eqref{rstr} satisfies 
the axioms of a norm. 
Indeed, for $(r,s)(t,\cdot)=(0,0)$, 
we can take $h(t,\cdot)=0$, which 
implies that (see \eqref{fi})
$f_i(t,\cdot)=0$, $i=1,\dots,6$, and 
therefore $\|(r,s)\|_{(u,v)}=0$.
Taking into account that \eqref{h-t} 
is a linear equation, we have that 
$\lambda h$, $\lambda\in\R\setminus\{0\}$, 
is a solution of the equation
\[
(\lambda h)_t+uv(\lambda h)_x
=(u_xv+uv_x)(\lambda h)
+u(\lambda s)+v(\lambda r).
\]
Therefore
$\|(\lambda r,\lambda s)\|_{(u,v)}
=|\lambda| \|(r,s)\|_{(u,v)}$, for all 
$\lambda\in\R\setminus\{0\}$.
Finally, using the linearity of the 
operator $\mathcal{I}_\alpha$ 
and that $(h_1+h_2)$ satisfies
\[
(h_1+h_2)_t+uv(h_1+h_2)_x
=(u_xv+uv_x)(h_1+h_2)
+u(s_1+s_2)+v(r_1+r_2),
\]
provided that
\begin{equation*}
\begin{split}
	&(h_1)_t+uv(h_1)_x
	=(u_xv+uv_x)h_1
	+us_1+vr_1,\\
	&(h_2)_t+uv(h_2)_x
	=(u_xv+uv_x)h_2
	+us_2+vr_2, 
\end{split}
\end{equation*}
we verify the triangle inequality: 
$\|(r_1+r_2,s_1+s_2)\|_{(u,v)}
\leq\|(r_1,s_1)\|_{(u,v)}
+\|(r_2,s_2)\|_{(u,v)}$.
\end{remark}

For smooth solutions $(u,v)$, the norm 
of the associated tangent vectors,
as introduced in \eqref{rstr}, 
satisfies the following estimate.

\begin{theorem}\label{Lipn}
Consider a smooth solution $(u,v)(t,x)$, 
$t\in[-T,T]$, of the Cauchy 
problem \eqref{t-c-N-n}--\eqref{iid}.  
Then the following uniform estimate holds:
\begin{equation}\label{ae1}
	\|(r,s)(t)\|_{(u,v)(t,\cdot)}\leq C\,
	\|(r,s)(0)\|_{(u_0,v_0)(\cdot)},
\end{equation}
for some constant 
$C=C(T,E_{u_0},E_{v_0},H_0,\alpha)>0$.
\end{theorem}

\begin{proof}
To establish \eqref{ae1}, it is 
enough to verify that
(here and below, $C$ denotes 
a positive constant depending on 
$E_{u_0},E_{v_0},H_0$, and $\alpha$)
\begin{equation}\label{Gr}
	\frac{d}{dt}
	\sum\limits_{i=1}^{6}
	\mathcal{I}_\alpha
	\left(|f_i(t,\cdot)|\right)
	\leq
	C\sum\limits_{i=1}^{6}
	\mathcal{I}_\alpha
	\left(|f_i(t,\cdot)|\right).
\end{equation}

First, we observe that 
(cf.\,\,\cite[Lemma 3.1]{CCCS18}; 
throughout, we omit the arguments of $f_i(t,x)$)
\begin{equation}\label{cl}
	\begin{split}
		\frac{d}{dt}\mathcal{I}_\alpha(|f_i|)
		&=\int_{-\infty}^\infty
		\left[
		\left(
		|f_i|e^{-\alpha|x|}
		\right)_{\!t}
		+
		\left(
		uv|f_i|e^{-\alpha|x|}
		\right)_{\!x}\,
		\right]\,dx\\
		&=\int_{-\infty}^\infty
		\left[\sign(f_i)
		\left(
		(f_i)_t
		+(uvf_i)_x
		\right)e^{-\alpha|x|}
		-\sign(x)
		uv|f_i|e^{-\alpha|x|}
		\right]\,dx\\
		&\leq\mathcal{I}_\alpha
		(|(f_i)_t+(uvf_i)_x|)
		+E_{u_0}^{1/2}E_{v_0}^{1/2}
		\mathcal{I}_\alpha(|f_i|),
		\quad i=1,\dots,6,
	\end{split}
\end{equation}
where we have used Sobolev inequality 
\eqref{Sob} for $u$ and $v$, as well as \eqref{consq}.
Combining \eqref{Gr} and \eqref{cl}, we conclude 
that it suffices to establish 
the following inequalities:
\begin{equation}\label{clin}
	\mathcal{I}_\alpha
	(|(f_i)_t+(uvf_i)_x|)
	\leq
	C\sum\limits_{j=1}^{6}
	\mathcal{I}_\alpha
	\left(|f_j|\right),\quad
	i=1,\dots,6.
\end{equation}
In what follows, we will use the following estimates
(see \cite[Section 5]{HQ21}, 
\cite[Section 3.1]{KR25} and recall \eqref{consq}):
\begin{equation}\label{P1-e}
	\begin{split}
		\|P_1(t,\cdot)\|_{L^p},\,\,
		\|\px P_1(t,\cdot)\|_{L^p}\leq
		\frac{1}{2}\left\|e^{-|\cdot|}\right\|_{L^p}
		\left\|u^2v+uu_xv_x+\frac{1}{2}vu_x^2
		\right\|_{L^1}
		\leq C_pE_{u_0}E_{v_0}^{1/2},
	\end{split}
\end{equation}
and
\begin{equation}
	\label{P2-e}
	\begin{split}
		\|P_2(t,\cdot)\|_{L^p},\,\,
		\|\px P_2(t,\cdot)\|_{L^p}\leq
		\frac{1}{4}\left\|e^{-|\cdot|}\right\|_{L^p}
		\left\|u_x^2v_x\right\|_{L^1}
		\leq C_pK_{u_0},
	\end{split}
\end{equation}
for some constant $C_p>0$, $p\in [1,\infty]$.

\medskip

\textbf{Step 1, $i=1$ in \eqref{clin}.}
Using \eqref{tcNnd} and \eqref{h-t}, we 
have (recall \eqref{D} and \eqref{fi}):
\begin{equation}\label{f1cl}
	(f_1)_t+(uvf_1)_x
	=(hD)_t+(uvhD)_x=f_{1,1}+hf_{1,2}+hf_{1,3},
\end{equation}
where
\begin{equation}\label{f1,i}
	\begin{split}
		&f_{1,1}=\left(u(s+v_xh)+v(r+u_xh)\right)D,\\
		&f_{1,2}=2u_x\left(1+v_x^2\right)\left(
		u^2v+\frac{1}{2}v-P_1-\px P_2
		\right),\\
		&f_{1,3}=2\left(1+u_x^2\right)v_x\left(
		uv^2+\frac{1}{2}u-S_1-\px S_2
		\right).
	\end{split}
\end{equation}
Using the Sobolev inequality \eqref{Sob} 
for $u$ and $v$, we obtain (recall \eqref{fi})
\begin{equation}\label{f1,1}
	|f_{1,1}|\leq C\left(|f_2|+|f_3|\right).
\end{equation}
Applying again the Sobolev inequality for 
the terms involving $u^2v$, $uv^2$, $u$ and $v$,
taking into account \eqref{P1-e}, \eqref{P2-e} with 
$p=\infty$, \eqref{Cin} for $v_x$ and $u_x$, 
we arrive at (see \eqref{D})
\begin{equation}\label{f1,2}
	|f_{1,2}|,|f_{1,3}|\leq CD.
\end{equation}
Combining \eqref{f1cl}, \eqref{f1,1} 
and \eqref{f1,2}, we conclude that 
(recall $f_1=hD$ and $|h|D=|hD|$)
\begin{equation}\label{f1est}
	\left|(f_1)_t+(uvf_1)_x\right|
	\leq C\left(|f_1|+|f_2|+|f_3|\right),
\end{equation}
which, in view of \eqref{mon}, implies
\eqref{clin} for $i=1$.

\medskip

\textbf{Step 2, $i=2$ in \eqref{clin}.}
Using \eqref{tcNnd}, \eqref{rs-t} 
and \eqref{h-t}, we have (see \eqref{f1,i}):
\begin{equation*}
	\begin{split}
		(f_2)_t+(uvf_2)_x
		&=((r+u_xh)D)_t+(uv(r+u_xh)D)_x\\
		&=\left(\left(
		u^2v+\frac{1}{2}u_x^2v+uu_xv_x
		-P_1-\px P_2
		\right)\!h-P_4-\px P_3\right)\!D\\
		&\quad+(r+u_xh)(f_{1,2}+f_{1,3}).
	\end{split}
\end{equation*}
Invoking \eqref{f1,2}, the Sobolev inequality 
\eqref{Sob} for the terms containing $u^2 v$, 
and \eqref{P1-e} with $p=\infty$, 
yields the following estimate:
\begin{equation}
	\label{f2c2}
	\left|
	(f_2)_t+(uvf_2)_x
	\right|
	\leq C(|f_1|+|f_2|)
	+\left|
	\left(
	\frac{1}{2}u_x^2v+uu_xv_x
	\right)\!h-\px P_3
	\right|D
	+|P_4|D.
\end{equation}

\textbf{Step 2.1.} Let us estimate $|P_4|D$.
Notice that (see \eqref{P34};
cf.\,\,$I_{22}$ in 
\cite[Lemma	3.1, Item 2.]{CCCS18})
\begin{equation}\label{P4s}
	P_4=\frac{1}{2}(1-\px^2)^{-1}
	\left(
	2u_xv_x(r_x+u_{xx}h)
	+u_x^2v_xh_x+u_x^2(s_x+v_{xx}h)
	-\left(u_x^2v_xh\right)_x
	\right).
\end{equation}
Integrating by parts, we obtain
\begin{equation}\label{P4s1}
	\left|(1-\px^2)^{-1}
	\left(\left(u_x^2v_xh\right)_x\right)
	\right|
	\leq(1-\px^2)^{-1}\left(|f_1|\right).
\end{equation}
Recalling the definition of $f_5$, we obtain 
the following inequality 
from \eqref{P4s} and \eqref{P4s1}:
\begin{equation}\label{P4s2}
	|P_4|\leq
	\left|(1-\px^2)^{-1}
	\left(
	2u_xv_x(r_x+u_{xx}h)
	+u_x^2v_xh_x
	\right)
	\right|
	+(1-\px^2)^{-1}
	\left(|f_5|+|f_1|\right).
\end{equation}
Observe that 
\begin{equation}\label{P4s3}
	\begin{split}
		2u_xv_x(r_x+u_{xx}h)
		+u_x^2v_xh_x
		=&\left(2u_x(r_x+u_{xx}h)
		+2v_x\frac{1+u_x^2}{1+v_x^2}
		(s_x+v_{xx}h)
		+(1+u_x^2)h_x
		\right)\!v_x\\
		&\qquad\qquad-2v_x^2\frac{1+u_x^2}{1+v_x^2}
		(s_x+v_{xx}h)-v_xh_x.
	\end{split}
\end{equation}
Using \eqref{Cin} for $v_x$ and that 
$v_x^2/(1+v_x^2)\leq 1$ in \eqref{P4s3},
we conclude from \eqref{P4s2} that
\begin{equation}\label{P4s4}
	|P_4|\leq
	(1-\px^2)^{-1}\left(
	|f_6|+3|f_5|+|f_1|
	\right)
	+\left|
	(1-\px^2)^{-1}\left(
	v_xh_x
	\right)
	\right|.
\end{equation}
To estimate the last term in \eqref{P4s4}, we 
integrate by parts and use 
the following inequalities:
\begin{equation}\label{P4s5}
	\begin{split}
		\left|
		(1-\px^2)^{-1}\left(
		v_xh_x
		\right)
		\right|
		&=\frac{1}{2}\left|\int_{-\infty}^\infty
		e^{-|x-y|}\left(\sign(x-y)
		v_yh+v_{yy}h\right)\,dy
		\right|\\
		&\leq(1-\px^2)^{-1}(|f_1|)
		+\frac{1}{2}\left|\int_{-\infty}^\infty
		e^{-|x-y|}v_{yy}h\,dy\right|\\
		&\leq(1-\px^2)^{-1}(|f_1|+|f_5|)
		+\frac{1}{2}\left|\int_{-\infty}^\infty
		e^{-|x-y|}s_y\,dy\right|.
	\end{split}
\end{equation}
Integrating by parts, we estimate 
the last term in \eqref{P4s5} as follows:
\begin{equation}\label{P4s6}
	\frac{1}{2}\left|\int_{-\infty}^\infty
	e^{-|x-y|}s_y\,dy\right|
	\leq\frac{1}{2}\int_{-\infty}^\infty
	e^{-|x-y|}|s|\,dy
	\leq(1-\px^2)^{-1}(|f_1|+|f_3|).
\end{equation}
Combining \eqref{P4s4}, \eqref{P4s5} 
and \eqref{P4s6}, we arrive at
\begin{equation}\label{P4D}
	|P_4|D\leq C(1-\px^2)^{-1}
	(|f_1|+|f_3|+|f_5|+|f_6|)\cdot D.
\end{equation}

\textbf{Step 2.2.}
Here we estimate 
$\left|
\left(
\frac{1}{2}u_x^2v+uu_xv_x
\right)h-\px P_3
\right|D$ in \eqref{f2c2}.
Observe that (cf.\,\,$I_{23}$
in \cite[Lemma 3.1, Item 2]{CCCS18}):
\begin{equation}\label{P3s1}
	\left(
	\frac{1}{2}u_x^2v+uu_xv_x
	\right)h
	=-\frac{1}{2}\left(
	\int_x^\infty-\int_{-\infty}^x
	\right)
	\py
	\left(
	e^{-|x-y|}
	\left(
	\frac{1}{2}u_y^2v+uu_yv_y
	\right)\!h
	\right)dy
	=I_{2,1}+I_{2,2},
\end{equation}
where
\begin{equation*}
	\begin{split}
		&I_{2,1}=-\frac{1}{2}\left(
		\int_x^\infty-\int_{-\infty}^x
		\right)
		e^{-|x-y|}
		\left[
		\left(
		u_yu_{yy}v+
		\frac{3}{2}u_y^2v_y+uu_{yy}v_y+uu_yv_{yy}
		\right)\!h
		\right.\\
		&\left.
		\qquad\qquad\qquad\qquad\qquad\qquad
		\qquad\quad+
		\left(
		\frac{1}{2}u_y^2v+uu_yv_y
		\right)\!h_y
		\right]dy,
		\\
		&I_{2,2}=\frac{1}{2}
		\int_{-\infty}^\infty
		e^{-|x-y|}
		\left(
		\frac{1}{2}u_y^2v+uu_yv_y
		\right)\!h\,dy.
	\end{split}
\end{equation*}
Using Sobolev inequality \eqref{Sob} for $u$ 
and $v$ in $I_{2,2}$, we conclude from \eqref{P3s1} that
\begin{equation}\label{P3s}
	\left|
	\left(
	\frac{1}{2}u_x^2v+uu_xv_x
	\right)h-\px P_3
	\right|\leq
	|I_{2,1}-\px P_3|+|I_{2,2}|\leq
	C(1-\px^2)^{-1}(|f_1|)
	+|I_{2,1}-\px P_3|.
\end{equation}
Recalling \eqref{P34}, direct calculations show that
\begin{equation}\label{P3sd}
	I_{2,1}-\px P_3=-\frac{1}{2}\sum\limits_{i=3}^{9}I_{2,i},
\end{equation}
where
\begin{equation}\label{I23-9}
	\begin{split}
		&I_{2,3}=\left(
		\int_x^\infty-\int_{-\infty}^x
		\right)
		e^{-|x-y|}
		\left(
		u_yvr_y+u_yu_{yy}vh+\frac{1}{2}u_y^2vh_y
		\right)dy,\\
		&I_{2,4}=\left(
		\int_x^\infty-\int_{-\infty}^x
		\right)
		e^{-|x-y|}
		\left(uv_y(r_y+u_{yy}h)\right)\,dy,\\
		&I_{2,5}=\left(
		\int_x^\infty-\int_{-\infty}^x
		\right)
		e^{-|x-y|}
		\left(u_yv_y(r+u_yh)\right)\,dy,\\
		&I_{2,6}=\left(
		\int_x^\infty-\int_{-\infty}^x
		\right)
		e^{-|x-y|}
		\left(2uvr+u^2s\right)\,dy,\\
		&I_{2,7}=\left(
		\int_x^\infty-\int_{-\infty}^x
		\right)
		e^{-|x-y|}
		\left(uu_y(s_y+v_{yy}h)\right)\,dy,\\
		&I_{2,8}=\frac{1}{2}\left(
		\int_x^\infty-\int_{-\infty}^x
		\right)
		e^{-|x-y|}
		\left(u_y^2(s+v_yh)\right)dy,\\
		&I_{2,9}=
		\left(
		\int_x^\infty-\int_{-\infty}^x
		\right)
		e^{-|x-y|}
		\left(
		uu_yv_yh_y
		\right)\,dy.
	\end{split}
\end{equation}
Applying the Sobolev inequality to $u$ and 
using \eqref{Cin} for $u_y$ and $v_y$, we obtain
\begin{equation}\label{I24,7}
	|I_{2,4}|\leq C(1-\px^2)^{-1}|f_4|,\quad
	|I_{2,7}|\leq C(1-\px^2)^{-1}|f_5|.
\end{equation}
Using \eqref{Cin} for $u_y$, $v_y$, we obtain
\begin{equation}\label{I25,8}
	|I_{2,5}|\leq C(1-\px^2)^{-1}|f_2|,\quad
	|I_{2,8}|\leq C(1-\px^2)^{-1}|f_3|.
\end{equation}
Observing that $r=(r+u_yh)-u_yh$ and 
$s=(s+v_yh)-v_yh$, we can estimate $I_{2,6}$ as follows:
\begin{equation}\label{I2,6}
	|I_{2,6}|\leq C(1-\px^2)^{-1}(|f_1|+|f_2|+|f_3|).
\end{equation}
It remains to estimate $I_{2,3}$ and $I_{2,9}$.
Notice that (see $I_{2,3}$ in \eqref{I23-9})
\begin{equation*}
	\begin{split}
		u_yvr_y+u_yu_{yy}vh+\frac{1}{2}u_y^2vh_y&=
		\frac{v}{1+v_y^2}
		\left(u_y(1+v_y^2)(r_y+u_{yy}h)
		+\frac{1}{2}u_y^2(1+v_y^2)h_y
		\right)\\
		&=\frac{v}{1+v_y^2}\left(
		f_6-v_y(1+u_y^2)(s_y+v_{yy}h)
		\right)
		-\frac{v}{2}h_y.
	\end{split}
\end{equation*}
Observing that $\frac{1}{1+v_y^2}$ 
and $\frac{|v_y|}{1+v_y^2}$ are bounded by $1$, 
and integrating by parts in the 
term $\frac{v}{2}h_y$, we derive the following 
estimate for $I_{2,3}$:
\begin{equation}\label{I2,3}
	\begin{split}
		|I_{2,3}|&\leq
		C(1-\px^2)^{-1}(|f_1|+|f_5|+|f_6|)
		+|vh|\\
		&\leq
		C(1-\px^2)^{-1}(|f_1|+|f_5|+|f_6|)
		+C|h|.
	\end{split}
\end{equation}
To estimate $I_{2,9}$, we notice that
\begin{equation*}
	uu_yv_yh_y=\frac{uu_yv_y}
	{(1+u_y^2)(1+v_y^2)}f_6
	-\frac{2uu_y^2v_y}{1+u_y^2}(r_y+u_{yy}h)
	-\frac{2uu_yv_y^2}{1+v_y^2}(s_y+v_{yy}h),
\end{equation*}
which implies that
\begin{equation}\label{I2,9}
	|I_{2,9}|\leq
	C(1-\px^2)^{-1}(|f_4|+|f_5|+|f_6|).
\end{equation}

Finally, combining \eqref{P3s}, \eqref{P3sd},
\eqref{I24,7}, \eqref{I25,8}, \eqref{I2,6},
\eqref{I2,3} and \eqref{I2,9}, we conclude that
\begin{equation}\label{P3D}
	\left|
	\left(
	\frac{1}{2}u_x^2v+uu_xv_x
	\right)h-\px P_3
	\right|D\leq
	C(1-\px^2)^{-1}\left(
	\sum\limits_{i=1}^6|f_i|
	\right)\cdot D+C|f_1|.
\end{equation}
Recalling that $\mathcal{I}_\alpha$ is 
linear, and using \eqref{ID} and \eqref{mon}, 
we obtain from \eqref{f2c2}, \eqref{P4D}, 
and \eqref{P3D} the estimate \eqref{clin} 
for $i=2$.

\medskip

\textbf{Step 3, $i=3$ in \eqref{clin}.}
The proof of \eqref{clin} for $f_3$ is 
analogous to the case $i=2$; see Step 2 above.

\medskip

\textbf{Step 4, $i=4$ in \eqref{clin}.}
Direct calculations show that \eqref{tcNnd}, 
\eqref{rsx-t}, and \eqref{h-t} yield
\begin{equation}\label{f4c1}
	\begin{split}
		(f_4)_t+(uvf_4)_x
		&=\left((r_x+u_{xx}h)
		\left(1+v_x^2\right)\right)_t
		+\left(uv(r_x+u_{xx}h)
		\left(1+v_x^2\right)\right)_x\\
		&=\left(
		2uvr+u^2(s+v_xh)+h(2uu_xv-\px P_1-P_2)
		-\frac{1}{2}u_x^2s-P_3-\px P_4
		\right)\left(1+v_x^2\right)\\
		&\quad+2v_x(r_x+u_{xx}h)
		\left(uv^2+\frac{u}{2}-S_1-\px S_2\right).
	\end{split}
\end{equation}
Taking into account that (here we apply the Sobolev 
inequality \eqref{Sob} to $u$ and $v$, 
and \eqref{Cin} to $u_x$ and $v_x$)
\begin{equation*}
	\begin{split}
		&|2uvr(1+v_x^2)|\leq
		|2uv(r+u_xh)(1+v_x^2)|+|2uvu_x(1+v_x^2)h|
		\leq C(|f_1|+|f_2|),\\
		&|u^2(s+v_xh)(1+v_x^2)|\leq
		C|f_3|,
	\end{split}
\end{equation*}
and (see \eqref{P1-e} for $p=\infty$)
\begin{equation*}
	\begin{split}
		&|h(2uu_xv-\px P_1-P_2)(1+v_x^2)|
		\leq C|f_1|,\\
		&|2v_x(r_x+u_{xx}h)
		\left(uv^2+\frac{u}{2}-S_1-\px S_2\right)|
		\leq C|f_4|,
	\end{split}
\end{equation*}
we have from \eqref{f4c1} that
\begin{equation}\label{f4c2}
	|(f_4)_t+(uvf_4)_x|\leq
	C\sum\limits_{i=1}^4|f_i|
	+\left(\left|
	\frac{1}{2}u_x^2s-\px P_4
	\right|
	+|P_3|\right)
	\left(1+v_x^2\right).
\end{equation}

\medskip

\textbf{Step 4.1.}
Here we estimate $|P_3|\left(1+v_x^2\right)$.
Using that (see \eqref{P34})
\begin{equation*}
	\begin{split}
		&u_xv_xr=u_xv_x(r+u_xh)-u_x^2v_xh,\quad
		2uvr=2uv(r+u_xh)-2uu_xvh,\\
		& \frac{1}{2}u_x^2s=
		\frac{1}{2}u_x^2(s+v_xh)
		-\frac{1}{2}u_x^2v_xh,\quad
		u^2s=u^2(s+v_xh)-u^2v_xh,
	\end{split}
\end{equation*}
and applying the Sobolev inequality \eqref{Sob} 
for $u,v$ and \eqref{Cin} for 
$u_x,v_x$, we obtain the following inequality:
\begin{equation}\label{P3e1}
	|P_3|\leq C(1-\px^2)^{-1}
	\left(|f_1|+|f_2|+|f_3|\right)
	+\left|
	(1-\px^2)^{-1}
	(u_xvr_x+uv_xr_x+uu_xs_x)
	\right|.
\end{equation}
Considering that
\begin{equation}\label{P3e2}
	\begin{split}
		u_xvr_x+uv_xr_x+uu_xs_x
		=\,&u_xv(r_x+u_{xx}h)
		-u_xu_{xx}vh
		+uv_x(r_x+u_{xx}h)
		-uu_{xx}v_xh\\
		&+uu_x(s_x+v_{xx}h)
		-uu_xv_{xx}h,
	\end{split}
\end{equation}
we conclude from \eqref{P3e1} that 
(the terms $|uv_x(r_x+u_{xx}h)|$ and
$|uu_x(s_x+v_{xx}h)|$ in \eqref{P3e2} are 
estimated by $C|f_4|$ and $C|f_5|$, 
respectively, see \eqref{fi})
\begin{equation}\label{P3e3}
	|P_3|\leq C(1-\px^2)^{-1}
	\left(
	\sum\limits_{i=1}^5|f_i|
	\right)
	+\left|
	(1-\px^2)^{-1}
	\left(
	f_{4,1}
	\right)
	\right|,
\end{equation}
where
\begin{equation*}
	f_{4,1}=u_xv(r_x+u_{xx}h)
	-u_xu_{xx}vh-uu_{xx}v_xh-uu_xv_{xx}h.
\end{equation*}
Observe that
\begin{equation}\label{f4,1}
	f_{4,1}=u_xv(r_x+u_{xx}h)
	-\frac{1}{2}\left(u_x^2\right)_xvh
	-(uu_xv_xh)_x
	+u_x^2v_xh+uu_xv_xh_x.
\end{equation}
Integrating by parts the terms 
$(1-\px^2)^{-1}\left(
\frac{1}{2}\left(u_x^2\right)_xvh
\right)$
and
$(1-\px^2)^{-1}\left(
(uu_xv_xh)_x
\right)$, and
estimating $|u_x^2v_xh|$ by $C|f_1|$,
see \eqref{f4,1}, we deduce 
from \eqref{P3e3} that
\begin{equation}
	\label{P3e4}
	|P_3|\leq C(1-\px^2)^{-1}
	\left(
	\sum\limits_{i=1}^5|f_i|
	\right)
	+\left|
	(1-\px^2)^{-1}
	\left(
	u_xv(r_x+u_{xx}h)
	+\frac{1}{2}u_x^2vh_x
	+uu_xv_xh_x
	\right)
	\right|.
\end{equation}
Taking into account that the terms
$(1-\px^2)^{-1}
\left(
u_xv(r_x+u_{xx}h)
+\frac{1}{2}u_x^2 v h_x
\right)$
and
$(1-\px^2)^{-1}(u u_x v_x h_x)$
can be estimated exactly as in the cases 
of $I_{2,3}$ and $I_{2,9}$,
respectively (see \eqref{I23-9}, 
\eqref{I2,3}, and \eqref{I2,9}, Step~2.2),
we infer from \eqref{P3e4} that 
(noticing that no term of the form $C|h|$ 
appears after integrating by parts
in 
$(1-\px^2)^{-1}
\left(
u_xv(r_x+u_{xx}h)
+\frac{1}{2}u_x^2 v h_x
\right)$, cf.\,\eqref{I2,3})
\begin{equation}\label{|P3|}
	|P_3|\leq
	C(1-\px^2)^{-1}
	\left(
	\sum_{i=1}^6 |f_i|
	\right),
\end{equation}
and therefore,
\begin{equation}\label{P3eD}
	|P_3|\left(1+v_x^2\right)
	\leq |P_3|\, D
	\leq 
	C(1-\px^2)^{-1}
	\left(
	\sum_{i=1}^6 |f_i|
	\right) D.
\end{equation}

\medskip

\textbf{Step 4.2.}
Let us estimate 
$
\left|
\frac{1}{2}u_x^2s-\px P_4
\right|
\left(1+v_x^2\right)
$.
Using that (cf.\,\,\eqref{P3s1})
\begin{equation*}
	\begin{split}
		u_x^2s&=
		\frac{1}{2}\left(
		\int_{-\infty}^x-
		\int_x^\infty
		\right)
		\py\left(
		e^{-|x-y|}u_y^2s
		\right)dy\\
		&=
		(1-\px^2)^{-1}\left(u_x^2s\right)
		+\frac{1}{2}\left(
		\int_{-\infty}^x-
		\int_x^\infty
		\right)
		e^{-|x-y|}\left(
		2u_yu_{yy}s+u_y^2s_y
		\right)dy,
	\end{split}
\end{equation*}
we obtain (here we use that 
$u_x^2s=u_x^2(s+v_xh)-u_x^2v_xh$)
\begin{equation}\label{P4e1}
	\left|
	\frac{1}{2}u_x^2s-\px P_4
	\right|
	\leq C(1-\px^2)^{-1}(|f_1|+|f_3|)
	+\frac{1}{2}|I_{4,1}|,
\end{equation}
where
\begin{equation*}
	I_{4,1}=\left(
	\int_{-\infty}^x-
	\int_x^\infty
	\right)e^{-|x-y|}
	\left(
	u_yu_{yy}s-u_yv_yr_y
	\right)\,dy.
\end{equation*}
Notice that
\begin{equation*}
	u_yu_{yy}s-u_yv_yr_y=
	\frac{1}{2}(u_y^2)_y(s+v_yh)
	-u_yv_y(r_y+u_{yy}h).
\end{equation*}
Integrating by parts in the term 
$\frac{1}{2}(u_y^2)_y(s+v_y h)$, we arrive at
\begin{equation}\label{I4,1s}
	I_{4,1}=
	u_x^2(s+v_xh)-
	I_{4,2}-I_{4,3},
\end{equation}
with
\begin{equation}\label{I4,23}
	\begin{split}
		&I_{4,2}=
		\frac{1}{2}
		\int_{-\infty}^\infty
		e^{-|x-y|}u_y^2(s+v_yh)\,dy+
		\frac{1}{2}\left(
		\int_{-\infty}^x-
		\int_x^\infty
		\right)e^{-|x-y|}
		\left(
		u_y^2(s_y+v_{yy}h)
		\right)dy,\\
		&I_{4,3}=
		\left(
		\int_{-\infty}^x-
		\int_x^\infty
		\right)e^{-|x-y|}
		\left(
		u_yv_y(r_y+u_{yy}h)
		+\frac{1}{2}u_y^2v_yh_y
		\right).
	\end{split}
\end{equation}
Combining \eqref{P4e1}, \eqref{I4,1s}, and
\eqref{I4,23}, we obtain
\begin{equation}\label{P4e2}
	\left|
	\frac{1}{2}u_x^2s-\px P_4
	\right|
	\leq C(1-\px^2)^{-1}(|f_1|+|f_3|+|f_5|)
	+u_x^2|s+v_xh|+
	|I_{4,3}|.
\end{equation}
Applying \eqref{P4s3} for $I_{4,3}$, we conclude that
\begin{equation}\label{I4,3}
	|I_{4,3}|\leq C(1-\px^2)^{-1}
	(|f_5|+|f_6|)+\left|
	(1-\px^2)^{-1}(v_xh_x)
	\right|.
\end{equation}
Given that
\begin{equation*}
	v_xh_x=\frac{v_xf_6}{(1+u_x^2)(1+v_x^2)}
	-\frac{2u_x}{1+u_x^2}v_x(r_x+u_{xx}h)
	-\frac{2v_x^2}{1+v_x^2}(s_y+v_{yy}h),
\end{equation*}
we have from \eqref{I4,3} that
\begin{equation}\label{I4,3e}
	|I_{4,3}|\leq C(1-\px^2)^{-1}
	(|f_4|+|f_5|+|f_6|).
\end{equation}
Combining \eqref{P4e2} and \eqref{I4,3e}, we obtain
\begin{equation}\label{P4eD}
	\left|
	\frac{1}{2}u_x^2s-\px P_4
	\right|\left(1+v_x^2\right)
	\leq C(1-\px^2)^{-1}
	\left(
	\sum\limits_{i=1,i\neq2}^6
	|f_i|
	\right) D
	+|f_3|.
\end{equation}
Finally, recalling that $\mathcal{I}_\alpha$ 
is linear and using \eqref{ID} and \eqref{mon}, 
we obtain \eqref{clin} for $i=4$ from 
\eqref{f4c2}, \eqref{P3eD}, and \eqref{P4eD}.

\medskip

\textbf{Step 5, $i=5$ in \eqref{clin}.}
The proof of \eqref{clin} for $i=5$ is 
analogous to the case $i=4$, see Step 4.

\medskip

\textbf{Step 6, $i=6$ in \eqref{clin}.}
Taking into account that
$f_6=2u_xf_4+2v_xf_5+h_xD$, we have
\begin{equation}\label{f6c1}
	\begin{split}
		(f_6)_t+(uvf_6)_x
		=&\,2u_{tx}f_4+2uu_{xx}vf_4
		+2v_{tx}f_5+2uvv_{xx}f_5\\
		&+2u_x\left(
		(f_{4})_t+(uvf_4)_x
		\right)
		+2v_x\left(
		(f_{5})_t+(uvf_5)_x
		\right)
		+(h_xD)_t+(uvh_xD)_x.
	\end{split}
\end{equation}

Using that (see \eqref{h-tx})
\begin{equation*}
	\begin{split}
		(h_xD)_t+(uvh_xD)_x
		=\left(
		2u_xv_xh+v(r_x+u_{xx}h)
		+u(s_x+v_{xx}h)+u_xs+v_xr
		\right)D+f_{6,1},
	\end{split}
\end{equation*}
with
\begin{equation}\label{f6,1}
	\begin{split}
		f_{6,1}=&\,2u_x\left(1+v_x^2\right)h_x
		\left(
		u^2v+\frac{v}{2}-P_1-\px P_2
		\right)
		+2v_x\left(1+u_x^2\right)h_x
		\left(
		uv^2+\frac{u}{2}-S_1-\px S_2
		\right),
	\end{split}
\end{equation}
we obtain from \eqref{f6c1} and \eqref{tcNnd} that
\begin{equation}\label{f6c2}
	\begin{split}
		(f_6)_t+(uvf_6)_x
		=&\,2\left(
		u^2v+\frac{v}{2}-P_1-\px P_2
		\right)f_4
		+2\left(
		uv^2+\frac{u}{2}-S_1-\px S_2
		\right)f_5\\
		&+2u_x\left(
		(f_{4})_t+(uvf_4)_x
		\right)
		+2v_x\left(
		(f_{5})_t+(uvf_5)_x
		\right)\\
		&+\left(
		2u_xv_xh+u_xs+v_xr
		\right)D+f_{6,1}.
	\end{split}
\end{equation}
Notice that (cf.\,\,\eqref{f4c1})
\begin{equation}
	\label{f5c1}
	\begin{split}
		(f_5)_t+(uvf_5)_x
		&=\left((s_x+v_{xx}h)
		\left(1+u_x^2\right)\right)_t
		+\left(uv(s_x+v_{xx}h)
		\left(1+u_x^2\right)\right)_x\\
		&=\left(
		2uvs+v^2(r+u_xh)+h(2uvv_x-\px S_1-S_2)
		-\frac{1}{2}v_x^2r-S_3-\px S_4
		\right)\left(1+u_x^2\right)\\
		&\quad\qquad\qquad+2u_x(s_x+v_{xx}h)
		\left(u^2v+\frac{v}{2}-P_1-\px P_2\right).
	\end{split}
\end{equation}
Taking into account that (see \eqref{f6,1})
\begin{equation*}
	\begin{split}
		&u_x\left(1+v_x^2\right)h_x
		=\frac{u_x}{1+u_x^2}f_6
		-\frac{2u_x^2}{1+u_x^2}
		\left(1+v_x^2\right)
		(r_x+u_{xx}h)
		-2u_xv_x(s_x+v_{xx}h),\\
		&v_x\left(1+u_x^2\right)h_x
		=\frac{v_x}{1+v_x^2}f_6
		-\frac{2v_x^2}{1+v_x^2}
		\left(1+u_x^2\right)
		(s_x+v_{xx}h)
		-2u_xv_x(r_x+u_{xx}h),
	\end{split}
\end{equation*}
we obtain (see \eqref{f4c1} and \eqref{f5c1})
\begin{equation}
\label{f45c1}
	\begin{split}
		&2u_x\left(
		(f_{4})_t+(uvf_4)_x
		\right)
		+2v_x\left(
		(f_{5})_t+(uvf_5)_x
		\right)+f_{6,1}\\
		&\quad=2u_x\left(
		2uvr+u^2(s+v_xh)+h(2uu_xv-\px P_1-P_2)
		-\frac{1}{2}u_x^2s-P_3-\px P_4
		\right)\left(1+v_x^2\right)\\
		&\qquad
		+2v_x\left(
		2uvs+v^2(r+u_xh)+h(2uvv_x-\px S_1-S_2)
		-\frac{1}{2}v_x^2r-S_3-\px S_4
		\right)\left(1+u_x^2\right)\\
		&\qquad
		+\left(\frac{2u_x}{1+u_x^2}f_6
		-\frac{4u_x^2}{1+u_x^2}
		\left(1+v_x^2\right)
		(r_x+u_{xx}h)\right)
		\left(
		u^2v+\frac{v}{2}-P_1-\px P_2
		\right)\\
		&\qquad
		+\left(
		\frac{2v_x}{1+v_x^2}f_6
		-\frac{4v_x^2}{1+v_x^2}
		\left(1+u_x^2\right)
		(s_x+v_{xx}h)
		\right)
		\left(
		uv^2+\frac{u}{2}-S_1-\px S_2
		\right).
	\end{split}
\end{equation}
Combining \eqref{f45c1} and \eqref{f6c2}, 
we arrive at the following inequality
(here we use \eqref{P1-e} with $p=\infty$, the 
Sobolev inequality \eqref{Sob} and \eqref{Cin} 
for $u_x$ and $v_x$):
\begin{equation}\label{f6c3}
	\left|(f_6)_t+(uvf_6)_x\right|
	\leq C\sum\limits_{i=1}^6
	|f_i|+|f_{6,2}|,
\end{equation}
where
\begin{equation}\label{f6,2}
	\begin{split}
		f_{6,2}&=
		(2u_xv_xh+u_xs+v_xr)D-2u_x\left(
		\frac{1}{2}u_x^2s+P_3+\px P_4
		\right)\left(1+v_x^2\right)\\
		&\quad-2v_x\left(
		\frac{1}{2}v_x^2r+S_3+\px S_4
		\right)\left(1+u_x^2\right)\\
		&=2u_xv_xD-2u_x(P_3+\px P_4)
		\left(1+v_x^2\right)
		-2v_x(S_3+\px S_4)\left(1+u_x^2\right)
		+u_xs\left(1+v_x^2\right)
		+v_xr\left(1+u_x^2\right)\\
		&=u_x\left(1+v_x^2\right)(s+v_xh)
		+v_x\left(1+u_x^2\right)(r+u_xh)
		-2u_x\left(1+v_x^2\right)P_3
		-2v_x\left(1+u_x^2\right)S_3+f_{6,3},
	\end{split}
\end{equation}
where
\begin{equation}\label{f6,3}
	f_{6,3}=u_x\left(1+v_x^2\right)
	\left(u_x^2v_xh-2\px P_4\right)
	+v_x\left(1+u_x^2\right)
	\left(u_xv_x^2h-2\px S_4\right).
\end{equation}
We obtain from \eqref{f6c3}, \eqref{f6,2}, 
and \eqref{|P3|} that (the estimate 
for $|S_3|$ is the same as \eqref{|P3|})
\begin{equation}\label{f6c4}
	\left|(f_6)_t+(uvf_6)_x\right|
	\leq C\sum\limits_{i=1}^6
	|f_i|
	+C(1-\px^2)^{-1}\left(
	\sum\limits_{i=1}^6
	|f_i|\right) D
	+|f_{6,3}|.
\end{equation}
Using that
\begin{equation*}
	\begin{split}
		u_x^2v_xh&=
		\frac{1}{2}\left(
		\int_{-\infty}^x-
		\int_x^\infty
		\right)
		\py\left(
		e^{-|x-y|}u_y^2v_yh
		\right)dy\\
		&=(1-\px^2)^{-1}\left(u_y^2v_yh\right)
		+\frac{1}{2}\left(
		\int_{-\infty}^x-
		\int_x^\infty
		\right)
		e^{-|x-y|}\left(
		2u_yu_{yy}v_yh+u_y^2v_{yy}h+u_y^2v_yh_y
		\right)dy,
	\end{split}
\end{equation*}
we obtain (recall \eqref{I4,23})
\begin{equation*}
	\begin{split}
		u_x^2v_xh-2\px P_4=
		(1-\px^2)^{-1}\left(u_y^2v_yh\right)
		+\frac{1}{2}\left(
		\int_{-\infty}^x-
		\int_x^\infty
		\right)
		e^{-|x-y|}\left(
		u_y^2(s_y+v_{yy}h)
		\right)dy+I_{4,3}.
	\end{split}
\end{equation*}
Applying \eqref{I4,3e}, we arrive at
\begin{equation*}
	\left|u_x^2v_xh-2\px P_4\right|
	\leq
	C(1-\px^2)^{-1}(|f_1|+|f_4|+|f_5|+|f_6|).
\end{equation*}
Arguing similarly for 
$u_xv_x^2h-2\px S_4$, we conclude from 
\eqref{f6,3} and \eqref{f6c4} that
\begin{equation*}
	\left|(f_6)_t+(uvf_6)_x\right|
	\leq C\sum\limits_{i=1}^6
	|f_i|
	+C(1-\px^2)^{-1}\left(
	\sum\limits_{i=1}^6
	|f_i|\right) D,
\end{equation*}
which, in view of \eqref{mon} 
and \eqref{ID}, yields \eqref{clin} for $i=6$.
\end{proof}

\subsection{Tangent vectors in the transformed variables}

Following the approach of Section \ref{Tvss}, 
we define the seminorm
of the tangent vector associated with 
the following one-parameter family of perturbed solutions of 
the ODE system \eqref{ODE} (cf.\,\eqref{rs}
and see \eqref{Rtr} below):
\begin{equation}\label{Ubf-ve}
	\mathbf{U}^\ve(t,\xi)=\mathbf{U}(t,\xi)
	+\ve\mathbf{R}(t,\xi)+\so(\ve),
	\quad \mathbf{R}=\left.\pve\mathbf{U}^\ve\right|_{\ve=0},
\end{equation}
where we use the following notations
\begin{equation}\label{URbf}
\begin{split}
	&\mathbf{U}^\ve(t,\xi)=
	\left(U^\ve,V^\ve,W^\ve,Z^\ve,q^\ve\right)(t,\xi),
	\quad
	\mathbf{R}(t,\xi)=(R,S,A,B,Q)(t,\xi),\\
	&\mathbf{U}(t,\xi)
	=\left.\mathbf{U}^\ve(t,\xi)\right|_{\ve=0}
	=(U,V,W,Z,q)(t,\xi).
\end{split}
\end{equation}
We assume that 
$\mathbf{U}^\ve$ and $\mathbf{R}$ satisfy the following 
properties, consistent with $\mathbf{U}^\theta$ 
in Definition \ref{PrpT}:
\begin{itemize}
	\item $\mathbf{U}^\ve\in C\left([-T,T],
	\mathcal{C}^k_{\Omega_{\mathrm{tcN}}^1}\right)$,
	$k\geq6$, is a (global) solution of \eqref{ODE}
	for all fixed $\ve$;

	\item $\mathbf{R}\in C\left([-T,T],
	\mathcal{C}^k_{E_0^1}\right)$, 
	$k\geq6$, see \eqref{OmL1tcN};

	\item both $W^\ve$ and $Z^\ve$ satisfy \eqref{p-ins} for all
	$t,\xi\in[-T,T]\times\R$ and all fixed $\ve$.
\end{itemize}
Notice that since $W^\ve$ and $Z^\ve$ 
satisfy \eqref{p-ins}, we conclude that 
\begin{equation*}
	y^\ve(t,\cdot) \text{ is strictly monotone } \forall
	t\in[-T,T] \text{ and } \forall\ve,
\end{equation*}
where $y^\ve$ can be expressed in terms 
of $\mathbf{U}^\ve$ by (recall \eqref{pxiy})
\begin{equation}\label{y^ve}
	\begin{split}
		&y^\ve(t,\xi)=y_0(\xi)+\int_{-\infty}^\xi\left(
		\left(q^\ve\cos^2\frac{W^\ve}{2}\cos^2\frac{Z^\ve}{2}
		\right)(t,\xi')-y_0^\prime(\xi^\prime)
		\right)d\xi',\\
		&y(t,\xi)=\left.y^\ve(t,\xi)\right|_{\ve=0}
		=y_0(\xi)+\int_{-\infty}^\xi\left(\left(
		q\cos^2\frac{W}{2}\cos^2\frac{Z}{2}\right)(t,\xi')
		-y_0^\prime(\xi^\prime)\right)d\xi'.
	\end{split}
\end{equation}
Here we have used that $y(t,\xi)=y_0(\xi)+\so(1)$ as $|\xi|\to\infty$ for all $t$,
see \eqref{char1}.
Observe that the strict monotonicity of $y(t,\cdot)$ 
allows us to perform the change of variables 
$x=y(t,\xi)$ in the integrals appearing in 
\eqref{rstr}, see 
\eqref{J-alph-i}--\eqref{J-alph} below.

Analogously to the horizontal shift $h$ (see \eqref{xve}),
we define the perturbation of $\xi$ as follows:
\begin{equation}\label{xi-p}
	\xi^\ve(\xi)=\xi+\ve\eta(\xi)+\so(\ve),
	\quad\text{where }x^\ve=y^\ve(t,\xi^\ve),
	\quad \eta=\left.\pve \xi^\ve\right|_{\ve=0}.
\end{equation}
Observe that \eqref{char1} yields
$\pt x^\ve=\left(u^\ve v^\ve\right)(x^\ve)$,
which is consistent with \eqref{ch1}.
Now define $H(t,\xi)$ as follows 
(cf.\,\,\eqref{xve} and recall \eqref{y^ve}):
\begin{equation}\label{Hve}
	H(t,\xi)=\left.\frac{d}{d\ve}
	y^\ve(t,\xi^\ve)\right|_{\ve=0}
	=z(t,\xi)+\eta(\xi) y_\xi(t,\xi),
\end{equation}
where $\eta$ is given in \eqref{xi-p}, while 
$z$ can be found from \eqref{y^ve} and \eqref{URbf}:
\begin{equation}\label{zve}
	\begin{split}
		z(t,\xi)&=\left.\pve y^\ve(t,\xi)\right|_{\ve=0}\\
		&=\int_{-\infty}^\xi
		\left(Q\cos^2\frac{W}{2}\cos^2\frac{Z}{2}
		-\frac{q}{2}A\sin W\cos^2\frac{Z}{2}
		-\frac{q}{2}B\cos^2\frac{W}{2}\sin Z
		\right)(t,\xi')\,d\xi'.
	\end{split}
\end{equation}
Notice that since $\mathbf{R}(t,\cdot)\in E_0$, 
see \eqref{E-0}, we have that $Q\in L^1(\R)$ 
and the integral in \eqref{zve} is finite.
Recalling \eqref{xve} and that $x=y(t,\xi)$, we have
\begin{equation*}
	H(t,\xi)=h(t,y(t,\xi)),\quad \xi\in\R.
\end{equation*}

Next, we derive the expressions for the 
derivatives of the functions in \eqref{fi} 
with respect to $\ve$, expressed in 
terms of the variables $\mathbf{U}$ and $\mathbf{R}$.
Recalling  \eqref{xi-p} and \eqref{Hve}, we have
\begin{equation}\label{xve-tr}
	\left.\frac{d}{d\ve}
	x^\ve\right|_{\ve=0}=
	\left.\frac{d}{d\ve}\left(
	y^\ve(t,\xi^\ve)\right)\right|_{\ve=0}=
	z(t,\xi)+\eta(\xi)y_\xi(t,\xi).
\end{equation}
Then \eqref{uvdef}, \eqref{Ubf-ve}, 
and \eqref{xi-p} imply that
(see also \eqref{URbf})
\begin{equation}\label{uve-tr}
	\begin{split}
		&\left.\frac{d}{d\ve}\left(
		u^\ve(t,x^\ve)\right)\right|_{\ve=0}=
		\left.\frac{d}{d\ve}\left(
		u^\ve(t,y^\ve(t,\xi^\ve))\right)\right|_{\ve=0}=
		\left.\frac{d}{d\ve}\left(
		U^\ve(t,\xi^\ve)\right)\right|_{\ve=0}=
		R(t,\xi)+\eta(\xi)U_\xi(t,\xi),\\
		&\left.\frac{d}{d\ve}\left(
		v^\ve(t,x^\ve)\right)\right|_{\ve=0}=
		S(t,\xi)+\eta(\xi)V_\xi(t,\xi),
	\end{split}
\end{equation}
where, in the last equation, we apply the 
same reasoning as in the 
derivation of the first equation.
Using \eqref{pxuv}, \eqref{Ubf-ve}, 
and \eqref{xi-p}, we obtain
\begin{equation}\label{tan-tr}
\begin{split}
	&\left.\frac{d}{d\ve}\left(
	\arctan\left(u^\ve_x(t,x^\ve)\right)\right)\right|_{\ve=0}=
	\frac{1}{2}\left.\frac{d}{d\ve}
	W^\ve(t,\xi^\ve)\right|_{\ve=0}=
	\frac{1}{2}\left(
	A(t,\xi)+\eta(\xi)W_\xi(t,\xi)
	\right),\\
	&\left.\frac{d}{d\ve}\left(
	\arctan\left(v^\ve_x(t,x^\ve)\right)\right)
	\right|_{\ve=0}=
	\frac{1}{2}\left(
	B(t,\xi)+\eta(\xi)Z_\xi(t,\xi)
	\right).
\end{split}
\end{equation}
Finally, using the identity (see \eqref{q})
\begin{equation*}
	D^\ve(x^\ve)\, y^\ve_\xi(t,\xi^\ve)
	=q^\ve(t,\xi^\ve),
\end{equation*}
we deduce from \eqref{Ubf-ve} and \eqref{xi-p} that 
(here we use $y_\xi\neq0$)
\begin{equation}\label{meas-tr}
	\begin{split}
		\left.\frac{d}{d\ve}\left(
		D^\ve(x^\ve)\frac{dx^\ve}{dx}\right)\right|_{\ve=0}&=
		\left.\frac{d}{d\ve}\left(D^\ve(x^\ve)
		\frac{y_\xi^\ve(t,\xi^\ve)d\xi^\ve}
		{y_\xi(t,\xi)d\xi}\right)\right|_{\ve=0}\\
		&=y_\xi^{-1}(t,\xi)\left.\frac{d}{d\ve}\left(
		q^\ve(t,\xi^\ve)+\ve\eta^\prime(\xi)
		q^\ve(t,\xi^\ve)+\so(\ve)\right)\right|_{\ve=0}\\
		&=y_\xi^{-1}(t,\xi)\left(
		Q(t,\xi)+\eta(\xi)q_\xi(t,\xi)+\eta^\prime(\xi)q(t,\xi)
		\right).
	\end{split}
\end{equation}

Performing the change of variables $x=y(t,\xi)$ and 
using \eqref{xve-tr}, we obtain the following 
expression for $\mathcal{I}_\alpha(|f_1(t,\cdot)|)$ 
(recall \eqref{I-op} and \eqref{fi}):
\begin{equation}\label{If1ph1}
	\begin{split}
		\mathcal{I}_\alpha(|f_1(t,\cdot)|)&=
		\int_{-\infty}^{\infty}e^{-\alpha|x|}
		\left(\left.\frac{d}{d\ve}
		x^\ve\right|_{\ve=0}\right)D(t,x)\,dx\\
		&=\int_{-\infty}^{\infty}e^{-\alpha|y(t,\xi)|}
		\left(z(t,\xi)+\eta(\xi)y_\xi(t,\xi)
		\right)D(t,y(t,\xi))y_\xi(t,\xi)\,d\xi\\
		&=\int_{-\infty}^{\infty}e^{-\alpha|y(t,\xi)|}
		\left(z(t,\xi)+\eta(\xi)y_\xi(t,\xi)
		\right)q(t,\xi)\,d\xi,
	\end{split}
\end{equation}
where the final equality uses \eqref{q}.
Using \eqref{uve-tr}, \eqref{tan-tr}, \eqref{meas-tr}, 
and similar reasoning as in \eqref{If1ph1}, we 
obtain the following expressions for the 
integrals in the right-hand side of \eqref{rstr}:
\begin{equation}\label{J-alph-i}
	\int_{-\infty}^\infty
	e^{-\alpha|y(t,\xi)|}
	\left|\phi_{i}(t,\xi)\right|\,d\xi=
	\mathcal{I}_\alpha(|f_i(t,\cdot)|),
	\quad i=1,\dots,6,
\end{equation}
where
\begin{equation}\label{J-alph}
	\begin{split}
		&\phi_{1}(t,\xi)=\left(
		z(t,\xi)+\eta(\xi)y_\xi(t,\xi)
		\right)q(t,\xi),\\
		&\phi_{2}(t,\xi)=\left(
		R(t,\xi)+\eta(\xi)U_\xi(t,\xi)
		\right)q(t,\xi),\\
		&\phi_{3}(t,\xi)=\left(
		S(t,\xi)+\eta(\xi)V_\xi(t,\xi)
		\right)q(t,\xi),\\
		&\phi_{4}(t,\xi)=\frac{1}{2}\left(
		A(t,\xi)+\eta(\xi)W_\xi(t,\xi)
		\right)q(t,\xi),\\
		&\phi_{5}(t,\xi)=\frac{1}{2}\left(
		B(t,\xi)+\eta(\xi)Z_\xi(t,\xi)
		\right)q(t,\xi),\\
		&\phi_{6}(t,\xi)=
		Q(t,\xi)+\eta(\xi)q_\xi(t,\xi)
		+\eta^\prime(\xi)q(t,\xi),
	\end{split}
\end{equation}
with $y$ and $z$ given in 
\eqref{y^ve} and \eqref{zve}, respectively, 
and an arbitrary bounded smooth function $\eta(\xi)$.
Using \eqref{J-alph}, we define the
\textit{seminorm} of the tangent vector
$\mathbf{R}$ as per \eqref{rstr}:
\begin{equation}\label{Rtr}
	\|\mathbf{R}(t,\cdot)\|_{\mathbf{U}(t,\cdot)}=
	\inf\limits_{\eta\in C^{\infty}}
	\sum\limits_{i=1}^{6}
	\int_{-\infty}^\infty
	e^{-\alpha|y(t,\xi)|}
	\left|\phi_{i}(t,\xi)\right|\,d\xi.
\end{equation}
	
Notice that \eqref{Rtr} can vanish for nonzero tangent vectors $\mathbf{R}$.
To obtain positive definiteness, we must pass to the quotient space
under relabeling (cf.~\cite[Section 3]{HR07}
and see Proposition \ref{infr} below).
The corresponding equivalence classes are defined below.
\begin{definition}[Relabeling]\label{Drelab}
	Consider $\mathbf U,\overline{\mathbf U}\in \mathcal{C}^k_{\Omega_{\mathrm{tcN}}^1}$, $k\in\N$.
		We say that $\mathbf U\sim \overline{\mathbf U}$
	if they satisfy the following relations:
	\begin{equation}\label{relab}
	\begin{split}
		&(\overline U,\overline V,\overline W,\overline Z)(\sigma(\xi))
		=(U,V,W,Z)(\xi),\\
		&\overline q(\sigma(\xi))\sigma_\xi(\xi)=q(\xi),
	\end{split}
	\end{equation}
	where 
	$\mathbf U=(U,V,W,Z,q)$, $\overline{\mathbf U}
	=(\overline U,\overline V,\overline W,\overline Z,\overline q)$,
	and $\sigma$ belongs to the group
	$G^k$ of increasing bi-Lipschitz
	homeomorphisms $\sigma:\R\to\R$ such that
	\begin{equation*}
	\begin{split}
		&\sigma-\operatorname{id}\in C^k(\R),\\
		&\sigma_\xi-1\in L^1(\R),\quad
		0<c_\sigma\leq\sigma_\xi(\xi)\leq C_\sigma\quad
		\hbox{for a.e. }\,\xi\in\R.
	\end{split}
	\end{equation*}
	In the case $\mathbf U,\overline{\mathbf U}\in \Omega_{\mathrm{tcN}}^1$,
	the definition is similar, but with
	$\sigma-\mathrm{id}$ from $W^{1,\infty}(\R)$.
\end{definition}

Note that the density transformation
law $\overline q(\sigma(\xi))\sigma_\xi(\xi)=q(\xi)$ preserves
the compatibility relations \eqref{ccUV}.
Moreover, using that
\begin{equation*}
	\|\bar{q}-1\|_{L^1}=
	\int_{-\infty}^\infty|\bar{q}(\sigma)-1|\,d\sigma
	=\int_{-\infty}^\infty|\bar{q}(\sigma(\xi))-1|
	\sigma_\xi(\xi)\,d\xi
	\leq\|q-1\|_{L^1}+\|\sigma_\xi-1\|_{L^1},
\end{equation*}
as well as (here we take $s=\sigma(\xi)$)
\begin{equation*}
	\|(\sigma^{-1})_s-1\|_{L^1}=
	\int_{-\infty}^\infty
	\left|\left(\sigma^{-1}(s)\right)_s-1\right|\,ds
	=\int_{-\infty}^\infty
	\left|\frac{1}{\sigma_\xi(\xi)}-1\right|\sigma_\xi(\xi)\,d\xi
	=\|\sigma_\xi-1\|_{L^1},
\end{equation*}
we conclude that relabeling preserves $\Omega_{\mathrm{tcN}}^1$.

The next proposition establishes a family of
tangent vectors $\mathbf{R}$ having zero cost
with respect to the
seminorm \eqref{Rtr}.
\begin{proposition}\label{infr}
	A tangent direction $\mathbf{R}$ has zero
	seminorm $\|\mathbf{R}\|_{\mathbf{U}}$ if
	$\mathbf{U}$ and $\mathbf{U}^\ve$ belong to the same equivalence
	class under the relabeling
	$\sigma^{-1}(\xi)=\xi+\ve\rho(\xi)$,
$\rho\in C^\infty_0(\R)$:
\begin{equation}\label{infrelab}
\begin{split}
	&(U^\ve,V^\ve,W^\ve,Z^\ve)(\xi)
	=(U,V,W,Z)\left(\tau^\ve(\xi)\right),\\
	&q^\varepsilon
	=q\left(\tau^\ve(\xi)\right)\cdot\tau^\ve_\xi(\xi).
\end{split}
\end{equation}
\end{proposition}
\begin{proof}
First, consider $\mathbf{U}^\ve$
defined in \eqref{infrelab}.
Direct computations show that (recall \eqref{Ubf-ve}--\eqref{URbf})
\begin{equation}\label{inftan}
	(R,S,A,B,Q)(\xi)
	=\left(\rho U_\xi,\rho V_\xi,\rho W_\xi,\rho Z_\xi,
	(\rho q)_\xi\right)(\xi).
\end{equation}
	Combining \eqref{pxiy} and \eqref{zve}, we obtain
	the following expression for $z(\xi)$:
\begin{equation}\label{infz}
\begin{split}
	z(\xi)&=\int_{-\infty}^\xi
	\left((\rho q)_\xi\cos^2\frac{W}{2}\cos^2\frac{Z}{2}
	-\frac{\rho q}{2}W_\xi\sin W\cos^2\frac{Z}{2}
	-\frac{\rho q}{2}Z_\xi\cos^2\frac{W}{2}\sin Z
	\right)(\xi^\prime)\,d\xi^\prime\\
	&=\int_{-\infty}^\xi(\rho y_\xi)_{\xi}(\xi^\prime)\,d\xi^\prime
	=(\rho y_\xi)(\xi),
\end{split}
\end{equation}
where we have used that $\rho(\xi)\to0$ as $\xi\to-\infty$.
Taking $\eta(\xi)=-\rho(\xi)$ in \eqref{J-alph} gives,
together with \eqref{inftan} and \eqref{infz},
that $\phi_i\equiv0$ for all $i=1,\dots,6$.
\end{proof}

Motivated by the preceding discussion, we
introduce the notion of length of a piecewise regular path 
$\mathbf{U}^\theta$,
see Definition \ref{prp}.
This applies to any piecewise regular path connecting 
$\mathbf{U}^0,\mathbf{U}^1\in
\mathcal{C}^k_{\Omega_{\mathrm{tcN}}^1}$,
$k\geq6$.
We stress that the construction does not 
involve the ODE system \eqref{ODE}: 
$\theta$ parametrizes the path, and any 
time variable $t$---if present---is simply 
held fixed.

\begin{definition}[Length of a path]\label{lPR}
Consider a piecewise regular path $\mathbf{U}^\theta$, see Definition \eqref{prp}.
Then we define the length of $\mathbf{U}^\theta$ 
as an integral of the norm of its tangent vector. 
Namely,
\begin{equation}\label{npath}
	\left\|\mathbf{U}^\theta\right\|_{\mathcal{L}}
	=\int_0^1
	\left\|\frac{d}{d\theta}\mathbf{U}^\theta\right\|
	_{\mathbf{U}^\theta}d\theta
	=\int_0^1
	\left\|\mathbf{R}^\theta\right\|
	_{\mathbf{U}^\theta}d\theta,
\end{equation}
where $\mathbf{R}^\theta
=\frac{d}{d\theta}\mathbf{U}^\theta$, 
$\theta\in[0,1]\setminus\{\theta_i^*\}_{i=0}^{N^*}$,
see \eqref{Uprp}, and
\begin{equation}\label{Rthg}
\begin{split}
	&\left\|\mathbf{R}^\theta\right\|_{\mathbf{U}^\theta}
	=\inf\limits_{\eta\in C^{\infty}}
	\sum\limits_{i=1}^{6}
	\int_{-\infty}^\infty
	e^{-\alpha|y^\theta(\xi)|}
	\left|\phi_{i}^\theta(\xi)\right|\,d\xi,\quad
	\theta\in[0,1]\setminus\{\theta_i^*\}_{i=0}^{N^*}.
\end{split}
\end{equation}
	Here the functions $y^\theta$ and 
$\phi_{i}^\theta$, $i=1,\dots,6$, are given in
\eqref{y^ve} and \eqref{J-alph}, respectively, with
the functions
\[
\mathbf{U}^\theta=
\left(
U^\theta,V^\theta,W^\theta,Z^\theta,q^\theta\right)
\quad \text{and}\quad
\mathbf{R}^\theta=
\left(
R^\theta,S^\theta,A^\theta,B^\theta,Q^\theta\right),
\]
replacing the functions
\[
\mathbf{U}(t,\xi)=(U,V,W,Z,q)(t,\xi)\quad
\text{and}\quad
\mathbf{R}(t,\xi)=(R,S,A,B,Q)(t,\xi),
\]
respectively.
\end{definition}

In Theorem \ref{egammat} below, we aim to 
establish a uniform estimate for the 
length of the path $\mathbf{U}^\theta(t)$ for all $t$ 
under the ODE system \eqref{ODE}.  
If $y_\xi(t,\xi)\neq 0$ for every 
$\xi\in\R$, then Theorem \ref{egammat} follows 
directly from Theorem \ref{Lipn} 
after performing the change of variables 
$x=y(t,\xi)$ in \eqref{Rtr}.  
The main difficulty, therefore, is to handle 
the singular points $(t,\xi)$ at 
which $y_\xi(t,\xi)=0$, because the 
corresponding solution $(u,v)(t,x)$ is no 
longer smooth at those $x$.

Since the singular points are isolated 
and finite in number for the class of 
initial paths $\mathbf{U}^\theta_0$ under 
consideration, we show that they do 
not contribute to the estimate of 
$\frac{d}{dt}\|\mathbf{R}(t,\cdot)\|_{\mathbf{U}(t,\cdot)}$.  
More precisely, we perturb $\phi_i$ in a small 
neighborhood of each critical 
point so that the ratio $\frac{\phi_i(t,\xi)}{y_\xi(t,\xi)}$ 
vanishes as $\xi$ approaches the critical value (see \eqref{Rve1d} below).  
This allows us to treat the time derivatives of 
the integrals on the right-hand side of \eqref{Rtr} 
away from any neighborhood of these singular points.  
After performing the change of variables $x=y(t,\xi)$, 
the functions $u$ and $v$ 
in these integrals are of class $C^k$, and we 
can use the same estimates as in the 
proof of Theorem \ref{Lipn}.
Here we impose the stronger regularity condition
$k\geq13$. It guarantees that, at each singular point,
$\pxi^i y(t,\xi)\neq0$ for some $i\in\{2,\dots,13\}$; see
items (1)--(3) after \eqref{phima} and \eqref{Rve1d}.

\begin{theorem}\label{egammat}
Let $\mathbf{U}^\theta_0$
be a piecewise
regular path under the ODE system \eqref{ODE} 
for $t\in[-T,T]$, $T>0$, with regularity
$k\geq13$ (see Definition \ref{PrpT}).
Denote by $\mathbf{U}^\theta(t)$ 
the evolution of $\mathbf{U}^\theta_0$ under 
the ODE system \eqref{ODE}. 
Then we have the following uniform bound for the length 
of $\mathbf{U}^\theta(t)$, see \eqref{npath}:
\begin{equation}\label{pgest}
	\left\|\mathbf{U}^\theta(t)\right\|_{\mathcal{L}}\leq
	C\left\|\mathbf{U}^\theta_0\right\|_{\mathcal{L}},\quad
	\text{for all }\,t\in[-T,T],
\end{equation}
for some $C=C(T,K,\alpha)>0$, where $K>0$ 
is defined as in \eqref{Unccl}.
\end{theorem}

\begin{proof}
Recalling \eqref{npath},
we conclude that to prove \eqref{pgest}
it is enough to establish that
\begin{equation}\label{GaeU}
	\frac{d}{dt}\left\|\mathbf{R}^\theta(t,\cdot)\right\|
	_{\mathbf{U}^\theta(t,\cdot)}\leq
	C\left\|\mathbf{R}^\theta(t,\cdot)\right\|
	_{\mathbf{U}^\theta(t,\cdot)},
	\quad\text{for all }\,
	\theta\in[0,1]\setminus\{\theta_i^*\}_{i=0}^{N^*},
\end{equation}
for any $t\in[-T,T]$ and some $C=C(T,K,\alpha)>0$.

Fix arbitrary $t\in[-T,T]$ and 
$\theta\in[0,1]\setminus\{\theta_i^*\}_{i=0}^{N^*}$.
Suppose that $(t,\xi)\not\in\left(
\Gamma^{W^\theta}\cup\Gamma^{Z^\theta}\right)$,
see \eqref{GammaWZ}, for all $\xi\in\R$.
In this case, we perform the change of variables
$x=y^\theta(t,\xi)$ in the integrals on the right-hand 
side of \eqref{Rtr}, where $y^\theta$ is defined 
as in \eqref{y^ve}, but with $W^\theta$, $Z^\theta$, 
and $q^\theta$ replacing $W$, $Z$, and $q$, respectively.
This transformation reduces the proof of 
\eqref{GaeU} to that of \eqref{Gr}, which was 
already established in Theorem \ref{Lipn}.

Now assume that there exists $\xi\in\R$ 
such that $(t,\xi)\in\left(
\Gamma^{W^\theta}\cup\Gamma^{Z^\theta}\right)$.
In view of \eqref{p-ins} for $W^\theta$ and $Z^\theta$,
there exists only a finite number of $\xi_i$, $i=1,\dots,M$, 
such that $(t,\xi_i)\in\left(
\Gamma^{W^\theta}\cup\Gamma^{Z^\theta}\right)$, 
see Remark \ref{PrpT-r}.
It is crucial for the subsequent analysis that 
there is no ``plateau'' in 
$\Gamma^{W^\theta}\cup\Gamma^{Z^\theta}$, i.e., 
there is no interval $[\tilde{\xi}_{1},\tilde{\xi}_{2}]$,
$\tilde{\xi}_{1}<\tilde{\xi}_{2}$, such that
$(t,\xi)\in\Gamma^{W^\theta}\cup\Gamma^{Z^\theta}$
for all $\xi\in[\tilde{\xi}_{1},\tilde{\xi}_{2}]$.
In what follows, we establish the following estimates
(we omit the superscript $\theta$ for simplicity):
\begin{equation}\label{phima}
	\frac{d}{dt}\int_{-\infty}^\infty
	e^{-\alpha|y(t,\xi)|}
	\left|\phi_{i}(t,\xi)\right|\,d\xi
	\leq C\sum\limits_{j=1}^6
	\int_{-\infty}^\infty
	e^{-\alpha|y(t,\xi)|}
	\left|\phi_j(t,\xi)\right|\,d\xi,
	\quad i=1,\dots,6.
\end{equation}
We present the details for $i=1$; the cases
$i=2,\dots,6$ follow analogously.
For simplicity of the arguments and 
notation, we assume $M=1$, i.e., there is only one 
point $(t,\xi_1)\in\Gamma^{W^\theta}\cup\Gamma^{Z^\theta}$ 
for the fixed $t$ and $\theta$.
The analysis for the general case $M\in\N$ 
proceeds along the same lines.

Let $\tilde{\phi}_1^{\ve_0}(t,\xi)$
be defined as an $L^1(\R)$-small perturbation of $\phi_1(t,\xi)$
in \eqref{J-alph}.
Namely, for any $\ve_0>0$, we consider a function 
$\tilde{\phi}_1^{\ve_0}(t,\xi)$ such that
\begin{enumerate}
	\item $\tilde{\phi}_1^{\ve_0}(t,\xi)=\phi_1(t,\xi)$
	for $|\xi-\xi_1|\geq\ve_0$,

	\item $\int_{\xi_1-\ve_0}^{\xi_1+\ve_0}
	\left|\tilde{\phi}_1^{\ve_0}(t,\xi)\right|
	d\xi+\left|\frac{d}{dt}\int_{\xi_1-\ve_0}^{\xi_1+\ve_0}
	\left|\tilde{\phi}_1^{\ve_0}(t,\xi)\right|
	\,d\xi\right|=\mathcal{O}(\ve_0)$, as $\ve_0\to0$,

	\item $\left.\pxi^n\tilde{\phi}_1^{\ve_0}(t,\xi)
	\right|_{\xi=\xi_1}=0$ for all $n=0,\dots,13$.
\end{enumerate}
We restrict the third condition to $n\leq 13$, 
because there exists $i\in\{2,\dots,13\}$
such that $\pxi^i y(t,\xi_1)\neq0$; see \eqref{Tayy2} and
Remark \ref{multyz} (see
also the fractions in \eqref{Rve1d} below).
Taking into account that
$\phi_1(t,\xi)=\tilde{\phi}_1^{\ve_0}(t,\xi)$ for 
$|\xi-\xi_1|\geq\ve_0$, we have
\begin{equation*}
\begin{split}
	\int_{-\infty}^\infty
	e^{-\alpha|y(t,\xi)|}
	\left|\phi_1(t,\xi)\right|\,d\xi
	=&\int_{-\infty}^\infty
	e^{-\alpha|y(t,\xi)|}
	\left|\tilde{\phi}_1^{\ve_0}(t,\xi)\right|\,d\xi\\
	&+\int_{\xi_1-\ve_0}^{\xi_1+\ve_0}
	e^{-\alpha|y(t,\xi)|}\left(
	\left|\phi_1(t,\xi)\right|
	-\left|\tilde{\phi}_1^{\ve_0}(t,\xi)\right|
	\right)\,d\xi.
\end{split}
\end{equation*}
Then using the uniform in $\ve_0$ bound given in
item (2) above
we obtain the following 
expansion for any $0<\ve_1\ll\ve_0$:
\begin{equation}\label{phi1te}
	\begin{split}
		\frac{d}{dt}\int_{-\infty}^\infty
		e^{-\alpha|y(t,\xi)|}
		\left|\phi_1(t,\xi)\right|\,d\xi&=
		\frac{d}{dt}\int_{-\infty}^\infty
		e^{-\alpha|y(t,\xi)|}
		\left|\tilde{\phi}_1^{\ve_0}(t,\xi)\right|\,d\xi
		+\mathcal{O}(\ve_0)\\
		&=
		\frac{d}{dt}\left(
		\int_{-\infty}^{\xi_1-\ve_1}+
		\int_{\xi_1+\ve_1}^{\infty}
		\right)
		e^{-\alpha|y(t,\xi)|}
		\left|\tilde{\phi}_1^{\ve_0}(t,\xi)\right|\,d\xi
		+\mathcal{O}(\ve_0)+\mathcal{O}(\ve_1)\\
		&=:\frac{d}{dt}\tilde{I}_1
		+\mathcal{O}(\ve_0)+\mathcal{O}(\ve_1).
	\end{split}
\end{equation}
Changing the variables $x=y(t,\xi)$ 
on the right-hand side of 
\eqref{phi1te}, we arrive at (recall \eqref{fi})
\begin{equation}\label{I1-til-a}
	\frac{d}{dt}\tilde{I}_1
	=\frac{d}{dt}\left(
	\int_{-\infty}^{x_1^-}+
	\int_{x_1^+}^{\infty}\right)
	e^{-\alpha|x|}
	\left|\tilde{f}_1^{\ve_0}(t,x)\right|\,dx,
\end{equation}
where
\begin{equation}\label{xpm}
	x_1^\pm=y(t,\xi_1\pm\ve_1),
	\quad
	\tilde{f}_1^{\ve_0}(t,x)
	=\frac{\tilde{\phi}_1^{\ve_0}(t,\xi)}{y_\xi(t,\xi)},
	\quad 
	x=y(t,\xi),
	\quad \xi\in\R\setminus(\xi_1-\ve_1,\xi_1+\ve_1).
\end{equation}
Taking into account that
\begin{equation*}
	\left(
	\int_{-\infty}^{x_1^-}+
	\int_{x_1^+}^{\infty}\right)
	\left(
	e^{-\alpha|x|}(uv)(t,x)
	\left|\tilde{f}_1^{\ve_0}(t,x)\right|\right)_x
	\,dx=-R_{\ve_1}(t),
\end{equation*}
where
\begin{equation}\label{Rve1d-a}
	R_{\ve_1}(t)=
	e^{-\alpha|x_1^+|}(uv)(t,x_1^+)
	\left|\tilde{f}_1^{\ve_0}(t,x_1^+)\right|
	-e^{-\alpha|x_1^-|}(uv)(t,x_1^-)
	\left|\tilde{f}_1^{\ve_0}(t,x_1^-)\right|,
\end{equation}
we obtain from \eqref{I1-til-a} the 
following expression for
$\frac{d}{dt}\tilde{I}_1$:
\begin{equation}\label{I1-til}
	\begin{split}
		\frac{d}{dt}\tilde{I}_1
		=\left(
		\int_{-\infty}^{x_1^-}+
		\int_{x_1^+}^{\infty}\right)
		\left[\left(
		e^{-\alpha|x|}
		\left|\tilde{f}_1^{\ve_0}(t,x)\right|\right)_t
		+\left(
		e^{-\alpha|x|}(uv)(t,x)
		\left|\tilde{f}_1^{\ve_0}(t,x)\right|\right)_x
		\right]\,dx+R_{\ve_1}(t).
	\end{split}
\end{equation}
Observing that $R_{\ve_1}(t)$ from \eqref{Rve1d-a} 
can be written in the transformed variables as 
follows (see \eqref{uvdef} and \eqref{xpm};  we omit 
the dependence on $t$):
\begin{equation}\label{Rve1d}
	\begin{split}
		R_{\ve_1}=e^{-\alpha|y(\xi_1+\ve_1)|}(UV)(\xi_1+\ve_1)
		\frac{\left|\tilde{\phi}_1^{\ve_0}(\xi_1+\ve_1)\right|}
		{y_\xi(\xi_1+\ve_1)}
		-e^{-\alpha|y(\xi_1-\ve_1)|}(UV)(\xi_1-\ve_1)
		\frac{\left|\tilde{\phi}_1^{\ve_0}(\xi_1-\ve_1)\right|}
		{y_\xi(\xi_1-\ve_1)},
	\end{split}
\end{equation}
and recalling 
$\left.\pxi^n\tilde{\phi}_1^{\ve_0}(t,\xi)
\right|_{\xi=\xi_1}=0$ for $n=0,\dots,13$
(see item (3) above) and Remark \ref{multyz}, 
we conclude from \eqref{Rve1d} that
\begin{equation}\label{Rve1e}
	R_{\ve_1}(t)=\mathcal{O}(\ve_1).
\end{equation}

Using the estimate \eqref{f1est} in \eqref{I1-til},
we obtain that
\begin{equation}\label{I1-til1}
	\begin{split}
		\frac{d}{dt}\tilde{I}_1
		&\leq C\sum\limits_{i=1}^{3}
		\left(
		\int_{-\infty}^{x_1^-}+
		\int_{x_1^+}^{\infty}\right)
		e^{-\alpha|x|}
		\left|\tilde{f}_i^{\ve_0}(t,x)\right|\,dx
		+R_{\ve_1}(t)\\
		&\leq C\sum\limits_{i=1}^{3}
		\int_{-\infty}^{\infty}e^{-\alpha|x|}
		\left|\tilde{f}_i^{\ve_0}(t,x)\right|\,dx
		+R_{\ve_1}(t)\\
		&=C\sum\limits_{i=1}^{3}
		\int_{-\infty}^{\infty}e^{-\alpha|y(t,\xi)|}
		\left|\tilde{\phi}_i^{\ve_0}(t,\xi)\right|\,d\xi
		+R_{\ve_1}(t),
	\end{split}
\end{equation}
where in the last equality we have 
changed the variables $x=y(t,\xi)$.
Here $\tilde{\phi}^{\ve_0}_i$, $i=1,2,3$, are defined 
as $L^1(\R)$-small perturbations of $\phi_i(t,\xi)$,
satisfying assumptions similar to those imposed on
$\tilde{\phi}^{\ve_0}_1$
above. Using item (2) above for $\tilde{\phi}^{\ve_0}_i$, 
$i=1,2,3$, we conclude from \eqref{I1-til1} that
\begin{equation}\label{I1-til2}
	\frac{d}{dt}\tilde{I}_1\leq
	C\sum\limits_{i=1}^{3}
	\int_{-\infty}^{\infty}e^{-\alpha|y(t,\xi)|}
	\left|\phi_i(t,\xi)\right|\,d\xi
	+\mathcal{O}\left(\ve_0\right)
	+R_{\ve_1}(t).
\end{equation}
Finally, combining \eqref{phi1te}, \eqref{I1-til2}, and
\eqref{Rve1e}, we arrive at \eqref{phima} for $i=1$.
\end{proof}

\subsection{Geodesic distance in $\Omega_{\mathrm{tcN}}^1$}\label{GdO}

In this section we introduce 
a new metric on the set 
$\Omega_{\mathrm{tcN}}^1$
defined in \eqref{tcNOm}.
To this end, we first fix arbitrary $\mathbf{U}^0,\mathbf{U}^1\in
\mathcal{C}^k_{\Omega_{\mathrm{tcN}}^1}$
and consider paths $\mathbf{U}^\theta$ as 
in Definition \ref{lPR} that connect 
$\mathbf{U}^0$ and $\mathbf{U}^1$.  
We additionally assume that the energies \eqref{consqU} of 
$\mathbf{U}^\theta$ are uniformly 
bounded by some constant $K>0$ for all 
$\theta\in[0,1]$.  Taking the infimum of the 
lengths $\left\|\mathbf{U}^\theta\right\|_{\mathcal{L}}$ over all 
such paths yields the geodesic distance 
$d_{\Omega_{\mathrm{tcN}}^1,K}(\cdot,\cdot)$
between the sufficiently 
regular functions $\mathbf{U}^0$ and $\mathbf{U}^1$ 
(see Definition \ref{Gdist} below). 
We emphasize that this definition is independent 
of the ODE system \eqref{ODE}; time $t$ does not 
appear, or, if present in the 
functions, is regarded as a fixed parameter.

By approximating $\mathbf{U}^\theta$ 
with regular paths governed by 
\eqref{ODE} and applying Theorem 4.8, we conclude that 
$d_{\Omega_{\mathrm{tcN}}^1,K}(\cdot,\cdot)$ satisfies the Lipschitz property on 
$\mathcal{C}^k_{\Omega_{\mathrm{tcN}}^1}$,
$k\geq13$, as stated in 
Theorem \ref{LRO}. Finally, using a simple completion 
argument, we extend this metric to the 
entire set $\Omega_{\mathrm{tcN}}^1$ (see Proposition \ref{COm} 
and the upcoming Definition \ref{OGdist}).

\begin{definition}[Geodesic distance in 
$\mathcal{C}^k_{\Omega_{\mathrm{tcN}}^1}$]
\label{Gdist}
Consider
$
\mathbf{U}^0,\mathbf{U}^1\in
\mathcal{C}^k_{\Omega_{\mathrm{tcN}}^1},
$
$k\geq 13$,
where $\mathcal{C}^k_{\Omega_{\mathrm{tcN}}^1}$
is defined in \eqref{CktcN},
and a constant $K>0$ such that
\begin{equation}
	\label{K-unb}
	(1+C_0)\max\left\{E_u^0(0),E_v^0(0),
	H^0(0),E_u^1(0),E_v^1(0),
	H^1(0)\right\}\leq K,\quad C_0>0,
\end{equation}
where $(E_u^0(0),E_v^0(0),H^0(0))$ and 
$(E_u^1(0),E_v^1(0),H^1(0))$ are defined 
as in \eqref{consqU} with, 
respectively, $\mathbf{U}^0$ and 
$\mathbf{U}^1$ instead of $\mathbf{U}$.
Then we define a geodesic distance 
$d_{\Omega_{\mathrm{tcN}}^1,K}(\cdot,\cdot)$ between 
$\mathbf{U}^0$ and $\mathbf{U}^1$ as 
the infimum over all piecewise regular paths $\mathbf{U}^\theta$ satisfying Definition \ref{prp} with $k\geq 13$ and
connecting $\mathbf{U}^0$ and $\mathbf{U}^1$:
\begin{equation}
	\label{d-Om-r}
	d_{\Omega_{\mathrm{tcN}}^1,K}
	\left(\mathbf{U}^0,\mathbf{U}^1\right)
	=\inf\limits_{\mathbf{U}^\theta}
	\left\{\left\|\mathbf{U}^\theta\right\|_{\mathcal{L}}:
	\sup\limits_{\theta\in[0,1]}\max\left\{
	E_u^\theta(0),E_v^\theta(0),H^\theta(0)\right\}
	\leq K\right\},
\end{equation}
where the length $\left\|\mathbf{U}^\theta\right\|_{\mathcal{L}}$ 
of the path $\mathbf{U}^\theta$
is given in Definition \ref{lPR},
while $E_u^\theta(0)$, $E_v^\theta(0)$, and
$H^\theta(0)$ are defined as in \eqref{consqU} with
$\mathbf{U}^\theta$ 
instead of $\mathbf{U}$.
\end{definition}

We make two observations about Definition \ref{Gdist}.
First, $d_{\Omega_{\mathrm{tcN}}^1,K}(\cdot,\cdot)$, defined in
\eqref{d-Om-r}, is a pseudometric: it does not distinguish
points in the same equivalence class under
Definition \ref{Drelab}.
Consequently, the corresponding metric is defined on the quotient
$\mathcal{C}^k_{\Omega_{\mathrm{tcN}}^1}/\mathord{\sim}$.

Second, strictly speaking,
we can define \eqref{d-Om-r} only if there exists at least one
piecewise regular path $\mathbf{U}^\theta$ connecting
$
\mathbf{U}^0,\mathbf{U}^1\in
\mathcal{C}^k_{\Omega_{\mathrm{tcN}}^1}.
$
Note that the straight-line path
$\mathbf{U}^{\theta}=(1-\theta)\mathbf{U}^0+\theta\mathbf{U}^1$,
$\theta\in[0,1]$, does not, in general, satisfy the compatibility
conditions \eqref{ccUV} for $\theta\neq0,1$.
If the endpoints $\mathbf{U}^i$, $i=0,1$, are such that
$(U^i,V^i,W^i,Z^i)(\xi)$ are compactly supported and neither
$W^i(\xi)$ nor $Z^i(\xi)$ attains the value $\pi$,
we can construct a path as a map of the straight-line 
path 
\[(u^\theta,v^\theta)(\xi)=(1-\theta)(u^0,v^0)(\xi)+\theta(u^1,v^1)(\xi),
\]
in the Eulerian variables.
Since $(u^i,v^i)$, $i=0,1$, correspond to $\mathbf{U}^i$ and are
compactly supported, we have
$(u^\theta,v^\theta)\in\Sigma$ for all $\theta\in[0,1]$, although
$\Sigma$ is not a vector space (recall \eqref{Sig}).
Moreover, since $W^i,Z^i\neq\pi$, the functions $(u^i,v^i)$ are
regular and exhibit no energy concentration.
For general endpoints, we define the metric
by completion arguments, as detailed in Proposition \ref{COm}
and Definition \ref{OGdist} below.
	
Finally, the metric 
$d_{\Omega_{\mathrm{tcN}}^1,K}
\left(\cdot,\cdot\right)$ depends on the energy constant $K$.
Thus, we always choose $K>0$ large enough to satisfy \eqref{K-unb}
for the endpoints under consideration.

Theorems \ref{dpath} and \ref{egammat} allow 
us to establish the fundamental 
Lipschitz property of the 
pseudometric $d_{\Omega_{\mathrm{tcN}}^1,K}(\cdot,\cdot)$
from Definition \ref{Gdist}.

\begin{theorem}\label{LRO}
Consider initial data $\mathbf{U}^0_0,\mathbf{U}^1_0\in
\mathcal{C}^k_{\Omega_{\mathrm{tcN}}^1}$, 
$k\geq 13$,
where $\mathcal{C}^k_{\Omega_{\mathrm{tcN}}^1}$
is defined in \eqref{CktcN}, 
and the corresponding global solutions 
$\mathbf{U}^0(t),\mathbf{U}^1(t)\in
\mathcal{C}^k_{\Omega_{\mathrm{tcN}}^1}$
of the ODE system \eqref{ODE}  given in Theorem \ref{Ckgl}.
Then the geodesic distance 
$d_{\Omega_{\mathrm{tcN}}^1,K}(\cdot,\cdot)$
satisfies the Lipschitz property
\begin{equation}\label{LUt}
	d_{\Omega_{\mathrm{tcN}}^1,K}\left(
	\mathbf{U}^0(t),\mathbf{U}^1(t)\right)
	\leq C d_{\Omega_{\mathrm{tcN}}^1,K}
	\left(\mathbf{U}_0^0,\mathbf{U}^1_0\right),
	\quad C=C(T,K,\alpha)>0,
\end{equation}
for all $t\in[-T,T]$, where $K$ is as in
Definition \ref{Gdist}.
\end{theorem}

\begin{proof}
By Definition \ref{Gdist}, for any $\ve>0$ there exists a 
piecewise regular initial path
$\mathbf{U}_0^\theta$ (see Definition \eqref{prp}) 
connecting $\mathbf{U}_0^0$ and $\mathbf{U}_0^1$ such that
(recall \eqref{npath})
\begin{equation*}
	\left|
	d_{\Omega_{\mathrm{tcN}}^1,K}\left(\mathbf{U}_0^0,\mathbf{U}_0^1\right)
	-\left\|\mathbf{U}_0^\theta\right\|_{\mathcal{L}}
	\right|<\ve,
\end{equation*}
where $\|\cdot\|_{\mathcal{L}}$ 
is defined in Definition \ref{lPR}.
Then, applying Theorem \ref{dpath} to every path
$\mathbf{U}_0^\theta\in\mathcal{P}^{k,1}_{\mathrm{tcN}}
[\theta_i^*,\theta_{i+1}^*]$,
$i=0,\dots,N^*-1$, we obtain a piecewise regular initial path
$\widehat{\mathbf{U}}_0^\theta$ under 
the ODE system \eqref{ODE} for which 
(recall that, with time regarded as fixed, both 
$\left(\mathbf{U}_0^\theta,\frac{d}{d\theta}\mathbf{U}_0^\theta\right)$ 
and $\left(\widehat{\mathbf{U}}_0^\theta, 
\frac{d}{d\theta}\widehat{\mathbf{U}}_0^\theta\right)$ 
are elements of 
${\mathcal{P}^{k,1}_{\mathrm{tcN}}}[\theta_i^*,\theta_{i+1}^*]$)
\begin{equation}\label{1ehU}
	\left\|\left(
	\widehat{\mathbf{U}}_0^\theta-\mathbf{U}_0^\theta,\,
	\frac{d}{d\theta}\left(\widehat{\mathbf{U}}_0^\theta
	-\mathbf{U}_0^\theta\right)
	\right)
	\right\|_{\mathcal{P}^{k,1}_{\mathrm{tcN}}[\theta_i^*,\theta_{i+1}^*]}
	<\ve,\quad 
	\widehat{\mathbf{U}}_0^{\theta_i^*}=\mathbf{U}_0^{\theta_i^*},
	\quad i=0,\dots,N^*-1,
\end{equation}
where we have used that Theorem \ref{dpath} preserves the endpoints,
see item (1).
Estimate \eqref{1ehU} shows that the functions
inside the integrals in \eqref{Rtr}--\eqref{npath} 
are uniformly approximated by
$\left(\widehat{\mathbf{U}}_0^\theta,
\frac{d}{d\theta}\widehat{\mathbf{U}}_0^\theta\right)$
on every interval $[\theta_i^*,\theta_{i+1}^*]$, $i=0,\dots,N^*$,
which gives the following estimate:
\begin{equation}\label{1hdOp}
	\left|
	d_{\Omega_{\mathrm{tcN}}^1,K}\left(\mathbf{U}_0^0,\mathbf{U}_0^1\right)
	-\left\|\widehat{\mathbf{U}}_0^\theta\right\|_{\mathcal{L}}
	\right|<C\ve,\quad C>0,\quad\text{and}\quad
	\widehat{\mathbf{U}}_0^i=\mathbf{U}_0^i,\,\, i=0,1.
\end{equation}

Applying Theorem \ref{egammat} to the piecewise regular 
path $\widehat{\mathbf{U}}_0^\theta$ under the ODE system \eqref{ODE},
we obtain
\begin{equation}\label{hegammat}
	\left\|\widehat{\mathbf{U}}^\theta(t)\right\|
	_{\mathcal{L}}\leq
	C(T,K,\alpha)\left\|\widehat{\mathbf{U}}^\theta_0\right\|
	_{\mathcal{L}},\quad
	\text{for all }\,t\in[-T,T],
\end{equation}
where $\widehat{\mathbf{U}}^\theta(t)$ is the evolution of
$\widehat{\mathbf{U}}^\theta_0$.
Let $\mathbf{U}^\theta(t)$ be 
the evolution of the initial path 
$\mathbf{U}_0^\theta$. Then the uniform bound given 
in item (3b) of Theorem \ref{dpath} implies that
\begin{equation}\label{inUh}
	\left\|\widehat{\mathbf{U}}^\theta(t)
	\right\|_\mathcal{L}\geq
	\left\|\mathbf{U}^\theta(t)\right\|_\mathcal{L}
	-C\ve\geq d_{\Omega_{\mathrm{tcN}}^1,K}
	\left(\mathbf{U}^0(t),\mathbf{U}^1(t)\right)
	-C\ve,\quad C>0.
\end{equation}
Combining \eqref{hegammat}, \eqref{1hdOp}, 
and \eqref{inUh}, we arrive at \eqref{LUt}.
\end{proof}

To establish the Lipschitz property \eqref{LUt}
for arbitrary initial data in $\Omega_{\mathrm{tcN}}$,
we extend 
$d_{\Omega_{\mathrm{tcN}}^1,K}(\cdot,\cdot)$, introduced 
in Definition \ref{Gdist}, from 
$\mathcal{C}^k_{\Omega_{\mathrm{tcN}}^1}$, $k\geq 13$, to 
the entire set $\Omega_{\mathrm{tcN}}^1$
via a completion argument.
To this end, we show that every element
$\mathbf{U}\in\Omega_{\mathrm{tcN}}^1$
admits a sequence
$\{\mathbf{U}_n\}_{n=1}^\infty\subset
\mathcal{C}^k_{\Omega_{\mathrm{tcN}}^1}$, 
that represents this element and is Cauchy in the geodesic distance.
We construct $\mathbf{U}_n$ so that
(i) the corresponding $(U_n,V_n,W_n,Z_n)(\xi)$ are compactly supported and
(ii) $W_n$ and $Z_n$ do not attain the singular value $\pi$
for all $\xi\in\R$.
The following proposition provides the details of these arguments.

\begin{proposition}\label{COm}
	For any 
	\begin{equation*}
		\mathbf{U}\in\Omega_{\mathrm{tcN}}^1,\quad
		\mathbf{U}=(U,V,W,Z,q),
	\end{equation*}
	there exists a sequence 
	\begin{equation*}
		\{\mathbf{U}_n\}_{n=1}^\infty\subset
		\mathcal{C}^k_{\Omega_{\mathrm{tcN}}^1},\quad
		\mathbf{U}_n=(U_n,V_n,W_n,Z_n,q_n),\quad k\geq13,
	\end{equation*}
	such that
	\begin{enumerate}
	\item $\mathbf{U}_n$ is a Cauchy sequence
	for the geodesic distance:
	$d_{\Omega_{\mathrm{tcN}}^1,K}(\mathbf{U}_n,\mathbf{U}_m)\to0$,
	$n,m\to\infty$;
	\item $\mathbf{U}_n$ converges to $\mathbf{U}$ in
	the following sense:
	\begin{equation}\label{apprU}
		\begin{split}
			&U_n\to U,\quad V_n\to V,
			\quad\text{in }W^{1,1}(\R)\cap H^1(\R)\cap W^{1,4}(\R),\\
			&W_n\to W,\quad Z_n\to Z,
			\quad\text{in }L^1(\R)\cap L^2(\R),\\
			&q_n\to q,\quad\text{in }L^1(\R)\cap L^2(\R);
		\end{split}
	\end{equation}
	\item $W_n$, $Z_n$, and $q_n$ satisfy the following bounds:
	\begin{equation}\label{Linf-W_n}
		\|W_n\|_{L^\infty(\R)}<\pi,\quad
		\|Z_n\|_{L^\infty(\R)}<\pi,\quad
		\|q_n\|_{L^\infty(\R)}<C,\quad
		\text{for all }\,n\in\N,
	\end{equation}
	with a constant $C>0$ independent of $n$;
	\item the total support of $\mathbf{U}_n$
	satisfies the following property:
	\begin{equation}\label{suppUn}
		(U_n,V_n,W_n,Z_n)(\xi)=(0,0,0,0),\quad
		|\xi|\geq R_n+2;\quad
		q_n(\xi)=1,\quad |\xi|\geq R_n,\quad n\in\N,
	\end{equation}
	with some $R_n\to\infty$.
	\end{enumerate}
Moreover, if $\{\mathbf U_n\}_{n=1}^\infty$ and
$\{\widehat{\mathbf U}_n\}_{n=1}^\infty$ are two
such sequences representing the same $\mathbf U$, then
\[d_{\Omega_{\mathrm{tcN}}^1,K}
\left(\mathbf U_n,\widehat{\mathbf U}_n\right)\to0.\]
\end{proposition}
\begin{proof}
	\textbf{Step 1.}
Let us show that for any
$\mathbf{U}\in\Omega_{\mathrm{tcN}}^1$ there exists a sequence  
$\{\mathbf{U}_n\}_{n\in\N}\subset
\mathcal{C}^k_{\Omega_{\mathrm{tcN}}^1}$ 
satisfying items (2), (3), and (4) of the proposition.
	The main challenge is to prove that $\mathbf{U}_n$ satisfies the
compatibility conditions \eqref{ccUV} for all $n\in\N$.

Since $W$ and $Z$ in $\mathbf{U}$ can be discontinuous, we may
always choose representatives satisfying
$W(\xi),Z(\xi)\in(-\pi,\pi]$ for all $\xi\in\R$; see
Remark \ref{WZval}.
Standard approximation arguments then imply that, for any
$R_n\to\infty$,
there exist functions $W_n^0$, $Z_n^0$, and
$q_n^0$ in $C^{k-1}(\R)$ such that
\begin{enumerate}[(i)]
	\item $W_n^0(\xi)=Z_n^0(\xi)=0$ and $q_n^0(\xi)=1$ for all
	$|\xi|\geq R_n$;
	\item $\|W_n^0-W\|_{L^1(\R)\cap L^2(\R)}
	+\|Z_n^0-Z\|_{L^1(\R)\cap L^2(\R)}
	+\|q_n^0-q\|_{L^1(\R)}\to0$;
	\item $\|W_n^0\|_{C(\R)}<\pi$ and 
	$\|Z_n^0\|_{C(\R)}<\pi$ for all $n\in\N$.
\end{enumerate}
The last property can be ensured by multiplying the approximants
by $(1-\ve_n)$ where $\ve_n\to0$ as $n\to\infty$ sufficiently fast.

\medskip 

\textbf{Step 1.1.}
Introduce the constants $c_{i,n}$, $i=1,2$, as follows (see \eqref{ccUV}):
\begin{equation}\label{ci}
\begin{split}
	&c_{1,n}=\frac{1}{2}\int_{-\infty}^\infty
	\left(q_n^0\sin W_n^0\cos^2\frac{Z_n^0}{2}\right)(\xi)\,d\xi,\\
	&c_{2,n}=\frac{1}{2}\int_{-\infty}^\infty
	\left(q_n^0\cos^2\frac{W_n^0}{2}\sin Z_n^0\right)(\xi)\,d\xi.
\end{split}
\end{equation}
Let us show that 
\begin{equation}\label{lim-ci}
	|c_{i,n}|\to0,\quad i=1,2,\quad n\to\infty.
\end{equation}
Since $\mathbf{U}$ satisfies \eqref{ccUV} and
$U\in H^1(\R)$,
\begin{equation}\label{lccUV}
	\frac{1}{2}\left(\int_{-\infty}^{-R_n}+\int_{R_n}^\infty\right)
	\left(q\sin W\cos^2\frac{Z}{2}\right)\to0,\quad n\to\infty,
\end{equation}
for any $R_n\to\infty$.
Combining \eqref{ci} and \eqref{lccUV}, we obtain
the following estimate (see item (1) above):
\begin{equation}\label{est-ci}
\begin{split}
	2|c_{1,n}|&=\left|
	\int_{-\infty}^\infty
	\left(q_n^0\sin W_n^0\cos^2\frac{Z_n^0}{2}\right)(\xi)\,d\xi
	-\int_{-\infty}^\infty
	\left(q\sin W\cos^2\frac{Z}{2}\right)(\xi)\,d\xi\right|\\
	&=\left|
	\int_{-R_n}^{R_n}
	\left(q_n^0\sin W_n^0\cos^2\frac{Z_n^0}{2}
	-q\sin W\cos^2\frac{Z}{2}
	\right)(\xi)\,d\xi\right|+\so(1)
	=\left|I^0_n\right|+\so(1),
\end{split}
\end{equation}
as $n\to\infty$. The integral $I^0_n$ can be estimated as follows
(here, $L^{p}=L^p(\R)$, $p=1,\infty$):
\begin{equation}\label{I_n^0-est}
\begin{split}
	\left|I^0_n\right|&\leq\int_{-R_n}^{R_n}
	\left|\left(q_n^0-q\right)
	\sin W_n^0\cos^2\frac{Z_n^0}{2}\right|(\xi)\,d\xi
	+\|q\|_{L^\infty}
	\int_{-R_n}^{R_n}\left|\sin W_n^0\cos^2\frac{Z_n^0}{2}
	-\sin W\cos^2\frac{Z}{2}\right|(\xi)\,d\xi\\
	&\leq\left\|q_n^0-q\right\|_{L^1}
	+\|q\|_{L^\infty}\int_{-R_n}^{R_n}\left|W_n^0-W\right|(\xi)\,d\xi
	+\|q\|_{L^\infty}\int_{-R_n}^{R_n}
	\left|\cos^2\frac{Z_n^0}{2}-\cos^2\frac{Z}{2}\right|(\xi)\,d\xi.
\end{split}
\end{equation}
Observe that (recall \eqref{sinx}; cf.~\eqref{I_n} above)
\begin{equation}\label{cosZ_n^0}
	\int_{-R_n}^{R_n}
	\left|\cos^2\frac{Z_n^0}{2}-\cos^2\frac{Z}{2}\right|(\xi)\,d\xi
	\leq\frac{1}{8}\int_{-R_n}^{R_n}
	\left|Z_n^0-Z\right|(\xi)\left|Z_n^0+Z\right|(\xi)\,d\xi
	\leq C\|Z_n^0-Z\|_{L^2(\R)},
\end{equation}
where in the last estimate we have used the Cauchy-Schwarz inequality.
Combining \eqref{I_n^0-est}, \eqref{cosZ_n^0}, and that
$\|W_n^0-W\|_{L^1(\R)}+\|q_n^0-q\|_{L^1(\R)}\to0$,
we conclude that $|I_n^0|\to0$, as $n\to\infty$, which, 
together with \eqref{est-ci}, yields \eqref{lim-ci} for $i=1$.
Arguing similarly for $i=2$, we obtain \eqref{lim-ci}.
 
\medskip

\textbf{Step 1.2.}
Using arguments similar to those in the proof of Lemma \ref{L2},
we ``perturb'' $W_n^0$ and $Z_n^0$ to define
functions $U_n$ and $V_n$ in $H^1(\R)$ through the
compatibility conditions \eqref{ccUV}.
Introduce the following smooth bump functions of fixed shape,
beyond the total support where $(W_n^0,Z_n^0,q_n^0)=(0,0,1)$:
\begin{equation}\label{psi_n}
\begin{split}
	&\psi_n(\xi)=\psi(\xi-R_n-1),\quad
	\chi_n(\xi)=\chi(\xi-R_n-1),\quad n\in\N,\\
	&\text{with }\,\psi,\chi\in C^\infty(\R)
	\,\text{ such that }\,
	\psi(\xi)=\chi(\xi)=0,\,\,\,|\xi|\geq1,
\end{split}
\end{equation}
and normalize them by
\begin{equation}\label{nzpsi_n}
	\int_{\R}\psi_n\,d\xi
	=\int_{\R}\chi_n\,d\xi=2.
\end{equation}
Define (cf.~\eqref{WZ-pert}, \eqref{q-pert})
\begin{equation}\label{W_n}
	W_n(\xi;a)=W_n^0+a\psi_n(\xi),\quad
	Z_n(\xi;b)=Z_n^0+b\chi_n(\xi),\quad
	q_n(\xi)=q_n^0(\xi),
\end{equation}
where the constants $a,b\in\R$,
$(a,b)=(a_n,b_n)$, are determined from the following
two-dimensional equation (cf.~\eqref{Hnuab} and \eqref{ccUV}):
\begin{equation}\label{H0ab}
H_n^0(a,b)=\frac{1}{2}\begin{pmatrix}
	\int_{-\infty}^\infty q_n(\xi)\sin W_n(\xi;a)
	\cos^2\frac{Z_n(\xi;b)}{2}\,d\xi\\
	\int_{-\infty}^\infty q_n(\xi)\cos^2\frac{W_n(\xi;a)}{2}
	\sin Z_n(\xi;b)\,d\xi
\end{pmatrix}=
\begin{pmatrix}
	0\\0
\end{pmatrix}.
\end{equation}
Taking into account that
(see \eqref{nzpsi_n}, \eqref{psi_n}, and item (i) above)
\begin{equation*}
	\det\left.D_{(a,b)}H_n^0(a,b)\right|_{(a,b)=(0,0)}
	=\det\frac{1}{2}\begin{pmatrix}
		\int_{-\infty}^\infty\psi_n(\xi)\,d\xi& 0\\
		0 &\int_{-\infty}^\infty\chi_n(\xi)\,d\xi
	\end{pmatrix}=1,
\end{equation*}
we conclude that there exists a unique $(a,b)=(a_n,b_n)$ such that
$H_n^0(a_n,b_n)=0$.
Recalling \eqref{lim-ci}, we have that
$H_n^0(0,0)=(c_{1,n},c_{2,n})\to(0,0)$ and therefore
$(a_n,b_n)\to(0,0)$ as well.
Thus, using item (ii) and \eqref{W_n}, we have that
(cf.~\eqref{apprU}; here $L^p=L^p(\R)$, $p=1,2$)
\begin{equation}\label{lim-W_n}
	\|W_n-W\|_{L^1\cap L^2}+\|Z_n-Z\|_{L^1\cap L^2}+
	\|q_n-q\|_{L^1\cap L^2}\to0,\quad n\to\infty.
\end{equation}

Define the compactly supported functions $U_n, V_n$
from the compatibility conditions as follows (see \eqref{W_n}):
\begin{equation*}
\begin{split}
	&U_n(\xi)=\frac{1}{2}\int_{-\infty}^\xi
	\left(q_n\sin W_n\cos^2\frac{Z_n}{2}\right)(\eta)\,d\eta,\\
	&V_n(\xi)=\frac{1}{2}\int_{-\infty}^\xi
	\left(q_n\cos^2\frac{W_n}{2}\sin Z_n\right)(\eta)\,d\eta.\\
\end{split}
\end{equation*}
By construction, the sequence $(U_n,V_n,W_n,Z_n,q_n)$
satisfies the compatibility conditions \eqref{ccUV} for all $n\in\N$.
Moreover, their supports satisfy \eqref{suppUn}.

\medskip

\textbf{Step 1.3.}
It remains to prove that $U_n\to U$ and $V_n\to V$
in $W^{1,1}\cap H^1(\R)\cap W^{1,4}(\R)$.
The limit \eqref{lim-W_n} immediately implies that
$\pxi U_n\to\pxi U$ and $\pxi V_n\to\pxi V$ in $L^1(\R)$.
Recalling that $W,Z,q\in L^\infty(\R)$ and \eqref{Linf-W_n},
we have these limits in $L^2(\R)\cap L^4(\R)$ as well.

Applying the one-dimensional Poincar\'e inequality anchored at the
right endpoint, we obtain the following estimates:
\begin{equation}\label{estUn}
\begin{split}
	\|U_n-U\|_{L^1(-R_n-2,R_n+2)}
	\leq&\,2(R_n+2)|U_n-U|(R_n+2)\\
	&+2(R_n+2)\|\pxi U_n-\pxi U\|_{L^1(-R_n-2,R_n+2)},
\end{split}
\end{equation}
and
\begin{equation}\label{estVn}
\begin{split}
	\|V_n-V\|_{L^1(-R_n-2,R_n+2)}
	\leq&\,2(R_n+2)|V_n-V|(R_n+2)\\
	&+2(R_n+2)\|\pxi V_n-\pxi V\|_{L^1(-R_n-2,R_n+2)},
\end{split}
\end{equation}
Using that $|U_n-U|+|V_n-V|$ belongs to $L^1(\R)$, we conclude that
there exists a sequence $\{R_n\}_{n=1}^\infty$, $R_n\to\infty$,
such that
\begin{equation}\label{lim-Un}
\left.\xi\left(|U_n-U|+|V_n-V|\right)(\xi)\right|_{\xi=(R_n+2)}\to0,\quad n\to\infty.
\end{equation}
Also, we can choose $q_n^0=q_n$ such that
\[
\|q_n-q\|_{(L^1\cap L^2)(-R_n,R_n)}\leq R_n^{-p},
\quad\text{with any }\,p>0,
\]
which, together with \eqref{suppUn}, \eqref{estUn}, \eqref{estVn}, and
\eqref{lim-Un}, yields that
\[
	\|U_n-U\|_{L^1(\R)}+\|V_n-V\|_{L^1(\R)}\to0,\quad n\to\infty,
	\]
	and we have established items (2) and (3).

\medskip

\textbf{Step 2.}
Our goal is to show that 
$\{\mathbf{U}_n\}_{n\in\N}$ is a Cauchy 
sequence with respect to the metric 
$d_{\Omega_{\mathrm{tcN}}^1,K}(\cdot,\cdot)$.
Note that the main difficulty here is that we cannot consider
a straight-line path
\begin{equation*}
	\widetilde{\mathbf{U}}_{n,m}^\theta
	=(1-\theta)\mathbf{U}_n+\theta\mathbf{U}_m,
	\quad \theta\in[0,1],
\end{equation*}
since $\widetilde{\mathbf{U}}_{n,m}^\theta$ does not satisfy the compatibility
condition \eqref{ccUV} for $\theta\neq0,1$.
Therefore, we introduce the following ``perturbed'' straight-line
paths for $W_n$, $Z_n$, and $q_n$ (cf.~\eqref{WZ-pert}
and Remark \ref{altp} below):
\begin{equation}\label{thetaW}
\begin{split}
	&W^\theta_{n,m}(\xi;a)=(1-\theta)W_n(\xi)
	+\theta W_m(\xi)+a\psi_{n,m}(\xi),\\
	&Z^\theta_{n,m}(\xi;b)=(1-\theta)Z_n(\xi)
	+\theta Z_m(\xi)+b\chi_{n,m}(\xi),\\
	&q^\theta_{n,m}(\xi)=(1-\theta)q_n(\xi)+\theta q_m(\xi),
\end{split}
\end{equation}
where the fixed-shape translated bumps
$\psi_{n,m}$ and $\chi_{n,m}$ are defined as follows
(see \eqref{psi_n}--\eqref{nzpsi_n}):
\begin{equation}\label{psi_nm}
	\psi_{n,m}(\xi)=\psi(\xi-\max\{R_n,R_m\}-1),\quad
	\chi_{n,m}(\xi)=\chi(\xi-\max\{R_n,R_m\}-1),\quad
	n,m\in\N.
\end{equation}
The unknown functions $a=a_{n,m}(\theta)$ and $b=b_{n,m}(\theta)$
are determined through the following system of equations (cf.~\eqref{H0ab}):
\begin{equation}\label{H_nab}
	H_{n,m}(a,b,\theta)=\frac{1}{2}
	\begin{pmatrix}
	\int_{-\infty}^\infty q^\theta_{n,m}(\xi)
	\sin W^\theta_{n,m}(\xi;a)
	\cos^2\frac{Z^\theta_{n,m}(\xi;b)}{2}\,d\xi\\
	\int_{-\infty}^\infty q^\theta_{n,m}(\xi)
	\cos^2\frac{W^\theta_{n,m}(\xi;a)}{2}
	\sin Z^\theta_{n,m}(\xi;b)\,d\xi
	\end{pmatrix}=
	\begin{pmatrix}
		0\\0
	\end{pmatrix}.
\end{equation}
Taking into account that
\begin{equation}\label{DH_nab}
	\det\left.D_{(a,b)}H_{n,m}(a,b,\theta)\right|_{(a,b,\theta)=(0,0,\theta)}
	=\det\frac{1}{2}\begin{pmatrix}
		\int_{-\infty}^\infty\psi_{n,m}(\xi)\,d\xi& 0\\
		0 &\int_{-\infty}^\infty\chi_{n,m}(\xi)\,d\xi
	\end{pmatrix}=1,\quad\theta\in[0,1],
\end{equation}
we conclude that we can define the functions $a_{n,m}(\theta)$
and $b_{n,m}(\theta)$ near $(a,b)=0$.
To obtain $a_{n,m},b_{n,m}\in C^{k-1}[0,1]$, it suffices to
prove that the integrals in \eqref{H_nab} are small
for $a=b=0$ uniformly in $\theta\in[0,1]$.

Using that
\[
\int_{-\infty}^\infty\left(q_n\sin W_n\cos^2\frac{Z_n}{2}
\right)(\xi)\,d\xi=0,\quad n\in\N,
\]
we have the following equation for the first
component of \eqref{H_nab}:
\begin{equation*}
\begin{split}
	2H_{n,m}(0,0,\theta)&=
	2H_{n,m}(0,0,\theta)
	-(1-\theta)\int_{-\infty}^\infty\left(q_n\sin W_n\cos^2\frac{Z_n}{2}
	\right)(\xi)\,d\xi\\
	&\quad-\theta\int_{-\infty}^\infty\left(q_m\sin W_m\cos^2\frac{Z_m}{2}
	\right)(\xi)\,d\xi
	=(1-\theta)I_{n,m}^{(1)}(\theta)+\theta I_{n,m}^{(2)}(\theta),
\end{split}
\end{equation*}
with
\begin{equation}\label{Inmi}
\begin{split}
	&I_{n,m}^{(1)}(\theta)=\int_{-\infty}^\infty q_n\left(
	\sin\left((1-\theta)W_n+\theta W_m\right)
	\cos^2\frac{(1-\theta)Z_n+\theta Z_m}{2}
	-\sin W_n\cos^2\frac{Z_n}{2}
	\right)(\xi)\,d\xi,\\
	&I_{n,m}^{(2)}(\theta)=\int_{-\infty}^\infty q_m\left(
	\sin\left((1-\theta)W_n+\theta W_m\right)
	\cos^2\frac{(1-\theta)Z_n+\theta Z_m}{2}
	-\sin W_m\cos^2\frac{Z_m}{2}
	\right)(\xi)\,d\xi.
\end{split}
\end{equation}
Recalling that $\|q_n\|_{L^\infty(\R)}<C$ for all $n\in\N$
and $\|W_n-W\|_{L^1(\R)}+\|Z_n-Z\|_{L^1(\R)}\to0$, we conclude that
\begin{equation}\label{limInmi}
	I_{n,m}^{(i)}(\theta)\to0,\quad i=1,2,\quad n,m\to\infty,
\end{equation}
uniformly in $\theta\in[0,1]$.
Thus, we can define the desired functions
$a_{n,m},b_{n,m}\in C^{k-1}[0,1]$.
Moreover, applying similar arguments for
$\partial_\theta H_{n,m}(0,0,\theta)$, we obtain that
\begin{equation}\label{anmC1}
	\|a_{n,m}\|_{C^1[0,1]}+\|b_{n,m}\|_{C^1[0,1]}\to0,\quad n,m\to\infty.
\end{equation}

Finally, in view of \eqref{H_nab}, we can define 
the functions $U^\theta_{n,m}$ and $V^\theta_{n,m}$
from the compatibility condition as follows:
\begin{equation}\label{Uthetanm}
\begin{split}
	&U_{n,m}^\theta(\xi)=\frac{1}{2}\int_{-\infty}^\xi
	\left(q_{n,m}^\theta\sin W_{n,m}^\theta
	\cos^2\frac{Z_{n,m}^\theta}{2}\right)(\eta)\,d\eta,\\
	&V_{n,m}^\theta(\xi)=\frac{1}{2}\int_{-\infty}^\xi
	\left(q_{n,m}^\theta\cos^2\frac{W_{n,m}^\theta}{2}
	\sin Z_{n,m}^\theta\right)(\eta)\,d\eta,\\
\end{split}
\end{equation}
This defines the path
\begin{equation}\label{theta-U-nm}
	\mathbf{U}_{n,m}^\theta(\xi)
	=\left(U_{n,m}^\theta,V_{n,m}^\theta,
	W_{n,m}^\theta,Z_{n,m}^\theta,q_{n,m}^\theta\right)(\xi),
	\quad (\xi,\theta)\in\R\times[0,1],
\end{equation}
between $\mathbf{U}_n$ and
$\mathbf{U}_m$; it is compactly supported
for all $n,m,\theta$.

\medskip

\textbf{Step 3.}
Using Definition~\ref{Gdist}, we obtain
\begin{equation}\label{infund}
	d_{\Omega_{\mathrm{tcN}}^1,K}(\mathbf{U}_n,\mathbf{U}_m)
	\leq
	\int_0^1
	\inf_{\eta\in C^\infty(\R)}
	\sum_{i=1}^{6}
	\int_{-\infty}^{\infty}
	e^{-\alpha|y_{n,m}^\theta(\xi)|}\,
	\left|\bigl(\phi_i^\theta\bigr)_{n,m}(\xi)\right|
	\,d\xi\,d\theta,
\end{equation}
where $y_{n,m}^\theta$ and 
$\bigl(\phi_i^\theta\bigr)_{n,m}$, $i=1,\dots,6$, 
are defined by \eqref{y^ve} and \eqref{J-alph} with
$\mathbf{U}_{n,m}^\theta$ given in \eqref{theta-U-nm}
and (recall \eqref{URbf})
\begin{equation}\label{Rbftheta}
	\mathbf{R}_{n,m}^\theta(\xi)=
	\partial_{\theta}\mathbf{U}_{n,m}^\theta(\xi)
	=\left(R_{n,m}^\theta,S_{n,m}^\theta,
	A_{n,m}^\theta,B_{n,m}^\theta,Q_{n,m}^\theta\right)(\xi),
\end{equation}
in place of 
$\mathbf{U}$ and $\mathbf{R}$, respectively.
Taking $\eta\equiv0$ in \eqref{infund} 
and using \eqref{J-alph},
we obtain that
\begin{equation}\label{D-Om-est}
	d_{\Omega_{\mathrm{tcN}}}(\mathbf{U}_n,\mathbf{U}_m)
	\leq\int_0^1\int_{-\infty}^\infty
	e^{-\alpha|y_{n,m}^\theta(\xi)|}
	\Phi_{n,m}^\theta(\xi)\,d\xi,
\end{equation}
with (see \eqref{Rbftheta} and \eqref{zve})
\begin{equation*}
	\Phi_{n,m}^\theta(\xi)
	=\left(\left|z_{n,m}^\theta\right|
	+\left|R_{n,m}^\theta\right|
	+\left|S_{n,m}^\theta\right|
	+\frac{1}{2}\left|A_{n,m}^\theta\right|
	+\frac{1}{2}\left|B_{n,m}^\theta\right|\right)(\xi)
	q_{n,m}^\theta(\xi)
	+\left|Q_{n,m}^\theta\right|(\xi).
\end{equation*}

Using \eqref{pxiy} for $y_{n,m}^\theta$, we have that
(we omit the argument $\xi$)
\begin{equation*}
	\pxi y_{n,m}^\theta-1
	=q_{n,m}^\theta-1
	-q_{n,m}^\theta\sin^2\frac{W_{n,m}^\theta}{2}
	-q_{n,m}^\theta\sin^2\frac{Z_{n,m}^\theta}{2}
	+q_{n,m}^\theta\sin^2\frac{W_{n,m}^\theta}{2}
	\sin^2\frac{Z_{n,m}^\theta}{2},
\end{equation*}
and therefore (see \eqref{sinx} and \eqref{apprU})
\begin{equation}\label{ynmth-est}
	|y_{n,m}^\theta(\xi)-\xi|\leq
	\|q_{n,m}^\theta-1\|_{L^1(\R)}
	+C\left(\left\|W_{n,m}^\theta\right\|_{L^2(\R)}^2
	+\left\|Z_{n,m}^\theta\right\|_{L^2(\R)}^2\right)
	\leq C,
\end{equation}
with some $C=C(K)>0$, which does not depend on $n,m,\theta$.
Then, taking into account that (recall \eqref{thetaW})
\begin{equation}\label{pthetaW}
\begin{split}
	&A_{n,m}^\theta(\xi)=\partial_\theta W_{n,m}^\theta(\xi)=
	(W_m-W_n)(\xi)+a_{n,m}^\prime(\theta)\psi_{n,m}(\xi),\\
	&B_{n,m}^\theta(\xi)=\partial_\theta Z_{n,m}^\theta(\xi)=
	(Z_m-Z_n)(\xi)+b_{n,m}^\prime(\theta)\chi_{n,m}(\xi),\\
	&Q_{n,m}^\theta(\xi)=\partial_\theta q_{n,m}^\theta(\xi)=
	(q_m-q_n)(\xi),
\end{split}
\end{equation}
and using \eqref{apprU}, \eqref{Linf-W_n}, \eqref{anmC1},
and \eqref{ynmth-est},
we obtain from \eqref{D-Om-est}
the following estimate:
\begin{equation}\label{2D-Om-est}
	d_{\Omega_{\mathrm{tcN}}^1}(\mathbf{U}_n,\mathbf{U}_m)
	\leq C(\alpha,K)\int_0^1\int_{-\infty}^\infty e^{-\alpha|\xi|}
	\left(\left|z_{n,m}^\theta\right|
	+\left|R_{n,m}^\theta\right|
	+\left|S_{n,m}^\theta\right|\right)(\xi)\,d\xi\,d\theta
	+\so(1),
\end{equation}
where $\so(1)\to0$ as $n,m\to\infty$.
Recalling \eqref{zve} for $z_{n,m}^\theta$ we have that (see \eqref{pthetaW})
\begin{equation}\label{zth-est}
\begin{split}
	\left|z_{n,m}^\theta(\xi)\right|&\leq\|q_m-q_n\|_{L^1(\R)}
	+C\int_{-\infty}^\infty
	\left|(W_m-W_n)(\xi)+a_{n,m}^\prime(\theta)\psi_{n,m}(\xi)\right|
	\left|W_{n,m}^\theta(\xi)\right|d\xi\\
	&\quad+C\int_{-\infty}^\infty
	\left|(Z_m-Z_n)(\xi)+b_{n,m}^\prime(\theta)\chi_{n,m}(\xi)\right|
	\left|Z_{n,m}^\theta(\xi)\right|d\xi\\
	&\leq C\left(\|q_m-q_n\|_{L^1(\R)}+
	\|W_m-W_n\|_{L^2(\R)}+\|Z_m-Z_n\|_{L^2(\R)}\right)\\
	&\quad+C|a_{n,m}^\prime(\theta)|
	\int_{-\infty}^\infty \left|\psi_{n,m}(\xi)\right|d\xi
	+C|b_{n,m}^\prime(\theta)|
	\int_{-\infty}^\infty \left|\chi_{n,m}(\xi)\right|d\xi,\quad C=C(K)>0,
\end{split}
\end{equation}
where in the last estimate we have applied the Cauchy-Schwarz inequality
as well as the uniform in $n,m,\theta$ bounds on
$\|W_{n,m}^\theta\|_{L^2(\R)\cap L^\infty(\R)}$ and
$\|Z_{n,m}^\theta\|_{L^2(\R)\cap L^\infty(\R)}$.
Using \eqref{psi_nm} (see also \eqref{psi_n}) as well as convergences in \eqref{apprU} and \eqref{anmC1}, we conclude from 
\eqref{zth-est} and \eqref{2D-Om-est} that
\begin{equation}\label{3D-Om-est}
	d_{\Omega_{\mathrm{tcN}}}(\mathbf{U}_n,\mathbf{U}_m)
	\leq C(\alpha,K)\int_0^1\int_{-\infty}^\infty e^{-\alpha|\xi|}
	\left(\left|R_{n,m}^\theta\right|
	+\left|S_{n,m}^\theta\right|\right)(\xi)\,d\xi\,d\theta
	+\so(1),
\end{equation}
with $\so(1)\to0$ as $n,m\to\infty$.

Finally, we estimate
$\int_{-\infty}^\infty\left|R_{n,m}^\theta(\xi)\right|$.
Employing \eqref{Uthetanm} and \eqref{pthetaW}, we have
the following inequalities:
\begin{equation}\label{Rthest}
\begin{split}
	\left|R_{n,m}^\theta(\xi)\right|&=
	\left|\partial_\theta U_{n,m}^\theta(\xi)\right|
	=\frac{1}{2}\left|\int_{-\infty}^\xi\left(Q_{n,m}^\theta
	\sin W_{n,m}^\theta\cos^2\frac{Z_{n,m}^\theta}{2}
	+q_{n,m}^\theta A_{n,m}^\theta
	\cos W_{n,m}^\theta\cos^2\frac{Z_{n,m}^\theta}{2}\right.\right.\\
	&\left.\left.\qquad\qquad\qquad\qquad\qquad\qquad
	+\frac{q_{n,m}^\theta}{2} B_{n,m}^\theta
	\sin W_{n,m}^\theta\sin Z_{n,m}^\theta
	\right)(\eta)\,d\eta\right|\\
	&\leq\frac{1}{2}\|q_m-q_n\|_{L^1(\R)}
	+C\|Z_m-Z_n\|_{L^2(\R)}\\
	&\quad+C\left(|b^\prime_{n,m}(\theta)|\cdot\|\chi_{n,m}\|_{L^1(\R)}
	+|a^\prime_{n,m}(\theta)|\cdot\|\psi_{n,m}\|_{L^1(\R)}\right)
	+\frac{1}{2}I_{n,m}^{\theta}(\xi),
\end{split}
\end{equation}
where $C=C(K)>0$ and
\begin{equation}\label{Inmth}
	I_{n,m}^\theta(\xi)=
	\left|\int_{-\infty}^\xi
	\left(q_{n,m}^\theta\cos W_{n,m}^\theta\cos^2\frac{Z_{n,m}^\theta}{2}
	\right)(\eta)\left(W_m-W_n\right)(\eta)\,d\eta\right|.
\end{equation}
Combining \eqref{apprU}, \eqref{psi_nm}, \eqref{anmC1}, and \eqref{Rthest},
we have that
\begin{equation}\label{Rnmthest}
	\left|R_{n,m}^\theta(\xi)\right|=
	\frac{1}{2}I_{n,m}^\theta(\xi)+\so(1),
\end{equation}
where $\so(1)\to0$ uniformly in $\theta$ as $n,m\to\infty$.
Finally, using that $W_n$ is a Cauchy sequence in $L^1(\R)$,
we conclude that $\left\|R_{n,m}^\theta(\xi)\right\|_{L^\infty(\R)}=\so(1)$.
Applying similar arguments for $\left|S_{n,m}^\theta(\xi)\right|$,
we conclude that $\mathbf{U}_n$ is Cauchy
in $d_{\Omega_{\mathrm{tcN}}^1,K}(\cdot,\cdot)$.

\medskip

\textbf{Step 4.}
Consider two sequences $\{\mathbf U_n\}_{n=1}^\infty$ and
$\{\widehat{\mathbf U}_n\}_{n=1}^\infty$
representing the same $\mathbf U$.
Employing the same estimates as in step 3 above 
with $\widehat{\mathbf U}_n$ in place of $\mathbf U_m$,
we conclude that
$d_{\Omega_{\mathrm{tcN}}^1,K}
(\mathbf U_n,\widehat{\mathbf U}_n)\to0$.
\end{proof}
\begin{remark}\label{compL1}
		The $L^1$ convergence of $W_n$ and $Z_n$ to $W$ and $Z$,
		respectively, is essential in \eqref{3D-Om-est},
	where we prove that
	$\left|R_{n,m}^\theta(\xi)\right|+\left|S_{n,m}^\theta(\xi)\right|\to0$, as $n,m\to\infty$.
	Note that we can choose $W_n^0$ and $Z_n^0$
		arbitrarily close in $L^2$ to $W$ and $Z$ on $(-R_n,R_n)$,
	for any sequence $R_n\to\infty$, and thus have \eqref{lim-ci}
	and \eqref{limInmi} from, respectively, \eqref{I_n^0-est}
	and \eqref{Inmi}.
\end{remark}

\begin{remark}[Step 2]\label{altp}
	An alternative path connecting $\mathbf{U}_n$ and $\mathbf{U}_m$
	can be defined as a direct transform of the straight-line path
	$\left(u_{n,m}^\theta,v_{n,m}^\theta\right)(x)=
	(1-\theta)\left(u_n,v_n\right)+\theta(u_m,v_m)(x)$
	in the Eulerian variables.
	Note that neither $W_n$ nor $Z_n$ attains $\pi$;
	the corresponding functions $u_n$ and $v_n$ are regular and
	compactly supported, and they exhibit no concentration of energy.
	However, 
	the Eulerian straight-line path does not treat a general
	element of $\Sigma$: from
	$u_x,v_x\in L^2$ and $u_xv_x\in L^2$ one cannot infer the separate
	$L^4$ bounds needed for the terms $f_4,f_5,f_6$ in \eqref{fi} with $h=0$
	and therefore the $W^{1,1}$ assumption does not allow us to close
	these estimates.
	One would also need to prove that $(u_n,v_n)$ is Cauchy
	in $\Sigma$, which is also challenging because energy may concentrate
	and the derivative of the corresponding characteristic may vanish
	in the limit $n\to\infty$.
\end{remark}

Introduce
\[
X=\mathcal C^k_{\Omega_{\mathrm{tcN}}^1}/\mathord{\sim},\quad k\geq13,
\quad\text{and }\,
\overline X\,\text{ as the completion of }\,
(X,d_{\Omega_{\mathrm{tcN}}^1,K}(\cdot,\cdot)),
\]
where $\sim$ is the relabeling equivalence relation described
in Definition \ref{Drelab}.
Applying Proposition \ref{COm}, we 
can identify each
$\mathbf U\in\Omega_{\mathrm{tcN}}^1/\mathord{\sim}$
with any sequence $\{\mathbf U_n\}_{n=1}^\infty$ constructed
in Proposition \ref{COm} representing this element:
\[
\iota(\mathbf U)=[\{\mathbf U_n\}_{n=1}^\infty]\in\overline X.
\]
Thus, the metric space
$\left(\mathcal{C}^k_{\Omega_{\mathrm{tcN}}^1},
d_{\Omega_{\mathrm{tcN}}^1,K}\right)$, $k\geq13$,
admits a completion containing the image
$\iota(\Omega_{\mathrm{tcN}}^1/\mathord{\sim})$. 
This allows us to extend the distance between any two elements 
of this set by continuity, as detailed below.

\begin{definition}[Geodesic distance in 
	$\Omega_{\mathrm{tcN}}^1/\mathord{\sim}$]
\label{OGdist}
We define $d_{\Omega_{\mathrm{tcN}}^1,K}(\cdot,\cdot)$
on $\Omega_{\mathrm{tcN}}^1/\mathord{\sim}$
as the metric induced by the
completion of the metric space
$\left(\mathcal{C}^k_{\Omega_{\mathrm{tcN}}^1}/\mathord{\sim},
d_{\Omega_{\mathrm{tcN}}^1,K}\right)$, $k\geq 13$,
where $d_{\Omega_{\mathrm{tcN}}^1,K}(\cdot,\cdot)$
is the geodesic distance introduced in
Definition~\ref{Gdist} and $\sim$ is the relabeling relation 
given in Definition \ref{Drelab}.
\end{definition}

Now the goal is to show that the metric
$d_{\Omega_{\mathrm{tcN}}^1,K}(\cdot,\cdot)$
defined on the whole $\Omega_{\mathrm{tcN}}^1/\mathord{\sim}$
satisfies the Lipschitz property \eqref{LUt}.
We begin with the following technical lemma,
	which ensures that
the continuous flow on the completion is exactly the ODE 
flow from Theorem~\ref{gwp}.
\begin{lemma}\label{Lstab}
Consider arbitrary initial data
$\mathbf U, \widehat{\mathbf U}\in\Omega_{\mathrm{tcN}}^1$
as well as the following metric (cf.~\eqref{apprU}):
\begin{equation}\label{dstar}
\begin{split}
	d_*\left(\mathbf U,\widehat{\mathbf U}\right)
	=\norm{U-\hat U}_{L^1(\R)}
	+\norm{V-\hat V}_{L^1(\R)}
	+\norm{W-\hat W}_{L^1(\R)}
	+\norm{Z-\hat Z}_{L^1(\R)}
	+\norm{q-\hat q}_{L^1(\R)}.
\end{split}
\end{equation}
Then the evolutions of $\mathbf U_0$ and $\widehat{\mathbf U}_0$
under the ODE system \eqref{ODE} satisfy the following stability estimate:
\begin{equation}\label{stabil}
\sup_{t\in[-T,T]}
d_*\left(\mathbf U(t),\widehat{\mathbf U}(t)\right)
\leq C(T,K)
d_*\left(\mathbf U,\widehat{\mathbf U}\right).
\end{equation}
\end{lemma}
\begin{proof}
	First, the $L^1$ convergences in \eqref{apprU}, together with
	the common $L^\infty$ bounds, imply convergence in the
	$L^2\cap L^4$ norms. By \cite[Lemma 4.1]{KR25},
	\begin{equation*}
		\mathcal{E}_{\mathbf{U}}(t,\xi,\eta),\,
		\mathcal{E}_{\widehat{\mathbf{U}}}(t,\xi,\eta)\leq
		\exp\left(\frac{q^-}{4}(R_2^2-|\xi|)\right),
	\end{equation*}
	where $\mathcal{E}_{\mathbf{U}}$ and
	$\mathcal{E}_{\widehat{\mathbf{U}}}$ are defined as in \eqref{E},
	using the functions from ${\mathbf{U}}$ and
	$\widehat{\mathbf{U}}$, respectively.
	Applying Young's inequality to all differences of
	the nonlocal terms $P_j$, $\px P_j$, $S_j$, and $\px S_j$, $j=1,2$,
	appearing in \eqref{ODE}, we can bound their $L^1(\R)$ norms by
	$d_*\left(\mathbf U(t),\widehat{\mathbf U}(t)\right)$
	(cf.~\cite[Proposition 4.2]{KR25}).
	Thus, considering the differences of the equations in \eqref{ODE}
	(cf.~the proof of \cite[Lemma 4.1]{HQ21}),
	we eventually obtain that
	\begin{equation*}
	\begin{split}
		\frac{d}{dt}d_*\left(\mathbf U(t),\widehat{\mathbf U}(t)\right)
		\leq C(K)d_*\left(\mathbf U(t),\widehat{\mathbf U}(t)\right),
	\end{split}
	\end{equation*}
	which yields \eqref{stabil} by Gronwall's inequality.
\end{proof}

Consider any $\mathbf U\in\Omega_{\mathrm{tcN}}^1$
and a regular representative sequence
$\{\mathbf U_n\}_{n=1}^\infty$
constructed in Proposition \ref{COm}.
Lemma \ref{Lstab} yields that
\begin{equation*}
	d_*\left(\mathbf U_n(t),\mathbf U(t)\right)
	\leq C(T,K)d_*\left(\mathbf U_n,\mathbf U\right),
\end{equation*}
and therefore $\mathbf U_n(t)\to\mathbf U(t)$ with respect to
$d_*(\cdot,\cdot)$.
Taking into account that (see \eqref{ccUV})
\begin{equation*}
\begin{split}
	\norm{U_n(t)- U(t)}_{*}
	+\norm{V_n(t)-V(t)}_{*}\leq C(T,K)
	\sum\limits_{p=1,2,4}d_*^{1/p}(\mathbf U_n(t),\mathbf U(t)),
\end{split}
\end{equation*}
where 
\[
\norm{U}_{*}=\norm{U}_{W^{1,1}(\R)\cap H^1(\R)\cap W^{1,4}(\R)},
\]
we conclude that $\{\mathbf U_n(t)\}_{n=1}^\infty$
is a representative of $\mathbf U(t)$,
that is
\begin{equation*}
	\iota(\mathbf U(t))=[\{\mathbf U_n(t)\}_{n=1}^\infty].
\end{equation*}
This, together with Theorem~\ref{LRO}, Proposition~\ref{COm},
and Definition~\ref{OGdist}
implies that the
Lipschitz property \eqref{LUt} holds for arbitrary initial data
$\mathbf{U}_0^0,\mathbf{U}_0^1\in
\Omega_{\mathrm{tcN}}^1/\mathord{\sim}$.

\subsection{Lipschitz metric in the original variables}
\label{LipDS}

In this final section we define a Lipschitz metric on 
$\mathcal{D}$, see \eqref{D-set}, by means of the metric 
$d_{\Omega_{\mathrm{tcN}}^1,K}(\cdot,\cdot)$ introduced for 
the transformed (ODE system) variables in Section \ref{GdO}.  
This will allow us to establish Theorem \ref{ThLip}.  
The definition of the new metric is as follows:

\begin{definition}[Lipschitz metric on $\mathcal{D}$]\label{Ldist}
Let $\mathbf{u},\hat{\mathbf{u}}
\in\mathcal{D}$, where 
\[
\mathbf{u}=\left(u,v,\mu;D_W,D_Z\right),
\qquad
\hat{\mathbf{u}}=\left(\hat{u},\hat{v},\hat{\mu};
\hat{D}_W,\hat{D}_Z\right),
\]
and $u,v,\hat{u},\hat{v}\in W^{1,1}(\R)$.
Consider $\mathbf{U}^0,\mathbf{U}^1\in\Omega_{\mathrm{tcN}}^1/\mathord{\sim}$ 
as the corresponding elements associated 
with $\mathbf{u}$ and $\hat{\mathbf{u}}$ via \eqref{id2}.  
We define a metric $d_{\mathcal{D},K}(\cdot,\cdot)$ 
on $\mathcal{D}$ by means of the 
geodesic distance $d_{\Omega_{\mathrm{tcN}}^1,K}(\cdot,\cdot)$ 
introduced in Definition~\ref{OGdist} as follows:
\begin{equation}\label{Lmds}	
	d_{\mathcal{D},K}\left(\mathbf{u},\hat{\mathbf{u}}\right)
	:=d_{\Omega_{\mathrm{tcN}}^1,K}\left(\mathbf{U}^0,\mathbf{U}^1\right).
\end{equation}
\end{definition}

We now show that the metric $d_{\mathcal{D},K}(\cdot,\cdot)$ 
introduced in Definition \ref{Ldist} ensures 
the Lipschitz continuity of global solutions to 
the two-component Novikov system.  
The central difficulty is that, in general, 
there is no bijection between the Euler 
variables and the Bressan-Constantin 
variables for $t\neq 0$ (see Figure \ref{inv-d-t}).  
Thus, we must explicitly verify that the distance 
between the Eulerian solutions 
$\mathbf{u}(t)$ and $\hat{\mathbf{u}}(t)$ 
coincides with the distance between 
the corresponding transformed solutions 
$\mathbf{U}^0(t)$ and $\mathbf{U}^1(t)$.  

Equivalently, since the distance between 
$\mathbf{u}(t)$ and $\hat{\mathbf{u}}(t)$ is, by 
Definition \ref{Ldist} and 
\eqref{Du(t)d}, equal to the geodesic distance between 
the associated transformed variables 
$\widetilde{\mathbf{U}}^0$ 
and $\widetilde{\mathbf{U}}^1$, 
it remains to prove that
\[
d_{\Omega_{\mathrm{tcN}}^1,K}\left(\widetilde{\mathbf{U}}^0,
\widetilde{\mathbf{U}}^1\right)
=
d_{\Omega_{\mathrm{tcN}}^1,K}\left(\mathbf{U}^0(t),
\mathbf{U}^1(t)\right),
\]
see \eqref{UU-teq} below.

\begin{figure}
\centering{\includegraphics[scale=0.5]{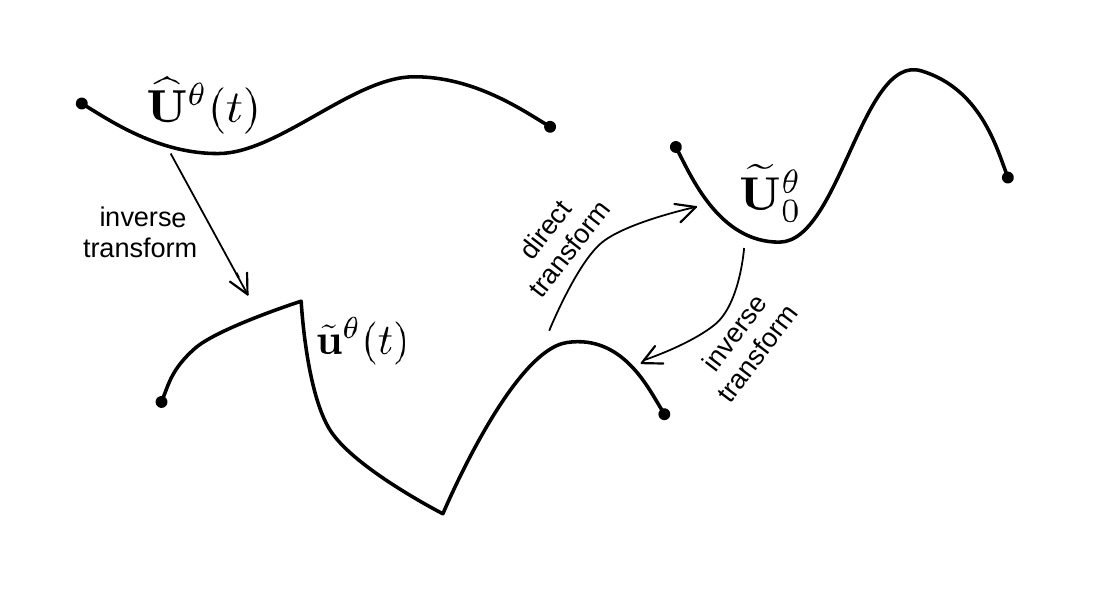}}
\caption{Mapping the path $\widehat{\mathbf{U}}^\theta(t)$ 
in the transformed variables to a path 
$\tilde{\mathbf{u}}^\theta(t)$ in the Euler 
variables ($t$ is considered fixed here), and 
then lifting it back to the Bressan-Constantin 
variables, produces the path $\widetilde{\mathbf{U}}^\theta_0$ 
(note that $\tilde{\mathbf{u}}^\theta(t)$ is a piecewise regular 
path in the sense of \cite[Definition 2]{BC17A}).
Because the transformation is not bijective, one 
generally has $\widetilde{\mathbf{U}}^\theta_0
\neq \widehat{\mathbf{U}}^\theta(t)$ for
$t\neq0$ and $\theta\in[0,1]$.
Nevertheless, since both 
$\widehat{\mathbf{U}}^\theta(t)$ 
and $\widetilde{\mathbf{U}}^\theta_0$ correspond 
to the same path $\tilde{\mathbf{u}}^\theta(t)$ in 
the Euler variables, they belong to the same equivalence
class (see Definition \ref{Drelab}), and
their lengths $\|\cdot\|_{\mathcal{L}}$ 
coincide (see Definition \ref{lPR}).}
\label{path-m}
\end{figure}

To this end, using \eqref{1hdOp}, we approximate 
$d_{\Omega_{\mathrm{tcN}}^1,K}\!\left(\mathbf{U}^0(t),\mathbf{U}^1(t)\right)$ 
by the length of a regular path 
$\widehat{\mathbf{U}}^\theta(t)$ under the ODE flow such that
\begin{equation}\label{appr-1}
	\left\|\widehat{\mathbf{U}}^\theta(t)\right\|_{\mathcal{L}}
	\approx 
	d_{\Omega_{\mathrm{tcN}}^1,K}\left(\mathbf{U}^0(t),\mathbf{U}^1(t)\right).
\end{equation}
Mapping $\widehat{\mathbf{U}}^\theta(t)$ back 
to the Euler variables yields a 
path $\tilde{\mathbf{u}}^\theta(t)$.  
When this path is transformed again to the 
Bressan-Constantin variables, we 
generally obtain a different path, denoted 
$\widetilde{\mathbf{U}}_0^\theta$ (see 
Figure \ref{path-m}), satisfying
\begin{equation}\label{appr-2}
	\left\|\widetilde{\mathbf{U}}_0^\theta\right\|_{\mathcal{L}}
	\approx 
	d_{\Omega_{\mathrm{tcN}}^1,K}\!\left(\widetilde{\mathbf{U}}^0,
	\widetilde{\mathbf{U}}^1\right).
\end{equation}

The crucial observation is that both paths, 
$\widetilde{\mathbf{U}}_0^\theta$ and 
$\widehat{\mathbf{U}}^\theta(t)$, represent 
the same intermediate states 
$\tilde{\mathbf{u}}^\theta(t)$ in the Euler 
variables and therefore belong to the same equivalence
class (see Definition \ref{Drelab}). Consequently,
$$
\left\|\widetilde{\mathbf{U}}_0^\theta\right\|_{\mathcal{L}}
=\left\|\widehat{\mathbf{U}}^\theta(t)\right\|_{\mathcal{L}}.
$$
Combining this fact with 
\eqref{appr-1} and \eqref{appr-2} establishes 
\eqref{UU-teq}.

Now, let us provide a detailed 
proof of Theorem \ref{ThLip}.

\begin{proof}[Proof of Theorem \ref{ThLip}]
It suffices to establish \eqref{Lipgs} for initial data 
$\mathbf{u}_0$, $\hat{\mathbf{u}}_0$ 
whose corresponding transformed data 
$\mathbf{U}_0^0,\mathbf{U}_0^1$ lie in 
$\mathcal{C}^k_{\Omega_{\mathrm{tcN}}^1}$, $k\geq13$
(recall that $q_0^0=q_0^1=1$, see \eqref{id2} for 
the ODE initial data). Then, using completion 
arguments as in Section \ref{GdO}, 
we obtain \eqref{Lipgs} for any 
$\mathbf{u}_0,\hat{\mathbf{u}}_0\in\mathcal{D}$.
Consider the global conservative solutions
\begin{equation}\label{eq:u-tildeu}
	\begin{split}
		\mathbf{u}(t)=
		\left(u(t),v(t),\mu_{(t)};
		D_W(t),D_Z(t)\right),\quad
		\hat{\mathbf{u}}(t)=
		\left(\hat{u}(t),\hat{v}(t),\hat{\mu}_{(t)};
		\hat{D}_W(t),\hat{D}_Z(t)\right),
	\end{split}
\end{equation}
given by Theorem \ref{Thm}, 
which correspond to the global solutions
$\mathbf{U}^0(t),\mathbf{U}^1(t)
\in\mathcal{C}^k_{\Omega_{\mathrm{tcN}}^1}$
of the associated ODE system \eqref{ODE}--\eqref{id2}. 
By definition \eqref{Lmds}, we have
\begin{equation}\label{Du(t)d}
	d_{\mathcal{D},K}\left(\mathbf{u}(t),
	\hat{\mathbf{u}}(t)\right)
	=d_{\Omega_{\mathrm{tcN}}^1,K}\left(
	\widetilde{\mathbf{U}}^0,
	\widetilde{\mathbf{U}}^1\right),
\end{equation}
where 
\begin{equation}\label{wtilU}
	\widetilde{\mathbf{U}}^i(\sigma)=
	\left(
	\tilde{U}^i,\tilde{V}^i,\tilde{W}^i,\tilde{Z}^i,1
	\right)(\sigma),\quad i=0,1,
\end{equation}
are defined as in \eqref{id2}, but now using the data at time $t$. 
More precisely, given \eqref{eq:u-tildeu}, the data 
$\widetilde{\mathbf{U}}^0$ and $\widetilde{\mathbf{U}}^1$ 
are obtained from $\mathbf{u}(t)$ and 
$\hat{\mathbf{u}}(t)$, respectively, in the same way as in \eqref{id2}, 
with $(u_0,v_0,\mu_0;D_{W,0},D_{Z,0})$ replaced 
by the corresponding time-dependent data; compare 
with Figure \ref{inv-d-t}, where 
$\mathbf{U}_0^0$, $\mathbf{U}^0(t)$, and $\widetilde{\mathbf{U}}^0$
play the roles of $\mathbf{U}_0$, $\mathbf{U}(t)$, and
$\widetilde{\mathbf{U}}$, respectively.)
Here, the new relabeling parameter $\sigma=\sigma(\xi)$ 
is given by \eqref{sigd}.

In general, 
$\widetilde{\mathbf{U}}^0\neq\mathbf{U}^0(t)$ and
$\widetilde{\mathbf{U}}^1\neq\mathbf{U}^1(t)$, since, 
for example, $\tilde{q}^0=1$ and $\tilde{q}^1=1$ 
in $\widetilde{\mathbf{U}}^0$ and 
$\widetilde{\mathbf{U}}^1$, respectively (see \eqref{wtilU}). 
Thus, we must show that
\begin{equation}\label{UU-teq}
	d_{\Omega_{\mathrm{tcN}}^1,K}\left(
	\widetilde{\mathbf{U}}^0,
	\widetilde{\mathbf{U}}^1\right)
	=d_{\Omega_{\mathrm{tcN}}^1,K}\left(
	\mathbf{U}^0(t),\mathbf{U}^1(t)\right).
\end{equation}

Following the proof of Theorem \ref{LRO}, we obtain
a piecewise regular path
$\widehat{\mathbf{U}}^\theta(t)$, $\theta\in[0,1]$, 
under the ODE system \eqref{ODE} such that its 
length approximates the geodesic distance between
$\mathbf{U}^0(t)$ and $\mathbf{U}^1(t)$ 
(see \eqref{1hdOp}), and the endpoints of the 
path $\widehat{\mathbf{U}}^\theta(t)$ 
coincide with
$\mathbf{U}^0(t)$ 
and $\mathbf{U}^1(t)$ (see \eqref{1ehU} and \eqref{1hdOp}).
More specifically, for any $\ve>0$ 
there exists a regular path 
$\widehat{\mathbf{U}}^\theta(t)$ 
under the ODE system \eqref{ODE} 
such that (recall \eqref{npath})
\begin{equation}\label{appep}
\begin{split}
	&\left|d_{\Omega_{\mathrm{tcN}}^1,K}\left(
	\mathbf{U}^0(t),\mathbf{U}^1(t)\right)
	-\left\|\widehat{\mathbf{U}}^\theta(t)\right\|
	_{\mathcal{L}}\right|<\ve,\quad
	\widehat{\mathbf{U}}^i(t)=\mathbf{U}^i(t),\quad i=0,1.
\end{split}
\end{equation}

We map $\widehat{\mathbf{U}}^\theta(t)$ into 
the Euler variables for each fixed 
$\theta$ using \eqref{uvdef}--\eqref{mut} and \eqref{DN-a}.  
This yields the path
\begin{equation}\label{t-bfu}
	\tilde{\mathbf{u}}^\theta(t)=
	\left(\tilde{u}^\theta(t),\tilde{v}^\theta(t),
	\tilde{\mu}_{(t)}^\theta;
	\tilde{D}^\theta_W(t),\tilde{D}^\theta_Z(t)\right),
\end{equation}
in the original variables (paths of this type are referred to as 
piecewise regular in \cite[Definition~2]{BC17A}).  
We then transform $\tilde{\mathbf{u}}^\theta(t)$ back to the 
Bressan-Constantin variables using \eqref{id2}, 
with the right-hand side of \eqref{t-bfu} 
in place of $(u_0,v_0,\mu_0;D_{W,0},D_{Z,0})$, 
thereby obtaining the path 
$\widetilde{\mathbf{U}}_0^\theta$ 
(see Figure \ref{path-m}).

Since $\widehat{\mathbf{U}}^\theta(t)$ is 
a regular path under the ODE system \ref{ODE}, 
the characteristic $y^\theta(t,\cdot)$ 
is strictly monotone for all 
$\theta\neq\theta_i$, $i=1,\dots,N$
(recall Definition \ref{PrpT}).
Thus, we may perform the change of variables 
$x=y^\theta(t,\xi)$ (with $t$ fixed) in the integrals 
appearing in the definition of the path length of
$\widehat{\mathbf{U}}^\theta(t)$ (see Definition \ref{lPR}).
Taking into account that both
$\widetilde{\mathbf{U}}^\theta_0$
and $\widehat{\mathbf{U}}^\theta(t)$
are equal to the same path
$\tilde{\mathbf{u}}^\theta(t)$
in the Euler variables,
we conclude that their lengths are the same:
\begin{equation}\label{pthi1}
	\left\|\widetilde{\mathbf{U}}^\theta_0\right\|
	_{\mathcal{L}}
	=\left\|\widehat{\mathbf{U}}^\theta(t)\right\|
	_{\mathcal{L}}.
\end{equation}
Moreover, since (see Remark \ref{RemKR} 
and \cite[Section 7]{KR25})
\begin{equation*}
	\left(
	U^i,V^i,W^i,Z^i
	\right)(t,\xi)=
	\left(
	\tilde{U}^i,\tilde{V}^i,\tilde{W}^i,\tilde{Z}^i
	\right)(\sigma(\xi)),\quad i=0,1,
\end{equation*}
where $\sigma(\xi)$ is given by \eqref{sigd}, we 
conclude from \eqref{appep} that the endpoints 
$\theta=0$ and $\theta=1$
of $\widetilde{\mathbf{U}}^\theta_0$
are $\widetilde{\mathbf{U}}^0$ and
$\widetilde{\mathbf{U}}^1$, that is
\begin{equation}\label{pthi2}
	\widetilde{\mathbf{U}}^i_0
	=\widetilde{\mathbf{U}}^i,
	\quad i=0,1.
\end{equation}
We have from \eqref{pthi2} and Definition \ref{Gdist} 
of the geodesic distance in $\Omega$ that
\begin{equation}\label{1in01}
	d_{\Omega_{\mathrm{tcN}}^1,K}\left(
	\widetilde{\mathbf{U}}^0,\widetilde{\mathbf{U}}^1\right)
	\leq
	\left\|\widetilde{\mathbf{U}}^\theta_0\right\|_{\mathcal{L}}
	+C\ve,
	\quad C>0.
\end{equation}
The first inequality in \eqref{appep} and
\eqref{pthi1} imply the following estimate:
\begin{equation}\label{0in1}
	d_{\Omega_{\mathrm{tcN}}^1,K}\left(
	\mathbf{U}^0(t),\mathbf{U}^1(t)\right)>
	\left\|\widetilde{\mathbf{U}}^\theta_0\right\|
	_{\mathcal{L}}-\ve.
\end{equation}
Then combining \eqref{0in1} 
and \eqref{1in01}, we obtain that
\begin{equation}\label{in01a}
	d_{\Omega_{\mathrm{tcN}}^1,K}\left(
	\mathbf{U}^0(t),\mathbf{U}^1(t)\right)>
	d_{\Omega_{\mathrm{tcN}}^1,K}\left(
	\widetilde{\mathbf{U}}^0,\widetilde{\mathbf{U}}^1\right)
	-C\ve,\quad C>0.
\end{equation}
Since $\ve>0$ is arbitrary, we 
have from \eqref{in01a} that
\begin{equation}\label{in01}
	d_{\Omega_{\mathrm{tcN}}^1,K}\left(
	\widetilde{\mathbf{U}}^0,\widetilde{\mathbf{U}}^1\right)
	\leq
	d_{\Omega_{\mathrm{tcN}}^1,K}\left(
	\mathbf{U}^0(t),\mathbf{U}^1(t)\right).
\end{equation}

The reverse inequality is obtained 
by an analogous construction.  
In this case, we build a sufficiently regular path 
$\widetilde{\mathbf{U}}_0^\theta$ 
under the ODE system \eqref{ODE} that 
approximates the geodesic distance between 
$\widetilde{\mathbf{U}}^0$ and 
$\widetilde{\mathbf{U}}^1$, as in \eqref{appep}.  
In particular, we have (cf.~the first inequality in \eqref{in01a})
\begin{equation}\label{2in01}
	d_{\Omega_{\mathrm{tcN}}^1,K}\left(
	\widetilde{\mathbf{U}}^0,\widetilde{\mathbf{U}}^1\right)
	> \left\|\widetilde{\mathbf{U}}_0^\theta\right\|_{\mathcal{L}}
	-\ve.
\end{equation}

Consider the backward evolution 
$\widetilde{\mathbf{U}}^\theta(-t)$ of the path 
$\widetilde{\mathbf{U}}_0^\theta$ for the time $t$.  
Mapping this path first into the Eulerian 
variables and then back into the 
Bressan--Constantin variables (cf.~the transform 
$\widehat{\mathbf{U}}^\theta(t)\mapsto
\tilde{\mathbf{u}}^\theta(t)\mapsto
\widetilde{\mathbf{U}}_0^\theta$ described above and illustrated in 
Figure \ref{path-m}), we obtain an initial regular path 
$\widehat{\mathbf{U}}_0^\theta$ for the ODE system \eqref{ODE}.  
The endpoints of $\widehat{\mathbf{U}}_0^\theta$ at $\theta=0,1$ 
coincide with the initial data
$\mathbf{U}_0^0$ and $\mathbf{U}_0^1$.

Evolving $\widehat{\mathbf{U}}_0^\theta$ forward in time 
by $t$ produces a path $\widehat{\mathbf{U}}^\theta(t)$ whose length 
$\|\cdot\|_{\mathcal{L}}$ coincides with that of the original path 
$\widetilde{\mathbf{U}}_0^\theta$, whose 
endpoints are
$\mathbf{U}^0(t)$ and $\mathbf{U}^1(t)$.  
Therefore (cf.~\eqref{pthi1} and \eqref{1in01}),
\begin{equation}\label{1pthi1}
\begin{split}
	&\left\|\widehat{\mathbf{U}}^\theta(t)\right\|_{\mathcal{L}}
	=\left\|\widetilde{\mathbf{U}}_0^\theta\right\|_{\mathcal{L}},\\[0.2cm]
	&d_{\Omega_{\mathrm{tcN}}^1,K}\!\left(\mathbf{U}^0(t),\mathbf{U}^1(t)\right)
	\leq 
	\left\|\widehat{\mathbf{U}}^\theta(t)\right\|_{\mathcal{L}}
	+ C\ve,
	\quad C>0.
\end{split}
\end{equation}

Combining \eqref{2in01} with \eqref{1pthi1}, 
and arguing as in \eqref{in01a}, we conclude that
\[
d_{\Omega_{\mathrm{tcN}}^1,K}\!\left(\mathbf{U}^0(t),\mathbf{U}^1(t)\right)
\leq 
d_{\Omega_{\mathrm{tcN}}^1,K}\!\left(\widetilde{\mathbf{U}}^0,
\widetilde{\mathbf{U}}^1\right),
\]
which, together with 
\eqref{in01}, yields \eqref{Lipgs}.
\end{proof}
\medskip

\textbf{Data availability statement.}
Data availability is not applicable to this article as no new data were created or analysed in this study.

\textbf{Conflict of interest.}
The authors declare no conflict of interest.


\end{document}